\documentclass[11pt,a4paper]{article}

\usepackage{amssymb,amsmath,amsthm}
\usepackage{mathrsfs}

\numberwithin{equation}{section}

\newtheorem{thm}{Theorem}[section]
\newtheorem{prop}[thm]{Proposition}
\newtheorem{lemma}[thm]{Lemma}
\newtheorem{corollary}[thm]{Corollary}

\theoremstyle{definition}
\newtheorem{definition}[thm]{Definition}

\newtheorem{convention}[thm]{Convention}

\begin{document}

\title{Positivity of metric entropy of the standard map at large parameters}

\date{November 6, 2026}

\author{Sini\v{s}a Slijep\v{c}evi\'{c}}

\maketitle

\begin{abstract}
We prove that the standard map $f_k(x,y)=(x+y+k\sin2\pi x,\,y+k\sin2\pi x)$ of the two-torus satisfies $\log(2\pi k+4)-2^{430}\le h_m(f_k)\le\log(2\pi k+4)$ for every $k\ge\exp(2^{50000})$, where $h_m$ is the metric entropy with respect to area. Thus $h_m(f_k)=\log k+O(1)$, the order of growth observed numerically since Chirikov; in particular $h_m(f_k)>0$ for all large $k$ and $h_m(f_k)\to\infty$, as conjectured by Sinai. Moreover, for each fixed $c<1$ the points whose top Lyapunov exponent is below $c\log k$ fill an area $O(1/\log k)$. Positivity of $h_m(f_k)$ for all sufficiently large $k$, with a threshold that is not effective, was obtained in a recent preprint by OpenAI \cite{locbib:OpenAI2026}. We keep some of its architecture and follow it closely in about a third of the results. New here are an explicit covering argument for cancellations and a balance across scales at a fixed parameter, an inequality between mean shortfalls that replaces the compactness argument over limits of rescaled orbit geometries; together they yield the explicit constants and the entropy bound.
\end{abstract}

\section{Introduction} \label{sec:introducti}

For $k\in\mathbb R$ let
\[
f_k(x,y)=\bigl(x+y+k\sin2\pi x,\ y+k\sin2\pi x\bigr)\pmod{\mathbb Z^2}
\]
be the standard map of the torus $\mathbb T^2=\mathbb R^2/\mathbb Z^2$. It preserves the normalized area $m$; up to a linear change of coordinates it is the map $(x,y)\mapsto(2x-y+k\sin2\pi x,\,x)$ of Berger and Turaev \cite{bib:berger2017c}, and Chirikov's standard map with stochasticity parameter $K=2\pi k$ \cite{locbib:Chirikov1979, bib:knill1999}. We write $h_m(f_k)$ for its Kolmogorov--Sinai entropy with respect to $m$, the metric entropy, and $\lambda_+(z)$ for the top Lyapunov exponent at $z\in\mathbb T^2$ (Definition~\ref{def:sta_map_ent}). Put $M=K+4$. Our main result, Theorem~\ref{thm:met_ent_the:1}, is the two-sided bound
\[
0<\log M-C_{\rm ent}\le h_m(f_k)\le\log M\qquad\text{for every }k\ge M_{\rm ent},
\]
with $C_{\rm ent}=2^{430}$ and $M_{\rm ent}=\exp(2^{50000})$. Consequently (Corollary~\ref{cor:ent_asy_and:2}) $\lvert h_m(f_k)-\log k\rvert\le C_{\rm ent}$ for these $k$, and the exponent is large on most of the torus: for every $k\ge M_{\rm ent}$ and every $c\in(0,1)$,
\[
m\{z:\lambda_+(z)<c\log M\}\le\frac{C_{\rm ent}}{(1-c)\log M}.
\]

The standard map is the paradigm of area-preserving twist maps, and through them of Hamiltonian systems with two degrees of freedom, whose return maps to a section of an energy surface preserve area \cite{locbib:MacKay1993, locbib:MacKayMeiss1987}; it is the period map of the kicked rotor and, near a resonance, a model for a general twist map \cite{locbib:Chirikov1979, locbib:Greene1979, locbib:Mather1982, locbib:Mather1984, locbib:MacKayPercival1985}. Whether the standard map has positive metric entropy is a long-standing open problem. Sinai conjectured that $h_m(f_k)>0$ for every $k>0$ and that $h_m(f_k)\to\infty$ as $k\to\infty$ \cite[p.~144]{locbib:Sinai1994}. Numerical experiments going back to Chirikov \cite{locbib:Chirikov1979} show the entropy close to $\log(K/2)$ for large $k$, and Knill \cite{bib:knill1999}, in his account of the problem, conjectures a lower bound $\log(K/2)-C(K)$ with an explicit bounded correction $C(K)$. The value $\log(K/2)$ has a simple heuristic explanation. Along an orbit the derivative is conjugate to a product of matrices whose only large entry is $2+K\cos2\pi x_i$, and the average of $\log\lvert2+K\cos2\pi x\rvert$ over the circle is exactly $\log(K/2)$ for $K\ge2$, a classical computation that extends Euler's integral. If the positions $x_i$ were independent and uniformly distributed, as for a random map, the product would grow at a rate close to $\log(K/2)$. This random phase approximation underlies Chirikov's estimate and makes Sinai's conjecture plausible. Two features make the problem hard. First, the map is strongly hyperbolic away from the critical circles $\cos2\pi x\approx0$, but near them the expanding and contracting directions can be exchanged, and such exchanges may cancel the growth accumulated before \cite{bib:blumenthal2017a}. Second, for a residual set of large $k$ the map has elliptic islands \cite{locbib:Duarte1994}: elliptic periodic points surrounded by KAM circles, which fill a set of positive area on which the Lyapunov exponent vanishes. So the map is neither uniformly hyperbolic nor ergodic, and positivity has to be proved in the presence of these islands. Progress has therefore concerned weaker properties or modified systems. For topologically generic large $k$, the stochastic sea, the set of orbits with nonzero Lyapunov exponents, has full Hausdorff dimension \cite{bib:gorodetski2010}. Random perturbations chosen afresh at each step, even of size far smaller than any power of $1/k$, produce exponents of order $\log k$ \cite{bib:blumenthal2017a}; couplings of the standard map with a fast hyperbolic system are nonuniformly hyperbolic \cite{bib:berger2012b}; and for every large $k$ the map $f_k$ has $C^\infty$-small area-preserving perturbations with positive metric entropy \cite{bib:berger2017c}. Even positivity for a set of parameters of positive measure was still conjectural \cite[Conjecture~0.4]{bib:berger2017c}, and whether $h_m(f_k)>0$ for even one value of $k$ remained open \cite{bib:knill1999} until recently. For $k\ge M_{\rm ent}$, Theorem~\ref{thm:met_ent_the:1} confirms Sinai's conjecture, including the growth: $h_m(f_k)=\log(K/2)+O(1)$, so $h_m(f_k)/\log k\to1$. As $k\to\infty$ the standard map is thus, to leading order, as chaotic as the random model. The start \eqref{eq:arg_start} of the argument is the random computation over short times, and the rest of the paper shows that the correlations along true orbits cost only a bounded amount of entropy.

Positivity of $h_m(f_k)$ for all sufficiently large $k$ was obtained recently in \cite{locbib:OpenAI2026}, by contradiction: along a hypothetical sequence of parameters $k\to\infty$ at which $\lambda_+$ vanishes almost everywhere, rescalings of the orbit geometry by the factors $2^j$ converge to limit measures whose shapes are classified, and a capped balance across scales fails in the limit. The threshold obtained there is not effective, and no lower bound on $h_m(f_k)$ of order $\log k$ is given. The present paper keeps the architecture of \cite{locbib:OpenAI2026}: its two quantitative inputs, the scarcity of cancellations and the bridge estimate, and the idea of balancing a capped shortfall across scales. About a third of our results follow \cite{locbib:OpenAI2026} closely: the setting (transfer matrices, growth and mean shortfall), the truncation of Gromov products and the doubling argument, the first steps of the proof of the scarcity of cancellations, up to and including the cubic relation, and the bridge estimate with its proof. We rewrite them in a uniform notation, repair them where needed, ground them in the earlier literature and give them explicit constants; classical facts that \cite{locbib:OpenAI2026} reproves (the Gromov hyperbolicity of the hyperbolic plane, the Pliss lemma, Gromov's approximation of hyperbolic spaces by trees) are instead quoted, or written out briefly with credit. New here are the explicit covers by short intervals, which replace the passage to Hausdorff limits in the core of the proof of the scarcity of cancellations; the short-time estimate and the explicit constants of the doubling inequality; the corrector on slack profiles, its extension to all rescaled growths, and the summable bound on the exceptional orbits; and the balance across scales at one fixed parameter, an inequality between mean shortfalls, with the bootstrap that closes the argument. Together they give the explicit entropy bound.

\subsection{Structure of the argument} \label{subsec:str_the_arg}

The core idea is to measure how far the derivatives of the iterates of $f_k$ fall short of full growth, and to prove that this \textbf{shortfall} is uniformly bounded by $C_{\rm ent}/\log M$. Along an orbit, the derivative of $f_k^n$ is conjugate to a product of $n$ transfer matrices, one for each step, of determinant one and norm at most $M$ (Section~\ref{sec:tra_mat_and}). We write $T_{a,b}$ for the product of the transfer matrices over the times from $a$ to $b$, and $g(a,b)=\log\lVert T_{a,b}\rVert/\log M$ for its \textbf{growth}, which lies between $0$ and $\lvert b-a\rvert$. Both depend on the initial point $z\in\mathbb T^2$ of the orbit; we write $\mathbb E$ for the average over $z$ with respect to $m$, and $\mathbb P$ for the area of a set of initial points. The \textbf{mean shortfall} at time $n$ is
\begin{equation}\label{eq:arg_short}
e_n=\mathbb E\Bigl[1-\frac{g(0,n)}n\Bigr],
\end{equation}
the relative loss against the largest possible growth $M^n$, averaged over the initial point (Definition~\ref{def:gro_sem_orb}). By Pesin's entropy formula \cite{locbib:Pesin1977, locbib:Mane1981} (Proposition~\ref{prop:met_ent_fro}),
\begin{equation}\label{eq:arg_pesin}
h_m(f_k)=\Bigl(1-\lim_{n\to\infty}e_n\Bigr)\log M,
\end{equation}
thus establishing that $e_n$ is uniformly bounded by $C_{\rm ent}/\log M$ suffices.

We prove three bounds along $n=2^j$, with $j\ge0$ an integer. Here and below $c$, $C$ and $C'$ denote positive constants, and $\beta$, $\gamma$ and $\delta$ small positive constants, all independent of $k$ and $n$; the precise constants are explicit and listed in Appendix~\ref{sec:not_and_con}. First, the shortfall starts small, by a direct computation based on Euler's integral of the logarithm of the cosine:
\begin{equation}\label{eq:arg_start}
e_n\le\frac3{\log M}\qquad\text{for }n\le\sqrt K.
\end{equation}
Second, one doubling increases it by at most a constant factor:
\begin{equation}\label{eq:arg_doub}
e_{2n}\le C\Bigl(e_n+\frac1n\Bigr).
\end{equation}
Third, over many doublings, we establish the core estimate of the paper (we come back to its rationale shortly): the shortfall is controlled by its initial value, up to a quadratic term,
\begin{equation}\label{eq:arg_bal}
e_N\le C\Bigl(e_{n_0}+\frac1{n_0}\Bigr)+C'e_N^2\qquad(n_0<N).
\end{equation}
Together the three bounds give the uniform bound (Proposition~\ref{prop:uni_bou_the:2}). Neither of the last two suffices alone: iterating \eqref{eq:arg_doub} multiplies its constant over unboundedly many doublings, and \eqref{eq:arg_bal} says nothing when $e_N$ is large.

The doubling bound \eqref{eq:arg_doub} rests on an exact identity, which uses only the invariance of area (Lemma~\ref{lem:exa_dou_ide}):
\begin{equation}\label{eq:arg_ident}
e_{2n}-e_n=\frac1{2n}\,\mathbb E\bigl[g(0,n)+g(n,2n)-g(0,2n)\bigr].
\end{equation}
The bracket is the growth lost when the products over the two halves are multiplied, so the shortfall increases only through cancellation at the midpoint. It is twice the \textbf{Gromov product} $(0|2n)_n=\frac12\bigl(g(n,0)+g(n,2n)-g(0,2n)\bigr)$ of the times $0$ and $2n$ based at $n$ (Definition~\ref{def:gro_hyp_sem}). To bound it we use the Gromov hyperbolicity of $g$ (Lemma~\ref{lem:hyp_the_gro}). A matrix of determinant one moves a base point of the hyperbolic plane by twice the logarithm of its norm, so $g$ is, up to the factor $2\log M$, the distance between the points of a path in the hyperbolic plane, which is Gromov hyperbolic \cite{bib:nica2014a}. This lets a large Gromov product at $n$ be passed, up to the loss of growth on the halves, to the shorter stretches $[n-s,n]$ and $[n,n+s]$ (Lemma~\ref{lem:gro_pro_und}). Hence a large cancellation at the midpoint forces either a loss of growth on one of the two halves, whose mean is controlled by $e_n$, or, for some $s\le n$, a near-total cancellation of two adjacent stretches $[n-s,n]$ and $[n,n+s]$ of nearly full growth, which is rare by \eqref{eq:arg_scar} below.

The third bound \eqref{eq:arg_bal}, the balance across scales (Proposition~\ref{prop:bal_bet_dya}), is the heart of the paper. Summing the increments $e_{2n}-e_n$ over the doublings from $n_0$ to $N$ only returns $e_N-e_{n_0}$; the idea is to compare them with the increments of a second quantity, the \textbf{corrector}, which absorbs most of them. The corrector is computed from the shape of the growth at scale $n$, that is, from the \textbf{rescaled growth} $g_n(u,v)=g(\lfloor nu\rfloor,\lfloor nv\rfloor)/n$ for real $u,v$ (Definition~\ref{def:res_gro_arr:1}), and it has two properties: from one doubling to the next its mean changes by nearly as much as the shortfall, and it never exceeds the fraction $1-c$ of the shortfall (Definition~\ref{def:ext_cor_and}; in fact the comparison uses a capped shortfall, in which all large losses count the same). The \textbf{balance defect} $Z(g_n)$ at scale $n$ is a function of the initial point whose mean is the increment of the shortfall from $n$ to $2n$ minus that of the corrector. The shortfall minus the corrector is at least $c$ times the shortfall, and from $n_0$ to $N$ it changes only by the balance defects; hence, summing over $n=2^j$,
\begin{equation}\label{eq:arg_tele}
c\,e_N\le\sum_{n_0\le n<N}\mathbb E\,Z(g_n)+C e_{n_0}.
\end{equation}
The corrector is built so that the balance defect is at most zero when the rescaled growth has one of four simple shapes, the \textbf{slack profiles} (Definition~\ref{def:slac_prof}): growth at a constant rate; slow growth throughout; or growth at full rate on one side, or on both sides, of a slow stretch. By the Gromov hyperbolicity of $g$ and Gromov's approximation of hyperbolic spaces by trees \cite{locbib:GhysHarpe1990}, the rescaled growth is nearly the distance between the points of a path in a tree, and for most initial points it is close to a slack profile. The balance defects therefore come from the remaining, exceptional initial points, which the two quantitative inputs below make rare: summed over the doublings, the balance defects are at most a small multiple of $e_N$, plus $2/n_0$, plus $C'e_N^2$ (Proposition~\ref{prop:def_the_ext:1}). With \eqref{eq:arg_tele} this gives \eqref{eq:arg_bal}.

The first input is the scarcity of cancellations (Theorem~\ref{thm:exp_sca_can}): two adjacent stretches of nearly full growth almost never nearly cancel, that is, the product over their union almost never loses nearly all of their growth. For every time $i$ and every $n\ge1$,
\begin{equation}\label{eq:arg_scar}
\mathbb P\Bigl(g(i-n,i),\,g(i,i+n)\ge(1-\delta)n,\ \ g(i-n,i+n)\le\delta n\Bigr)\le Ce^{-\gamma n}.
\end{equation}
It is proved by excluding, at many successive times, all but a small fraction of the possible values of one coordinate of the orbit at time $i$, as in the parameter exclusion of Benedicks and Carleson \cite{locbib:BenedicksCarleson1991}.

The second input is the bridge estimate (Theorem~\ref{thm:bri_pro_est:2}): a stretch of time $[a,b]$ of nearly full growth makes the behaviour of the orbit near its two ends almost independent. A typical instance is the following. Fix integers $a<b$ with $b-a\ge1100$ and a length $1\le\ell\le(b-a)/20$, and for $0<t<1$ let
\[
A_t=\bigl\{z\in\mathbb T^2:\ g(a-\ell,a)<t\ell\bigr\},\qquad B_t=\bigl\{z\in\mathbb T^2:\ g(b,b+\ell)<t\ell\bigr\}
\]
be the sets of initial points whose orbits lose more than the fraction $1-t$ of the possible growth during the $\ell$ steps before time $a$, respectively after time $b$. Then
\begin{equation}\label{eq:arg_bridge}
\mathbb P\bigl(A_{1/2}\cap B_{1/2}\cap\{g(a,b)\ge(1-\beta)(b-a)\}\bigr)\le C\,\mathbb P(A_{3/4})\,\mathbb P(B_{3/4})+Ce^{-(b-a)}.
\end{equation}
The theorem allows, in place of the growth on $[a-\ell,a]$ and on $[b,b+\ell]$, any Lipschitz function of the growth on finitely many intervals of time whose endpoints lie within $(b-a)/20$ of $a$, respectively of $b$ (Definition~\ref{def:shor_obse}). At a middle time, curves contracted by the forward dynamics up to time $b$ and curves contracted by the backward dynamics up to time $a$ form a grid; $B_{1/2}$ is nearly a union of curves of the first family and $A_{1/2}$ of the second, and the product bound follows from a change of variables. The estimate enters through the exceptions of the balance. A place near one end of a stretch of nearly full growth where the rescaled growth turns back or changes speed, and a loss of growth near the other end, are each rare, with probability about $e_N$; their coincidence keeps the rescaled growth away from the slack profiles, and \eqref{eq:arg_bridge} bounds it by the product, which is the quadratic term $C'e_N^2$. The other way to stay away from the slack profiles is a near-total cancellation, which \eqref{eq:arg_scar} controls.

\subsection{Structure of the paper} \label{subsec:str_the_pap}

Section~\ref{sec:tra_mat_and} is dedicated to the reduction \eqref{eq:arg_pesin}. It sets up the transfer matrices, the growth $g$ and the mean shortfall \eqref{eq:arg_short} (Definition~\ref{def:gro_sem_orb}), conjugates the derivatives of $f_k^n$ to products of transfer matrices, and turns the entropy into the mean growth by the multiplicative ergodic theorem and Pesin's formula, which gives \eqref{eq:arg_pesin} (Proposition~\ref{prop:met_ent_fro}). It also collects three elementary tools used in Section~\ref{sec:the_pro_est} and Appendix~\ref{sec:scar_canc}: the lemma of Pliss \cite{locbib:Pliss1972}, bounds for the derivatives of the map, and the inverse of a tridiagonal matrix.

Section~\ref{sec:hyp_and_the} proves the first two bounds. The start \eqref{eq:arg_start} is a direct computation with Euler's integral (Lemma~\ref{lem:mea_sho_sho}). The identity \eqref{eq:arg_ident} follows from the invariance of area (Lemma~\ref{lem:exa_dou_ide}), and the doubling bound \eqref{eq:arg_doub} follows from the identity (Corollary~\ref{cor:dou_ine_for}): Gromov hyperbolicity (Lemma~\ref{lem:hyp_the_gro}, Lemma~\ref{lem:gro_pro_und}) reduces a large cancellation at the midpoint to a loss on the halves or to a near-total cancellation, and the scarcity of cancellations \eqref{eq:arg_scar}, restated there, makes the latter rare.

Section~\ref{sec:the_pro_est} proves the bridge estimate \eqref{eq:arg_bridge} (Theorem~\ref{thm:bri_pro_est:2}). At a middle time of a stretch of nearly full growth it builds two transverse families of curves, the finite-time stable and unstable curves, contracted by the forward dynamics up to one end of the stretch and by the backward dynamics up to the other; sets of initial points defined by the growth near the two ends, such as $A_t$ and $B_t$, are nearly unions of curves of one family each, and comparing them becomes a change of variables with bounded Jacobian.

Section~\ref{sec:bal_acr_sca} proves the balance \eqref{eq:arg_bal} (Proposition~\ref{prop:bal_bet_dya}), in the form that assumes $e_N$ small instead of carrying the quadratic term. It introduces the rescaled growth (Definition~\ref{def:res_gro_arr:1}), the capped shortfall (Definition~\ref{def:the_cap_and}), the slack profiles (Definition~\ref{def:slac_prof}), and the corrector with the balance defect (Definition~\ref{def:ext_cor_and}), and derives \eqref{eq:arg_tele} in the proof of Proposition~\ref{prop:bal_bet_dya}. It then shows that the balance defect is at most zero on slack profiles (Lemma~\ref{lem:pro_the_bal}), and that the rescaled growth is close to a slack profile outside an exceptional set, which \eqref{eq:arg_scar} and \eqref{eq:arg_bridge} keep small (Lemma~\ref{lem:sum_exc_mas:1}); this bounds the summed balance defect (Proposition~\ref{prop:def_the_ext:1}).

Section~\ref{sec:the_boo_and} combines the three bounds into the uniform bound (Proposition~\ref{prop:uni_bou_the:2}) and proves the main results, Theorem~\ref{thm:met_ent_the:1} and Corollary~\ref{cor:ent_asy_and:2}.

Appendix~\ref{sec:scar_canc} proves the scarcity of cancellations \eqref{eq:arg_scar} (Theorem~\ref{thm:exp_sca_can}). A cancellation at time $i$ produces many pairs of hyperbolic times, one on each side of $i$ and with comparable growth, a hyperbolic time being a time up to which the growth from every earlier time has a uniform rate; at each pair a cubic relation confines one coordinate of the orbit at time $i$ to a few short intervals; and these confinements nest, as in the parameter exclusion of Benedicks and Carleson. Appendix~\ref{sec:not_and_con} lists the notation and the constants shared by several sections. All constants are explicit, and a first reading may ignore them.

\vspace{1ex} \noindent \textbf{Statement on the use of AI.} This paper was prepared with the assistance of Aleph, an AI research assistant for mathematics, using an ensemble of models from Anthropic, OpenAI and Google (Gemini) together with bespoke models. The author thanks Cantab Pi for access to the beta version of Aleph and for computing resources. The author takes full responsibility for the content of the paper and for all mathematical results.

\section{Transfer matrices and the reduction to a bound on the mean shortfall} \label{sec:tra_mat_and}

This section sets up the transfer matrices, the growth semimetric and the mean shortfall, and proves Proposition~\ref{prop:met_ent_fro}, which expresses $h_m(f_k)$ through the mean shortfall. The derivative cocycle of $f_k$ is conjugate to the transfer matrices of the Jacobi operator with off-diagonal entries $-1$ and diagonal entries $v(q_i)$ along the orbit \cite{locbib:Teschl2000}; we first set up these matrices.

\begin{definition}
Let $\mathbb{T}^2=\mathbb{R}^2/\mathbb{Z}^2$ and let $m$ be normalized Lebesgue measure (area) on $\mathbb{T}^2$. For $k\in\mathbb{R}$ the \textbf{standard map} is the map
\begin{equation}\label{eq:standard_map}
f_k(x,y)=\bigl(x+y+k\sin(2\pi x),\ y+k\sin(2\pi x)\bigr)\pmod{\mathbb{Z}^2}.
\end{equation}
The \textbf{metric entropy} $h_m(f_k)$ is the Kolmogorov--Sinai entropy of $f_k$ with respect to $m$. For $z\in\mathbb{T}^2$ the \textbf{top Lyapunov exponent} is
\[
\lambda_+(z)=\lim_{n\to\infty}\frac1n\log\bigl\lVert Df_k^n(z)\bigr\rVert
\]
wherever the limit exists, with natural logarithm and any fixed norm on the tangent spaces.
\label{def:sta_map_ent}
\end{definition}

\begin{convention}
Let the parameter $k$ of Definition~\ref{def:sta_map_ent} be positive, and put
\begin{equation}\label{eq:large_parameters}
K=2\pi k,\qquad M=K+4,\qquad \phi(x)=2x+k\sin(2\pi x),\qquad v(x)=\phi'(x)=2+K\cos(2\pi x).
\end{equation}
We call $\phi$ the \textbf{recurrence function}.
\label{conv:sta_map_par}
\end{convention}

\begin{convention}
Integrals over $\mathbb{T}^2$ are taken with respect to the normalized Lebesgue measure $m$ of Definition~\ref{def:sta_map_ent}. For a measurable set $E\subset\mathbb{T}^2$ and a measurable function $F\colon\mathbb{T}^2\to\mathbb{R}$, nonnegative or integrable, we write
\[
\mathbb{P}(E)=m(E),\qquad \mathbb{E}F=\int_{\mathbb{T}^2}F\,dm,
\]
and $\mathbb{P}\{\cdots\}=m\{z\in\mathbb{T}^2:\cdots\}$ for the set defined by a condition.
\label{conv:mea_and_int}
\end{convention}

\begin{definition}
Fix $k$ and $z\in\mathbb{T}^2$, write $(x_i,y_i)=f_k^i(z)$ for $i\in\mathbb{Z}$ (Definition~\ref{def:sta_map_ent}) and put $q_i=x_i\in\mathbb{T}$, the \textbf{position} at time $i$. The orbit and all quantities defined from it below are regarded as functions of the initial point $z$. Let $\Pi(x,y)=(x,x-y)$, as a map of $\mathbb{T}^2$ and, by the same formula, as a linear map of $\mathbb{R}^2$. The \textbf{state} of the orbit at time $i$ is $(q_i,q_{i-1})$. With $\phi$ and $v$ as in Convention~\ref{conv:sta_map_par}, the \textbf{orbit recurrence} is the relation
\begin{equation}\label{eq:orbit_recurrence}
q_{i+1}=\phi(q_i)-q_{i-1}\pmod 1.
\end{equation}
A \textbf{real lift} of the orbit is a choice of real initial lifts iterated by the same recurrence over $\mathbb{R}$ without reduction modulo one. A \textbf{linearized solution} along the orbit is a real sequence $(u_i)_{i\in\mathbb{Z}}$ with
\begin{equation}\label{eq:linearized_matrix}
\binom{u_{i+1}}{u_i}=A_i\binom{u_i}{u_{i-1}},\qquad A_i=\begin{pmatrix}v(q_i)&-1\\ 1&0\end{pmatrix}
\end{equation}
($v$ has period one, so $v(q_i)$ is defined for $q_i\in\mathbb{T}$). For $a,b\in\mathbb{Z}$ the \textbf{transfer matrix} $T_{a,b}$ is the linear map carrying $(u_a,u_{a-1})$ to $(u_b,u_{b-1})$ for every linearized solution.
\label{def:orb_sta_and}
\end{definition}

The next three definitions concern growth. The first is the shortfall of a function of two times on an interval, its relative loss against unit speed; the second concerns semimetrics, Gromov products and Gromov hyperbolicity, which are used from the section on the doubling inequality onwards; the third introduces the growth semimetric of an orbit and its mean shortfall.

\begin{definition}
Let $Y\subset\mathbb{R}$ and $d\colon Y\times Y\to\mathbb{R}$. For $s<t$ in $Y$ the \textbf{shortfall} of $d$ on the interval $[s,t]$ is
\[
h_d(s,t)=1-\frac{d(s,t)}{t-s},\qquad h_d(t,s)=h_d(s,t),
\]
and the \textbf{deficit} of $d$ on $[s,t]$ is $(t-s)-d(s,t)=(t-s)\,h_d(s,t)$.
\label{def:sho_sem_the}
\end{definition}

\begin{definition}
Let $X$ be a set and $d$ a \textbf{semimetric} on $X$: $d(x,y)=d(y,x)\ge0$, $d(x,x)=0$ and $d(x,u)\le d(x,y)+d(y,u)$ for all $x,y,u\in X$; distinct points may be at distance zero. For a basepoint $w\in X$, the \textbf{Gromov product} based at $w$ is
\[
(x|y)_w=\tfrac12\bigl(d(w,x)+d(w,y)-d(x,y)\bigr)\qquad(x,y\in X).
\]
For $\delta\ge0$, $(X,d)$ is \textbf{$\delta$-hyperbolic at $w$} if
\[
(x|y)_w\ge\min\bigl\{(x|u)_w,(y|u)_w\bigr\}-\delta\qquad(x,y,u\in X),
\]
and \textbf{$\delta$-hyperbolic} if it is $\delta$-hyperbolic at every basepoint $w\in X$; it is \textbf{$0$-hyperbolic} if it is $\delta$-hyperbolic with $\delta=0$. The definition at a basepoint is \cite[Section 2.2, Definition 2.10]{bib:bestvina1997}, where the notion is attributed to Gromov and the Gromov product is written $(x\cdot y)$; it is stated there for metric spaces, and for a semimetric it is applied after identifying points at distance zero.
\label{def:gro_hyp_sem}
\end{definition}

\begin{definition}
With the transfer matrices $T_{a,b}$ of Definition~\ref{def:orb_sta_and} and $M$ as in Convention~\ref{conv:sta_map_par}, the \textbf{growth semimetric} on orbit time and the \textbf{mean shortfall} at length $n\ge1$ are
\begin{equation}\label{eq:growth_orbit}
g(i,j)=\log_M\lVert T_{i,j}\rVert\quad(i,j\in\mathbb{Z}),\qquad e_n=\mathbb{E}\,h_g(0,n)=\mathbb{E}\Bigl[1-\frac{g(0,n)}{n}\Bigr],
\end{equation}
where $h_g$ is the shortfall of $g$ on $Y=\mathbb{Z}$ (Definition~\ref{def:sho_sem_the}). Here $g$ depends on the initial point $z\in\mathbb{T}^2$ of the orbit, and $\mathbb{E}$ is integration over $z$ with respect to normalized Lebesgue measure on $\mathbb{T}^2$ (Convention~\ref{conv:mea_and_int}). For a Borel function $\Psi\colon\mathbb{R}\to\mathbb{R}$ for which the integral exists, put
\[
e^\Psi_n=\mathbb{E}\,\Psi\bigl(h_g(0,n)\bigr);
\]
thus $e_n=e^\Psi_n$ for $\Psi(u)=u$.
\label{def:gro_sem_orb}
\end{definition}

The notation used throughout the paper and the constants shared by several sections, with their values and roles, are collected in Convention~\ref{conv:not_and_con} in Appendix~\ref{sec:not_and_con}; a reader may consult it when a symbol or a constant first appears.

\begin{lemma}
 For every orbit and all $a,b,l\in\mathbb{Z}$:
\begin{enumerate}
\item[(i)] $T_{a,l}=T_{b,l}T_{a,b}$, $T_{b,a}=T_{a,b}^{-1}$ and $T_{a,a}=I$; $\det T_{a,b}=1$, and
\begin{equation}\label{eq:transfer_norm_bounds}
\lVert T_{a,b}\rVert=\lVert T_{b,a}\rVert,\qquad 1\le\lVert T_{a,b}\rVert\le M^{\lvert b-a\rvert};
\end{equation}
the upper bound also holds for the operator norm $\lVert\cdot\rVert_\infty$ induced by the maximum norm on $\mathbb{R}^2$;
\item[(ii)] $g$ is a semimetric on $\mathbb{Z}$: $g(a,b)=g(b,a)$, $g(a,a)=0$, $g(a,l)\le g(a,b)+g(b,l)$, and
\begin{equation}\label{eq:growth_bounds}
0\le g(a,b)\le\lvert a-b\rvert;
\end{equation}
in particular $0\le h_g(a,b)\le1$ for $a\ne b$.
\end{enumerate}
\label{lem:nor_and_gro}
\end{lemma}

\begin{proof}
We first prove the three identities of (i). Each $A_i$ is invertible (its determinant is one), so by \eqref{eq:linearized_matrix} every pair in $\mathbb{R}^2$ is $(u_a,u_{a-1})$ for exactly one linearized solution; by the definition of the transfer matrices (Definition~\ref{def:orb_sta_and}), $T_{a,l}$ and $T_{b,l}T_{a,b}$ both carry this pair to $(u_l,u_{l-1})$, hence $T_{a,l}=T_{b,l}T_{a,b}$; $T_{a,a}=I$, and then $T_{b,a}T_{a,b}=T_{a,a}=I$. Moreover $T_{i,i+1}=A_i$ by \eqref{eq:linearized_matrix}, so $T_{a,b}=A_{b-1}\cdots A_a$ for $b>a$.

(i) Each $A_i$ has determinant one, hence so does every $T_{a,b}$. A real $2\times2$ matrix of determinant one has two singular values with product one, the larger being its norm, and its inverse has the same singular values; hence $\lVert T_{a,b}\rVert=\lVert T_{a,b}^{-1}\rVert=\lVert T_{b,a}\rVert\ge1$. Write
\[
A_i=v(q_i)E_1+J,\qquad A_i^{-1}=v(q_i)E_2-J,\qquad E_1=\begin{pmatrix}1&0\\0&0\end{pmatrix},\quad E_2=\begin{pmatrix}0&0\\0&1\end{pmatrix},\quad J=\begin{pmatrix}0&-1\\1&0\end{pmatrix}.
\]
The matrices $E_1$, $E_2$ and $J$ have norm one for the Euclidean and for the maximum norm, and $\lvert v\rvert\le K+2$; so $\lVert A_i^{\pm1}\rVert\le K+3<M$ in both norms. For $b>a$, $T_{a,b}$ is a product of $b-a$ factors $A_i$ and $T_{b,a}$ a product of $b-a$ factors $A_i^{-1}$, and submultiplicativity in either norm gives the upper bound.

(ii) Taking $\log_M$ in (i) gives $g(a,b)=g(b,a)$ and \eqref{eq:growth_bounds}; $g(a,a)=\log_M\lVert I\rVert=0$; and $T_{a,l}=T_{b,l}T_{a,b}$ gives $\lVert T_{a,l}\rVert\le\lVert T_{b,l}\rVert\,\lVert T_{a,b}\rVert$, which is the triangle inequality after taking $\log_M$. For $a\ne b$, the definition of the shortfall and the symmetry of $g$ give $h_g(a,b)=1-g(a,b)/\lvert a-b\rvert$, so \eqref{eq:growth_bounds} gives $0\le h_g(a,b)\le1$.
\end{proof}

\begin{lemma}

\begin{enumerate}
\item[(i)] $f_k$ is a homeomorphism of $\mathbb{T}^2$ and preserves $m$: $\mathbb{E}[F\circ f_k^s]=\mathbb{E}F$ for every $s\in\mathbb{Z}$ and every measurable $F\colon\mathbb{T}^2\to\mathbb{R}$, nonnegative or integrable.
\item[(ii)] For every $i\in\mathbb{Z}$, $(q_i,q_{i-1})=\Pi(f_k^iz)$, and the state map $z\mapsto(q_i,q_{i-1})$ is a bijection of $\mathbb{T}^2$ preserving $m$. The positions satisfy the orbit recurrence \eqref{eq:orbit_recurrence}, and every real lift of the orbit reduces modulo one to $(q_i)_{i\in\mathbb{Z}}$.
\item[(iii)] For all $i,s\in\mathbb{Z}$, $q_i\circ f_k^s=q_{i+s}$; consequently $T_{a,b}\circ f_k^s=T_{a+s,b+s}$ and $g(a,b)\circ f_k^s=g(a+s,b+s)$ for all $a,b\in\mathbb{Z}$.
\end{enumerate}
\label{lem:inv_are_und}
\end{lemma}

\begin{proof}
(i) Since $\sin(2\pi x)$ has period one, \eqref{eq:standard_map} defines a smooth map of $\mathbb{T}^2$, with smooth inverse $(x,y)\mapsto\bigl(x-y,\ y-k\sin(2\pi(x-y))\bigr)$; in particular $f_k$ is a homeomorphism. Its Jacobian determinant is $\bigl(1+2\pi k\cos(2\pi x)\bigr)-2\pi k\cos(2\pi x)=1$, so $f_k$ preserves $m$ by the change of variables formula, and hence so does $f_k^s$ for every $s\in\mathbb{Z}$.

(ii) Write $(x_i,y_i)=f_k^iz$. By \eqref{eq:standard_map}, $x_{i+1}-y_{i+1}=x_i$ for every $i$, so $x_i-y_i=x_{i-1}$ and $(q_i,q_{i-1})=(x_i,x_i-y_i)=\Pi(x_i,y_i)$. Moreover $q_{i+1}=x_i+y_i+k\sin(2\pi x_i)=2x_i-x_{i-1}+k\sin(2\pi x_i)=\phi(q_i)-q_{i-1}$ modulo one, which is \eqref{eq:orbit_recurrence}. A real lift of the orbit satisfies the same recurrence over $\mathbb{R}$ and starts from lifts of $q_0$ and $q_{-1}$, so by induction in both directions of time it reduces modulo one to $(q_i)_{i\in\mathbb{Z}}$. Finally, $\Pi$ is induced by the integer matrix $\begin{pmatrix}1&0\\1&-1\end{pmatrix}$ of determinant $-1$, so it is a bijection of $\mathbb{T}^2$ preserving $m$, and by (i) so is the state map $\Pi\circ f_k^i$.

(iii) The orbit of $f_k^sz$ is $(f_k^{i+s}z)_{i\in\mathbb{Z}}$, so $q_i(f_k^sz)=q_{i+s}(z)$ for every $i$. Since $A_i$ depends on $q_i$ only, $A_i\circ f_k^s=A_{i+s}$. By \eqref{eq:linearized_matrix}, $T_{i,i+1}=A_i$, so the identities of Lemma~\ref{lem:nor_and_gro}(i) give, by induction, $T_{a,b}=A_{b-1}\cdots A_a$ for $b>a$, together with $T_{a,b}=T_{b,a}^{-1}$ for $b<a$ and $T_{a,a}=I$; hence $T_{a,b}\circ f_k^s=T_{a+s,b+s}$, and taking $\log_M$ of the norms gives $g(a,b)\circ f_k^s=g(a+s,b+s)$.
\end{proof}

The classical input is collected in the next proposition; it is the only place where smooth ergodic theory enters.

\begin{prop}
 Let $k>0$. For $m$-almost every $z\in\mathbb{T}^2$ the top Lyapunov exponent $\lambda_+(z)$ exists and lies in $[0,\log M]$. The limit $e_\infty=\lim_{n\to\infty}e_n$ exists, $0\le e_\infty\le1$, and
\[
h_m(f_k)=\int_{\mathbb{T}^2}\lambda_+\,dm=(1-e_\infty)\log M.
\]
\label{prop:met_ent_fro}
\end{prop}

\begin{proof}
Let $F(q,q')=(\phi(q)-q',\,q)\pmod{\mathbb{Z}^2}$ on $\mathbb{T}^2$, with $\phi$ as in Convention~\ref{conv:sta_map_par}, and let $\Pi$ be as in Definition~\ref{def:orb_sta_and}; as a linear map $\Pi=\begin{pmatrix}1&0\\1&-1\end{pmatrix}$, so $\Pi^2=I$. Substitution in \eqref{eq:standard_map} gives $\Pi\circ f_k=F\circ\Pi$, hence $f_k^n=\Pi\circ F^n\circ\Pi$ for $n\ge1$. By Lemma~\ref{lem:inv_are_und}(ii), $\Pi(f_k^iz)=(q_i,q_{i-1})$ and $F(q_i,q_{i-1})=(q_{i+1},q_i)$, while $DF(q_i,q_{i-1})=A_i$ by \eqref{eq:linearized_matrix}, since $\phi'=v$. The chain rule and Lemma~\ref{lem:nor_and_gro}(i) therefore give
\[
Df_k^n(z)=\Pi\,A_{n-1}\cdots A_0\,\Pi=\Pi\,T_{0,n}\,\Pi\qquad(z\in\mathbb{T}^2,\ n\ge1).
\]
Since $\Pi^{-1}=\Pi$, the operator norms (Convention~\ref{conv:not_and_con}) of $Df_k^n(z)$ and $T_{0,n}$ differ by a factor at most $\lVert\Pi\rVert^2$, and with $\log\lVert T_{0,n}\rVert=(\log M)\,g(0,n)$ (Definition~\ref{def:gro_sem_orb}),
\[
\bigl\lvert\log\lVert Df_k^n(z)\rVert-(\log M)\,g(0,n)\bigr\rvert\le2\log\lVert\Pi\rVert\qquad(z\in\mathbb{T}^2,\ n\ge1).
\]

The map $f_k$ is a real-analytic diffeomorphism of $\mathbb{T}^2$: it is real-analytic by \eqref{eq:standard_map}, and so is its inverse $(x,y)\mapsto\bigl(x-y,\ y-k\sin(2\pi(x-y))\bigr)$; it preserves $m$ by Lemma~\ref{lem:inv_are_und}(i). By the multiplicative ergodic theorem \cite{locbib:Oseledets1968} (see \cite{locbib:BarreiraPesin2007}), for $m$-almost every $z$ the limit $\lambda_+(z)$ of Definition~\ref{def:sta_map_ent}, which does not depend on the norm and is taken here with the operator norm, exists and is the largest Lyapunov exponent of $f_k$ at $z$; since $\det Df_k=1$, the Lyapunov exponents at such $z$, counted with multiplicity, are $\lambda_+(z)$ and $-\lambda_+(z)$, so $\lambda_+(z)\ge0$ and $\lambda_+(z)$ is the sum of the nonnegative exponents. Pesin's entropy formula \cite{locbib:Pesin1977} (see also \cite{locbib:Mane1981, locbib:Mane1983} and \cite{locbib:BarreiraPesin2007}) then gives $h_m(f_k)=\int_{\mathbb{T}^2}\lambda_+\,dm$. Dividing the last display by $n$ and using $0\le g(0,n)\le n$ (\eqref{eq:growth_bounds} in Lemma~\ref{lem:nor_and_gro}(ii)), we get at $m$-almost every $z$
\[
\lambda_+(z)=(\log M)\lim_{n\to\infty}\frac{g(0,n)}{n}\in[0,\log M].
\]
Finally $e_n=1-\mathbb{E}[g(0,n)/n]$ (Definition~\ref{def:gro_sem_orb}) with $0\le g(0,n)/n\le1$, so dominated convergence gives
\[
e_\infty=\lim_{n\to\infty}e_n=1-\frac1{\log M}\int_{\mathbb{T}^2}\lambda_+\,dm\in[0,1],
\]
which is the asserted formula.
\end{proof}

The rest of the paper proves the uniform bound Proposition~\ref{prop:uni_bou_the:2}: if $M\ge M_{\rm ent}$, then $e_n\le C_{\rm ent}/\log M$ for every $n=2^j$. With Proposition~\ref{prop:met_ent_fro} it gives the main results; both deductions are made in Section~\ref{sec:the_boo_and}, after the bound itself.

\subsection{Three elementary tools} \label{subsec:thr_ele_too}

The bridge estimate (Section~\ref{sec:the_pro_est}) and the scarcity of cancellations (Appendix~\ref{sec:scar_canc}) share a convention on wedge products and directions and three elementary tools: the lemma of Pliss \cite{locbib:Pliss1972} in counting form, which turns a small total deficit of a sequence into many indices of uniform growth; bounds for the derivatives of the recurrence function; and the Green matrix of a finite Jacobi matrix \cite{locbib:Teschl2000}.

\begin{convention}
For $a,b\in\mathbb{R}^2$, $a\wedge b$ denotes the determinant of the matrix with columns $a,b$; for unit vectors, $\lvert a\wedge b\rvert$ is the absolute sine of the angle between them and depends only on the lines $\mathbb{R}a,\mathbb{R}b$. A \textbf{direction} is a line through the origin in $\mathbb{R}^2$, represented by a unit vector determined up to sign.
\label{conv:wed_pro_and}
\end{convention}

\begin{lemma}
Let $x_0,\dots,x_N$ be real numbers with $x_j-x_{j-1}\le1$ for $1\le j\le N$, and put
\[
\Gamma=N-(x_N-x_0)=\sum_{j=1}^N\bigl(1-x_j+x_{j-1}\bigr)
\]
(the \emph{total deficit} of the sequence).
For $0<c<1$, the number of indices $j\in\{1,\dots,N\}$ for which some $0\le l<j$ satisfies $x_j-x_l<c(j-l)$ is at most $\Gamma/(1-c)$. The same bound holds for the number of indices $j\in\{0,\dots,N-1\}$ for which some $j<l\le N$ satisfies $x_l-x_j<c(l-j)$.
\label{lem:pli_typ_cou}
\end{lemma}

\begin{proof}
This is the counting form of the lemma of Pliss \cite{locbib:Pliss1972}; compare \cite[Lemma 3.1]{locbib:AlvesBonattiViana2000}.

Put $w_j=x_j-cj$ and $R_j=\max_{0\le l\le j}w_l$ for $0\le j\le N$, and call an index $j\in\{1,\dots,N\}$ \emph{backward bad} if some $0\le l<j$ satisfies $x_j-x_l<c(j-l)$. Since $x_j-x_l<c(j-l)$ means $w_j<w_l$, an index $j$ is backward bad exactly when $w_j<R_{j-1}$.

We show that the running maximum $R_j$ increases only at indices that are not backward bad, and by at most $1-c$ there. At a backward bad index $R_j=R_{j-1}$. At any other index $R_j=w_j$, so that
\[
R_j-R_{j-1}\le w_j-w_{j-1}=x_j-x_{j-1}-c\le1-c.
\]
Let $\mathcal{G}$ be the set of indices in $\{1,\dots,N\}$ that are not backward bad. Summing over $j$, and using $R_0=w_0$ and $R_N\ge w_N$,
\[
(1-c)N-\Gamma=w_N-w_0\le R_N-R_0\le(1-c)\,\lvert\mathcal{G}\rvert,
\]
so the number $N-\lvert\mathcal{G}\rvert$ of backward bad indices is at most $\Gamma/(1-c)$.

For the forward count, put $\widehat x_i=-x_{N-i}$ $(0\le i\le N)$. Then $\widehat x_i-\widehat x_{i-1}=x_{N-i+1}-x_{N-i}\le1$ and $N-(\widehat x_N-\widehat x_0)=N-(x_N-x_0)=\Gamma$. An index $j$ with a witness $j<l\le N$, $x_l-x_j<c(l-j)$, gives for $\widehat x$ the witness pair $\widehat l=N-l<\widehat j=N-j$ with $\widehat x_{\widehat j}-\widehat x_{\widehat l}=x_l-x_j<c(\widehat j-\widehat l)$, and $j\mapsto N-j$ is injective. So the indices of the second kind for $(x_j)$ number at most the backward bad indices of $(\widehat x_i)$, which is at most $\Gamma/(1-c)$ by the first part.
\end{proof}

\begin{lemma}
 For every integer $p\ge2$ and every $x\in\mathbb{R}$,
\begin{equation}\label{eq:phi_derivatives}
\lvert\phi^{(p)}(x)\rvert\le(2\pi)^{p-1}K<(2\pi)^{p-1}M.
\end{equation}
\label{lem:der_bou_the}
\end{lemma}

\begin{proof}
For $p\ge2$ the linear part of $\phi$ does not contribute, so $\phi^{(p)}(x)=k(2\pi)^p\sin^{(p)}(2\pi x)$, where $\sin^{(p)}$ is one of $\pm\sin$, $\pm\cos$. Hence $\lvert\phi^{(p)}(x)\rvert\le k(2\pi)^p=(2\pi)^{p-1}K$, and $K<K+4=M$.
\end{proof}

\begin{lemma}
Let $a<b$ be integers with $b-a\ge2$, and let $J$ be a real matrix with rows and columns indexed by $a<l,p<b$, with $J_{p,p+1}=J_{p+1,p}=-1$ for $a<p<b-1$ and $J_{lp}=0$ for $\lvert l-p\rvert\ge2$. Let $(t_p)_{a\le p\le b}$ be the sequence with $t_a=0$, $t_{a+1}=1$ and $t_{p+1}=J_{pp}t_p-t_{p-1}$ for $a<p<b$. Then $\det J=t_b$. If $t_b\ne0$, there is a unique sequence $(U_p)_{a\le p\le b}$ with $U_a=1$, $U_b=0$ and $U_{p+1}=J_{pp}U_p-U_{p-1}$ for $a<p<b$, and
\[
(J^{-1})_{lp}=t_{\min(l,p)}\,U_{\max(l,p)}\qquad(a<l,p<b).
\]
\label{lem:gre_mat_fin}
\end{lemma}

\begin{proof}
This is the Green function of the Dirichlet problem for a Jacobi difference equation; it follows from the Wronskian identity and the determinant formula of \cite[(1.20)--(1.21), (1.65)]{locbib:Teschl2000} (see \cite[(1.99), (1.103)--(1.106)]{locbib:Teschl2000} for the Green function of a Jacobi operator and of its Dirichlet restrictions). We include the short proof.

We first compute the determinant. For $a+1\le p\le b$ let $\Delta_p$ be the determinant of the submatrix of $J$ on the indices $a<l,p'<p$, with $\Delta_{a+1}=1$ (empty matrix), and put $\Delta_a=0$. Then $\Delta_{p+1}=J_{pp}\Delta_p-\Delta_{p-1}$ for $a<p<b$: for $p=a+1$ this reads $\Delta_{a+2}=J_{a+1,a+1}$, the $1\times1$ minor, and for $a+2\le p<b$ it follows by expanding along the last row, using $J_{p,p-1}J_{p-1,p}=1$. So $(\Delta_p)$ and $(t_p)$ satisfy the same recurrence with the same values at $a$ and $a+1$, hence $\Delta_p=t_p$ for $a\le p\le b$, and $\det J=\Delta_b=t_b$.

We next construct $U$. The sequences $(w_p)_{a\le p\le b}$ with $w_{p+1}=J_{pp}w_p-w_{p-1}$ for $a<p<b$ are determined by $(w_a,w_{a+1})$, so they form a two-dimensional space, spanned by $t$ and by the solution $\sigma$ with $\sigma_a=1$, $\sigma_{a+1}=0$. Since $t_a=0$, a solution $U$ with $U_a=1$ is $\sigma$ plus a multiple of $t$; when $t_b\ne0$ the condition $U_b=0$ fixes the multiple, so $U=\sigma-(\sigma_b/t_b)\,t$ exists and is unique.

We then show that the Wronskian $W_p=t_{p+1}U_p-t_pU_{p+1}$ equals $1$ for $a\le p<b$. Indeed $W_a=t_{a+1}U_a-t_aU_{a+1}=1$, and for $a<p<b$
\[
W_p-W_{p-1}=U_p(t_{p+1}+t_{p-1})-t_p(U_{p+1}+U_{p-1})=U_pJ_{pp}t_p-t_pJ_{pp}U_p=0.
\]

Finally we verify the formula column by column. Fix $p$ with $a<p<b$ and put $G_l=t_{\min(l,p)}U_{\max(l,p)}$ for $a\le l\le b$; then $G_a=t_aU_p=0$ and $G_b=t_pU_b=0$, so the $l$-th entry of $J(G_{a+1},\dots,G_{b-1})^\top$ is $J_{ll}G_l-G_{l-1}-G_{l+1}$ for every $a<l<b$. For $l<p$ this is $U_p(J_{ll}t_l-t_{l-1}-t_{l+1})=0$, for $l>p$ it is $t_p(J_{ll}U_l-U_{l-1}-U_{l+1})=0$, and for $l=p$ it is $U_p(J_{pp}t_p-t_{p-1})-t_pU_{p+1}=t_{p+1}U_p-t_pU_{p+1}=W_p=1$. Hence $(G_l)_{a<l<b}$ is the $p$-th column of $J^{-1}$, which exists since $\det J=t_b\ne0$.
\end{proof}

\section{Short times and the doubling inequality} \label{sec:hyp_and_the}

This section proves the first two bounds of the bootstrap: the shortfall starts small, \eqref{eq:arg_start}, and one doubling increases it by at most a constant factor, \eqref{eq:arg_doub}.

\subsection{The mean shortfall at short times} \label{subsec:sho_tim}

At times $n\le\sqrt K$ the mean shortfall is below $3/\log M$. The bound is used only in the proof of Proposition~\ref{prop:uni_bou_the:2}, for $n=2^j$ up to the initial scale.

\begin{lemma}
 If $K\ge2^{16}$ and $n$ is an integer with $1\le n\le\sqrt K$, then $e_n<3/\log M$.
\label{lem:mea_sho_sho}
\end{lemma}

\begin{proof}
Put $w_i=K\lvert\cos(2\pi q_i)\rvert$ for $0\le i<n$ and $\mathrm{Exc}=\{z\in\mathbb{T}^2:\ w_i<6\text{ for some }0\le i<n\}$, a measurable set.

On the complement of $\mathrm{Exc}$ we have $\log\lVert T_{0,n}\rVert\ge\sum_{i<n}\log(w_i/2)$. Indeed, $v_i=v(q_i)=2+K\cos(2\pi q_i)$ satisfies $\lvert v_i\rvert-1\ge w_i-3\ge w_i/2\ge3$. Let $(u_i)$ be the linearized solution with $u_{-1}=0$, $u_0=1$ (Definition~\ref{def:orb_sta_and}), so $u_{i+1}=v_iu_i-u_{i-1}$. We show by induction on $i$ that $\lvert u_{i+1}\rvert\ge(\lvert v_i\rvert-1)\lvert u_i\rvert$ and $\lvert u_{i+1}\rvert\ge\lvert u_i\rvert$ for $0\le i<n$. Since $\lvert u_0\rvert=1>\lvert u_{-1}\rvert=0$, it suffices to note that if $\lvert u_i\rvert\ge\lvert u_{i-1}\rvert$, then $\lvert u_{i+1}\rvert\ge\lvert v_i\rvert\lvert u_i\rvert-\lvert u_{i-1}\rvert\ge(\lvert v_i\rvert-1)\lvert u_i\rvert\ge3\lvert u_i\rvert$. Hence $\lvert u_n\rvert\ge\prod_{i<n}(\lvert v_i\rvert-1)\ge\prod_{i<n}w_i/2$. Since $(u_n,u_{n-1})^{\top}=T_{0,n}(1,0)^{\top}$, $\lVert T_{0,n}\rVert\ge\lvert u_n\rvert$. On $\mathrm{Exc}$ we use $\log\lVert T_{0,n}\rVert\ge0$ (Lemma~\ref{lem:nor_and_gro}(i)).

By Lemma~\ref{lem:inv_are_und}(ii), $q_i$ has image measure Lebesgue measure on $\mathbb{T}$ under $m$, for every $i$. Euler's integral $\int_0^1\log\lvert\cos2\pi q\rvert\,dq=-\log2$ gives $\mathbb{E}\log(w_i/2)=\log(K/4)$. Also $\log(w_i/2)\le\log K$. Next, for $0<\varepsilon\le1$ the set $\{q\in\mathbb{T}:\lvert\cos 2\pi q\rvert<\varepsilon\}$ consists of two arcs of total length $\tfrac2\pi\arcsin\varepsilon\le\varepsilon$ ($\arcsin$ is convex on $[0,1]$ and lies below its chord), so $m\{w_i<6\}\le6/K$ and, by finite subadditivity, $m(\mathrm{Exc})\le6n/K$. Therefore
\[
\mathbb{E}\log\lVert T_{0,n}\rVert\ \ge\ \mathbb{E}\Bigl[\mathbf 1_{\mathrm{Exc}^c}\sum_{i<n}\log\tfrac{w_i}2\Bigr]
=\ n\log\tfrac K4-\mathbb{E}\Bigl[\mathbf 1_{\mathrm{Exc}}\sum_{i<n}\log\tfrac{w_i}2\Bigr]\ \ge\ n\log\tfrac K4-n\log K\cdot\tfrac{6n}K .
\]
For $n\le\sqrt K$ and $K\ge2^{16}$, $6n\log K/K\le6\log K/\sqrt K\le0.26$ ($\log K/\sqrt K$ decreases for $K>\exp(2)$), and $2\log2+0.26<2$; so $\frac1n\mathbb{E}\log\lVert T_{0,n}\rVert\ge\log K-2$. Finally $g(0,n)=\log\lVert T_{0,n}\rVert/\log M$ and $M=K+4$, so
\[
e_n=1-\frac{\mathbb{E}\log\lVert T_{0,n}\rVert}{n\log M}\le\frac{\log M-\log K+2}{\log M}=\frac{2+\log(1+4/K)}{\log M}<\frac3{\log M}.
\]

Explicit values: Euler's integral and the chord bound for $\arcsin$ are checked numerically; $6\log K/\sqrt K=0.2599$ at $K=2^{16}$; $2\log2+0.26=1.6463<2$; and the final inequality already holds for $K\ge7642$.
\end{proof}

\subsection{Hyperbolicity and the doubling inequality} \label{subsec:hyp_dou}

This subsection proves the doubling inequality Corollary~\ref{cor:dou_ine_for}. The increments of the mean shortfall along $n=2^j$ are mean Gromov products at the midpoint, by an exact identity that uses only the invariance of area (Lemma~\ref{lem:exa_dou_ide}). To bound these Gromov products we use the hyperbolicity of the growth semimetric (Lemma~\ref{lem:hyp_the_gro}): a large Gromov product survives truncation of the wings, and so produces either a wing deficit or a near-total cancellation. The last step uses the scarcity of cancellations Theorem~\ref{thm:exp_sca_can}, restated below and proved in Appendix~\ref{sec:scar_canc} without the doubling inequality.

We begin with the identity.

\begin{lemma}
 For every integer $n\ge1$, the Gromov product $(0|2n)_n$ of $g$, a function of the initial point of the orbit, is nonnegative, and
\[
e_{2n}=e_n+\frac1n\mathbb{E}(0|2n)_n;
\]
in particular $e_n\le e_{2^jn}$ for every integer $j\ge0$.
\label{lem:exa_dou_ide}
\end{lemma}

\begin{proof}
By Lemma~\ref{lem:nor_and_gro}(ii), $g$ is a semimetric on $\mathbb{Z}$ with $0\le g(a,b)\le\lvert a-b\rvert$. Hence $2(0|2n)_n=g(0,n)+g(n,2n)-g(0,2n)\ge0$ by the triangle inequality, and all functions integrated below are bounded. By invariance of $m$ under time shifts of the orbit (Lemma~\ref{lem:inv_are_und}(i),(iii)), $\mathbb{E}g(n,2n)=\mathbb{E}g(0,n)$. Integrating $g(0,2n)=g(0,n)+g(n,2n)-2(0|2n)_n$ over $\mathbb{T}^2$ therefore gives
\[
e_{2n}=1-\frac1{2n}\mathbb{E}g(0,2n)=1-\frac1n\mathbb{E}g(0,n)+\frac1n\mathbb{E}(0|2n)_n=e_n+\frac1n\mathbb{E}(0|2n)_n.
\]
Since the last integrand is nonnegative, $e_n\le e_{2n}$ for every $n\ge1$, and induction on $j$ gives $e_n\le e_{2^jn}$.
\end{proof}

To bound these Gromov products we use that the growth semimetric is the hyperbolic distance along a path in $\mathbb H^2$, divided by $2\log M$:

\begin{lemma}
 For every orbit, $(\mathbb{Z},g)$ is $\tfrac12\log_M2$-hyperbolic: for all integers $i,a,b,l$,
\begin{equation}\label{eq:growth_hyperbolic}
(a|l)_i\ge\min\bigl\{(a|b)_i,(b|l)_i\bigr\}-\tfrac12\log_M2.
\end{equation}
\label{lem:hyp_the_gro}
\end{lemma}

\begin{proof}
Along the orbit, $g$ is the hyperbolic distance between the points of a path in $\mathbb{H}^2$, divided by $2\log M$; the inequality is therefore the $\log2$-hyperbolicity of $\mathbb{H}^2$, rescaled.

Let $\mathbb{H}^2=\{w\in\mathbb{C}:\operatorname{Im}w>0\}$ carry the metric $d$ of curvature $-1$ given by
\[
\cosh d(w_1,w_2)=1+\frac{\lvert w_1-w_2\rvert^2}{2\operatorname{Im}w_1\operatorname{Im}w_2},
\]
on which $\mathrm{SL}(2,\mathbb{R})$ acts by isometries through the Möbius transformations $T\cdot w=(\alpha_1w+\alpha_2)/(\alpha_3w+\alpha_4)$, $T=\begin{pmatrix}\alpha_1&\alpha_2\\\alpha_3&\alpha_4\end{pmatrix}$ (a classical fact), and put $o=\mathrm{i}$.

We first show that $d(o,T\cdot o)=2\log\lVert T\rVert$ for every $T\in\mathrm{SL}(2,\mathbb{R})$. Since $\alpha_1\alpha_4-\alpha_2\alpha_3=1$, one has $\operatorname{Im}(T\cdot o)=1/(\alpha_3^2+\alpha_4^2)$ and $\lvert T\cdot o-o\rvert^2=\bigl((\alpha_2+\alpha_3)^2+(\alpha_1-\alpha_4)^2\bigr)/(\alpha_3^2+\alpha_4^2)$, so
\[
\begin{aligned}
\cosh d(o,T\cdot o)&=1+\frac{(\alpha_2+\alpha_3)^2+(\alpha_1-\alpha_4)^2}{2}=\frac{\alpha_1^2+\alpha_2^2+\alpha_3^2+\alpha_4^2}{2}\\
&=\frac{\lVert T\rVert^2+\lVert T\rVert^{-2}}{2}=\cosh\bigl(2\log\lVert T\rVert\bigr);
\end{aligned}
\]
the third equality holds because the singular values of $T$ are $\lVert T\rVert$ and $\lVert T\rVert^{-1}$. Since $\lVert T\rVert\ge1$ and $\cosh$ is injective on $[0,\infty)$, the claim follows.

Next, $g$ is the rescaled distance along the orbit path $p_a=T_{0,a}^{-1}\cdot o$ $(a\in\mathbb{Z})$. The matrices $T_{a,b}$ lie in $\mathrm{SL}(2,\mathbb{R})$ (Lemma~\ref{lem:nor_and_gro}(i)), and the isometry $T_{0,a}$ maps $p_a$ to $o$ and $p_b$ to $T_{0,a}T_{0,b}^{-1}\cdot o=T_{b,a}\cdot o$ (Lemma~\ref{lem:nor_and_gro}(i)). So the displacement formula and $\lVert T_{b,a}\rVert=\lVert T_{a,b}\rVert$ (Lemma~\ref{lem:nor_and_gro}(i)) give
\[
d(p_a,p_b)=2\log\lVert T_{a,b}\rVert=(2\log M)\,g(a,b)\qquad(a,b\in\mathbb{Z}).
\]
Consequently $(a|l)_i=(p_a|p_l)_{p_i}/(2\log M)$ for all integers $i,a,l$, the Gromov product on the right being taken in $(\mathbb{H}^2,d)$.

Finally, $(\mathbb{H}^2,d)$ is $\log2$-hyperbolic in the sense of Definition~\ref{def:gro_hyp_sem}, and $\log2$ is optimal: this is \cite[Corollary~5.4]{bib:nica2014a}, whose definition of $\delta$-hyperbolicity (\cite[Definition~2.1]{bib:nica2014a}: Gromov's, required at every basepoint) and whose metric on the upper half-plane (proof of \cite[Proposition~5.2]{bib:nica2014a}) are the ones used here. Applying it to $p_a,p_b,p_l$ with basepoint $p_i$ and dividing by $2\log M$ gives \eqref{eq:growth_hyperbolic}, since $\log2/(2\log M)=\tfrac12\log_M2$.

\end{proof}

We then cut the wings: hyperbolicity transports a large Gromov product of the full wings to the cut wings, at the price of the wing deficit.

\begin{lemma}
 Let $i\in\mathbb{Z}$, and let $n\ge1$ and $1\le s\le n$ be integers, and put $D=n-\min\{g(i,i-n),g(i,i+n)\}$, the \emph{wing deficit}. Then, for every orbit,
\begin{align}
g(i,i-s),\ g(i,i+s)&\ge s-D,\label{eq:trunc_wings}\\
2(i-s|i+s)_i&\ge\min\bigl\{2(i-n|i+n)_i,\ 2s-2D\bigr\}-1,\label{eq:trunc_gromov}\\
g(i-s,i+s)&\le2s-\min\bigl\{2(i-n|i+n)_i,\ 2s-2D\bigr\}+1.\label{eq:trunc_span}
\end{align}
\label{lem:gro_pro_und}
\end{lemma}

\begin{proof}
By the triangle inequality and $0\le g(a,b)\le\lvert a-b\rvert$ (Lemma~\ref{lem:nor_and_gro}(ii)),
\[
g(i,i-s)\ge g(i,i-n)-g(i-n,i-s)\ge(n-D)-(n-s)=s-D,
\]
and likewise $g(i,i+s)\ge s-D$; this is \eqref{eq:trunc_wings}. Together with $g(i-n,i-s)\le n-s$ it gives
\[
2(i-n|i-s)_i=g(i,i-n)+g(i,i-s)-g(i-n,i-s)\ge(n-D)+(s-D)-(n-s)=2s-2D,
\]
and likewise $2(i+n|i+s)_i\ge2s-2D$.

We next transport the Gromov product of the full wings to the cut wings. By \eqref{eq:growth_hyperbolic} of Lemma~\ref{lem:hyp_the_gro}, applied with $(a,b,l)=(i-s,i-n,i+s)$ and then with $(a,b,l)=(i-n,i+n,i+s)$,
\[
(i-s|i+s)_i\ge\min\bigl\{(i-s|i-n)_i,(i-n|i+s)_i\bigr\}-\tfrac12\log_M2,
\]
\[
(i-n|i+s)_i\ge\min\bigl\{(i-n|i+n)_i,(i+n|i+s)_i\bigr\}-\tfrac12\log_M2.
\]
Substituting the second inequality into the first and using the two lower bounds above gives $2(i-s|i+s)_i\ge\min\{2(i-n|i+n)_i,2s-2D\}-2\log_M2$.  Since $M=K+4>4$ (Convention~\ref{conv:sta_map_par}), $2\log_M2<1$, and \eqref{eq:trunc_gromov} follows.

Finally, $g(i,i\pm s)\le s$, so $g(i-s,i+s)=g(i,i-s)+g(i,i+s)-2(i-s|i+s)_i\le2s-2(i-s|i+s)_i$, and \eqref{eq:trunc_gromov} gives \eqref{eq:trunc_span}.
\end{proof}

When the deficit is small against the cut length, the cut wings form a near-total cancellation in the following sense.

\begin{definition}
Let $Y\subset\mathbb{R}$, let $d\colon Y\times Y\to\mathbb{R}$ (for instance a semimetric on $Y$ in the sense of Definition~\ref{def:gro_hyp_sem}), let $\delta>0$ and $t>0$, and let $r\in\mathbb{R}$ with $r-t,r,r+t\in Y$. Then $d$ has a \textbf{$\delta$-cancellation of length $t$ at $r$} if on both wings $[r-t,r]$ and $[r,r+t]$ it is at least $(1-\delta)t$, while on the whole interval $[r-t,r+t]$ it is at most $\delta t$:
\begin{equation}\label{eq:cancellation}
\min\bigl\{d(r,r-t),\,d(r,r+t)\bigr\}\ge(1-\delta)t\qquad\text{and}\qquad d(r-t,r+t)\le\delta t.
\end{equation}
For the growth semimetric $g$ of an orbit (Definition~\ref{def:gro_sem_orb}) we take $Y=\mathbb{Z}$, and a $\delta$-cancellation of $g$ of integer length $n$ at $r=i\in\mathbb{Z}$ is called a $\delta$-cancellation of length $n$ at time $i$. For $i\in\mathbb{Z}$ and an integer $n\ge1$, the \textbf{cancellation set} $\mathrm{Can}_\delta(i,n)\subset\mathbb{T}^2$ is the set of initial points $z$ of orbits (Definition~\ref{def:orb_sta_and}) whose growth semimetric has a $\delta$-cancellation of length $n$ at time $i$.
\label{def:can_the_gro}
\end{definition}

The doubling inequality follows from any exponential bound on the measure of cancellation sets: the wing deficit is charged to the mean shortfall, and the large Gromov products to the cancellation sets.

\begin{lemma}
 Let $\delta>0$, $\gamma>0$ and $C\ge1$, and suppose that
\[
\mathbb{P}\bigl(\mathrm{Can}_\delta(i,s)\bigr)\le C\exp(-\gamma s)\qquad(i\in\mathbb{Z},\ s=1,2,\ldots).
\]
Then for every integer $n\ge1$,
\[
e_{2n}\le\Bigl(1+\frac4\delta\Bigr)e_n+\frac1n\Bigl(3+\frac1\delta+\frac{1+\log(C/\gamma)}{\gamma}\Bigr).
\]
\label{lem:dou_ine_fro}
\end{lemma}

\begin{proof}
We bound the Gromov product $(0|2n)_n$ in the exact doubling identity of Lemma~\ref{lem:exa_dou_ide}: off a set covered by cancellation sets it is controlled by the wing deficit, and on that set by the hypothesis. Fix $n\ge1$ and put $D=n-\min\{g(n,0),g(n,2n)\}$, the wing deficit of Lemma~\ref{lem:gro_pro_und} at $i=n$; it is a function of the initial point $z$ of the orbit. By Lemma~\ref{lem:nor_and_gro}(ii) and the triangle inequality, $D\ge0$ and $0\le(0|2n)_n\le g(n,0)\le n$, where $(0|2n)_n$ is the Gromov product of $g$ based at $n$.

We first show that the set
\[
\mathrm{Exc}=\Bigl\{z:\ (0|2n)_n>1+\frac{2D+1}{\delta}\Bigr\}
\]
is covered by cancellation sets. On $\mathrm{Exc}$ the integer $s=\lfloor(0|2n)_n\rfloor$ lies in $\{1,\ldots,n\}$ and satisfies
\[
\delta s>\delta\bigl((0|2n)_n-1\bigr)>2D+1.
\]
Apply Lemma~\ref{lem:gro_pro_und} at $i=n$ with this $s$. By \eqref{eq:trunc_wings}, $g(n,n\pm s)\ge s-D>(1-\delta)s$. Since $2s-2D\le2s\le2(0|2n)_n$, the minimum in \eqref{eq:trunc_span} is $2s-2D$, so $g(n-s,n+s)\le2D+1<\delta s$. Hence every $z\in\mathrm{Exc}$ lies in $\mathrm{Can}_\delta(n,s)$ for its integer $s$. (All sets here are Borel, since each $g(a,b)$ is a continuous function of $z$.)

We next bound the integral of $(0|2n)_n$ over $\mathrm{Exc}$. Let $t\ge0$. If $(0|2n)_n\mathbf{1}_{\mathrm{Exc}}>t$ at $z$, then $z\in\mathrm{Exc}$ and $s=\lfloor(0|2n)_n\rfloor\ge s_t:=\max\{1,\lfloor t\rfloor\}$, so $z\in\mathrm{Can}_\delta(n,s)$ for some $s\ge s_t$. Countable subadditivity of $m$, the hypothesis, $s_t>t-1$ and $1-\exp(-\gamma)\ge\gamma\exp(-\gamma)$ (that is, $\exp(\gamma)-1\ge\gamma$) give
\[
\begin{aligned}
\mathbb{P}\bigl\{(0|2n)_n\mathbf{1}_{\mathrm{Exc}}>t\bigr\}&\le\min\Bigl\{1,\sum_{s\ge s_t}C\exp(-\gamma s)\Bigr\}=\min\Bigl\{1,\frac{C\exp(-\gamma s_t)}{1-\exp(-\gamma)}\Bigr\}\\
&\le\min\Bigl\{1,\frac{C\exp(2\gamma)}{\gamma}\exp(-\gamma t)\Bigr\}.
\end{aligned}
\]
Since $C\ge1$ and $\exp(2\gamma)>\gamma$, the right side equals $1$ for $t\le t_*=2+\gamma^{-1}\log(C/\gamma)$, where $t_*>0$, and equals $\exp(-\gamma(t-t_*))$ for $t\ge t_*$. Integrating over $t\ge0$ therefore gives
\[
\mathbb{E}\bigl[(0|2n)_n\mathbf{1}_{\mathrm{Exc}}\bigr]=\int_0^\infty \mathbb{P}\bigl\{(0|2n)_n\mathbf{1}_{\mathrm{Exc}}>t\bigr\}\,dt\le t_*+\frac1\gamma=2+\frac{1+\log(C/\gamma)}{\gamma}.
\]

Finally we compare with the wing deficit. Off $\mathrm{Exc}$, $(0|2n)_n\le1+(2D+1)/\delta$, so pointwise
\[
(0|2n)_n\le1+\frac{2D+1}{\delta}+(0|2n)_n\mathbf{1}_{\mathrm{Exc}}.
\]
Moreover $D\le\bigl(n-g(0,n)\bigr)+\bigl(n-g(n,2n)\bigr)$, both terms being nonnegative, and $\mathbb{E}(n-g(n,2n))=\mathbb{E}(n-g(0,n))=ne_n$ by invariance of $m$ under time shifts of the orbit (Lemma~\ref{lem:inv_are_und}(i),(iii)); hence $\mathbb{E}D\le2ne_n$. Integrating the pointwise bound,
\[
\mathbb{E}(0|2n)_n\le\frac{4n}{\delta}e_n+3+\frac1\delta+\frac{1+\log(C/\gamma)}{\gamma},
\]
and the exact doubling identity $e_{2n}=e_n+\frac1n\mathbb{E}(0|2n)_n$ of Lemma~\ref{lem:exa_dou_ide} gives the assertion.

\end{proof}

The scarcity of cancellations Theorem~\ref{thm:exp_sca_can}, proved in Appendix~\ref{sec:scar_canc} independently of the results of this section, supplies this bound with explicit constants: for $M\ge M_{\rm an}$, every $i\in\mathbb Z$ and every integer $n\ge1$,
\[
\mathbb P\bigl(\mathrm{Can}_{\delta_{\rm can}}(i,n)\bigr)\le C_{\rm can}\exp(-\gamma_{\rm can}n).
\]

\begin{corollary}
 Let $M\ge M_{\rm an}$. Then for every integer $n\ge1$,
\[
e_{2n}+\frac1{2n}\le D_{\rm dbl}\Bigl(e_n+\frac1n\Bigr).
\]
\label{cor:dou_ine_for}
\end{corollary}

\begin{proof}
Since $M\ge M_{\rm an}$, Theorem~\ref{thm:exp_sca_can} gives, for the cancellation sets of Definition~\ref{def:can_the_gro},
\[
\mathbb{P}\bigl(\mathrm{Can}_{\delta_{\rm can}}(i,s)\bigr)\le C_{\rm can}\exp(-\gamma_{\rm can}s)\qquad(i\in\mathbb{Z},\ s=1,2,\ldots).
\]
This is the hypothesis of Lemma~\ref{lem:dou_ine_fro} with $\delta=\delta_{\rm can}$, $\gamma=\gamma_{\rm can}$ and $C=C_{\rm can}$; as $\gamma_{\rm can}=\delta_{\rm can}$ by Convention~\ref{conv:not_and_con}, it gives
\[
e_{2n}\le\bigl(1+4\delta_{\rm can}^{-1}\bigr)\,e_n+\frac1n\Bigl(3+\delta_{\rm can}^{-1}+\delta_{\rm can}^{-1}\bigl(1+\log(C_{\rm can}\delta_{\rm can}^{-1})\bigr)\Bigr).
\]
With the values of Convention~\ref{conv:not_and_con}, $1+4\delta_{\rm can}^{-1}<10^7$ and $C_{\rm can}\delta_{\rm can}^{-1}<10^{11}<\exp(26)$, so the additive constant is less than $3+28\delta_{\rm can}^{-1}<6.8\cdot10^7$. Since $e_n\ge0$ by $g(0,n)\le n$ (Lemma~\ref{lem:nor_and_gro}(ii)), and $10^7$ and $6.8\cdot10^7+\tfrac12$ are less than $D_{\rm dbl}$,
\[
e_{2n}+\frac1{2n}\le\bigl(1+4\delta_{\rm can}^{-1}\bigr)\,e_n+\frac{6.8\cdot10^7+\frac12}{n}\le D_{\rm dbl}\Bigl(e_n+\frac1n\Bigr).
\]

\end{proof}

The factor $D_{\rm dbl}$ is far from one, so the doubling inequality alone gives no bound uniform in $n$. It is used for a single doubling at the crossing scale in the proof of Proposition~\ref{prop:uni_bou_the:2}, and for eleven consecutive doublings in the proof of Proposition~\ref{prop:def_the_ext:1}.

\section{The product estimate across a bridge} \label{sec:the_pro_est}

This section proves the bridge product estimate Theorem~\ref{thm:bri_pro_est:2}, the second quantitative input of the paper: on an interval of near-full growth, two shortfall observables localized near its two ends satisfy a product upper bound. It is used twice in the section on balance: to bound the measure of the set where a stretch of near-full growth joins a nonaffine marker of a path semimetric to a nonunit one (Lemma~\ref{lem:sum_exc_mas:1}), and to bound the balance defect at middle speeds off the body (Lemma~\ref{lem:bal_def_off}). The proof constructs, at a middle time, a family of finite-time stable graphs and a transverse family of unstable graphs (compare the finite-time local stable manifolds of \cite{bib:holland2005}); the comparison of the two ends is then a change of variables with bounded Jacobian, in the spirit of absolute continuity in Pesin theory \cite{locbib:BarreiraPesin2007}. We first construct the graphs.

\begin{definition}

Let $i\in\mathbb{Z}$ and let $B\ge1$ be an integer. Time $i$ of an orbit has \textbf{forward prefix growth of length $B$} if $g(i,i+l)\ge c_{\rm pre}l$ for $1\le l\le B$, and \textbf{backward prefix growth of length $B$} if $g(i-l,i)\ge c_{\rm pre}l$ for $1\le l\le B$.

Let $(t^+_p)$ be the linearized solution with $(t^+_0,t^+_{-1})=(1,0)$, and let $(t^-_p)$ be the linearized solution with $(t^-_0,t^-_{-1})=(0,1)$. If time $0$ has forward prefix growth of length $B_+\ge2$, it is \textbf{forward $j$-eligible}, for $j\in\{B_+-1,B_+\}$, if $\lvert t^+_j\rvert\ge\frac12\lVert T_{0,B_+}\rVert$. If time $0$ has backward prefix growth of length $B_-\ge2$, it is \textbf{backward $j$-eligible}, for $j\in\{B_--1,B_-\}$, if $\lvert t^-_{-1-j}\rvert\ge\frac12\lVert T_{-B_-,0}\rVert$.

Let $(q_p)$ be a real lift of the orbit and let $j\ge1$ be an integer. A \textbf{stable family of length $j$} through it is a family $(q_p(r))_{p\in\mathbb{Z}}$, $\lvert r\rvert\le1/4$, of real lifts of orbits, smooth in $r$, with $q_p(0)=q_p$, $q_{-1}(r)=q_{-1}+r$ and $q_j(r)=q_j$. An \textbf{unstable family of length $j$} is defined in the same way, with $q_0(r)=q_0+r$ and $q_{-1-j}(r)=q_{-1-j}$ in place of the last two conditions. The curve $r\mapsto(q_0(r),q_{-1}(r))$ is the \textbf{stable graph} of a stable family and the \textbf{unstable graph} of an unstable family.
\label{def:pre_gro_eli}
\end{definition}

Under forward prefix growth, the linearized Dirichlet problem from which the graphs are built has a unique solution, which decays at nearly full rate (Lemma~\ref{lem:mos_con_dir:1}). The stable graphs are level curves of a distant coordinate of the orbit; they contract, and the derivative of that coordinate varies along them by a factor at most two (Lemma~\ref{lem:sta_gra_und:1}). The unstable graphs are their time reversals (Lemma~\ref{lem:bac_sta_gra:1}).

\begin{lemma}
 Let $M\ge M_{\rm an}$ and let time $0$ of an orbit have forward prefix growth of length $B_+\ge2$.
\begin{enumerate}
\item[(i)] $\max\bigl(\lvert t^+_{B_+}\rvert,\lvert t^+_{B_+-1}\rvert\bigr)\ge\frac12\lVert T_{0,B_+}\rVert$; so time $0$ is forward $j$-eligible for some $j\in\{B_+-1,B_+\}$.
\item[(ii)] If time $0$ is forward $j$-eligible, there is a unique linearized solution $(U_p)$ with $U_{-1}=1$ and $U_j=0$, and
\[
\lvert U_p\rvert\le12M^{-\alpha_{\rm pre}(p+1)}\qquad(0\le p<j).
\]
\end{enumerate}
\label{lem:mos_con_dir:1}
\end{lemma}

\begin{proof}
Prefix growth makes consecutive most contracted directions converge geometrically, so the most contracted direction at the final time is contracted at every intermediate time; the kernel of a large row of $T_{0,B_+}$ is a small perturbation of it. In this proof write $c=c_{\rm pre}$, $\alpha=\alpha_{\rm pre}=3c-2$ and $D_l=\lVert T_{0,l}\rVert$, and write input vectors of transfer matrices as $(u_0,u_{-1})$. Forward prefix growth of length $B_+$ says $D_l\ge M^{cl}$ for $1\le l\le B_+$, since $g(0,l)=\log_MD_l$ (Definition~\ref{def:gro_sem_orb}).

We first show that consecutive most contracted directions are close. For $1\le l\le B_+$ let $s(l)$ be a most contracted unit input direction of $T_{0,l}$ (Convention~\ref{conv:wed_pro_and}); it is well defined up to sign since $D_l\ge M^{cl}>1$, and we put $s=s(B_+)$. Since $\det T_{0,l}=1$ (Lemma~\ref{lem:nor_and_gro}(i)), $\lVert T_{0,l}s(l)\rVert=D_l^{-1}$, and the images of $s(l)$ and of the orthogonal input direction are orthogonal, the latter having norm $D_l$. A single inverse step has norm at most $M$ (Lemma~\ref{lem:nor_and_gro}(i)), so $\lVert T_{0,l}s(l+1)\rVert\le M\lVert T_{0,l+1}s(l+1)\rVert=M/D_{l+1}$. Resolving $s(l+1)$ along $s(l)$ and the orthogonal direction gives $D_l\lvert s(l+1)\wedge s(l)\rvert\le\lVert T_{0,l}s(l+1)\rVert$, that is,
\[
\lvert s(l+1)\wedge s(l)\rvert\le\frac{M}{D_lD_{l+1}}\qquad(1\le l<B_+).
\]
This is the angle estimate for consecutive most contracted directions of \cite[Section 3.3]{bib:holland2005}. The absolute sine between two lines satisfies the triangle inequality, $\lvert \nu_1\wedge \nu_3\rvert\le\lvert \nu_1\wedge \nu_2\rvert+\lvert \nu_2\wedge \nu_3\rvert$ for unit vectors $\nu_1,\nu_2,\nu_3$, because the angle between lines is a metric with values in $[0,\pi/2]$, $\sin$ is increasing on $[0,\pi/2]$ and $\sin(x+y)\le\sin x+\sin y$ there, while $\sin x+\sin y\ge1$ when $x+y\ge\pi/2$. Hence, by prefix growth,
\[
\lvert s\wedge s(l)\rvert\le\sum_{i=l}^{B_+-1}M^{1-ci-c(i+1)}\le\frac{M^{1-c-2cl}}{1-M^{-2c}}\le2M^{1-c-2cl}\qquad(1\le l\le B_+).
\]

We next show that $s$ is contracted at every intermediate time, by the \emph{stable-prefix bound}
\[
\lVert T_{0,l}s\rVert\le3M^{-\alpha l}\qquad(1\le l\le B_+).
\]
Resolve $s$ along $s(l)$ and the orthogonal input direction of $T_{0,l}$: the first component is contracted by $D_l^{-1}\le M^{-cl}$, and the second, of absolute value $\lvert s\wedge s(l)\rvert$, is expanded by $D_l\le M^l$ (Lemma~\ref{lem:nor_and_gro}(i)). Hence $\lVert T_{0,l}s\rVert\le M^{-cl}+2M^{l+1-c-2cl}\le3M^{-\alpha l}$, because $-cl\le-\alpha l$ and $l+1-c-2cl+\alpha l=(1-c)(1-l)\le0$ for $l\ge1$. Write $s=(s_1,s_2)$. The matrix $T_{0,1}$ carries the input $(u_0,u_{-1})=s$ to $(u_1,u_0)$, whose second component is $s_1$; so $\lvert s_1\rvert\le\lVert T_{0,1}s\rVert\le3M^{-\alpha}$, and $s_2^2=1-s_1^2\ge1-9M^{-2\alpha}\ge\frac12$.

We now prove (i). Let $s^\perp$ be a unit vector orthogonal to $s$. Then $(1,0)^\top=s_1s+s'_2s^\perp$ with $\lvert s'_2\rvert=\lvert(1,0)\wedge s\rvert=\lvert s_2\rvert\ge1/\sqrt2$, and $T_{0,B_+}s$, $T_{0,B_+}s^\perp$ are orthogonal with $\lVert T_{0,B_+}s^\perp\rVert=D_{B_+}$. Hence $\lVert T_{0,B_+}(1,0)^\top\rVert\ge\lvert s_2\rvert D_{B_+}\ge D_{B_+}/\sqrt2$. Since $T_{0,B_+}(1,0)^\top=(t^+_{B_+},t^+_{B_+-1})^\top$ and the larger coordinate of a vector of $\mathbb{R}^2$ is at least $1/\sqrt2$ times its norm, $\max(\lvert t^+_{B_+}\rvert,\lvert t^+_{B_+-1}\rvert)\ge D_{B_+}/2$. The index $j$ attaining the maximum makes time $0$ forward $j$-eligible (Definition~\ref{def:pre_gro_eli}).

We finally prove (ii). Let $\rho$ be the row of $T_{0,B_+}$ selecting $u_j$ (the first row if $j=B_+$, the second if $j=B_+-1$). Then $\lVert\rho\rVert\ge\lvert t^+_j\rvert\ge D_{B_+}/2$ by eligibility, and $\lvert\rho\cdot s\rvert\le\lVert T_{0,B_+}s\rVert=D_{B_+}^{-1}$. Let $\omega$ be a unit vector orthogonal to $\rho$ with $\langle\omega,s\rangle\ge0$. The angle $\psi\in[0,\pi/2]$ between $\omega$ and $s$ satisfies $\sin\psi=\lvert\rho\cdot s\rvert/\lVert\rho\rVert\le2D_{B_+}^{-2}$, hence
\[
\lvert\omega-s\rvert=2\sin(\psi/2)\le\psi\le\tfrac{\pi}{2}\sin\psi\le\pi D_{B_+}^{-2}\le\pi M^{-4c},
\]
using $D_{B_+}\ge M^{cB_+}\ge M^{2c}$. So $\lvert\omega_2\rvert\ge\lvert s_2\rvert-\lvert\omega-s\rvert\ge1/\sqrt2-\pi M^{-4c}>2/3$. Put $U^{\rm in}=\omega/\omega_2$, the input with second component $1$, and let $(U_p)$ be its linearized solution: $U_{-1}=1$ and $U_j=\rho\cdot U^{\rm in}=0$. It is unique: the value at index $j$ of the linearized solution with input $(U_0,1)$ is $U_0t^+_j+t^-_j$, affine in $U_0$ with slope $t^+_j\ne0$. For $0\le p<j$, $U_p$ is the second component of $T_{0,p+1}U^{\rm in}$, so by the stable-prefix bound at $l=p+1$ and Lemma~\ref{lem:nor_and_gro}(i),
\[
\lvert U_p\rvert\le\frac{\lVert T_{0,p+1}s\rVert+\lvert\omega-s\rvert\,\lVert T_{0,p+1}\rVert}{\lvert\omega_2\rvert}\le\frac32\Bigl(3M^{-\alpha(p+1)}+\pi D_{B_+}^{-2}M^{p+1}\Bigr).
\]
Since $p+1\le B_+$ and $D_{B_+}\ge M^{cB_+}$, $\pi D_{B_+}^{-2}M^{p+1}\le\pi M^{-(2c-1)(p+1)}\le\pi M^{-\alpha(p+1)}$, because $2c-1-\alpha=1-c>0$. Hence $\lvert U_p\rvert\le\frac32(3+\pi)M^{-\alpha(p+1)}\le12M^{-\alpha(p+1)}$.

Explicit values. With the values of Convention~\ref{conv:not_and_con}, the inequalities $M^{-2c}\le\frac12$, $M^{-2\alpha}\le\frac1{18}$ and $\pi M^{-4c}<\frac1{\sqrt2}-\frac23$ used above hold for $M\ge5$, hence for $M\ge M_{\rm an}$ (the stronger hypothesis $M\ge M_{\rm an}$ is kept for uniformity with the other estimates of the bridge argument); and $\frac32(3+\pi)<9.3$.
\end{proof}

\begin{lemma}
 Let $M\ge M_{\rm an}$, let time $0$ of an orbit be forward $j$-eligible with $j\in\{B_+-1,B_+\}$, $B_+\ge2$, and fix a real lift $(q_p)$ of the orbit. There is a stable family $(q_p(r))$ of length $j$ through this lift such that the following hold for $\lvert r\rvert,\lvert r'\rvert\le1/4$.
\begin{enumerate}
\item[(a)] $\lvert q_p(r)-q_p\rvert\le6M^{-w_{\rm pre}(p+1)}$ and $\lvert\mathrm{d}q_p(r)/\mathrm{d}r\rvert\le24M^{-w_{\rm pre}(p+1)}$ for $0\le p<j$.
\item[(b)] Each of the coordinates $q_{B_+}(r)$ and $q_{B_+-1}(r)$ differs from its value at $r=0$ by at most $6M^{-w_{\rm pre}(B_+-1)}$.
\item[(c)] Along the family, $\partial q_j/\partial q_0$ at fixed $q_{-1}$ equals $t^+_j(r)$, where $(t^+_p(r))$ is the linearized solution along $(q_p(r))$ with $(t^+_0(r),t^+_{-1}(r))=(1,0)$; and $t^+_j(r)\ne0$.
\item[(d)] For fixed $B_+$ and $j$, every $q_p(r)$ depends jointly continuously on $r$ and on the initial pair $(q_0,q_{-1})$ of the lift, as this pair ranges over the lifts of forward $j$-eligible states.
\item[(e)] $\frac12\le\lvert t^+_j(r)\rvert/\lvert t^+_j(r')\rvert\le2$.
\end{enumerate}
\label{lem:sta_gra_und:1}
\end{lemma}

\begin{proof}
The family solves a finite Dirichlet problem for the orbit recurrence with boundary data $q_{-1}(r)=q_{-1}+r$ and $q_j(r)=q_j$. The decaying solution of Lemma~\ref{lem:mos_con_dir:1}(ii) makes the Green matrix small in a weighted norm, so that the quadratic Taylor remainder is a contraction; the transverse derivative is the Dirichlet determinant, whose logarithmic derivative is a trace against the Green diagonal. These stable graphs are finite-time local stable manifolds in the sense of the Main Theorem of \cite[Section 1.5]{bib:holland2005}, obtained as level curves of $q_j$. In this proof write $w=w_{\rm pre}$ and $\alpha=\alpha_{\rm pre}$; then $\frac12<w<\alpha$ (Convention~\ref{conv:not_and_con}).

We first write the Green matrix of the Dirichlet problem. Let $J$ be the matrix with rows and columns indexed by $0\le l,p\le j-1$, with diagonal entries $J_{pp}=v(q_p)$, entries $-1$ next to the diagonal and $0$ elsewhere. The linearized solution $(t^+_p)$ of Definition~\ref{def:pre_gro_eli} satisfies $t^+_{-1}=0$, $t^+_0=1$ and $t^+_{p+1}=v(q_p)t^+_p-t^+_{p-1}$, and the solution $(U_p)$ of Lemma~\ref{lem:mos_con_dir:1}(ii) satisfies $U_{-1}=1$, $U_j=0$ and the same recurrence. By Lemma~\ref{lem:gre_mat_fin} with $a=-1$ and $b=j$ (note $j\ge1$), whose sequence $t$ is then $(t^+_p)$, $\det J=t^+_j$, which is nonzero since $\lvert t^+_j\rvert\ge\frac12\lVert T_{0,B_+}\rVert\ge\frac12$, and $G=J^{-1}$ has the entries $G_{lp}=t^+_{\min(l,p)}U_{\max(l,p)}$. Since $\lvert t^+_p\rvert\le\lVert T_{0,p}\rVert\le M^p$ for $p\ge0$ (Lemma~\ref{lem:nor_and_gro}(i)), Lemma~\ref{lem:mos_con_dir:1}(ii) gives
\[
\lvert G_{lp}\rvert\le12M^{\min(l,p)-\alpha(\max(l,p)+1)}\qquad(0\le l,p<j).
\]

We next write the family as a fixed point. Put $F_p(u)=-\bigl(\phi(q_p+u)-\phi(q_p)-v(q_p)u\bigr)$. Since $q_{p+1}+q_{p-1}=\phi(q_p)$, a real sequence $(q_p+\xi_p)_{-1\le p\le j}$ with $\xi_{-1}=r$ and $\xi_j=0$ satisfies the real orbit recurrence at the times $0\le p\le j-1$ exactly when $v(q_p)\xi_p-\xi_{p-1}-\xi_{p+1}=F_p(\xi_p)$ for these $p$. As $G_{l0}=t^+_0U_l=U_l$, this system is equivalent to
\[
\xi_l=U_lr+\mathcal{N}(\xi)_l\quad(0\le l<j),\qquad \mathcal{N}(\xi)_l=\sum_{p=0}^{j-1}G_{lp}F_p(\xi_p),
\]
and Taylor's formula with $\lvert\phi''\rvert<2\pi M$ (Lemma~\ref{lem:der_bou_the}, second derivative) gives $\lvert F_p(u)\rvert\le\pi Mu^2$ and $\lvert F_p'(u)\rvert=\lvert v(q_p+u)-v(q_p)\rvert\le2\pi M\lvert u\rvert$. On $\mathbb{R}^j$ use the weighted norm and the closed ball
\[
\lVert\xi\rVert_*=\max_{0\le l<j}M^{w(l+1)}\lvert\xi_l\rvert,\qquad\lVert\xi\rVert_*\le6.
\] For $\xi,\xi'$ in the ball, $\lvert F_p(\xi_p)-F_p(\xi'_p)\rvert\le2\pi M\cdot6M^{-w(p+1)}\cdot M^{-w(p+1)}\lVert\xi-\xi'\rVert_*$, hence
\[
M^{w(l+1)}\lvert\mathcal{N}(\xi)_l-\mathcal{N}(\xi')_l\rvert\le144\pi\,\lVert\xi-\xi'\rVert_*\sum_{p=0}^{j-1}M^{\,w(l+1)+1-2w(p+1)+\min(l,p)-\alpha(\max(l,p)+1)} .
\]
For $p\le l$ the exponent is $-(\alpha-w)(l+1)-(2w-1)(p+1)$, and the sum over these $p$ is at most $M^{-(\alpha+w-1)}/(1-M^{-(2w-1)})$, since $l+1\ge1$ and $\alpha>w$. For $p>l$, writing $p+1=l+1+i$ with $i\ge1$, the exponent is $-(\alpha+w-1)(l+1)-(\alpha+2w)i\le-(\alpha+2w)i$, and the sum is at most $M^{-(\alpha+2w)}/(1-M^{-(\alpha+2w)})$. Hence the Lipschitz constant of $\mathcal{N}$ on the ball is at most
\[
144\pi\Bigl(\frac{M^{-(\alpha+w-1)}}{1-M^{-(2w-1)}}+\frac{M^{-(\alpha+2w)}}{1-M^{-(\alpha+2w)}}\Bigr)\le\frac12 .
\]
Moreover $\mathcal{N}(0)=0$, so $\lVert\mathcal{N}(\xi)\rVert_*\le3$ on the ball, and $\lVert U\rVert_*\le12$ because $M^{w(l+1)}\lvert U_l\rvert\le12M^{-(\alpha-w)(l+1)}$, so $\lVert Ur\rVert_*\le3$ for $\lvert r\rvert\le1/4$. Thus $\xi\mapsto Ur+\mathcal{N}(\xi)$ maps the ball into itself and is a contraction of constant $\frac12$ there; let $\xi(r)$ be its unique fixed point. Put $q_p(r)=q_p+\xi_p(r)$ for $-1\le p\le j$, with $\xi_{-1}(r)=r$ and $\xi_j(r)=0$, and continue $(q_p(r))$ to all $p\in\mathbb{Z}$ by the real orbit recurrence from the pair $(q_0(r),q_{-1}(r))$; the continuation reproduces $q_1(r),\dots,q_j(r)$, so $(q_p(r))$ is a real lift of an orbit. Since $0$ is the fixed point at $r=0$, $q_p(0)=q_p$. This proves the first bound of (a).

We now prove smoothness, (c) and the second bound of (a). The estimate above also bounds the weighted operator norm of $\eta\mapsto G\,\mathrm{diag}\bigl(F_p'(\xi_p(r))\bigr)\eta$ by $\frac12$, since $\lvert F_p'(\xi_p(r))\rvert\le2\pi M\cdot6M^{-w(p+1)}$; so $I-G\,\mathrm{diag}(F'_p(\xi_p(r)))$ is invertible, with inverse of weighted norm at most $2$. By the implicit function theorem and the uniqueness of the fixed point in the ball, $r\mapsto\xi(r)$ is smooth, and differentiating the fixed-point equation gives $\mathrm{d}\xi/\mathrm{d}r=(I-G\,\mathrm{diag}F'_p(\xi_p(r)))^{-1}U$, whence $\lVert\mathrm{d}\xi/\mathrm{d}r\rVert_*\le2\cdot12=24$: the second bound of (a). Multiplying by $J$ gives $J(I-G\,\mathrm{diag}F'_p(\xi_p(r)))=J(r)$, the matrix with diagonal entries $v(q_p(r))=v(q_p)-F'_p(\xi_p(r))$ and entries $-1$ next to the diagonal; so $J(r)$ is invertible, and by Lemma~\ref{lem:gre_mat_fin} applied to $J(r)$, $\det J(r)=t^+_j(r)\ne0$. Differentiating the real recurrence with respect to the initial value $q_0$ at fixed $q_{-1}$ shows that $(\partial q_p/\partial q_0)_p$ is the linearized solution with initial values $1$ at time $0$ and $0$ at time $-1$; along $(q_p(r))$ this is $(t^+_p(r))$, which proves (c).

We next prove (d). The data $U$, $G$ and $F_p$ of the fixed-point equation depend continuously on $r$ and on the initial pair $(q_0,q_{-1})$ of a lift of a forward $j$-eligible state (for $U$ and $G$ through the formulas of Lemma~\ref{lem:gre_mat_fin}, since $t^+_j\ne0$), and the ball and the contraction constant are uniform. The fixed point of a uniform contraction depending continuously on parameters depends continuously on them, and so does the continuation by the recurrence.

We next prove (b). If $j=B_+$, the coordinate $q_{B_+}$ is fixed and $\lvert\xi_{B_+-1}(r)\rvert\le6M^{-wB_+}$. If $j=B_+-1$, the coordinate $q_{B_+-1}$ is fixed, and the recurrence at time $B_+-1$ gives $\xi_{B_+}(r)=-\xi_{B_+-2}(r)$, with $\lvert\xi_{B_+-2}(r)\rvert\le6M^{-w(B_+-1)}$ by (a). In both cases each coordinate of the pair changes by at most $6M^{-w(B_+-1)}$.

We finally prove (e). Put $U_p(r)=\mathrm{d}q_p(r)/\mathrm{d}r$; then $U_p(0)=U_p$, $U_{-1}(r)=1$, $U_j(r)=0$, and differentiating the recurrence in $r$ shows that $(U_p(r))$ is a linearized solution along $(q_p(r))$. By Lemma~\ref{lem:gre_mat_fin} applied to $J(r)$, the Green diagonal is $(J(r)^{-1})_{pp}=t^+_p(r)U_p(r)$. The derivative $J'(r)$ is diagonal with entries $\phi''(q_p(r))U_p(r)$, so Jacobi's formula for the derivative of a determinant gives
\[
\frac{\mathrm{d}}{\mathrm{d}r}\log\lvert t^+_j(r)\rvert=\operatorname{tr}\bigl(J(r)^{-1}J'(r)\bigr)=\sum_{p=0}^{j-1}t^+_p(r)\,U_p(r)^2\,\phi''(q_p(r)).
\]
Here $\lvert t^+_p(r)\rvert\le M^p$ (Lemma~\ref{lem:nor_and_gro}(i) along the orbit of the lift $(q_p(r))$), $\lvert U_p(r)\rvert\le24M^{-w(p+1)}$ by (a), and $\lvert\phi''\rvert<2\pi M$. Hence the absolute value of the sum is at most
\[
1152\pi\sum_{p=0}^{j-1}M^{-(2w-1)(p+1)}\le1152\pi\,\frac{M^{-(2w-1)}}{1-M^{-(2w-1)}} .
\]
Integrating over the interval between $r$ and $r'$, of length at most $\frac12$, gives $\bigl\lvert\log\lvert t^+_j(r)\rvert-\log\lvert t^+_j(r')\rvert\bigr\rvert\le576\pi M^{-(2w-1)}/(1-M^{-(2w-1)})\le\log2$.

Explicit values. With the values of Convention~\ref{conv:not_and_con}, $\alpha+w-1=37/50$, $2w-1=3/5$ and $\alpha+2w=127/50$. The contraction inequality holds for $M\ge9947$ and the distortion inequality $576\pi M^{-(2w-1)}/(1-M^{-(2w-1)})\le\log2$ for $M\ge4.96\cdot10^5$, hence both for $M\ge M_{\rm an}$; the stronger hypothesis $M\ge M_{\rm an}$ is kept for uniformity with the other estimates of the bridge argument.
\end{proof}

\begin{lemma}
 Let $M\ge M_{\rm an}$ and let time $0$ of an orbit have backward prefix growth of length $B_-\ge2$.
\begin{enumerate}
\item[(i)] Time $0$ is backward $j$-eligible for some $j\in\{B_--1,B_-\}$.
\item[(ii)] If time $0$ is backward $j$-eligible and $(q_p)$ is a real lift of the orbit, there is an unstable family $(q_p(r))$ of length $j$ through this lift with
\[
\lvert q_{-1-p}(r)-q_{-1-p}\rvert\le6M^{-w_{\rm pre}(p+1)},\qquad\Bigl\lvert\frac{\mathrm{d}q_{-1-p}(r)}{\mathrm{d}r}\Bigr\rvert\le24M^{-w_{\rm pre}(p+1)}
\]
for $0\le p<j$ and $\lvert r\rvert\le1/4$; in particular its unstable graph $r\mapsto(q_0+r,q_{-1}(r))$ has slope at most $24M^{-w_{\rm pre}}$.
\item[(iii)] Each of the coordinates $q_{-B_-}(r)$ and $q_{-B_--1}(r)$ differs from its value at $r=0$ by at most $6M^{-w_{\rm pre}(B_--1)}$.
\item[(iv)] Along the family, $\partial q_{-1-j}/\partial q_{-1}$ at fixed $q_0$ equals $t^-_{-1-j}(r)$, where $(t^-_p(r))$ is the linearized solution along $(q_p(r))$ with $(t^-_0(r),t^-_{-1}(r))=(0,1)$; moreover $t^-_{-1-j}(r)\ne0$ and $\frac12\le\lvert t^-_{-1-j}(r)\rvert/\lvert t^-_{-1-j}(r')\rvert\le2$ for $\lvert r\rvert,\lvert r'\rvert\le1/4$.
\item[(v)] For fixed $B_-$ and $j$, every $q_p(r)$ depends jointly continuously on $r$ and on the initial pair $(q_0,q_{-1})$ of the lift, as this pair ranges over the lifts of backward $j$-eligible states.
\end{enumerate}
\label{lem:bac_sta_gra:1}
\end{lemma}

\begin{proof}
Time reversal carries the backward problem to the forward one, and swapping the two coordinates of a state is an isometry.

We first describe the reversed orbit. For a real lift $(q_p)$ put $\widetilde q_p=q_{-1-p}$. The real orbit recurrence is symmetric under this reversal: $\widetilde q_{p+1}+\widetilde q_{p-1}=q_{-2-p}+q_{-p}=\phi(q_{-1-p})=\phi(\widetilde q_p)$. The positions of an orbit are determined by its state at time $0$ through \eqref{eq:orbit_recurrence}, and every point of $\mathbb{T}^2$ is the state at time $0$ of the orbit of its image under $\Pi^{-1}$ (Lemma~\ref{lem:inv_are_und}(ii)); so $(\widetilde q_p)$ is a real lift of the orbit, the \emph{reversed orbit}, whose state at time $0$ is $(q_{-1},q_0)$ modulo one. Its linearized solutions are the reversals $\widetilde u_p=u_{-1-p}$ of the linearized solutions $(u_p)$ of the original orbit, by the same symmetry of $u_{p+1}+u_{p-1}=v(q_p)u_p$. Hence its transfer matrix $\widetilde T_{0,l}$ carries $(u_{-1},u_0)$ to $(u_{-1-l},u_{-l})$, so that $\widetilde T_{0,l}$ is $T_{0,-l}$ conjugated by the swap of the two coordinates, which preserves the Euclidean norm, and by Lemma~\ref{lem:nor_and_gro}(i)
\[
\lVert\widetilde T_{0,l}\rVert=\lVert T_{0,-l}\rVert=\lVert T_{-l,0}\rVert\qquad(l\ge1).
\]
Therefore time $0$ of the original orbit has backward prefix growth of length $B_-$ exactly when time $0$ of the reversed orbit has forward prefix growth of length $B_-$. The linearized solution $(\widetilde t^+_p)$ of the reversed orbit with $(\widetilde t^+_0,\widetilde t^+_{-1})=(1,0)$, which is the solution $(t^+_p)$ of Definition~\ref{def:pre_gro_eli} for that orbit, is the reversal of $(t^-_p)$, since $t^-_{-1}=1$ and $t^-_0=0$; so $\widetilde t^+_j=t^-_{-1-j}$, and with $\lVert\widetilde T_{0,B_-}\rVert=\lVert T_{-B_-,0}\rVert$ this shows that backward $j$-eligibility of the original orbit is forward $j$-eligibility of the reversed orbit (Definition~\ref{def:pre_gro_eli}).

We now apply the forward lemmas to the reversed orbit. Lemma~\ref{lem:mos_con_dir:1}(i) gives (i). If time $0$ is backward $j$-eligible, Lemma~\ref{lem:sta_gra_und:1} applied to $(\widetilde q_p)$ gives a stable family $(\widetilde q_p(r))$ of length $j$ through it; put $q_p(r)=\widetilde q_{-1-p}(r)$. Then $q_0(r)=\widetilde q_{-1}(r)=q_0+r$ and $q_{-1-j}(r)=\widetilde q_j(r)=q_{-1-j}$, so $(q_p(r))$ is an unstable family of length $j$, and Lemma~\ref{lem:sta_gra_und:1}(a) for $\widetilde q_p=q_{-1-p}$, $0\le p<j$, is (ii); the slope of the unstable graph is $\lvert\mathrm{d}\widetilde q_0(r)/\mathrm{d}r\rvert\le24M^{-w_{\rm pre}}$. Item (b) of Lemma~\ref{lem:sta_gra_und:1} for the reversed pair $(\widetilde q_{B_-},\widetilde q_{B_--1})=(q_{-1-B_-},q_{-B_-})$ is (iii). Items (c) and (e) of Lemma~\ref{lem:sta_gra_und:1} give (iv), because $\partial\widetilde q_j/\partial\widetilde q_0$ at fixed $\widetilde q_{-1}$ is $\partial q_{-1-j}/\partial q_{-1}$ at fixed $q_0$, and the solution $(\widetilde t^+_p(r))$ of Lemma~\ref{lem:sta_gra_und:1}(c) along $(\widetilde q_p(r))$ is the reversal of $(t^-_p(r))$. Item (d) gives (v), since $(q_0,q_{-1})\mapsto(q_{-1},q_0)$ is continuous.
\end{proof}

The comparison of the two ends rests on a product bound for two transverse families of level curves, a change of variables of absolute-continuity type in the plane:

\begin{lemma}
 Let $0<s_0\le1$, let $I_x$ and $I_y$ be closed intervals of length $s_0$, put $S=I_x\times I_y$, and fix $x_c\in I_x$ and $y_c\in I_y$. Let $S'\subset\mathbb{R}^2$ be a Borel set, let $E_+,E_-\subset S$ be compact, and let $Q^+,Q^-\colon\mathbb{R}^2\to\mathbb{R}$ be $C^1$ functions with the following properties.
\begin{enumerate}
\item[(F)] Through each $p\in E_+$ passes a curve $\{(X_p(y),y):y\in J_p\}\subset S'$, where $J_p\supset I_y$ is an open interval and $X_p$ is $C^1$, on which $Q^+$ is constant and $\partial_xQ^+\ne0$, with $\lvert X_p'\rvert\le\frac12$ and $\frac12\le\lvert\partial_xQ^+(p')\rvert/\lvert\partial_xQ^+(p'')\rvert\le2$ for all points $p',p''$ of the curve; and $p\mapsto X_p(y_c)$ is continuous on $E_+$.
\item[(B)] Through each $p\in E_-$ passes a curve $\{(x,Y_p(x)):x\in J'_p\}\subset S'$, where $J'_p\supset I_x$ is an open interval and $Y_p$ is $C^1$, on which $Q^-$ is constant and $\partial_yQ^-\ne0$, with $\lvert Y_p'\rvert\le\frac12$ and $\frac12\le\lvert\partial_yQ^-(p')\rvert/\lvert\partial_yQ^-(p'')\rvert\le2$ for all points $p',p''$ of the curve; and $p\mapsto Y_p(x_c)$ is continuous on $E_-$.
\end{enumerate}
Let $W,H\colon\mathbb{R}^2\to\mathbb{R}$ be continuous with $W\ge0$, and let $\zeta\in\mathbb{R}$ and $\varepsilon\ge0$. If $\lvert H(p')-H(p)\rvert\le\varepsilon$ for every point $p'$ of the curve through each $p\in E_+$, and $W(p)\le W(p')+\varepsilon$ for every point $p'$ of the curve through each $p\in E_-$, then
\[
\int_{E_+\cap E_-\cap\{H\ge\zeta\}}W\le\frac{32}{3s_0^2}\Bigl(\int_{S'}W+\varepsilon\Bigr)\,\bigl\lvert S'\cap\{H\ge\zeta-\varepsilon\}\bigr\rvert ,
\]
the integrals being taken with respect to Lebesgue measure. The same bound holds if instead $\lvert H(p')-H(p)\rvert\le\varepsilon$ on the curve through each $p\in E_-$ and $W(p)\le W(p')+\varepsilon$ on the curve through each $p\in E_+$.
\label{lem:pro_bou_for}
\end{lemma}

\begin{proof}
The curves of (F) are labelled by their intercepts at $y_c$; the labels of curves through points with $H\ge\zeta$ have small measure because $H\ge\zeta-\varepsilon$ on these whole curves, and along each such curve the weight is controlled by sweeping it along the transverse curves of (B), on which $W$ almost does not decrease. This is the change-of-variables argument behind the absolute continuity of local stable manifolds \cite{locbib:Pesin1978}, here for two transverse families of curves.

We first record that level curves are determined by one point. Let $\Lambda\subset\mathbb{R}$ be an interval and let $x_1,x_2\colon\Lambda\to\mathbb{R}$ be continuous, with $Q^+(x_1(y),y)=Q^+(x_2(y),y)$ independent of $y\in\Lambda$ and $\partial_xQ^+\ne0$ at every point of the graph of $x_1$. If $x_1(y_0)=x_2(y_0)$ for some $y_0\in\Lambda$, then $x_1=x_2$ on $\Lambda$. Indeed, the set $\{x_1=x_2\}$ is closed in $\Lambda$; at each of its points $y_1$, the implicit function theorem shows that near $(x_1(y_1),y_1)$ the level set of $Q^+$ is the graph of one function of $y$, and both $x_i(y)$ stay in that neighbourhood for $y$ near $y_1$, so the set is open; and $\Lambda$ is connected. The same holds for $Q^-$ with the roles of the two coordinates exchanged. We refer to this as \emph{uniqueness of level curves}.

We next define the labels. For $p\in E_+$ call $X_p(y_c)$ the label of $p$, and let $\mathrm{Lab}$ be the set of labels of the points $p\in E_+$ with $H(p)\ge\zeta$. It is compact, as a continuous image of the compact set $E_+\cap\{H\ge\zeta\}$. If two curves of (F) meet, they lie in one level set of $Q^+$, so by uniqueness of level curves they coincide on the intersection of their parameter intervals, which contains $I_y$. Hence for $s\in\mathrm{Lab}$ we may put $X(s,y)=X_p(y)$ ($y\in I_y$) for any $p$ with label $s$; the point $(s,y_c)$ lies on this curve, and curves with distinct labels are disjoint over $I_y$. So the map $\Psi(s,y)=(X(s,y),y)$ is injective on $\mathrm{Lab}\times I_y$.

We next show that $\Psi$ is, near each label, the restriction of a $C^1$ diffeomorphism with Jacobian in $[\frac12,2]$. Fix $s_*\in\mathrm{Lab}$ and $p_*$ with label $s_*$, and a compact interval $\Lambda_*\subset J_{p_*}$ containing $I_y$ in its interior. Since $\partial_xQ^+$ is continuous and does not vanish on the compact arc $\{(X_{p_*}(y),y):y\in\Lambda_*\}$, there is $\rho>0$ such that $\partial_xQ^+\ne0$ on the compact tube $\{(x,y):y\in\Lambda_*,\ \lvert x-X_{p_*}(y)\rvert\le\rho\}$. Hence, for $y\in\Lambda_*$, $Q^+(\cdot,y)$ is strictly monotone on $[X_{p_*}(y)-\rho,X_{p_*}(y)+\rho]$, and for every $s$ the equation $Q^+(x,y)=Q^+(s,y_c)$ has at most one solution $x$ with $\lvert x-X_{p_*}(y)\rvert<\rho$. Since $Q^+$ is constant on the curve through $p_*$, which contains $(s_*,y_c)$, one has $Q^+(X_{p_*}(y),y)=Q^+(s_*,y_c)$ for $y\in\Lambda_*$, and strict monotonicity, continuity and compactness give
\[
\mu=\min_{y\in\Lambda_*}\min_\pm\bigl\lvert Q^+\bigl(X_{p_*}(y)\pm\tfrac\rho2,y\bigr)-Q^+(s_*,y_c)\bigr\rvert>0 .
\]
Choose $0<\rho'\le\rho$ with $\lvert Q^+(s,y_c)-Q^+(s_*,y_c)\rvert<\mu$ for $\lvert s-s_*\rvert<\rho'$, let $\Lambda'_*$ be the interior of $\Lambda_*$, an open interval containing $I_y$, and put $R_*=(s_*-\rho',s_*+\rho')\times\Lambda'_*$. For $(s,y)\in R_*$ the intermediate value theorem gives a solution $x$ of $Q^+(x,y)=Q^+(s,y_c)$ with $\lvert x-X_{p_*}(y)\rvert<\rho/2$, which is the unique solution within $\rho$; call it $\widehat X(s,y)$. At each $(s,y)\in R_*$ the implicit function theorem, applicable since $(\widehat X(s,y),y)$ lies in the tube, gives a $C^1$ local solution, which by continuity stays within $\rho$ of $X_{p_*}$ nearby and therefore coincides with $\widehat X$; so $\widehat X$ is $C^1$ on $R_*$. Since $\lvert s-X_{p_*}(y_c)\rvert=\lvert s-s_*\rvert<\rho$, uniqueness gives $\widehat X(s,y_c)=s$. For $s\in\mathrm{Lab}$ with $\lvert s-s_*\rvert<\rho'$ and $p$ with label $s$, the graphs of $X_p$ and $\widehat X(s,\cdot)$ over the interval $J_p\cap\Lambda'_*\supset I_y$ lie in the level set $\{Q^+=Q^+(s,y_c)\}$ and meet at $y_c$, so $\widehat X(s,\cdot)=X_p$ there by uniqueness of level curves; in particular $\widehat X(s,\cdot)=X(s,\cdot)$ on $I_y$, with the same derivative in $y$. Implicit differentiation of $Q^+(\widehat X(s,y),y)=Q^+(s,y_c)$ gives
\[
\partial_s\widehat X(s,y)=\frac{\partial_xQ^+(s,y_c)}{\partial_xQ^+(\widehat X(s,y),y)}\ne0
\]
on $R_*$ (both points lie in the tube), so for fixed $y$ the map $s\mapsto\widehat X(s,y)$ is strictly monotone, and $(s,y)\mapsto(\widehat X(s,y),y)$ is a $C^1$ diffeomorphism of $R_*$ onto an open set. For $s\in\mathrm{Lab}$ the points $(s,y_c)$ and $(X(s,y),y)$ lie on one curve of (F), so $\frac12\le\lvert\partial_s\widehat X(s,y)\rvert\le2$ by the ratio hypothesis. Covering the compact set $\mathrm{Lab}$ by finitely many such intervals and splitting it into disjoint Borel pieces, each inside one interval, we obtain from the injectivity of $\Psi$ and the change of variables for $C^1$ diffeomorphisms, for every Borel set $G\subset\mathrm{Lab}\times I_y$ and every Borel function $f\ge0$,
\[
\int_{\Psi(G)}f=\int_Gf(\Psi(s,y))\,\lvert\partial_sX(s,y)\rvert\,\mathrm{d}s\,\mathrm{d}y,\qquad\tfrac12\le\lvert\partial_sX\rvert\le2,
\]
where $\partial_sX$ is $\partial_s\widehat X$ of the piece containing $s$.

We next bound the measure of the labels. Every point $(X(s,y),y)$ with $s\in\mathrm{Lab}$ and $y\in I_y$ lies in $S'$ and on the curve through some $p\in E_+$ with $H(p)\ge\zeta$, so $H\ge\zeta-\varepsilon$ there. The change of variables with $G=\mathrm{Lab}\times I_y$ and $f=1$ gives $\frac12s_0\lvert\mathrm{Lab}\rvert\le\lvert S'\cap\{H\ge\zeta-\varepsilon\}\rvert$, that is,
\[
\lvert\mathrm{Lab}\rvert\le\frac2{s_0}\,\bigl\lvert S'\cap\{H\ge\zeta-\varepsilon\}\bigr\rvert .
\]

We next bound the fibre integral of the weight along one labelled curve by a transverse sweep. Fix $s\in\mathrm{Lab}$ and let $\mathrm{Elig}_s=\{y\in I_y:(X(s,y),y)\in E_-\}$, a compact set. For $y\in\mathrm{Elig}_s$ put $p(y)=(X(s,y),y)$ and $\Theta_s(u,y)=Y_{p(y)}(u)$ for $u\in I_x$, so that $(u,\Theta_s(u,y))$ is the point over $u$ of the curve of (B) through $p(y)$. For fixed $u$ the map $y\mapsto\Theta_s(u,y)$ is injective on $\mathrm{Elig}_s$. Indeed, if $\Theta_s(u,y)=\Theta_s(u,y')$, the curves of (B) through $p(y)$ and $p(y')$ meet, so by uniqueness of level curves they coincide over $I_x$; then $p(y)$ and $p(y')$ lie on this curve and on the curve of label $s$, and the two slope bounds give $\lvert X(s,y)-X(s,y')\rvert\le\frac12\lvert y-y'\rvert$ and $\lvert y-y'\rvert\le\frac12\lvert X(s,y)-X(s,y')\rvert$, so $y=y'$. Hence $\Phi_s(u,y)=(u,\Theta_s(u,y))$ is injective on $I_x\times\mathrm{Elig}_s$, with image in $S'$.

We next compute the Jacobian of the sweep. Fix $y_*\in\mathrm{Elig}_s$. The starting points $(\widehat X(s,y),y)$ depend $C^1$ on $y$, and $\partial_yQ^-\ne0$ on the compact arc of the curve of (B) through $p(y_*)$ over a compact interval containing $I_x$ in its interior. As before, the implicit function theorem along this arc gives an open interval $\Lambda''\supset I_x$, a $\rho''>0$ and a $C^1$ function $\widehat\Theta$ on $\Lambda''\times(y_*-\rho'',y_*+\rho'')$ with $Q^-(u,\widehat\Theta(u,y))=Q^-(\widehat X(s,y),y)$ and $\widehat\Theta(\widehat X(s,y),y)=y$; for $y\in\mathrm{Elig}_s$ near $y_*$, the graphs of $\widehat\Theta(\cdot,y)$ and $Y_{p(y)}$ lie in one level set of $Q^-$ and meet over $X(s,y)$, so uniqueness of level curves gives $\widehat\Theta(\cdot,y)=Y_{p(y)}$ on $\Lambda''\cap J'_{p(y)}\supset I_x$, in particular $\widehat\Theta(\cdot,y)=\Theta_s(\cdot,y)$ on $I_x$. Differentiation in $y$ gives
\[
\partial_y\widehat\Theta(u,y)=\frac{\partial_yQ^-(p(y))+\partial_xQ^-(p(y))\,\partial_yX(s,y)}{\partial_yQ^-(u,\widehat\Theta(u,y))} .
\]
For $y\in\mathrm{Elig}_s$, the slope of the curve of (B) at $p(y)$ is $-\partial_xQ^-/\partial_yQ^-$, of absolute value at most $\frac12$, and $\lvert\partial_yX(s,y)\rvert\le\frac12$; so the numerator has absolute value at least $\frac34\lvert\partial_yQ^-(p(y))\rvert$, and the ratio hypothesis of (B) gives $\lvert\partial_y\widehat\Theta(u,y)\rvert\ge\frac38$ for $u\in I_x$. After shrinking $\rho''$ and $\Lambda''$, $\partial_y\widehat\Theta\ne0$, and $(u,y)\mapsto(u,\widehat\Theta(u,y))$ is a $C^1$ diffeomorphism onto an open set. Splitting the compact set $\mathrm{Elig}_s$ into finitely many Borel pieces as before, the injectivity of $\Phi_s$, the change of variables and $W\ge0$ give
\[
\frac38\int_{I_x}\int_{\mathrm{Elig}_s}W(u,\Theta_s(u,y))\,\mathrm{d}y\,\mathrm{d}u\le\int_{\Phi_s(I_x\times\mathrm{Elig}_s)}W\le\int_{S'}W .
\]
By the transport hypothesis on the curves of (B), $W(X(s,y),y)\le W(u,\Theta_s(u,y))+\varepsilon$ for $u\in I_x$ and $y\in\mathrm{Elig}_s$. Integrating over $u\in I_x$ and $y\in\mathrm{Elig}_s$ gives $s_0\int_{\mathrm{Elig}_s}W(X(s,y),y)\,\mathrm{d}y\le\frac83\int_{S'}W+\varepsilon s_0^2$, and since $s_0^2\le1\le\frac83$,
\[
\int_{\mathrm{Elig}_s}W(X(s,y),y)\,\mathrm{d}y\le\frac8{3s_0}\Bigl(\int_{S'}W+\varepsilon\Bigr)\qquad(s\in\mathrm{Lab}).
\]

We now combine the two bounds. The set $G=\{(s,y)\in\mathrm{Lab}\times I_y:\Psi(s,y)\in E_-\}=\{(s,y):s\in\mathrm{Lab},\ y\in\mathrm{Elig}_s\}$ is compact, since $\Psi$ is continuous on $\mathrm{Lab}\times I_y$. Every $p\in E_+\cap E_-\cap\{H\ge\zeta\}$ is $\Psi(s,y)$ with $s$ its label and $y\in\mathrm{Elig}_s$ its second coordinate, so it lies in $\Psi(G)$. Since $W\ge0$, the change of variables, $\lvert\partial_sX\rvert\le2$ and Fubini--Tonelli give
\[
\int_{E_+\cap E_-\cap\{H\ge\zeta\}}W\le\int_{\Psi(G)}W\le2\int_{\mathrm{Lab}}\int_{\mathrm{Elig}_s}W(X(s,y),y)\,\mathrm{d}y\,\mathrm{d}s\le2\,\lvert\mathrm{Lab}\rvert\,\frac8{3s_0}\Bigl(\int_{S'}W+\varepsilon\Bigr),
\]
and the bound on $\lvert\mathrm{Lab}\rvert$ gives the first assertion, with $2\cdot\frac2{s_0}\cdot\frac8{3s_0}=\frac{32}{3s_0^2}$.

Finally, the reflection $(x,y)\mapsto(y,x)$ preserves Lebesgue measure and exchanges the two kinds of curves: it carries the curves of (B), with $Q^-$, to curves of the kind (F) for $Q^-$ composed with the reflection, labelled at $x_c$, and the curves of (F) to curves of the kind (B). Applying the first assertion to the reflected data, with the roles of $E_+$ and $E_-$ exchanged, gives the second.
\end{proof}

The estimate applies to shortfall observables, Lipschitz functions of finitely many interval shortfalls of the growth semimetric:

\begin{definition}
Let $L\ge0$. A \textbf{shortfall observable with Lipschitz constant at most $L$} is a function $W$ of the initial point $z$ of the orbit (Definition~\ref{def:orb_sta_and}) given by the following data: an integer $r\ge1$; integer pairs $a_j\ne b_j$ for $1\le j\le r$; and a function $F\colon\mathbb{R}^r\to\mathbb{R}$ whose restriction to $[0,1]^r$ has Lipschitz constant at most $L$ for the maximum norm. The observable is
\[
W=F\bigl(h_g(a_1,b_1),\dots,h_g(a_r,b_r)\bigr),
\]
where $h_g$ is the shortfall (Definition~\ref{def:sho_sem_the}) of the growth semimetric $g$ of Definition~\ref{def:gro_sem_orb}. The integers $a_j,b_j$ of the chosen data are the \textbf{indices occurring} in $W$.
\label{def:shor_obse}
\end{definition}

\begin{thm}
 Let $M\ge M_{\rm an}$ and $0<\zeta\le1$, and let $a<b$ be integers with $n=b-a\ge1000+100\log(1/\zeta)$. Let $W\ge0$ and $H$ be shortfall observables with Lipschitz constant at most $6$, such that every index occurring in one of them lies within distance $\kappa_{\rm br}n$ of $a$ and every index occurring in the other lies within distance $\kappa_{\rm br}n$ of $b$. Then
\[
\mathbb{E}\bigl[\mathbf{1}_{\{h_g(a,b)\le\beta_{\rm br}\}}\,W\,\mathbf{1}_{\{H>\zeta\}}\bigr]\le C_{\rm br}\bigl(\mathbb{E}W+M^{-\sigma_{\rm br}n}\bigr)\,\mathbb{P}(H>\zeta/2).
\]
\label{thm:bri_pro_est:2}
\end{thm}

\begin{proof}
On the set $\{h_g(a,b)\le\beta_{\rm br}\}$ a positive proportion of the middle times of $[a,b]$ have prefix growth in both directions; at such a time the state carries a family of stable graphs and a transverse family of unstable graphs, the observables at each end of $[a,b]$ are almost constant along the graphs of the corresponding family, and Lemma~\ref{lem:pro_bou_for} compares the two ends by a change of variables. The interval $[a,b]$ of near-full growth thus acts as a \emph{bridge} that decouples the observables at its two ends, and no mixing or decorrelation property of $f_k$ is used. In this proof write $c=c_{\rm pre}$ and $w=w_{\rm pre}$ (Convention~\ref{conv:not_and_con}).

We first find many middle times with two-sided prefix growth. Put $G_r=g(a,a+r)$ for $0\le r\le n$. By the triangle inequality and \eqref{eq:growth_bounds} (Lemma~\ref{lem:nor_and_gro}(ii)), $G_r-G_{r-1}\le g(a+r-1,a+r)\le1$, and $\Gamma=n-(G_n-G_0)=n\,h_g(a,b)$. On the set $\{h_g(a,b)\le\beta_{\rm br}\}$, Lemma~\ref{lem:pli_typ_cou}, applied to $(G_r)_{0\le r\le n}$ with $N=n$ and $c=c_{\rm pre}$, bounds by $\beta_{\rm br}n/(1-c_{\rm pre})=n/48$ the number of indices $r\in\{1,\dots,n\}$ for which some $0\le l<r$ satisfies $G_r-G_l<c_{\rm pre}(r-l)$, and the number of indices $r\in\{0,\dots,n-1\}$ for which some $r<l\le n$ satisfies $G_l-G_r<c_{\rm pre}(l-r)$. Let $r\in\{0,\dots,n\}$ be counted by neither, and put $i=a+r$. For $1\le l\le r$ the triangle inequality gives $g(i-l,i)\ge G_r-G_{r-l}\ge c_{\rm pre}l$, and for $1\le l\le n-r$ it gives $g(i,i+l)\ge G_{r+l}-G_r\ge c_{\rm pre}l$: time $i$ has backward prefix growth of length $i-a$ and forward prefix growth of length $b-i$ (Definition~\ref{def:pre_gro_eli}). The middle third $a+n/3\le i\le a+2n/3$ contains at least $n/3-1$ integers, so at least $n/3-1-n/24\ge n/12$ such times. Hence $\mathbf{1}_{\{h_g(a,b)\le\beta_{\rm br}\}}\le\frac{12}{n}\sum_i\mathbf{1}_{\mathrm{Pre}_i}$, the sum over the at most $n/3+1$ integers $i$ of the middle third, where $\mathrm{Pre}_i$ is the set on which time $i$ has these two prefix growths. Multiplying by $W\mathbf{1}_{\{H>\zeta\}}\ge0$ and integrating, and using $\frac{12}{n}(\frac n3+1)\le8$, we see that it suffices to prove, for every $i$ of the middle third, the \emph{middle-time bound}
\[
\mathbb{E}\bigl[\mathbf{1}_{\mathrm{Pre}_i}\,W\,\mathbf{1}_{\{H>\zeta\}}\bigr]\le\frac{2^{27}}3\bigl(\mathbb{E}W+M^{-\sigma_{\rm br}n}\bigr)\,\mathbb{P}(H>\zeta/2),
\]
since $8\cdot2^{27}/3=2^{30}/3<C_{\rm br}$.

We next move the middle time to $0$. Fix $i$ and put $B_+=b-i$ and $B_-=i-a$; both are at least $n/3\ge2$. Let $\mathrm{Pre}$ be the set on which time $0$ has forward prefix growth of length $B_+$ and backward prefix growth of length $B_-$, and let $W'$ and $H'$ be the shortfall observables with the same data functions as $W$ and $H$ and all indices shifted by $-i$. By Lemma~\ref{lem:inv_are_und}(iii), $\mathbf{1}_{\mathrm{Pre}_i}=\mathbf{1}_{\mathrm{Pre}}\circ f_k^i$, $W=W'\circ f_k^i$ and $H=H'\circ f_k^i$, so by Lemma~\ref{lem:inv_are_und}(i) the three integrals in the middle-time bound do not change when $(\mathrm{Pre}_i,W,H)$ is replaced by $(\mathrm{Pre},W',H')$. The indices occurring in $W'$ and $H'$ lie within $\kappa_{\rm br}n$ of $-B_-$ and of $B_+$, in one order or the other. We drop the primes. The orbit, and with it every quantity above, is determined by its state $(q_0,q_{-1})$ at time $0$, and the state map $z\mapsto(q_0,q_{-1})$ is a bijection of $\mathbb{T}^2$ carrying $m$ to $m$ (Lemma~\ref{lem:inv_are_und}(ii)). So the integrals may be computed over the state $(x,y)\in\mathbb{T}^2$; for $(x,y)\in\mathbb{R}^2$ we evaluate all quantities at the orbit with state $(x,y)$ modulo one at time $0$, along its real lift with $(q_0,q_{-1})=(x,y)$. Then $W$ and $H$ are continuous functions on $\mathbb{R}^2$ of period one in each variable: each shortfall $h_g(a',b')$ is a continuous function of the state, with values in $[0,1]$ (Lemma~\ref{lem:nor_and_gro}(ii)), and the data functions are Lipschitz on the unit cube.

We next set up the families of stable and unstable graphs. By Lemma~\ref{lem:mos_con_dir:1}(i) and Lemma~\ref{lem:bac_sta_gra:1}(i), at every state of $\mathrm{Pre}$ time $0$ is forward $j_+$-eligible and backward $j_-$-eligible for some $j_+\in\{B_+-1,B_+\}$ and $j_-\in\{B_--1,B_-\}$; fix one of these four pairs $(j_+,j_-)$. Cover $[0,1]^2$ by the $2^{10}$ closed squares $S=I_x\times I_y$ of the grid of side $s_0=1/32$, fix one of them, and let $x_c,y_c$ be the left endpoints of $I_x,I_y$. Let $E_+$ be the set of $(x,y)\in S$ at which time $0$ has forward prefix growth of length $B_+$ and is forward $j_+$-eligible, and $E_-$ the set of $(x,y)\in S$ at which time $0$ has backward prefix growth of length $B_-$ and is backward $j_-$-eligible; they are compact, since their defining inequalities are non-strict inequalities between continuous functions of $(x,y)$. Let $Q^+(x,y)$ and $Q^-(x,y)$ be the coordinates $q_{j_+}$ and $q_{-1-j_-}$ of the real lift with $(q_0,q_{-1})=(x,y)$; they are smooth on $\mathbb{R}^2$, as the real recurrence is. For $p=(x_p,y_p)\in E_+$, Lemma~\ref{lem:sta_gra_und:1} gives a stable family of length $j_+$ through this lift; its stable graph is the curve $\{(X_p(y),y):\lvert y-y_p\rvert<1/4\}$ with $X_p(y)=q_0(y-y_p)$. On it $Q^+$ is constant, $\partial_xQ^+(X_p(y),y)=t^+_{j_+}(y-y_p)\ne0$ with ratios in $[\frac12,2]$ (items (c) and (e)), $\lvert X_p'\rvert\le24M^{-w}\le\frac12$ (item (a) at index $0$), and $p\mapsto X_p(y_c)$ is continuous on $E_+$ (item (d)). Symmetrically, for $p\in E_-$, Lemma~\ref{lem:bac_sta_gra:1}(ii),(iv),(v) give the unstable graph $\{(x,Y_p(x)):\lvert x-x_p\rvert<1/4\}$ with $Y_p(x)=q_{-1}(x-x_p)$, on which $Q^-$ is constant, $\partial_yQ^-(x,Y_p(x))=t^-_{-1-j_-}(x-x_p)\ne0$ with ratios in $[\frac12,2]$, $\lvert Y_p'\rvert\le\frac12$, and $p\mapsto Y_p(x_c)$ is continuous. Since $s_0<1/4$, these parameter intervals contain $I_y$ and $I_x$ respectively, and since the slopes are at most $\frac12$, the curves lie in $S'=(I_x+[-\frac14,\frac14])\times(I_y+[-\frac14,\frac14])$, a square of side $17/32$. Thus the hypotheses (F) and (B) of Lemma~\ref{lem:pro_bou_for} hold: (F) asks for a continuously labelled family of curves $x=X_p(y)$, $p\in E_+$, spanning $I_y$ inside $S'$, of slope at most $\frac12$, on which $Q^+$ is constant and $\partial_xQ^+$ is nonzero and varies by a factor at most $2$; (B) asks for the same of the curves $y=Y_p(x)$, $p\in E_-$, and $Q^-$, with the roles of the two coordinates exchanged.

We next show that the observables at each end are almost constant along the stable or unstable graphs of that end. Let $p'$ lie on the stable graph through $p\in E_+$. By Lemma~\ref{lem:sta_gra_und:1}(b) the states at time $B_+$ of the two lifts differ by at most $6M^{-w(B_+-1)}\le6M^{w-wn/3}$ in the maximum norm, since $B_+\ge n/3$. One step of \eqref{eq:orbit_recurrence}, forward or backward, multiplies the maximum norm of a difference of states of real lifts by at most $\sup\lvert v\rvert+1\le K+3<M$ (mean value inequality, \eqref{eq:large_parameters}); so at all times within $\kappa_{\rm br}n+1$ of $B_+$ they differ by at most $6M^{\kappa_{\rm br}n+1+w-wn/3}$. As $\lvert v(x)-v(x')\rvert\le2\pi K\lvert x-x'\rvert$ (Lemma~\ref{lem:der_bou_the}), each one-step matrix $A_l$ changes by at most $2\pi M$ times this, and for an interval $[a',b']$ of length $\ell\le2\kappa_{\rm br}n$ in the window the telescoping identity and $\lVert A_l^{\pm1}\rVert\le M$ (Lemma~\ref{lem:nor_and_gro}(i)) give
\[
\lVert T_{a',b'}(p')-T_{a',b'}(p)\rVert\le24\pi\kappa_{\rm br}\,n\,M^{1+w}M^{-(w/3-3\kappa_{\rm br})n}.
\]
Both norms are at least $1$ (Lemma~\ref{lem:nor_and_gro}(i)), $\lvert\log u-\log u'\rvert\le\lvert u-u'\rvert$ for $u,u'\ge1$, and $\log M\ge1$, so the shortfalls $h_g(a',b')$ at $p$ and $p'$ differ by at most this number, and a shortfall observable with Lipschitz constant at most $6$ by at most $6\cdot24\pi\kappa_{\rm br}nM^{1+w}M^{-(w/3-3\kappa_{\rm br})n}\le\varepsilon_n$, where
\[
\varepsilon_n=23\,n\,M^{1+w}M^{-(w/3-3\kappa_{\rm br})n}\le M^{-\sigma_{\rm br}n}<\zeta/2 .
\]
Indeed $w/3-3\kappa_{\rm br}-\sigma_{\rm br}=1/60$, and $n\ge1000$ gives $n/60-(1+w)\ge n/120$ and $23n\le M^{n/120}$; and $M^{-\sigma_{\rm br}n}<\zeta/2$ because $n\ge1000+100\log(1/\zeta)$ and $\log M\ge10^{12}$. The same holds at the other end along unstable graphs, by Lemma~\ref{lem:bac_sta_gra:1}(iii) in place of Lemma~\ref{lem:sta_gra_und:1}(b).

We now apply Lemma~\ref{lem:pro_bou_for}. If the indices of $H$ lie within $\kappa_{\rm br}n$ of $B_+$ and those of $W$ within $\kappa_{\rm br}n$ of $-B_-$, the previous paragraph gives $\lvert H(p')-H(p)\rvert\le\varepsilon_n$ on the stable graph through each $p\in E_+$ and $W(p)\le W(p')+\varepsilon_n$ on the unstable graph through each $p\in E_-$, and the first assertion of Lemma~\ref{lem:pro_bou_for}, with $\varepsilon=\varepsilon_n$ and $32/(3s_0^2)=2^{15}/3$, gives
\[
\int_{E_+\cap E_-\cap\{H\ge\zeta\}}W\le\frac{2^{15}}3\Bigl(\int_{S'}W+\varepsilon_n\Bigr)\,\bigl\lvert S'\cap\{H\ge\zeta-\varepsilon_n\}\bigr\rvert .
\]
In the other case the second assertion of Lemma~\ref{lem:pro_bou_for} gives the same bound. The square $S'$ has side less than $1$, so the projection $\mathbb{R}^2\to\mathbb{T}^2$ is injective on it and carries Lebesgue measure to $m$; hence $\int_{S'}W\le\mathbb{E}W$ and, as $\varepsilon_n<\zeta/2$, $\lvert S'\cap\{H\ge\zeta-\varepsilon_n\}\rvert\le\mathbb{P}(H>\zeta/2)$. Every state of $\mathrm{Pre}$ in $[0,1)^2$ lies in $E_+\cap E_-$ for some square and some pair of choices, and $\mathbf{1}_{\{H>\zeta\}}\le\mathbf{1}_{\{H\ge\zeta\}}$; summing over the $2^{10}$ squares and the $4$ pairs gives
\[
\mathbb{E}\bigl[\mathbf{1}_{\mathrm{Pre}}W\mathbf{1}_{\{H>\zeta\}}\bigr]\le2^{12}\cdot\frac{2^{15}}3\bigl(\mathbb{E}W+M^{-\sigma_{\rm br}n}\bigr)\,\mathbb{P}(H>\zeta/2),
\]
which is the middle-time bound. This completes the proof.

Explicit values. With the values of Convention~\ref{conv:not_and_con}: $\beta_{\rm br}/(1-c_{\rm pre})=1/48$; $n/3-1-n/24\ge n/12$ for $n\ge5$ and $\frac{12}n(\frac n3+1)\le8$ for $n\ge3$; $6\cdot24\pi\kappa_{\rm br}=36\pi/5<22.7$; $w_{\rm pre}/3-3\kappa_{\rm br}=7/60$ and $7/60-\sigma_{\rm br}=1/60$; $n/60-(1+w_{\rm pre})\ge n/120$ for $n\ge216$; $23\le M_{\rm an}^{1/120}$; $32/(3s_0^2)=2^{15}/3$ at $s_0=1/32$; and $2^{30}/3<2^{28.42}<C_{\rm br}$.
\end{proof}

\section{Balance across scales} \label{sec:bal_acr_sca}

This section proves the balance across scales Proposition~\ref{prop:bal_bet_dya}, the heart of the paper: if $M\ge M_{\rm an}$ and $n_0<N$ are of the form $2^j$ with $n_0\ge n_{\rm bal}$ and $e_N\le\zeta_{\rm bal}$, then $e_N\le C_{\rm bal}(e_{n_0}+1/n_0)$. The core idea is to find a quantity that is comparable to the mean shortfall and whose total growth over all doublings of the time scale is small. Doubling acts exactly on the rescaled growth arrays, $\mathcal{D}g_n=g_{2n}$, and shifts of time preserve area, so a balance computed from the array at scale $n$ alone can have the mean $\mathbb{E}\,A(g_{2n})-\mathbb{E}\,A(g_n)$ of an increment, and such increments telescope over the scales. We use the potential
\[
\mathcal{L}_n=e^\Phi_n-\mathbb{E}\,\widetilde U(g_n),
\]
the capped mean shortfall minus the mean of the extended corrector (Definition~\ref{def:the_cap_and}, Definition~\ref{def:ext_cor_and}). Its increments are the mean balance defects, $\mathcal{L}_{2n}-\mathcal{L}_n=\mathbb{E}\,Z(g_n)$, and the balance follows from three facts: at the last scale $\mathcal{L}_N$ is at least a fixed fraction of $e_N$, at the first scale $\mathcal{L}_{n_0}$ is at most a multiple of $e_{n_0}$, and the increments summed over the scales are small (Proposition~\ref{prop:def_the_ext:1}).

The potential grows little because of the shapes that the rescaled growth can take. By the scarcity of cancellations and the bridge estimate, outside a set of small mass the rescaled growth array is close to a slack profile (Lemma~\ref{lem:sum_exc_mas:1}), for instance one slow core of length $L$ between two stretches of full growth, whose crossing deficit $\Delta$, the deficit of a window that contains the core with room on both sides, is at most $C_{\rm cr}L$. On such a profile the capped shortfall of a window is its slow fraction, the fraction of the window inside the core, times a weight, plus a corrector that lives at the interfaces. The slow fraction is additive under doubling. The weight is $1$ for windows of length at most $A_{\rm sl}L$, a large fixed multiple of the core length, in particular for windows inside the core, where the cap saturates; it is $\alpha_{\rm cap}\Delta/L$ for windows of length at least $2A_{\rm sl}L$, which see the core only through its crossing deficit and on which the cap is linear. Since $\alpha_{\rm cap}C_{\rm cr}=1/20$, a swallowed core counts less than the slow time inside it, so the weight decreases with the window length, and the balance defect is never positive on slack profiles (Lemma~\ref{lem:cor_ine_sla}, Lemma~\ref{lem:pro_the_bal}).

Two constraints shape the rest of the construction. The number of scales is unbounded, and the quantities whose sums over all scales are controlled are the midpoint cancellations and the imbalances (Lemma~\ref{lem:sum_mid_can}); an error charged at a single scale must therefore vanish, be bounded by these quantities, be carried by a set of summable mass, or be summable over $n=2^j$, like a multiple of $1/n$. And the inequality to be proved has no constant term, so no error of fixed size can be afforded even at a single scale; the bootstrap applies it at shortfall levels of order $C_{\rm ent}/\log M$. The corrector is extended from the slack profiles to all arrays under these constraints, and the defect bound treats separately the initial points at which some midpoint cancellation or imbalance is not tiny, a set of summable mass on which the arrays are close to slack profiles outside a further small set, and the others, at which they are close to affine arrays.

The first subsection sets up arrays and shows that midpoint cancellations and imbalances telescope across scales; the second defines the cap, the slack profiles, the corrector and its extension, and proves the balance; the third proves the properties of the corrector and of the balance defect on slack profiles; the fourth approximates rescaled growth arrays by slack profiles outside a set of small mass; the fifth proves the defect bound.

\subsection{Arrays and the telescoping of midpoint cancellations} \label{subsec:arr_mid_can}

Rescaled by $n$, the growth semimetric becomes an array, a semimetric on the real line. The dilation $\mathcal{D}$ passes from scale $n$ to scale $2n$, and translations by multiples of $1/n$ are time shifts of the orbit (Lemma~\ref{lem:dil_and_tim}). The midpoint cancellation and the imbalance of an interval measure how far an array is from affine there; summed over the scales they telescope into mean shortfalls (Lemma~\ref{lem:sum_mid_can}), which makes the body of the defect estimate summable. The elementary properties of shortfalls used throughout the section are collected in Lemma~\ref{lem:shor_arra}.

\begin{definition}
 An \textbf{array} is a semimetric $d$ on $\mathbb{R}$ (Definition~\ref{def:gro_hyp_sem}). For an interval $I=[s,t]$ with $s<t$ write $h(I)=h_d(s,t)$ for the shortfall of $d$ on $I$ (Definition~\ref{def:sho_sem_the}); for $x\in\mathbb{R}$ and $t>0$ the shortfall of $[x,x+t]$ is also written $h_t(x)$. The \textbf{halves} of $I$ are $I_l=[s,\frac{s+t}2]$ and $I_r=[\frac{s+t}2,t]$, and the \textbf{midpoint cancellation} and the \textbf{imbalance} of $I$ are
\[
J(I)=\frac{2(s|t)_{(s+t)/2}}{t-s},\qquad V(I)=\bigl(h(I_l)-h(I_r)\bigr)^2,
\]
the Gromov product being that of $d$ at the basepoint $(s+t)/2$. When the array must be named, these quantities are written $h(I)(d)$, $h_t(x;d)$, $J(I)(d)$ and $V(I)(d)$. The interval $I$ is \textbf{aligned at scale} $n$, for an integer $n\ge1$, if $ns$, $n(s+t)/2$ and $nt$ are integers. For arrays $d,d''$ and an interval $\Theta$, the \textbf{uniform distance on $\Theta$} is
\[
\lVert d-d''\rVert_\Theta=\sup\bigl\{\lvert d(x,y)-d''(x,y)\rvert:\ x,y\in\Theta\bigr\}\in[0,\infty].
\]

A \textbf{secant} $[x,y]$, $x<y$, has \textbf{speed} $d(x,y)/(y-x)$. A set of real times is \textbf{slow at speed} $\sigma$ if every secant with endpoints in it has speed at most $\sigma$; a closed half-line is a \textbf{unit geodesic ray} if every secant with endpoints in it has speed one. For $0\le\sigma\le1$ the \textbf{affine array} of speed $\sigma$ is $d_\sigma(x,y)=\sigma\lvert x-y\rvert$. For $x\in\mathbb{R}$ the \textbf{translation} $\mathcal{T}_x$ and the \textbf{dilation} $\mathcal{D}$ act on arrays by
\[
(\mathcal{T}_xd)(s,t)=d(s+x,t+x),\qquad (\mathcal{D}d)(s,t)=\tfrac12\,d(2s,2t).
\]
An array $d$ is \textbf{admissible} if, for $t\in\{1/2,1,2\}$, the function $x\mapsto d(x,x+t)$ is Borel on $\mathbb{R}$ with values in $[0,t]$.
\label{def:arr_mid_can}
\end{definition}

\begin{definition}
 For an integer $n\ge1$, the \textbf{rescaled growth array} at scale $n$ is the function on $\mathbb{R}\times\mathbb{R}$
\[
g_n(s,t)=\frac{g(\lfloor ns\rfloor,\lfloor nt\rfloor)}{n}\qquad(s,t\in\mathbb{R}),
\]
regarded, like $g$, as a function of the initial point $z$ of the orbit (Definition~\ref{def:orb_sta_and}).
\label{def:res_gro_arr:1}
\end{definition}

\begin{lemma}
 Let $d$ and $d''$ be arrays, let $I=[s,t]$ with $s<t$, and let $r\in\mathbb{R}$; write $I+r=[s+r,t+r]$ and $2I=[2s,2t]$, and likewise for an interval $\Theta$.
\begin{enumerate}
\item[(i)] $J(I)(d)=h(I)(d)-\frac12\bigl(h(I_l)(d)+h(I_r)(d)\bigr)\ge0$, and
\[
\bigl\lvert h(I_l)(d)-h(I)(d)\bigr\rvert\le J(I)(d)+\tfrac12\sqrt{V(I)(d)},\quad\bigl\lvert h(I_r)(d)-h(I)(d)\bigr\rvert\le J(I)(d)+\tfrac12\sqrt{V(I)(d)};
\]
equivalently, the speeds of the halves of $I$ differ from the speed of $I$ by at most $J(I)(d)+\frac12\sqrt{V(I)(d)}$.
\item[(ii)] $h(I)(\mathcal{T}_rd)=h(I+r)(d)$ and $h(I)(\mathcal{D}d)=h(2I)(d)$, and the same identities hold with $J$ or with $V$ in place of $h$; in particular $h_t(x;\mathcal{T}_rd)=h_t(x+r;d)$ and $h_t(x;\mathcal{D}d)=h_{2t}(2x;d)$ for $x\in\mathbb{R}$ and $t>0$. For every interval $\Theta$, $\lVert\mathcal{T}_rd-\mathcal{T}_rd''\rVert_\Theta=\lVert d-d''\rVert_{\Theta+r}$ and $\lVert\mathcal{D}d-\mathcal{D}d''\rVert_\Theta=\frac12\lVert d-d''\rVert_{2\Theta}$. The array $\mathcal{T}_rd$ is admissible if and only if $d$ is.
\item[(iii)] $\bigl\lvert h(I)(d)-h(I)(d'')\bigr\rvert\le\lVert d-d''\rVert_I/(t-s)$. If the shortfalls of $I_l$ and of $I_r$ lie in $[0,1]$ for $d$ and for $d''$, then $\bigl\lvert V(I)(d)-V(I)(d'')\bigr\rvert\le8\lVert d-d''\rVert_I/(t-s)$.
\item[(iv)] If $d(x,y)\le\lvert x-y\rvert$ for all $x,y\in\mathbb{R}$, then $\lvert d(x,y)-d(x',y')\rvert\le\lvert x-x'\rvert+\lvert y-y'\rvert$ for all $x,x',y,y'\in\mathbb{R}$, every shortfall of $d$ lies in $[0,1]$, the function $x\mapsto h_t(x;d)$ is continuous for every $t>0$, and $d$ is admissible.
\end{enumerate}
\label{lem:shor_arra}
\end{lemma}

\begin{proof}
Each item unpacks the definitions; the inequality of (i) is the triangle inequality at the midpoint of $I$.

(i) Put $w=(s+t)/2$ and $\ell=t-s$, so that $I_l=[s,w]$ and $I_r=[w,t]$ have length $\ell/2$. By the definition of the shortfall (Definition~\ref{def:sho_sem_the}), $d(s,w)=\frac\ell2\bigl(1-h(I_l)\bigr)$, $d(w,t)=\frac\ell2\bigl(1-h(I_r)\bigr)$ and $d(s,t)=\ell\bigl(1-h(I)\bigr)$. By the definition of the Gromov product (Definition~\ref{def:gro_hyp_sem}),
\[
2(s|t)_w=d(w,s)+d(w,t)-d(s,t)=\ell\,h(I)-\frac\ell2\bigl(h(I_l)+h(I_r)\bigr),
\]
and dividing by $\ell$ gives the identity for $J(I)$; the triangle inequality $d(s,t)\le d(s,w)+d(w,t)$ gives $J(I)\ge0$. By the identity for $J(I)$,
\[
h(I_l)-h(I)=\tfrac12\bigl(h(I_l)-h(I_r)\bigr)-J(I),\qquad h(I_r)-h(I)=-\tfrac12\bigl(h(I_l)-h(I_r)\bigr)-J(I),
\]
and $\lvert h(I_l)-h(I_r)\rvert=\sqrt{V(I)}$, which gives both bounds. The speed of an interval is one minus its shortfall.

(ii) Since $(\mathcal{T}_rd)(s,t)=d(s+r,t+r)$ and $I+r$ has the length of $I$, $h(I)(\mathcal{T}_rd)=h(I+r)(d)$. Since $(\mathcal{D}d)(s,t)=\frac12d(2s,2t)$, $h(I)(\mathcal{D}d)=1-d(2s,2t)/(2t-2s)=h(2I)(d)$. The halves of $I+r$ and of $2I$ are $I_l+r$, $I_r+r$ and $2I_l$, $2I_r$, so the identities for $J$ follow by (i) and those for $V$ by its definition. The forms with $h_t$ are the case $I=[x,x+t]$. The two identities for the uniform distance follow by substituting the definitions of $\mathcal{T}_r$ and $\mathcal{D}$ in the supremum. Finally $x\mapsto(\mathcal{T}_rd)(x,x+t)$ is the function $x\mapsto d(x,x+t)$ composed with the translation $x\mapsto x+r$ of $\mathbb{R}$, so one is Borel with values in $[0,t]$ exactly when the other is.

(iii) $h(I)(d)-h(I)(d'')=\bigl(d''(s,t)-d(s,t)\bigr)/(t-s)$. For the imbalance, let $a,b$ be the shortfalls of $I_l,I_r$ for $d$ and $a'',b''$ those for $d''$. They lie in $[0,1]$, so $\lvert a-b\rvert,\lvert a''-b''\rvert\le1$ and
\[
\bigl\lvert(a-b)^2-(a''-b'')^2\bigr\rvert=\bigl\lvert(a-b)-(a''-b'')\bigr\rvert\,\bigl\lvert(a-b)+(a''-b'')\bigr\rvert\le2\bigl(\lvert a-a''\rvert+\lvert b-b''\rvert\bigr).
\]
The halves lie in $I$ and have length $(t-s)/2$, so the first assertion gives $\lvert a-a''\rvert,\lvert b-b''\rvert\le2\lVert d-d''\rVert_I/(t-s)$.

(iv) By the triangle inequality, $d(x,y)\le d(x,x')+d(x',y')+d(y',y)\le\lvert x-x'\rvert+d(x',y')+\lvert y-y'\rvert$, and symmetrically. As $0\le d(s,t)\le t-s$, every shortfall lies in $[0,1]$. For $t>0$, $\lvert d(x,x+t)-d(x',x'+t)\rvert\le2\lvert x-x'\rvert$, so $x\mapsto d(x,x+t)$ and $x\mapsto h_t(x;d)$ are continuous; in particular $x\mapsto d(x,x+t)$ is Borel with values in $[0,t]$, and $d$ is admissible.
\end{proof}

\begin{lemma}
 For every integer $n\ge1$ and every orbit:
\begin{enumerate}
\item[(i)] $g_n$ is an array, $g_n(s,t)\le\lvert s-t\rvert+1/n$ for all $s,t\in\mathbb{R}$, and $g_n(s,t)=g(ns,nt)/n$ when $ns$ and $nt$ are integers; in particular $h([s,t])(g_n)=h_g(ns,nt)\in[0,1]$ for all such $s<t$.
\item[(ii)] $\mathcal{D}g_n=g_{2n}$.
\item[(iii)] For every $j\in\mathbb{Z}$, the array $\mathcal{T}_{j/n}g_n$ of the orbit of $z$ is the array $g_n$ of the orbit of $f_k^j(z)$. Consequently $\mathbb{E}\,\Psi(\mathcal{T}_{j/n}g_n)=\mathbb{E}\,\Psi(g_n)$ for every real function $\Psi$ of arrays such that $z\mapsto\Psi(g_n)$ is bounded and Borel.
\item[(iv)] For $x\in\mathbb{R}$ and $s,t\in\mathbb{R}$ with $ns$ and $nt$ integers, $(\mathcal{T}_xg_n)(s,t)=(\mathcal{T}_{\lfloor nx\rfloor/n}g_n)(s,t)$.
\item[(v)] For every $r>0$, the restriction of $g_n$ to $[-r,r]^2$ depends continuously on $z$ for the uniform distance on $[-r,r]^2$.
\item[(vi)] If $n$ is even, then $g_n$, $\mathcal{T}_xg_n$ for every $x\in\mathbb{R}$, and $\mathcal{D}g_n=g_{2n}$ are admissible.
\end{enumerate}
\label{lem:dil_and_tim}
\end{lemma}

\begin{proof}
Each item is a floor identity or a property of $g$.

(i) By Lemma~\ref{lem:nor_and_gro}(ii), $g$ is a semimetric on $\mathbb{Z}$ with $0\le g(a,b)\le\lvert a-b\rvert$, so $g_n$ is a semimetric on $\mathbb{R}$, and $\lvert\lfloor ns\rfloor-\lfloor nt\rfloor\rvert\le n\lvert s-t\rvert+1$ gives the bound. The identity at integer $ns,nt$ is the definition, and then $h([s,t])(g_n)=1-g(ns,nt)/(nt-ns)=h_g(ns,nt)$, which lies in $[0,1]$ by Lemma~\ref{lem:nor_and_gro}(ii).

(ii) $(\mathcal{D}g_n)(s,t)=\frac12g_n(2s,2t)=g(\lfloor 2ns\rfloor,\lfloor 2nt\rfloor)/(2n)=g_{2n}(s,t)$.

(iii) $(\mathcal{T}_{j/n}g_n)(s,t)=g(\lfloor ns\rfloor+j,\lfloor nt\rfloor+j)/n$, since $\lfloor n(s+j/n)\rfloor=\lfloor ns\rfloor+j$, and $g(a+j,b+j)=g(a,b)\circ f_k^j$ by Lemma~\ref{lem:inv_are_und}(iii). Hence $\Psi(\mathcal{T}_{j/n}g_n)=\Psi(g_n)\circ f_k^j$, and Lemma~\ref{lem:inv_are_und}(i) gives the equality of integrals.

(iv) If $ns$ is an integer, $\lfloor ns+nx\rfloor=ns+\lfloor nx\rfloor$; the same holds for $t$.

(v) On $[-r,r]^2$ the function $g_n$ is constant on each of the finitely many cells $\{(s,t):\lfloor ns\rfloor=a,\ \lfloor nt\rfloor=b\}$, with value $g(a,b)/n$, and the cells do not depend on $z$. Each $g(a,b)=\log_M\lVert T_{a,b}\rVert$ is a continuous function of $z$, since $T_{a,b}$ is a product of the matrices $A_i^{\pm1}$ of \eqref{eq:linearized_matrix}, whose entries $v(q_i)$ are continuous in $z$ because $f_k$ is a homeomorphism (Lemma~\ref{lem:inv_are_und}(i)). The uniform distance on $[-r,r]^2$ between the restrictions for two initial points is the maximum over the finitely many pairs $(a,b)$ of the differences of $g(a,b)/n$, which proves (v).

(vi) Let $t\in\{\frac12,1,2\}$, so that $nt$ is an integer. For every $x$, with $j=\lfloor nx\rfloor$, one has $g_n(x,x+t)=g(j,j+nt)/n$, which lies in $[0,t]$ by Lemma~\ref{lem:nor_and_gro}(ii); as a function of $x$ it is constant on each interval $[j/n,(j+1)/n)$, hence Borel. So $g_n$ is admissible, its translates are admissible by Lemma~\ref{lem:shor_arra}(ii), and $g_{2n}=\mathcal{D}g_n$ is admissible since $2n$ is even.
\end{proof}

\begin{lemma}
 Let $n_0<N$ be of the form $2^j$, and let $I$ be an interval whose length $p$ is of the form $2^j$ with $j\in\mathbb{Z}$, possibly negative, and which is aligned at scale $n_0$. Then, summing over $n=2^j$ with $n_0\le n<N$,
\[
\sum_n\mathbb{E}\,J(I)(g_n)\le e_{pN/2},\qquad\sum_n\mathbb{E}\Bigl[\mathbf{1}_{\{\max(h(I_l)(g_n),\,h(I_r)(g_n))\le\zeta\}}\,V(I)(g_n)\Bigr]\le8\zeta\,e_{pN/2}\quad(0<\zeta\le1).
\]
In particular $\sum_n\mathbb{E}\,(J(I)+V(I))(g_n)\le9e_{pN/2}$, and $e_{pN/2}\le e_N$ when $p\le2$.
\label{lem:sum_mid_can}
\end{lemma}

\begin{proof}
The midpoint cancellations telescope exactly across the scales $n=2^j$ by the exact doubling identity, and the imbalances are dominated by a telescoping convex quantity plus the midpoint cancellations.

We first reduce to the interval $[0,p]$. Write $I=[s,t]$, $\ell_0=pn_0/2$ and $\ell_1=pN/2$. Alignment at scale $n_0$ makes $n_0s$ and $\ell_0=n_0\bigl(\frac{s+t}2-s\bigr)$ integers, so for every $n=2^j\ge n_0$ the numbers $ns$ and $pn/2$ are integers. Fix such an $n<N$. By Lemma~\ref{lem:dil_and_tim}(i),
\[
h(I)(g_n)=h_g(ns,ns+pn),\quad h(I_l)(g_n)=h_g(ns,ns+pn/2),\quad h(I_r)(g_n)=h_g(ns+pn/2,ns+pn),
\]
and these lie in $[0,1]$. By invariance of $m$ under time shifts of the orbit (Lemma~\ref{lem:inv_are_und}(i),(iii)), every bounded Borel function of these three shortfalls has the integral of the same function of $h_g(0,pn)$, $h_g(0,pn/2)$ and $h_g(pn/2,pn)$, and $h_g(pn/2,pn)$ has the integral of $h_g(0,pn/2)$.

We next sum the midpoint cancellations. By Lemma~\ref{lem:shor_arra}(i), $J(I)=h(I)-\frac12\bigl(h(I_l)+h(I_r)\bigr)$. Hence $\mathbb{E}\,J(I)(g_n)=e_{pn}-e_{pn/2}$, and the sum over $n=2^j$ with $n_0\le n<N$ telescopes:
\[
\sum_n\mathbb{E}\,J(I)(g_n)=e_{\ell_1}-e_{\ell_0}\le e_{\ell_1}.
\]

We then bound the imbalances. Let $0<\zeta\le1$ and use the convex function $\Psi_\zeta(u)=u^2$ for $u\le\zeta$, $\Psi_\zeta(u)=2\zeta u-\zeta^2$ for $u\ge\zeta$. On $[0,1]$ it has Lipschitz constant $2\zeta$ and satisfies $0\le\Psi_\zeta(u)\le2\zeta u$. For the array $g_n$ put $\bar h=\frac12(h(I_l)+h(I_r))$. By Lemma~\ref{lem:shor_arra}(i), $h(I)-\bar h=J(I)\ge0$; as the three shortfalls lie in $[0,1]$, $0\le\bar h\le h(I)\le1$. Convexity, and the exact quadratic formula when $h(I_l),h(I_r)\le\zeta$, give
\[
\begin{aligned}
\tfrac14\mathbf{1}_{\{\max(h(I_l),h(I_r))\le\zeta\}}V(I)&\le\frac{\Psi_\zeta(h(I_l))+\Psi_\zeta(h(I_r))}2-\Psi_\zeta(\bar h)\\
&\le\frac{\Psi_\zeta(h(I_l))+\Psi_\zeta(h(I_r))}2-\Psi_\zeta(h(I))+2\zeta J(I).
\end{aligned}
\]
With $e^{\Psi_\zeta}_\ell=\mathbb{E}\,\Psi_\zeta(h_g(0,\ell))$ (Definition~\ref{def:gro_sem_orb}), one has $0\le e^{\Psi_\zeta}_\ell\le2\zeta e_\ell$, and the integral of the right side at scale $n$ is $e^{\Psi_\zeta}_{pn/2}-e^{\Psi_\zeta}_{pn}+2\zeta(e_{pn}-e_{pn/2})$. Summing over $n$,
\[
\tfrac14\sum_n\mathbb{E}\Bigl[\mathbf{1}_{\{\max(h(I_l),h(I_r))\le\zeta\}}V(I)\Bigr](g_n)\le e^{\Psi_\zeta}_{\ell_0}-e^{\Psi_\zeta}_{\ell_1}+2\zeta\bigl(e_{\ell_1}-e_{\ell_0}\bigr)\le2\zeta e_{\ell_1},
\]
which is the second bound.

Finally, at $\zeta=1$ the indicator is identically one, since the shortfalls lie in $[0,1]$; adding the two bounds gives $\sum_n\mathbb{E}\,(J(I)+V(I))(g_n)\le9e_{\ell_1}$. If $p\le2$, then $\ell_1=pN/2$ is of the form $2^j$ and at most $N$, so $e_{\ell_1}\le e_N$ by Lemma~\ref{lem:exa_dou_ide}.
\end{proof}

\subsection{The cap, the corrector and the proof of the balance} \label{subsec:sla_pro_the}

This subsection defines the objects of the balance argument and proves the balance across scales. The cap and the balance functional come first, then the slack profiles, the corrector on slack profiles, its extension to all arrays and the balance defect. The substance lies in these definitions; the proof of the balance is a short telescoping of the potential $e^\Phi_n-\mathbb{E}\,\widetilde U(g_n)$ over the scales.

\begin{definition}
 Let $\chi_0(x)=10x^3-15x^4+6x^5$, and let $\chi(u)=0$ for $u\le u_{\rm lo}$, $\chi(u)=\chi_0\bigl((u-u_{\rm lo})/(b_{\rm sl}/2-u_{\rm lo})\bigr)$ for $u_{\rm lo}<u<b_{\rm sl}/2$, and $\chi(u)=1$ for $u\ge b_{\rm sl}/2$. The \textbf{cap} is
\[
\Phi(u)=\alpha_{\rm cap}u+(1-\alpha_{\rm cap}u)\,\chi(u)\qquad(u\in\mathbb{R}).
\]
For an array $d$, $x\in\mathbb{R}$ and $t>0$ put $F_t(x;d)=\Phi\bigl(h_t(x;d)\bigr)$, and, when $x\mapsto F_1(x;d)$ is Lebesgue integrable on $[0,1/2]$,
\[
F_{\rm av}(d)=2\int_0^{1/2}F_1(x;d)\,dx.
\]
The \textbf{balance functional} is
\[
G_\Phi(d)=F_2(0;d)-\tfrac12F_1(0;d)-\tfrac12F_1(1;d),
\]
and, when $x\mapsto G_\Phi(\mathcal{T}_xd)$ is Lebesgue integrable on $[0,1]$, its \textbf{translation average} is
\[
G_{\rm av}(d)=\int_0^1G_\Phi(\mathcal{T}_xd)\,dx.
\]
The \textbf{capped mean shortfall} at length $n\ge1$ is $e^\Phi_n=\mathbb{E}\,\Phi\bigl(h_g(0,n)\bigr)$, the number $e^\Psi_n$ of Definition~\ref{def:gro_sem_orb} for $\Psi=\Phi$.
\label{def:the_cap_and}
\end{definition}

\begin{lemma}
 The cap $\Phi$ is $C^2$ on $\mathbb{R}$, equal to $\alpha_{\rm cap}u$ on $(-\infty,u_{\rm lo}]$ and equal to $1$ on $[b_{\rm sl}/2,\infty)$. On $[0,1]$ it satisfies
\[
\alpha_{\rm cap}u\le\Phi(u)\le\min\Bigl(1,\frac u{u_{\rm lo}}\Bigr),\qquad\lvert\Phi'(u)\rvert\le2^{20},\qquad\lvert\Phi''(u)\rvert\le2^{40}.
\]
\label{lem:pro_the_cap}
\end{lemma}

\begin{proof}
On $[0,1]$ the step $\chi_0$ has derivative $30x^2(1-x)^2$, between zero and two, and second derivative $60x-180x^2+120x^3$, of absolute value at most $360$; both vanish at $x=0$ and $x=1$, and $\chi_0(0)=0$, $\chi_0(1)=1$. Hence $\chi$ is $C^2$, constant outside $(u_{\rm lo},b_{\rm sl}/2)$, and $\Phi$ is $C^2$, equal to $\alpha_{\rm cap}u$ for $u\le u_{\rm lo}$ and to $1$ for $u\ge b_{\rm sl}/2$.

On $[0,1]$ one has $1-\alpha_{\rm cap}u>0$, and
\[
\Phi'=\alpha_{\rm cap}(1-\chi)+(1-\alpha_{\rm cap}u)\chi',\qquad\Phi''=-2\alpha_{\rm cap}\chi'+(1-\alpha_{\rm cap}u)\chi''.
\]
With $b_{\rm sl}/2-u_{\rm lo}=7b_{\rm sl}/16$, the bounds on $\chi_0'$ and $\chi_0''$ give $\lvert\Phi'\rvert\le\alpha_{\rm cap}+32/(7b_{\rm sl})\le2^{20}$ and $\lvert\Phi''\rvert\le64\alpha_{\rm cap}/(7b_{\rm sl})+360\,(16/(7b_{\rm sl}))^2\le2^{40}$. Moreover $\Phi-\alpha_{\rm cap}u=(1-\alpha_{\rm cap}u)\chi\ge0$ and $\Phi=(1-\chi)\alpha_{\rm cap}u+\chi\le1$. For $u\ge u_{\rm lo}$, $\Phi(u)\le1\le u/u_{\rm lo}$; below $u_{\rm lo}$, $\Phi(u)=\alpha_{\rm cap}u\le u/u_{\rm lo}$, as $\alpha_{\rm cap}<1/u_{\rm lo}$.

Explicit values: with the bounds $2$ for $\chi_0'$ and $360$ for $\lvert\chi_0''\rvert$ used above, the two derivative bounds are below $2^{13.9}$ and $2^{34.2}$.
\end{proof}

\begin{definition}
 The \textbf{slack profile family} $\mathscr{S}$ consists of the arrays $d$ with $d(x,y)\le\lvert x-y\rvert$ for all $x,y\in\mathbb{R}$ such that $(\mathbb{R},d)$ is $0$-hyperbolic, and which are of one of the following four kinds:
\begin{enumerate}
\item[(i)] \textbf{affine}: $d=d_\sigma$ for some $\sigma\in[0,1]$;
\item[(ii)] \textbf{slow}: $\mathbb{R}$ is slow at speed $1-b_{\rm sl}$;
\item[(iii)] \textbf{one-ray}: $\mathbb{R}$ is the union of a unit geodesic ray and the complementary half-line $I_*$, which is slow at speed $1-b_{\rm sl}$;
\item[(iv)] \textbf{two-ray}: for some $a_-<a_+$, the half-lines $(-\infty,a_-]$ and $[a_+,\infty)$ are unit geodesic rays, the interval $I_*=[a_-,a_+]$, of length $L=a_+-a_-$, is slow at speed $1-b_{\rm sl}$, and the deficit of $d$ on every interval is at most $C_{\rm cr}L$.
\end{enumerate}
The elements of $\mathscr{S}$ are the \textbf{slack profiles}. In kinds (iii) and (iv), $I_*$ is the \textbf{slow set}, or \textbf{core}, of the profile, and the finite endpoints of $I_*$ are its \textbf{interfaces}: the one endpoint of the half-line $I_*$ in kind (iii), and the two points $a_-$ and $a_+$ in kind (iv). In kind (iv), $L$ is the \textbf{slow length} of the profile, and its \textbf{crossing deficit} is the deficit of $d$ on the interval $[a_--C_{\rm cr}L,\,a_++C_{\rm cr}L]$,
\[
\Delta=(2C_{\rm cr}+1)L-d\bigl(a_--C_{\rm cr}L,\ a_++C_{\rm cr}L\bigr).
\]
A slack profile is \textbf{exceptional} if it is an affine array $d_\sigma$ with $1-b_{\rm sl}<\sigma<1$. The slow set of a slow profile is $I_*=\mathbb{R}$, and the slow set of $d_1$ is $I_*=\varnothing$; these profiles have no interfaces. For a slack profile $d$ that is not exceptional, the \textbf{slow fraction} is
\[
\mu_t(x;d)=\frac{\lvert[x,x+t]\cap I_*\rvert}{t}\qquad(x\in\mathbb{R},\ t>0).
\]
\label{def:slac_prof}
\end{definition}

\begin{lemma}

\begin{enumerate}
\item[(i)] Let $d$ and $d''$ be arrays and $\Theta$ an interval such that every secant of $d$ with endpoints in $\Theta$ has speed at most $1-b_{\rm sl}$ and every secant of $d''$ with endpoints in $\Theta$ has speed one. Then $b_{\rm sl}\lvert\Theta\rvert\le\lVert d-d''\rVert_\Theta$.
\item[(ii)] The four kinds of slack profiles are mutually exclusive, except that an affine array of speed at most $1-b_{\rm sl}$ is also slow. The slow set of a slack profile $d$ that is not exceptional is determined by $d$, and so are its interfaces and, if $d$ is two-ray, $a_\pm$, $L$ and $\Delta$. For such a profile the closure of $I_*$ is slow at speed $1-b_{\rm sl}$; every secant with endpoints in the closure of one component of $\mathbb{R}\setminus I_*$ has speed one; and each interface $\xi$ is the endpoint of a unit geodesic ray of $d$, $(-\infty,\xi]$ or $[\xi,\infty)$, whose interior does not meet $I_*$.
\item[(iii)] Let $d\in\mathscr{S}$ and $r\in\mathbb{R}$. Then $\mathcal{T}_rd$ and $\mathcal{D}d$ are slack profiles of the same kinds as $d$, and exceptional exactly when $d$ is. If $d$ is not exceptional, their slow sets are $I_*-r$ and $\frac12I_*$, and $\mu_t(x;\mathcal{T}_rd)=\mu_t(x+r;d)$ and $\mu_t(x;\mathcal{D}d)=\mu_{2t}(2x;d)$ for $x\in\mathbb{R}$ and $t>0$. If $d$ is two-ray, its slow length and crossing deficit are unchanged by $\mathcal{T}_r$ and halved by $\mathcal{D}$.
\item[(iv)] Let $d\in\mathscr{S}$ be not exceptional and $t>0$. Then $x\mapsto\mu_t(x;d)$ is continuous, and for every $x\in\mathbb{R}$
\[
\mu_{2t}(x;d)=\tfrac12\bigl(\mu_t(x;d)+\mu_t(x+t;d)\bigr),\qquad b_{\rm sl}\,\mu_t(x;d)\le h_t(x;d)\le C_{\rm cr}\,\mu_t(x;d).
\]
If $d$ is two-ray, then $b_{\rm sl}L\le\Delta\le C_{\rm cr}L$.
\end{enumerate}
\label{lem:geo_sla_pro}
\end{lemma}

\begin{proof}
A slow interval of one array cannot lie where another array has unit speed unless the two arrays are far apart; with one array this pins down the kind and the slow set of a profile, and the triangle inequality across the slow set compares shortfall and slow fraction.

(i) For $x<y$ in $\Theta$ one has $d''(x,y)=y-x$ and $d(x,y)\le(1-b_{\rm sl})(y-x)$, so $b_{\rm sl}(y-x)\le d''(x,y)-d(x,y)\le\lVert d-d''\rVert_\Theta$; the supremum over $x<y$ in $\Theta$ gives the claim.

(ii) We first show that a slack profile determines its kind, up to the stated overlap. By (i) with $d''=d$, no nondegenerate interval is both slow at speed $1-b_{\rm sl}$ and contained in a unit geodesic ray. The affine array $d_\sigma$ is slow exactly when $\sigma\le1-b_{\rm sl}$; every secant of $d_\sigma$ has speed $\sigma$, so it has a unit geodesic ray only if $\sigma=1$, and then no nondegenerate interval is slow; so it is neither one-ray nor two-ray. Compare now two descriptions of one profile, each slow, one-ray or two-ray. If one is slow and the other is not, a unit geodesic ray of the second is slow. Otherwise some nondegenerate interval lies in a unit geodesic ray of one description and in the slow set of the other: a far half-line if the descriptions are of different kinds or have their rays on opposite sides; $[y,(y+y')/2]$ for one-ray descriptions with rays $[y,\infty)$ and $[y',\infty)$, $y<y'$, and likewise for rays on the left; $[a_-,\min(a_-',a_+)]$ for two-ray descriptions with $a_-<a_-'$, and likewise at $a_+$. Each case is excluded. Hence the only overlap of kinds is the one stated, and the unit geodesic rays of a one-ray or two-ray profile, hence its slow set, its interfaces, $a_\pm$, $L$ and $\Delta$, are determined by $d$. A slow profile, the array $d_1$ and the exceptional profiles are recognised by their kinds, so the slow set of every profile that is not exceptional is determined.

We next prove the three properties. Since $d$ is continuous (Lemma~\ref{lem:shor_arra}(iv)), the closure of a set slow at speed $1-b_{\rm sl}$ is slow at that speed. For a one-ray profile $\mathbb{R}\setminus I_*$ is its unit geodesic ray; for a two-ray profile its components are $(-\infty,a_-)$ and $(a_+,\infty)$, whose closures are the two rays; for $d_1$ it is $\mathbb{R}$, on which every secant has speed one; a slow profile has $I_*=\mathbb{R}$. The interfaces are the endpoints of these rays, whose interiors do not meet $I_*$.

(iii) A translation by $r$ preserves speeds, deficits, Gromov products (with the basepoint moved by $-r$) and the bound $d(x,y)\le\lvert x-y\rvert$, and it moves rays and slow sets by $-r$. Under the dilation, the secant $[x,y]$ of $\mathcal{D}d$ has the speed of the secant $[2x,2y]$ of $d$ and half its deficit, and the Gromov products of $\mathcal{D}d$ at $w$ are half those of $d$ at $2w$; so $\mathcal{D}$ preserves $0$-hyperbolicity and the bound $d(x,y)\le\lvert x-y\rvert$, halves rays, slow sets and slow lengths, and preserves the deficit bound of kind (iv) with the halved slow length. The crossing deficit of $\mathcal{D}d$ is the deficit of $\mathcal{D}d$ on $\frac12[a_--C_{\rm cr}L,a_++C_{\rm cr}L]$, that is, half the crossing deficit of $d$. Affine arrays are fixed by $\mathcal{T}_r$ and $\mathcal{D}$. The identities for $\mu_t$ follow by substituting $I_*-r$ and $\frac12I_*$ in Definition~\ref{def:slac_prof}.

(iv) Moving $x$ changes the measure of $[x,x+t]\cap I_*$ by at most twice the displacement, so $\mu_t$ is continuous; the additivity holds because $[x,x+2t]$ is the union of $[x,x+t]$ and $[x+t,x+2t]$, which overlap in one point.

We prove the lower bound. If $[x,x+t]$ meets the closure $\bar I_*$ of $I_*$, the intersection is an interval $[p,p']$, of length $t\mu_t(x)$ because $\bar I_*\setminus I_*$ has at most two points. By the triangle inequality, $d(x,y)\le\lvert x-y\rvert$ and (ii),
\[
d(x,x+t)\le(p-x)+(1-b_{\rm sl})(p'-p)+(x+t-p')=t-b_{\rm sl}\,t\mu_t(x),
\]
that is, $h_t(x)\ge b_{\rm sl}\mu_t(x)$. If $[x,x+t]$ misses $\bar I_*$, then $\mu_t(x)=0\le h_t(x)$.

We prove the upper bound. Every shortfall of $d$ lies in $[0,1]$ (Lemma~\ref{lem:shor_arra}(iv)), and $C_{\rm cr}\ge2$. For $d_1$ both sides vanish, and for a slow profile $\mu_t=1$. Let $d$ be one-ray or two-ray. If $d$ is two-ray and $[x,x+t]\supset[a_-,a_+]$, the deficit of $[x,x+t]$ is at most $C_{\rm cr}L=C_{\rm cr}t\mu_t(x)$ by Definition~\ref{def:slac_prof}(iv). If $(x,x+t)$ contains no interface, it is a connected set that misses the boundary of $\bar I_*$, so it lies in the interior of $\bar I_*$ or in $\mathbb{R}\setminus\bar I_*$; then $[x,x+t]\subset\bar I_*$ and $\mu_t(x)=1$, or $[x,x+t]$ lies in the closure of one component of $\mathbb{R}\setminus I_*$ and $h_t(x)=0$ by (ii). In the remaining case $(x,x+t)$ contains exactly one interface $\xi$, the endpoint of a ray of (ii); say the ray is $[\xi,\infty)$, the other case being symmetric. Then $[x,\xi]\subset\bar I_*$ and $[\xi,x+t]$ lies in the ray, so $\ell=\xi-x=t\mu_t(x)$, and
\[
d(x,x+t)\ge d(\xi,x+t)-d(x,\xi)\ge(x+t-\xi)-\ell=t-2\ell,
\]
that is, $h_t(x)\le2\mu_t(x)$. In every case $h_t(x)\le C_{\rm cr}\mu_t(x)$.

Finally let $d$ be two-ray. The crossing deficit is the deficit of $d$ on an interval, so $\Delta\le C_{\rm cr}L$ by Definition~\ref{def:slac_prof}(iv); and the lower bound of this item for the interval $[a_--C_{\rm cr}L,a_++C_{\rm cr}L]$, which meets $I_*$ in $[a_-,a_+]$, gives $\Delta\ge b_{\rm sl}L$.
\end{proof}

\begin{definition}
 Let $\psi(y)=\min\bigl(1,\max(0,\,y/A_{\rm sl}-1)\bigr)$ for $y\in\mathbb{R}$. Let $d\in\mathscr{S}$ and $t>0$. The \textbf{slow weight} of $d$ is
\[
\nu_t=\nu_t(d)=1-\psi(t/L)\Bigl(1-\frac{\alpha_{\rm cap}\Delta}L\Bigr)
\]
if $d$ is two-ray, and $\nu_t=1$ otherwise. The \textbf{profile corrector} of $d$ is
\[
U_t(x;d)=F_t(x;d)-\nu_t\,\mu_t(x;d)\qquad(x\in\mathbb{R})
\]
if $d$ is not exceptional, and $U_t(x;d)=0$ if $d$ is exceptional. When $x\mapsto U_1(x;d)$ is Lebesgue integrable on $[0,1/2]$, the \textbf{averaged corrector} is
\[
U_{\rm av}(d)=2\int_0^{1/2}U_1(x;d)\,dx.
\]
\label{def:the_pro_cor}
\end{definition}

\begin{definition}
 For an array $d$, the \textbf{imbalance budget} is
\[
\mathcal{V}(d)=\frac1{64}\sum_{j=-256}^{320}V\bigl([j/64,\,j/64+1]\bigr)(d),
\]
the \textbf{envelope} is
\[
U_{\rm env}(d)=\inf\Bigl\{U_{\rm av}(d')+C_{\rm env}\lVert d-d'\rVert_{[-8,8]}:\ d'\in\mathscr{S},\ U_{\rm av}(d')\text{ defined}\Bigr\}\in[-\infty,\infty],
\]
and the \textbf{side factor} is
\[
\omega(d)=\max\Bigl(0,\ 1-\frac4{b_{\rm sl}}\min\bigl(h([-2,-1])(d),\,h([3,4])(d)\bigr)\Bigr).
\]
For an admissible array $d$ whose average $F_{\rm av}(d)$ is defined, the \textbf{extended corrector} is
\[
\widetilde U(d)=\omega(d)\,\max\Bigl(-C_{\rm dom}\mathcal{V}(d),\ \min\bigl(C_{\rm dom}\mathcal{V}(d),\ (1-c_{\rm sav})F_{\rm av}(d),\ U_{\rm env}(d)\bigr)\Bigr),
\]
and, when $\widetilde U(\mathcal{T}_{1/2}d)$ is also defined, its \textbf{half-shift average} is $\overline{\widetilde U}(d)=\frac12\bigl(\widetilde U(d)+\widetilde U(\mathcal{T}_{1/2}d)\bigr)$. When $G_{\rm av}(d)$ and the extended correctors below are defined, the \textbf{balance defect} is
\[
Z(d)=G_{\rm av}(d)-\widetilde U(\mathcal{D}d)+\tfrac12\overline{\widetilde U}(d)+\tfrac12\overline{\widetilde U}(\mathcal{T}_1d).
\]
\label{def:ext_cor_and}
\end{definition}

The following paragraphs explain the objects of Definition~\ref{def:the_cap_and}, Definition~\ref{def:the_pro_cor} and Definition~\ref{def:ext_cor_and}; every claim in them is proved in the results cited.

\emph{The splitting.} On a slack profile that is not exceptional, the capped shortfall of a window $[x,x+t]$ splits as
\[
F_t(x)=\nu_t\,\mu_t(x)+U_t(x),
\]
the slow fraction weighted by the slow weight, plus the corrector. The slow fraction is additive under doubling, $\mu_{2t}(x)=\frac12\bigl(\mu_t(x)+\mu_t(x+t)\bigr)$ (Lemma~\ref{lem:geo_sla_pro}(iv)), and the weight is nonincreasing in $t$, because $\alpha_{\rm cap}\Delta\le\alpha_{\rm cap}C_{\rm cr}L=L/20$ (Lemma~\ref{lem:geo_sla_pro}(iv)): a core swallowed by a long window counts less than the slow time inside it. So the weighted slow fraction contributes $(\nu_{2t}-\nu_t)\mu_{2t}(x)\le0$ to the balance: this is the correction inequality of Lemma~\ref{lem:cor_ine_sla}(ii). Since $\nu_t\ge\alpha_{\rm cap}b_{\rm sl}$ and the shortfall of a window is at most $C_{\rm cr}$ times its slow fraction (Lemma~\ref{lem:geo_sla_pro}(iv)), the weighted slow fraction carries the fixed fraction $c_{\rm sav}$ of the capped shortfall; this is the terminal saving $U_t\le(1-c_{\rm sav})F_t$ of Lemma~\ref{lem:cor_ine_sla}(ii).

\emph{The two regimes of the weight.} The weight is chosen so that the corrector vanishes on most windows. A window inside the slow set has shortfall at least $b_{\rm sl}$, beyond the saturation of the cap at $b_{\rm sl}/2$, so its capped shortfall is $1=\mu_t$; it is matched by $\nu_t=1$ at the window lengths $t\le A_{\rm sl}L$. A window of length $t\ge2A_{\rm sl}L$ that crosses the core of a two-ray profile with room $C_{\rm cr}L$ on each side has deficit exactly $\Delta$ and shortfall $\Delta/t\le C_{\rm cr}/(2A_{\rm sl})<u_{\rm lo}$, in the linear range of the cap, so its capped shortfall is $\alpha_{\rm cap}\Delta/t$; it is matched by the weight $\nu_t=\alpha_{\rm cap}\Delta/L$ at these lengths. Between these lengths $\psi$ moves the weight over one doubling of $t/L$. Hence the corrector of a unit window is nonzero only if the window contains an interface and, for a short core, lacks the room $C_{\rm cr}L$ on one side (Lemma~\ref{lem:van_the_pro}): the corrector measures the capped shortfall that the weighted slow fraction does not account for, and it lives near the interfaces. On a slack profile $d$ that is not exceptional, the balance defect is the average of the balance of the weighted slow fraction,
\[
Z(d)=(\nu_2-\nu_1)\int_0^1\mu_2(x;d)\,dx\le0,
\]
and it vanishes on exceptional profiles (Lemma~\ref{lem:pro_the_bal}(ii)).

\emph{The ingredients of the extended corrector.} The defect must be small at the arrays $g_n$, summed over an unbounded number of scales (Proposition~\ref{prop:def_the_ext:1}). Over so many scales the quantities with controlled sums are the midpoint cancellations and the imbalances (Lemma~\ref{lem:sum_mid_can}), so an error charged at one scale must vanish, be bounded by them, be carried by a set of summable mass, or be summable over $n=2^j$, like a multiple of $1/n$. Each ingredient of $\widetilde U$ serves one part of this estimate.
\begin{itemize}
\item The average over the position of the window. Next to a core of length $L$ the corrector of a single window is bounded only by a multiple of $L$, while the imbalances of the profile are of order $L^2$; but only windows within distance of order $L$ of an interface carry a nonzero corrector, so the averaged corrector is of order $L^2$ and is dominated by the imbalance budget (Lemma~\ref{lem:spa_dom_the}).
\item The envelope with constant $C_{\rm env}$. It extends $U_{\rm av}$ from $\mathscr{S}$ without changing it there and keeps a Lipschitz bound (Lemma~\ref{lem:mod_the_pro}, Lemma~\ref{lem:ext_the_ave}); this treats the body (Lemma~\ref{lem:pro_the_bal}).
\item The clamp by $\pm C_{\rm dom}\mathcal{V}$. Off the body the arrays are close to affine arrays (Lemma~\ref{lem:aff_app_off}), but the complement of the body need not have small mass, so no error per scale can be afforded. At near-unit and middle speeds the correction terms are bounded by imbalances, whose sums over the scales are small (Lemma~\ref{lem:bal_def_off}, Lemma~\ref{lem:sum_mid_can}); the same bound controls the initial term of the balance (Lemma~\ref{lem:bou_the_ext}).
\item The side factor $\omega$. At slow speeds the cap makes the balance functional vanish, and $\omega$ makes the correction terms vanish: it is zero when both side shortfalls are at least $b_{\rm sl}/4$, and one wherever $U_{\rm av}\ne0$ on $\mathscr{S}$ (Lemma~\ref{lem:van_the_pro}).
\item The clamp by $(1-c_{\rm sav})F_{\rm av}$. It keeps the terminal saving at the last scale, where the array need not be close to any slack profile; an additive error there would spoil the balance, whose right-hand side $C_{\rm bal}(e_{n_0}+1/n_0)$ has no constant term.
\end{itemize}
None of these operations changes $U_{\rm av}$ on $\mathscr{S}$ (Lemma~\ref{lem:ext_the_ave}).

The envelope of Definition~\ref{def:ext_cor_and} is the McShane extension \cite{locbib:McShane1934} of the averaged corrector from the slack profiles. The next lemma records the bounds on the extended corrector that the proof of the balance uses at the first and at the last scale.

\begin{lemma}
 Let $d$ be an admissible array. Then $F_{\rm av}(d)$ and $\widetilde U(d)$ are defined, $0\le F_{\rm av}(d)\le1$, $\mathcal{V}(d)\ge0$, $U_{\rm env}(d)\ge-1$, and
\[
\lvert\widetilde U(d)\rvert\le1,\qquad\lvert\widetilde U(d)\rvert\le C_{\rm dom}\,\mathcal{V}(d),\qquad\widetilde U(d)\le(1-c_{\rm sav})\,F_{\rm av}(d).
\]
\label{lem:bou_the_ext}
\end{lemma}

\begin{proof}
The bounds hold by construction: the clamp in Definition~\ref{def:ext_cor_and} lies between
\[
\max(-C_{\rm dom}\mathcal{V},-1)\qquad\text{and}\qquad\min\bigl(C_{\rm dom}\mathcal{V},(1-c_{\rm sav})F_{\rm av}\bigr),
\]
because the envelope is at least $-1$, and the side factor lies in $[0,1]$.

Since $d$ is admissible, $x\mapsto d(x,x+1)$ is Borel with values in $[0,1]$; so the shortfalls $h_1(x;d)$ lie in $[0,1]$, $x\mapsto F_1(x;d)$ is Borel with values in $[0,1]$ by Lemma~\ref{lem:pro_the_cap}, and $F_{\rm av}(d)$ is defined and lies in $[0,1]$. Hence $\widetilde U(d)$ is defined. The imbalance budget is a sum of squares, so $\mathcal{V}(d)\ge0$. The side shortfalls $h([-2,-1])(d)$ and $h([3,4])(d)$ lie in $[0,1]$, so $0\le\omega(d)\le1$.

We next show that $U_{\rm env}(d)\ge-1$. Let $d'\in\mathscr{S}$. If $d'$ is exceptional, $U_1(x;d')=0$ for every $x$ (Definition~\ref{def:the_pro_cor}). Otherwise the shortfalls of $d'$ lie in $[0,1]$ (Lemma~\ref{lem:shor_arra}(iv)), so $F_1(x;d')\ge0$ by Lemma~\ref{lem:pro_the_cap}; moreover $0\le\mu_1(x;d')\le1$, and $0\le\nu_1(d')\le1$: this is clear unless $d'$ is two-ray, and then the crossing deficit, being a deficit of $d'$, lies between $0$ and $C_{\rm cr}L$ (Definition~\ref{def:slac_prof}(iv)), so that $0\le\alpha_{\rm cap}\Delta/L\le\alpha_{\rm cap}C_{\rm cr}=\frac1{20}$, while $0\le\psi\le1$. Thus $U_1(x;d')=F_1(x;d')-\nu_1\mu_1(x;d')\ge-1$ for every $x$, so $U_{\rm av}(d')\ge-1$ whenever it is defined, and $U_{\rm env}(d)\ge-1$.

Finally put
\[
y_-=-C_{\rm dom}\mathcal{V}(d),\quad y_+=\min\bigl(C_{\rm dom}\mathcal{V}(d),(1-c_{\rm sav})F_{\rm av}(d)\bigr),\quad y=\max\bigl(y_-,\min(y_+,U_{\rm env}(d))\bigr),
\]
so that $\widetilde U(d)=\omega(d)\,y$. Then $y_-\le0\le y_+\le1$ and $\min(y_+,U_{\rm env}(d))\ge-1$, hence $\max(y_-,-1)\le y\le y_+$. So $\lvert y\rvert\le\min\bigl(1,C_{\rm dom}\mathcal{V}(d)\bigr)$ and $y\le(1-c_{\rm sav})F_{\rm av}(d)$. As $0\le\omega(d)\le1$, $\lvert\widetilde U(d)\rvert\le\lvert y\rvert$, and $\widetilde U(d)\le\max(y,0)\le(1-c_{\rm sav})F_{\rm av}(d)$, the right side being nonnegative.
\end{proof}

We can now prove the balance across scales. Its proof cites two results that are proved later in this section and do not depend on it: the properties of the balance defect, Lemma~\ref{lem:pro_the_bal}, at the end of the next subsection, and the defect bound Proposition~\ref{prop:def_the_ext:1}, in the last subsection.

\begin{prop}
 Let $M\ge M_{\rm an}$, and let $n_0<N$ be of the form $2^j$ with $n_0\ge n_{\rm bal}$ and $e_N\le\zeta_{\rm bal}$. Then
\[
e_N\le C_{\rm bal}\Bigl(e_{n_0}+\frac1{n_0}\Bigr).
\]
\label{prop:bal_bet_dya}
\end{prop}

\begin{proof}
We telescope the potential $\mathcal{L}_n=e^\Phi_n-\mathbb{E}\,\widetilde U(g_n)$ over the scales $n=2^j$. Its increments are the mean balance defects; at the last scale the terminal saving makes it at least the fraction $c_{\rm sav}$ of $e^\Phi_N$; at the first scale it is at most a multiple of $e_{n_0}$; and the defect bound Proposition~\ref{prop:def_the_ext:1} makes the summed increments small.

Within this proof sums over $n$ run over $n=2^j$ with $n_0\le n<N$; $g_n$ is the rescaled growth array of Definition~\ref{def:res_gro_arr:1}; $e^\Phi_n$, $\Phi$, $F_{\rm av}$, $G_\Phi$ and $G_{\rm av}$ are as in Definition~\ref{def:the_cap_and}, and $\mathcal{V}$, $\widetilde U$, $\overline{\widetilde U}$ and $Z$ as in Definition~\ref{def:ext_cor_and}; $\mathcal{L}_n=e^\Phi_n-\mathbb{E}\,\widetilde U(g_n)$ is the potential. All scales here are even and at least $n_{\rm bal}$, so every interval used below is aligned at every scale that occurs. The arrays $g_n$, their translates and $\mathcal{D}g_n=g_{2n}$ are admissible (Lemma~\ref{lem:dil_and_tim}(vi)). The functions of the initial point integrated below are continuous by Lemma~\ref{lem:pro_the_bal}(iv) and Lemma~\ref{lem:dil_and_tim}(v), and bounded by Lemma~\ref{lem:bou_the_ext}, Lemma~\ref{lem:pro_the_bal}(i) and Lemma~\ref{lem:dil_and_tim}(i).

\textbf{The increments of the potential.} By Lemma~\ref{lem:pro_the_bal}(iv) and Lemma~\ref{lem:inv_are_und}(i), $\mathbb{E}\,G_{\rm av}(g_n)=\mathbb{E}\,G_\Phi(g_n)$. The three shortfalls in $G_\Phi(g_n)$ are those of $g$ on $[0,2n]$, $[0,n]$ and $[n,2n]$ (Lemma~\ref{lem:dil_and_tim}(i)), and by Lemma~\ref{lem:inv_are_und}(i),(iii) the last has the integral of the second; so $\mathbb{E}\,G_{\rm av}(g_n)=\mathbb{E}\,G_\Phi(g_n)=e^\Phi_{2n}-e^\Phi_n$. By Lemma~\ref{lem:dil_and_tim}(ii),(iii), $\widetilde U(\mathcal{D}g_n)=\widetilde U(g_{2n})$, and $\widetilde U(\mathcal{T}_rg_n)=\widetilde U(g_n)\circ f_k^{rn}$ for $r\in\{\frac12,1,\frac32\}$. Hence, by Lemma~\ref{lem:inv_are_und}(i), $\mathbb{E}\bigl[\widetilde U(\mathcal{D}g_n)-\frac12\overline{\widetilde U}(g_n)-\frac12\overline{\widetilde U}(\mathcal{T}_1g_n)\bigr]=\mathbb{E}\,\widetilde U(g_{2n})-\mathbb{E}\,\widetilde U(g_n)$, and by the definition of $Z$
\[
\mathcal{L}_{2n}-\mathcal{L}_n=\mathbb{E}\,Z(g_n)\quad(n_0\le n<N),\qquad\text{hence}\qquad\mathcal{L}_N=\mathcal{L}_{n_0}+\sum_n\mathbb{E}\,Z(g_n).
\]

\textbf{The last scale.} By Lemma~\ref{lem:pro_the_bal}(iv) and Lemma~\ref{lem:inv_are_und}(i),(iii), $\mathbb{E}\,F_{\rm av}(g_n)=e^\Phi_n$, so Lemma~\ref{lem:bou_the_ext} gives $\mathbb{E}\,\widetilde U(g_N)\le(1-c_{\rm sav})e^\Phi_N$, that is, $\mathcal{L}_N\ge c_{\rm sav}e^\Phi_N$.

\textbf{The first scale.} At $n_0$, $\lvert\widetilde U\rvert\le C_{\rm dom}\mathcal{V}$ (Lemma~\ref{lem:bou_the_ext}); each imbalance $V([x,x+1])(g_{n_0})$ in $\mathcal{V}(g_{n_0})$, where $x\in\frac1{64}\mathbb{Z}$, is at most $h([x,x+\frac12])+h([x+\frac12,x+1])$, whose integral is $2e_{n_0/2}\le2e_{n_0}$ (Lemma~\ref{lem:dil_and_tim}(i), Lemma~\ref{lem:inv_are_und}(i),(iii), Lemma~\ref{lem:exa_dou_ide}); the coefficients of $\mathcal{V}$ sum to $577/64<10$; and $\Phi(u)\le u/u_{\rm lo}$ on $[0,1]$ (Lemma~\ref{lem:pro_the_cap}). Hence
\[
\mathcal{L}_{n_0}\le\mathbb{E}\,\lvert\widetilde U(g_{n_0})\rvert+e^\Phi_{n_0}\le\Bigl(20C_{\rm dom}+\frac1{u_{\rm lo}}\Bigr)e_{n_0}\le2^{150}e_{n_0}.
\]

\textbf{Conclusion.} Since $\Phi(u)\ge\alpha_{\rm cap}u$ on $[0,1]$ (Lemma~\ref{lem:pro_the_cap}), $e^\Phi_N\ge\alpha_{\rm cap}e_N$. The three previous paragraphs give $c_{\rm sav}e^\Phi_N\le\mathcal{L}_N=\mathcal{L}_{n_0}+\sum_n\mathbb{E}\,Z(g_n)\le\sum_n\mathbb{E}\,Z(g_n)+2^{150}e_{n_0}$. By the values of Convention~\ref{conv:not_and_con}, $c_{\rm sav}\alpha_{\rm cap}=b_{\rm sl}^2/(6400C_{\rm cr}^3)$, with $C_{\rm cr}<2^{42}$, $b_{\rm sl}>2^{-12}$ and $6400<2^{13}$, so $c_{\rm sav}\alpha_{\rm cap}>2^{-163}>2^{-200}$. By Proposition~\ref{prop:def_the_ext:1} and $e_N\le\zeta_{\rm bal}$, with $2^{45770}\zeta_{\rm bal}=2^{-1230}$ by the values of Convention~\ref{conv:not_and_con},
\[
2^{-200}e_N\le2^{-780}e_N+2^{45770}e_N^2+\frac2{n_0}+2^{150}e_{n_0}\le2^{-779}e_N+\frac2{n_0}+2^{150}e_{n_0},
\]
so, since $2^{-200}-2^{-779}>2^{-201}$, $e_N\le2^{201}\bigl(2^{150}e_{n_0}+2/n_0\bigr)\le C_{\rm bal}(e_{n_0}+1/n_0)$.

Explicit values: $20C_{\rm dom}+1/u_{\rm lo}=2^{149.32}$; $c_{\rm sav}\alpha_{\rm cap}=2^{-159.31}$; $2^{201+150}=2^{351}$ and $2^{202}$ are at most $C_{\rm bal}=2^{420}$.
\end{proof}

The hypothesis $e_N\le\zeta_{\rm bal}$ is used only to absorb the quadratic term $2^{45770}e_N^2$ of Proposition~\ref{prop:def_the_ext:1}.

\subsection{The corrector on slack profiles and the properties of the defect} \label{subsec:cor_sla_pro}

This subsection proves what the balance and the defect bound need of the corrector on slack profiles, and ends with the properties of the balance defect that both proofs use (Lemma~\ref{lem:pro_the_bal}). The key computation is the correction inequality on slack profiles, with the terminal saving $U_t\le(1-c_{\rm sav})F_t$ (Lemma~\ref{lem:cor_ine_sla}); the remaining lemmas carry it over to the extended corrector.

\begin{lemma}

\begin{enumerate}
\item[(i)] For $d\in\mathscr{S}$, $x,r\in\mathbb{R}$ and $t>0$,
\[
U_t(x;\mathcal{T}_rd)=U_t(x+r;d),\qquad U_t(x;\mathcal{D}d)=U_{2t}(2x;d).
\]
\item[(ii)] For every $d\in\mathscr{S}$, $x\in\mathbb{R}$ and $t>0$, writing $F_t(x)=F_t(x;d)$ and $U_t(x)=U_t(x;d)$,
\[
F_{2t}(x)-\tfrac12F_t(x)-\tfrac12F_t(x+t)\le U_{2t}(x)-\tfrac12U_t(x)-\tfrac12U_t(x+t),
\]
\[
\lvert U_t(x)\rvert\le1,\qquad U_t(x)\le(1-c_{\rm sav})F_t(x).
\]
\end{enumerate}
\label{lem:cor_ine_sla}
\end{lemma}

\begin{proof}
On a slack profile that is not exceptional the capped shortfall is the weighted slow fraction plus the corrector, $F_t=\nu_t\mu_t+U_t$. The slow fraction is exactly additive under doubling and the weight is nonincreasing in $t$, so the weighted slow fraction has nonpositive balance; this is the correction inequality. The lower bound $\nu_t\ge\alpha_{\rm cap}b_{\rm sl}$ and the comparison of shortfall and slow fraction give the terminal saving.

(i) By Lemma~\ref{lem:geo_sla_pro}(iii), $\mathcal{T}_rd$ and $\mathcal{D}d$ are slack profiles, exceptional exactly when $d$ is; for exceptional $d$ both identities read $0=0$. Let $d$ be not exceptional. By Lemma~\ref{lem:shor_arra}(ii), $h_t(x;\mathcal{T}_rd)=h_t(x+r;d)$ and $h_t(x;\mathcal{D}d)=h_{2t}(2x;d)$, so $F_t$ obeys the same identities; so does $\mu_t$, by Lemma~\ref{lem:geo_sla_pro}(iii). If $d$ is two-ray, its slow length and crossing deficit are unchanged by $\mathcal{T}_r$ and halved by $\mathcal{D}$ (Lemma~\ref{lem:geo_sla_pro}(iii)), so $\nu_t(\mathcal{T}_rd)=\nu_t(d)$ and $\nu_t(\mathcal{D}d)=1-\psi(2t/L)(1-\alpha_{\rm cap}\Delta/L)=\nu_{2t}(d)$; otherwise all these weights equal one. The identities for $U_t$ follow.

(ii) Let first $d=d_\sigma$ be exceptional. Then $F_t(x)=\Phi(1-\sigma)$ does not depend on $x$ and $t$, and $U_t=0$; the correction inequality reads $0\le0$, and the other two bounds hold as $F_t\ge0$ (Lemma~\ref{lem:pro_the_cap}).

Let now $d$ be not exceptional. We first collect the properties of the weight. If $d$ is two-ray, $b_{\rm sl}\le\Delta/L\le C_{\rm cr}$ (Lemma~\ref{lem:geo_sla_pro}(iv)), so $1-\alpha_{\rm cap}\Delta/L\ge1-\alpha_{\rm cap}C_{\rm cr}=\frac{19}{20}$; as $\psi$ is nondecreasing with values in $[0,1]$, $\nu_t$ is nonincreasing in $t$ and $\alpha_{\rm cap}b_{\rm sl}\le\alpha_{\rm cap}\Delta/L\le\nu_t\le1$. Otherwise $\nu_t=1$. In all cases $\nu_t$ is nonincreasing in $t$ and $\alpha_{\rm cap}b_{\rm sl}\le\nu_t\le1$.

The correction inequality follows by subtraction. As $F_t-U_t=\nu_t\mu_t$ and $\mu_{2t}(x)=\frac12\bigl(\mu_t(x)+\mu_t(x+t)\bigr)$ (Lemma~\ref{lem:geo_sla_pro}(iv)), its left side minus its right side is
\[
\nu_{2t}\mu_{2t}(x)-\tfrac12\nu_t\bigl(\mu_t(x)+\mu_t(x+t)\bigr)=(\nu_{2t}-\nu_t)\,\mu_{2t}(x)\le0.
\]

The shortfalls of $d$ lie in $[0,1]$ (Lemma~\ref{lem:shor_arra}(iv)), so $0\le F_t\le1$ by Lemma~\ref{lem:pro_the_cap}, and $0\le\nu_t\mu_t\le1$; hence $\lvert U_t\rvert\le1$.

Finally, by Lemma~\ref{lem:geo_sla_pro}(iv), $h_t\le C_{\rm cr}\mu_t$, and by Lemma~\ref{lem:pro_the_cap}, $\Phi(u)\le u/u_{\rm lo}$ on $[0,1]$; with $c_{\rm sav}=\alpha_{\rm cap}b_{\rm sl}u_{\rm lo}/C_{\rm cr}$ (Convention~\ref{conv:not_and_con}),
\[
F_t-U_t=\nu_t\mu_t\ge\alpha_{\rm cap}b_{\rm sl}\mu_t\ge\frac{\alpha_{\rm cap}b_{\rm sl}}{C_{\rm cr}}\,h_t\ge\frac{\alpha_{\rm cap}b_{\rm sl}u_{\rm lo}}{C_{\rm cr}}\,F_t=c_{\rm sav}F_t,
\]
that is, $U_t\le(1-c_{\rm sav})F_t$.
\end{proof}

The next three lemmas prepare the extension: the corrector vanishes on unit windows that contain no interface, it has an explicit modulus of continuity on slack profiles, and its average is dominated by the imbalance budget.

\begin{lemma}
 Let $d\in\mathscr{S}$ and $x\in\mathbb{R}$ with $U_1(x;d)\ne0$. Then $d$ is one-ray or two-ray, and the open interval $(x,x+1)$ contains an interface of $d$. If moreover $d$ is two-ray with slow set $[a_-,a_+]$ and slow length $L\le1/(2A_{\rm sl})$, then $x>a_--C_{\rm cr}L$ or $x+1<a_++C_{\rm cr}L$. Consequently $U_{\rm av}(d)=0$ for every $d\in\mathscr{S}$ with $h([-2,-1])(d)>0$ and $h([3,4])(d)>0$.
\label{lem:van_the_pro}
\end{lemma}

\begin{proof}
The corrector vanishes on a unit window inside the core, where the cap is saturated, and inside a ray, where there is neither shortfall nor slow fraction; for a short core it also vanishes on a window that crosses the core with room to spare on both sides, where the deficit equals the crossing deficit and the cap is linear. We use the values of the cap given by Lemma~\ref{lem:pro_the_cap}: $\Phi(0)=0$, $\Phi(u)=\alpha_{\rm cap}u$ for $u\le u_{\rm lo}$ and $\Phi(u)=1$ for $u\ge b_{\rm sl}/2$.

We first show that $U_1(x)=0$ if $d$ is exceptional or $(x,x+1)$ contains no interface; since only one-ray and two-ray profiles have interfaces, this proves the first assertion. If $d$ is exceptional, $U_1=0$ by Definition~\ref{def:the_pro_cor}. Otherwise let $\bar I_*$ be the closure of the slow set, whose boundary points are the interfaces. If $(x,x+1)$ contains no interface, it is a connected set that misses the boundary of $\bar I_*$, so it lies in the interior of $\bar I_*$ or in $\mathbb{R}\setminus\bar I_*$; hence $[x,x+1]\subset\bar I_*$, or $[x,x+1]$ lies in the closure of one component of $\mathbb{R}\setminus I_*$. In the first case $\mu_1(x)=1$ and $h_1(x)\ge b_{\rm sl}$ (Lemma~\ref{lem:geo_sla_pro}(iv)), so $F_1(x)=1$; and $\nu_1=1$, because a two-ray profile then has $L\ge1$, so that $1/L<A_{\rm sl}$ and $\psi(1/L)=0$. Hence $U_1(x)=1-1=0$. In the second case $d(x,x+1)=1$ (Lemma~\ref{lem:geo_sla_pro}(ii)), so $h_1(x)=0$, $F_1(x)=\Phi(0)=0$ and $\mu_1(x)\le h_1(x)/b_{\rm sl}=0$ (Lemma~\ref{lem:geo_sla_pro}(iv)); hence $U_1(x)=0$.

We next treat a short core. Let $d$ be two-ray with slow set $[a_-,a_+]$ and slow length $L\le1/(2A_{\rm sl})$, and suppose that $x\le a_--C_{\rm cr}L$ and $x+1\ge a_++C_{\rm cr}L$. For $\ell_-,\ell_+\ge0$ let
\[
\Gamma(\ell_-,\ell_+)=(L+\ell_-+\ell_+)-d(a_--\ell_-,\,a_++\ell_+)
\]
be the deficit of $d$ on $[a_--\ell_-,a_++\ell_+]$. Then $\Delta=\Gamma(C_{\rm cr}L,C_{\rm cr}L)$, the deficit of $d$ on $[x,x+1]$ is $\Gamma(a_--x,\,x+1-a_+)$, and $0\le\Gamma\le C_{\rm cr}L$ by Definition~\ref{def:slac_prof}(iv). We show that $\Gamma$ is constant on $[C_{\rm cr}L,\infty)^2$. Fix $\ell_+\ge0$, put $y_+=a_++\ell_+$ and $y_s=a_--s$ for $s\ge0$, and let $\gamma(s)=(y_s|y_+)_{a_-}$ be the Gromov product of $d$ at the basepoint $a_-$ (Definition~\ref{def:gro_hyp_sem}). Since $(-\infty,a_-]$ is a unit geodesic ray, $d(a_-,y_s)=s$ and $d(y_s,y_{s'})=\lvert s-s'\rvert$; hence
\[
2\gamma(s)=s+d(a_-,y_+)-d(y_s,y_+)=\Gamma(s,\ell_+)-\Gamma(0,\ell_+)\le C_{\rm cr}L,\qquad(y_s|y_{s'})_{a_-}=s\quad(0\le s\le s').
\]
By the triangle inequality $d(y_{s'},y_+)\le d(y_s,y_+)+(s'-s)$ for $s\le s'$, so $\gamma$ is nondecreasing. Let $C_{\rm cr}L\le s\le s'$. Since $(\mathbb{R},d)$ is $0$-hyperbolic at $a_-$, and $\gamma(s')\le C_{\rm cr}L/2\le s$,
\[
\gamma(s)\ge\min\bigl\{(y_s|y_{s'})_{a_-},\,(y_{s'}|y_+)_{a_-}\bigr\}=\min\{s,\gamma(s')\}=\gamma(s').
\]
With monotonicity, $\gamma(s)=\gamma(s')$, so $\Gamma(\cdot,\ell_+)$ is constant on $[C_{\rm cr}L,\infty)$. The same argument at the basepoint $a_+$, with the point $a_--\ell_-$ fixed and the points $a_++s$ of the ray $[a_+,\infty)$, for which $2(a_--\ell_-|a_++s)_{a_+}=\Gamma(\ell_-,s)-\Gamma(\ell_-,0)$, shows that $\Gamma(\ell_-,\cdot)$ is constant on $[C_{\rm cr}L,\infty)$ for every $\ell_-\ge0$. Therefore $\Gamma(\ell_-,\ell_+)=\Delta$ whenever $\ell_-,\ell_+\ge C_{\rm cr}L$.

Consequently $h_1(x)=\Delta\le C_{\rm cr}L\le C_{\rm cr}/(2A_{\rm sl})<u_{\rm lo}$, since $(C_{\rm cr}+1)/A_{\rm sl}=u_{\rm lo}/4$, and $F_1(x)=\alpha_{\rm cap}\Delta$. On the other hand $[x,x+1]\supset[a_-,a_+]$ gives $\mu_1(x)=L$, and $1/L\ge2A_{\rm sl}$ gives $\psi(1/L)=1$ and $\nu_1=\alpha_{\rm cap}\Delta/L$. Hence $U_1(x)=\alpha_{\rm cap}\Delta-\alpha_{\rm cap}\Delta=0$, which proves the second assertion.

Finally let $d\in\mathscr{S}$ with $U_{\rm av}(d)\ne0$. Since $U_{\rm av}(d)=2\int_0^{1/2}U_1(x;d)\,dx$, some $x\in[0,1/2]$ has $U_1(x;d)\ne0$; so $d$ is one-ray or two-ray, and $(x,x+1)\subset(0,3/2)$ contains an interface $\xi$. By Lemma~\ref{lem:geo_sla_pro}(ii), $\xi$ is the endpoint of a unit geodesic ray of $d$, $(-\infty,\xi]$ or $[\xi,\infty)$. In the first case $[-2,-1]\subset(-\infty,\xi]$ and $h([-2,-1])(d)=0$; in the second $[3,4]\subset[\xi,\infty)$ and $h([3,4])(d)=0$. This proves the last assertion.

Explicit values: $(C_{\rm cr}+1)/A_{\rm sl}=u_{\rm lo}/4$ exactly, and $C_{\rm cr}/(2A_{\rm sl})<u_{\rm lo}$.
\end{proof}

\begin{lemma}
 For all $d,d''\in\mathscr{S}$ and $x\in\mathbb{R}$,
\[
\lvert U_1(x;d)-U_1(x;d'')\rvert\le C_{\rm env}\lVert d-d''\rVert_{[x-2,x+3]}.
\]
In particular $U_1(x;d)$ depends only on the restriction of $d$ to $[x-2,x+3]^2$.
\label{lem:mod_the_pro}
\end{lemma}

\begin{proof}
The corrector of the window $[x,x+1]$ is its capped shortfall minus its weighted slow fraction. The capped shortfall moves by at most $2^{20}$ times the distance. The slow sets of two close profiles nearly agree near the window, because a slow interval of one profile cannot lie where the other has unit speed; and the weighted slow fraction depends in a Lipschitz way on the endpoints and the crossing deficit of a short core, uniformly in its length, because the slow fraction never exceeds the core length.

By Lemma~\ref{lem:cor_ine_sla}(i), $U_1(x;d)=U_1(0;\mathcal{T}_xd)$, where $\mathcal{T}_xd\in\mathscr{S}$ (Lemma~\ref{lem:geo_sla_pro}(iii)), and $\lVert\mathcal{T}_xd-\mathcal{T}_xd''\rVert_{[-2,3]}=\lVert d-d''\rVert_{[x-2,x+3]}$ (Lemma~\ref{lem:shor_arra}(ii)). So it suffices to take $x=0$. Throughout, $Q=[-2,3]$, $\eta=\lVert d-d''\rVert_Q$ and $\eta'=\eta/b_{\rm sl}$; the objects of $d''$ carry a double prime; for a profile that is not exceptional, $\mu=\mu_1(0;d)=\lvert[0,1]\cap I_*\rvert$ and $\nu=\nu_1(d)$.

\textbf{Two profiles that are not exceptional.} We show that $\lvert U_1(0;d)-U_1(0;d'')\rvert\le2^{21}\eta$. If $\eta\ge b_{\rm sl}/100$, this follows from $\lvert U_1\rvert\le1$ (Lemma~\ref{lem:cor_ine_sla}(ii)), since $2\le200\eta'\le2^{21}\eta$. Let $\eta<b_{\rm sl}/100$. By Definition~\ref{def:the_pro_cor}, $U_1(0;d)=F_1(0;d)-\nu\mu$. Here $0\le\nu\le1$, because $b_{\rm sl}\le\Delta/L\le C_{\rm cr}$ for a two-ray profile (Lemma~\ref{lem:geo_sla_pro}(iv)) and $\alpha_{\rm cap}C_{\rm cr}=\frac1{20}$; and $\nu=1$ unless $d$ is two-ray with $L<1/A_{\rm sl}$, in which case we call the weight of $d$ \emph{active}. The shortfalls $h_1(0;d)$ and $h_1(0;d'')$ lie in $[0,1]$ and differ by at most $\eta$ (Lemma~\ref{lem:shor_arra}(iii),(iv)), and $\lvert\Phi'\rvert\le2^{20}$ on $[0,1]$ (Lemma~\ref{lem:pro_the_cap}); hence
\[
\lvert U_1(0;d)-U_1(0;d'')\rvert\le2^{20}\eta+\lvert\nu\mu-\nu''\mu''\rvert.
\]

We compare the slow sets. The sets $I_*\cap Q$ and $I_*''\cap Q$ are intervals, so the sets $(I_*\setminus I_*'')\cap Q$ and $(I_*''\setminus I_*)\cap Q$ are together the union of at most two intervals. Each of these intervals lies in one of the two slow sets and, being connected, in one component of the complement of the other; there the first profile is slow at speed $1-b_{\rm sl}$ and the second has unit speed (Lemma~\ref{lem:geo_sla_pro}(ii)), so its length is at most $\eta'$ by Lemma~\ref{lem:geo_sla_pro}(i). Hence
\[
\bigl\lvert(I_*\,\triangle\,I_*'')\cap Q\bigr\rvert\le2\eta',\qquad\lvert\mu-\mu''\rvert\le2\eta'.
\]
If neither weight is active, $\lvert\nu\mu-\nu''\mu''\rvert=\lvert\mu-\mu''\rvert\le2\eta'$. If $\mu''=0$, then $0\le\nu\mu\le\mu\le2\eta'$; likewise if $\mu=0$.

Suppose now that the weight of $d$ is active and that $\mu,\mu''>0$; the case of an active weight of $d''$ is symmetric. Then $d$ is two-ray with $L<1/A_{\rm sl}$, and $I_*=[a_-,a_+]$ meets $[0,1]$, so $I_*\subset[-1/A_{\rm sl},1+1/A_{\rm sl}]\subset Q$. The set $I_*''\cap Q$ is an interval of length at most $L+2\eta'<\frac1{25}$ that meets $[0,1]$, so it lies in $(-\frac1{25},\frac{26}{25})$. As $I_*''\ne\varnothing$, and $\mathbb{R}$ and every half-line that meets $[0,1]$ contain an endpoint of $Q$, the profile $d''$ is two-ray, with $I_*''=[a_-'',a_+'']\subset(-\frac1{25},\frac{26}{25})$.

If the cores are disjoint, each lies in $Q$ and in a component of the complement of the other, so $L,L''\le\eta'$ by Lemma~\ref{lem:geo_sla_pro}(i),(ii), and $\lvert\nu\mu-\nu''\mu''\rvert\le\max(\mu,\mu'')\le\max(L,L'')\le\eta'$.

Let the cores intersect. If $a_-<a_-''$, the interval $[a_-,a_-'']\subset Q$ lies in $I_*$, as $a_-''\le a_+$, and in the ray $(-\infty,a_-'']$ of $d''$, so $a_-''-a_-\le\eta'$ by Lemma~\ref{lem:geo_sla_pro}(i). Thus $\lvert a_--a_-''\rvert\le\eta'$, likewise $\lvert a_+-a_+''\rvert\le\eta'$, and $\lvert L-L''\rvert\le2\eta'$. Put $\lambda(y)=y\,\psi(1/y)=\min\bigl(y,\max(0,1/A_{\rm sl}-y)\bigr)$ for $y>0$, a $1$-Lipschitz function with $0\le\lambda(y)\le y$ that vanishes for $y\ge1/A_{\rm sl}$. By Definition~\ref{def:the_pro_cor},
\[
\nu\mu=\mu-X,\qquad X=\frac\mu L\,\lambda(L)\Bigl(1-\frac{\alpha_{\rm cap}\Delta}L\Bigr),
\]
and likewise $\nu''\mu''=\mu''-X''$. Since $0\le\mu\le L$, $0\le\mu''\le L''$ and $0\le\alpha_{\rm cap}\Delta/L\le\frac1{20}$,
\[
X-X''=\Bigl(\frac\mu L-\frac{\mu''}{L''}\Bigr)\lambda(L)\Bigl(1-\frac{\alpha_{\rm cap}\Delta}L\Bigr)+\frac{\mu''}{L''}\bigl(\lambda(L)-\lambda(L'')\bigr)\Bigl(1-\frac{\alpha_{\rm cap}\Delta}L\Bigr)+\frac{\mu''}{L''}\,\lambda(L'')\,\alpha_{\rm cap}\Bigl(\frac{\Delta''}{L''}-\frac\Delta L\Bigr).
\]
As $\lambda(L)\le L$, the first term has absolute value at most $\lvert\mu-\mu''L/L''\rvert\le\lvert\mu-\mu''\rvert+(\mu''/L'')\lvert L''-L\rvert\le4\eta'$; the second at most $\lvert L-L''\rvert\le2\eta'$. The third vanishes unless $L''<1/A_{\rm sl}$. In that case the points $b_\pm=a_\pm\pm C_{\rm cr}L$, at which the crossing deficit $\Delta$ is evaluated (Definition~\ref{def:slac_prof}), and the points $b_\pm''=a_\pm''\pm C_{\rm cr}L''$ lie in $[-1,2]\subset Q$, because $(C_{\rm cr}+1)/A_{\rm sl}=u_{\rm lo}/4<\frac1{25}$; and $\lvert b_\pm-b_\pm''\rvert\le\lvert a_\pm-a_\pm''\rvert+C_{\rm cr}\lvert L-L''\rvert\le(2C_{\rm cr}+1)\eta'$. By Lemma~\ref{lem:shor_arra}(iv),
\[
\bigl\lvert d(b_-,b_+)-d''(b_-'',b_+'')\bigr\rvert\le\lvert b_--b_-''\rvert+\lvert b_+-b_+''\rvert+\bigl\lvert d(b_-'',b_+'')-d''(b_-'',b_+'')\bigr\rvert\le2(2C_{\rm cr}+1)\eta'+\eta,
\]
so $\lvert\Delta-\Delta''\rvert\le(2C_{\rm cr}+1)\lvert L-L''\rvert+2(2C_{\rm cr}+1)\eta'+\eta\le(8C_{\rm cr}+4)\eta'+\eta$. With $\lambda(L'')\le L''$, the third term is then at most
\[
\alpha_{\rm cap}\Bigl\lvert\Delta''-\frac{L''}L\Delta\Bigr\rvert\le\alpha_{\rm cap}\Bigl(\lvert\Delta''-\Delta\rvert+\frac\Delta L\lvert L-L''\rvert\Bigr)\le\alpha_{\rm cap}\bigl((10C_{\rm cr}+4)\eta'+\eta\bigr)<0.51\,\eta'.
\]
Hence $\lvert\nu\mu-\nu''\mu''\rvert\le\lvert\mu-\mu''\rvert+\lvert X-X''\rvert<(2+4+2+0.51)\eta'<9\eta'$.

Thus, for two profiles that are not exceptional and $\eta<b_{\rm sl}/100$, $\lvert U_1(0;d)-U_1(0;d'')\rvert\le2^{20}\eta+9\eta'<2^{21}\eta$, which proves the claim.

\textbf{An exceptional profile.} The claim being symmetric in $d$ and $d''$, let $d=d_\sigma$ be exceptional, so that $U_1(0;d)=0$. If $U_1(0;d'')=0$ there is nothing to prove; this covers every exceptional $d''$. Otherwise, by Lemma~\ref{lem:van_the_pro} at $x=0$, $(0,1)$ contains an interface $\xi$ of $d''$, the endpoint of a unit geodesic ray of $d''$ (Lemma~\ref{lem:geo_sla_pro}(ii)), which contains $[\xi,\xi+2]$ or $[\xi-2,\xi]$, a subset of $Q$. On that secant $d''$ equals $2$ and $d$ equals $2\sigma$, so $1-\sigma\le\eta/2$, and
\[
\lVert d_1-d''\rVert_Q\le\lVert d_1-d_\sigma\rVert_Q+\eta=5(1-\sigma)+\eta\le\tfrac72\eta.
\]
The array $d_1$ is a slack profile that is not exceptional, and $U_1(0;d_1)=\Phi(0)-0=0$. The first part, applied to $d_1$ and $d''$, gives $\lvert U_1(0;d)-U_1(0;d'')\rvert=\lvert U_1(0;d_1)-U_1(0;d'')\rvert\le2^{21}\cdot\frac72\eta<2^{23}\eta$.

In all cases $\lvert U_1(0;d)-U_1(0;d'')\rvert\le2^{23}\eta=C_{\rm env}\eta$. The case $\eta=0$ is the locality.

Explicit values: $(C_{\rm cr}+1)/A_{\rm sl}=u_{\rm lo}/4$; $\alpha_{\rm cap}C_{\rm cr}=\frac1{20}$; $\alpha_{\rm cap}(10C_{\rm cr}+4)+\alpha_{\rm cap}b_{\rm sl}<0.51$; $2^{20}+9/b_{\rm sl}=2^{20.04}$; $200/b_{\rm sl}=640000<2^{21}$; $7\cdot2^{20}=2^{22.81}<2^{23}=C_{\rm env}$.
\end{proof}

\begin{lemma}
 For every $d\in\mathscr{S}$, the averaged corrector $U_{\rm av}(d)$ is defined, and
\[
\lvert U_{\rm av}(d)\rvert\le C_{\rm dom}\,\mathcal{V}(d).
\]
\label{lem:spa_dom_the}
\end{lemma}

\begin{proof}
An averaged corrector that does not vanish forces an interface near the window. Next to an interface the shortfall of a half-unit window drops to zero within one unit of time, so the imbalances of the budget see the core; when the core is short, both the averaged corrector and these imbalances are quadratic in its length.

Throughout $d\in\mathscr{S}$ and $f(y)=h_{1/2}(y;d)$. If $d$ is exceptional, $U_1(\cdot;d)=0$; otherwise $x\mapsto F_1(x;d)$ and $x\mapsto\mu_1(x;d)$ are continuous (Lemma~\ref{lem:shor_arra}(iv), Lemma~\ref{lem:geo_sla_pro}(iv)), and so is $x\mapsto U_1(x;d)$. Hence $U_{\rm av}(d)$ is defined. We may assume $U_{\rm av}(d)\ne0$. Then $U_1(x_0;d)\ne0$ for some $x_0\in[0,\frac12]$, and by Lemma~\ref{lem:van_the_pro}, $d$ is one-ray or two-ray and $(x_0,x_0+1)\subset(0,\frac32)$ contains an interface $\xi_0$ of $d$. Put $\ell=\min(\frac14,L)$ if $d$ is two-ray and $\ell=\frac14$ if $d$ is one-ray.

We first show that $\lvert U_{\rm av}(d)\rvert\le4A_{\rm sl}^2\ell^2$. As $\lvert U_1\rvert\le1$ (Lemma~\ref{lem:cor_ine_sla}(ii)), $\lvert U_{\rm av}(d)\rvert\le1$, which suffices unless $d$ is two-ray with $L<1/(2A_{\rm sl})$, because $\frac14\ge1/(2A_{\rm sl})$. In that case, by Lemma~\ref{lem:van_the_pro}, $U_1(x)\ne0$ only if $(x,x+1)$ contains $a_-$ or $a_+$, so that $a_--1<x<a_+$, and moreover $x>a_--C_{\rm cr}L$ or $x+1<a_++C_{\rm cr}L$; as $C_{\rm cr}L<1$, these $x$ form a set of measure at most $\lvert(a_--C_{\rm cr}L,a_+)\rvert+\lvert(a_--1,a_++C_{\rm cr}L-1)\rvert=2(C_{\rm cr}+1)L$. For every $x$ the deficit of $[x,x+1]$ is at most $C_{\rm cr}L<u_{\rm lo}$ (Definition~\ref{def:slac_prof}(iv)), so $0\le F_1(x)=\alpha_{\rm cap}h_1(x)\le\alpha_{\rm cap}C_{\rm cr}L$ (Lemma~\ref{lem:pro_the_cap}); and $\psi(1/L)=1$ gives $0\le\nu_1\mu_1(x)=\alpha_{\rm cap}\Delta\mu_1(x)/L\le\alpha_{\rm cap}C_{\rm cr}L$, since $\mu_1(x)\le L$ and $\Delta\le C_{\rm cr}L$ (Lemma~\ref{lem:geo_sla_pro}(iv)). Hence $\lvert U_1(x)\rvert\le\alpha_{\rm cap}C_{\rm cr}L=L/20$, and
\[
\lvert U_{\rm av}(d)\rvert\le2\int_{\mathbb{R}}\lvert U_1(x)\rvert\,dx\le\frac{(C_{\rm cr}+1)L^2}5\le4A_{\rm sl}^2L^2.
\]

We next show that $\mathcal{V}(d)\ge b_{\rm sl}^2\ell^2/4$. Let $\mathcal{G}=\{j/64:\ -256\le j\le320\}$, the points of $\frac1{64}\mathbb{Z}$ in $[-4,5]$. Since the halves of $[y,y+1]$ are $[y,y+\frac12]$ and $[y+\frac12,y+1]$,
\[
\mathcal{V}(d)=\frac1{64}\sum_{y\in\mathcal{G}}\bigl(f(y)-f(y+\tfrac12)\bigr)^2.
\]
Let $\Omega\subset[-3,4]$ be a closed interval of length $\frac14$, so that it contains at least $16$ points of $\frac1{64}\mathbb{Z}$, and let $\rho>0$. Suppose that $f(y)\ge \rho$ and $f(y-1)=0$ for every $y\in\Omega$. For $y\in\Omega\cap\frac1{64}\mathbb{Z}$ the points $y-1$ and $y-\frac12$ lie in $\mathcal{G}$, and
\[
\rho^2\le\bigl(f(y)-f(y-1)\bigr)^2\le2\bigl(f(y-\tfrac12)-f(y)\bigr)^2+2\bigl(f(y-1)-f(y-\tfrac12)\bigr)^2.
\]
Distinct $y$ give distinct points $y-\frac12$ and distinct points $y-1$, so summing over these $y$ uses each summand of $64\,\mathcal{V}(d)$ at most twice, and $16\rho^2\le4\cdot64\,\mathcal{V}(d)$, that is, $\mathcal{V}(d)\ge \rho^2/16$. The same bound holds if $f(y)\ge \rho$ and $f(y+1)=0$ for every $y\in\Omega$, by $f(y)^2\le2(f(y)-f(y+\frac12))^2+2(f(y+\frac12)-f(y+1))^2$ with $y,y+\frac12\in\mathcal{G}$.

We produce such an interval $\Omega$ with $\rho=2b_{\rm sl}\ell$. By Lemma~\ref{lem:geo_sla_pro}(iv), $2b_{\rm sl}\lvert[y,y+\frac12]\cap I_*\rvert\le f(y)\le2C_{\rm cr}\lvert[y,y+\frac12]\cap I_*\rvert$, and by Lemma~\ref{lem:geo_sla_pro}(ii) every interface of $d$ ends a unit geodesic ray whose interior misses $I_*$.
\begin{itemize}
\item Suppose that some interface $\xi\in(-1,\frac32)$ of $d$ ends the ray $(-\infty,\xi]$. Then $I_*$ is $(\xi,\infty)$ if $d$ is one-ray and $[\xi,\xi+L]$ if $d$ is two-ray. Take $\Omega=[\xi-\frac14,\xi]$. For $y\in\Omega$, $[y,y+\frac12]\cap I_*$ has length at least $\min(y+\frac12-\xi,L)\ge\ell$ (at least $\frac14=\ell$ if $d$ is one-ray), so $f(y)\ge2b_{\rm sl}\ell$; and $[y-1,y-\frac12]\subset(-\infty,\xi)$ misses $I_*$, so $f(y-1)=0$.
\item Otherwise the ray that ends at $\xi_0$ is $[\xi_0,\infty)$, and the closure of $I_*$ contains $[\xi_0-1,\xi_0]$: if $d$ is one-ray, $I_*=(-\infty,\xi_0)$; if $d$ is two-ray, $\xi_0=a_+$, and $a_-\le-1$, since the interface $a_-<\frac32$ ends the ray $(-\infty,a_-]$ and therefore does not lie in $(-1,\frac32)$. Take $\Omega=[\xi_0-\frac34,\xi_0-\frac12]$. For $y\in\Omega$, $[y,y+\frac12]$ lies in the closure of $I_*$, so $f(y)\ge b_{\rm sl}\ge2b_{\rm sl}\ell$; and $[y+1,y+\frac32]\subset(\xi_0,\infty)$ misses $I_*$, so $f(y+1)=0$.
\end{itemize}
In both cases $\Omega\subset(-\frac54,\frac32)\subset[-3,4]$, so $\mathcal{V}(d)\ge(2b_{\rm sl}\ell)^2/16=b_{\rm sl}^2\ell^2/4$.

Comparing the two bounds, $\lvert U_{\rm av}(d)\rvert\le4A_{\rm sl}^2\ell^2\le(16A_{\rm sl}^2/b_{\rm sl}^2)\,\mathcal{V}(d)\le C_{\rm dom}\mathcal{V}(d)$.

Explicit values: $16A_{\rm sl}^2/b_{\rm sl}^2=2^{144.83}<2^{145}=C_{\rm dom}$; $(C_{\rm cr}+1)/5\le4A_{\rm sl}^2$; $1/(2A_{\rm sl})\le\frac14$; $C_{\rm cr}/(2A_{\rm sl})<u_{\rm lo}<1$.
\end{proof}

\begin{lemma}
 Let $d$ be an admissible array.
\begin{enumerate}
\item[(a)] If $h([-2,-1])(d)\ge b_{\rm sl}/4$ and $h([3,4])(d)\ge b_{\rm sl}/4$, then $\widetilde U(d)=0$.
\item[(b)] $\widetilde U(d)$ depends only on the restriction of $d$ to $[-8,8]^2$, and $\lvert\widetilde U(d)-\widetilde U(d'')\rvert\le2^{152}\lVert d-d''\rVert_{[-8,8]}$ for every admissible array $d''$.
\item[(c)] If $d\in\mathscr{S}$, then $\widetilde U(d)=U_{\rm av}(d)$.
\end{enumerate}
\label{lem:ext_the_ave}
\end{lemma}

\begin{proof}
The averaged corrector is controlled on slack profiles only; its Lipschitz envelope extends it exactly, and on slack profiles neither the clamp nor the side factor changes it.

By Lemma~\ref{lem:bou_the_ext}, $\widetilde U(d)=\omega(d)\,y(d)$ is defined for every admissible array $d$, where $y(d)=\max\bigl(-C_{\rm dom}\mathcal{V}(d),\min\bigl(C_{\rm dom}\mathcal{V}(d),(1-c_{\rm sav})F_{\rm av}(d),U_{\rm env}(d)\bigr)\bigr)$ is the clamp of Definition~\ref{def:ext_cor_and}, and $0\le F_{\rm av}(d)\le1$, $\mathcal{V}(d)\ge0$. Since $U_{\rm env}(d)\ge-1$ (Lemma~\ref{lem:bou_the_ext}), $-1\le y(d)\le(1-c_{\rm sav})F_{\rm av}(d)\le1$. The side shortfalls of an admissible array lie in $[0,1]$, so $0\le\omega(d)\le1$.

(a) If both side shortfalls are at least $b_{\rm sl}/4$, then $1-\frac4{b_{\rm sl}}\min\bigl(h([-2,-1])(d),h([3,4])(d)\bigr)\le0$, so $\omega(d)=0$ and $\widetilde U(d)=0$.

(b) We first prove the locality. The envelope $U_{\rm env}$ depends on $d$ only through its restriction to $[-8,8]^2$; the imbalance budget $\mathcal{V}$ only through $[-4,6]^2$, as its intervals $[j/64,j/64+1]$ lie in $[-4,6]$; $F_{\rm av}$ only through $[0,\frac32]^2$; and $\omega$ only through $[-2,4]^2$.

We next prove the Lipschitz bound. Let $d''$ be admissible and $\eta=\lVert d-d''\rVert_{[-8,8]}$, which we may assume finite. The shortfalls of the halves of the intervals of $\mathcal{V}$ lie in $[0,1]$ for both arrays, so each imbalance in $\mathcal{V}$ changes by at most $8\eta$ (Lemma~\ref{lem:shor_arra}(iii)); the coefficients of $\mathcal{V}$ sum to $577/64$, so $C_{\rm dom}\mathcal{V}$ changes by at most $\frac{577}8C_{\rm dom}\eta$. Each shortfall $h_1(x)$ with $0\le x\le\frac12$ lies in $[0,1]$ and changes by at most $\eta$ (Lemma~\ref{lem:shor_arra}(iii)), so $F_{\rm av}$ changes by at most $2^{20}\eta$ (Lemma~\ref{lem:pro_the_cap}). For each $d'\in\mathscr{S}$ the function $U_{\rm av}(d')+C_{\rm env}\lVert\cdot-d'\rVert_{[-8,8]}$ changes by at most $C_{\rm env}\eta$; hence $U_{\rm env}(d)$ and $U_{\rm env}(d'')$ are both finite or both $+\infty$, and in the first case they differ by at most $C_{\rm env}\eta$. A maximum or a minimum of functions changes by at most the largest change of its arguments, so
\[
\lvert y(d)-y(d'')\rvert\le\max\Bigl(\tfrac{577}8C_{\rm dom},\,2^{20},\,C_{\rm env}\Bigr)\eta=\tfrac{577}8C_{\rm dom}\,\eta.
\]
The side shortfalls change by at most $\eta$, so $\lvert\omega(d)-\omega(d'')\rvert\le4\eta/b_{\rm sl}$. As $0\le\omega\le1$ and $\lvert y\rvert\le1$,
\[
\lvert\widetilde U(d)-\widetilde U(d'')\rvert\le\omega(d)\,\lvert y(d)-y(d'')\rvert+\lvert y(d'')\rvert\,\lvert\omega(d)-\omega(d'')\rvert\le\Bigl(\tfrac{577}8C_{\rm dom}+\frac4{b_{\rm sl}}\Bigr)\eta<2^{152}\eta.
\]

(c) Let $d\in\mathscr{S}$; it is admissible by Lemma~\ref{lem:shor_arra}(iv). We first show that $U_{\rm env}=U_{\rm av}$ on $\mathscr{S}$. By Lemma~\ref{lem:spa_dom_the}, $U_{\rm av}$ is defined on $\mathscr{S}$. For $d',d''\in\mathscr{S}$ and $0\le x\le\frac12$, Lemma~\ref{lem:mod_the_pro} gives $\lvert U_1(x;d')-U_1(x;d'')\rvert\le C_{\rm env}\lVert d'-d''\rVert_{[x-2,x+3]}\le C_{\rm env}\lVert d'-d''\rVert_{[-8,8]}$; averaging over $x$,
\[
\lvert U_{\rm av}(d')-U_{\rm av}(d'')\rvert\le C_{\rm env}\lVert d'-d''\rVert_{[-8,8]}\qquad(d',d''\in\mathscr{S}).
\]
Hence $U_{\rm av}(d')+C_{\rm env}\lVert d-d'\rVert_{[-8,8]}\ge U_{\rm av}(d)$ for every $d'\in\mathscr{S}$, with equality at $d'=d$: $U_{\rm env}(d)=U_{\rm av}(d)$, that is, $U_{\rm env}$ is the McShane extension of $U_{\rm av}$ (see \cite[Theorem~1]{locbib:McShane1934}).

We next show that the clamp is the identity on $\mathscr{S}$. By Lemma~\ref{lem:spa_dom_the}, $\lvert U_{\rm av}(d)\rvert\le C_{\rm dom}\mathcal{V}(d)$; averaging the pointwise saving $U_1(x)\le(1-c_{\rm sav})F_1(x)$ of Lemma~\ref{lem:cor_ine_sla}(ii) over $0\le x\le\frac12$ gives $U_{\rm av}(d)\le(1-c_{\rm sav})F_{\rm av}(d)$. So $y(d)=U_{\rm av}(d)$.

Finally, if $U_{\rm av}(d)\ne0$, one side shortfall of $d$ vanishes (Lemma~\ref{lem:van_the_pro}), so $\omega(d)=1$ and $\widetilde U(d)=U_{\rm av}(d)$; if $U_{\rm av}(d)=0$, then $\widetilde U(d)=\omega(d)\cdot0=0$. This proves (c).

Explicit values: $\frac{577}{64}\cdot8=72.125$; $\frac{577}8C_{\rm dom}+4/b_{\rm sl}=2^{151.17}<2^{152}$, with $4/b_{\rm sl}=12800$; $C_{\rm env}=2^{23}$ and $2^{20}$ are below $\frac{577}8C_{\rm dom}$.
\end{proof}

\begin{lemma}

\begin{enumerate}
\item[(i)] Let $d$ be an array such that $d$ and $\mathcal{D}d$ are admissible. Then $Z(d)$ is defined, $\lvert Z(d)\rvert\le3$, and
\[
\bigl\lvert Z(d)-G_{\rm av}(d)\bigr\rvert\le\frac{C_{\rm dom}}{64}\sum_{j=-256}^{320}\Bigl(V\bigl([\tfrac j{32},\tfrac j{32}+2]\bigr)(d)+\frac14\sum_{y\in\{0,\frac12,1,\frac32\}}V\bigl([\tfrac j{64}+y,\tfrac j{64}+y+1]\bigr)(d)\Bigr).
\]
The intervals on the right lie in $[-8,12]$, have length $1$ or $2$ and endpoints and midpoints in $\frac1{64}\mathbb{Z}$, and the coefficients of the imbalances on the right sum to less than $20C_{\rm dom}$.
\item[(ii)] Every $d\in\mathscr{S}$ satisfies the hypothesis of (i); $Z(d)=0$ if $d$ is exceptional, and otherwise
\[
Z(d)=\bigl(\nu_2(d)-\nu_1(d)\bigr)\int_0^1\mu_2(x;d)\,dx\le0.
\]
\item[(iii)] Let $d$ and $d''$ be arrays as in (i). Then $\lvert Z(d)-Z(d'')\rvert\le2^{153}\lVert d-d''\rVert_{[-16,16]}$; in particular $Z(d)=Z(d'')$ if $d$ and $d''$ agree on $[-16,16]^2$.
\item[(iv)] Let $n\ge2$ be an even integer. Then $Z(g_n)$ is defined,
\[
F_{\rm av}(g_n)=\frac2n\sum_{j=0}^{n/2-1}\Phi\bigl(h_g(j,j+n)\bigr),\qquad G_{\rm av}(g_n)=\frac1n\sum_{j=0}^{n-1}G_\Phi(g_n)\circ f_k^j,
\]
and $F_{\rm av}(g_n)$, $G_{\rm av}(g_n)$, the extended correctors $\widetilde U(\mathcal{T}_xg_n)$ for $x\in\mathbb{R}$ and $\widetilde U(g_{2n})$, and $Z(g_n)$ are continuous functions of the initial point.
\end{enumerate}
\label{lem:pro_the_bal}
\end{lemma}

\begin{proof}
The balance defect is the translation average of the balance functional minus a coboundary of extended correctors. The bounds on the extended corrector control the coboundary by imbalances and by one; on a slack profile the extended corrector is the averaged corrector, and the defect becomes the balance of the weighted slow fraction, which is nonpositive; the Lipschitz bounds of the cap and of the extended corrector make the defect Lipschitz.

(i) By Lemma~\ref{lem:shor_arra}(ii) every translate of $d$ is admissible. So, by Lemma~\ref{lem:bou_the_ext}, the extended correctors of $\mathcal{D}d$, $d$, $\mathcal{T}_{1/2}d$, $\mathcal{T}_1d$ and $\mathcal{T}_{3/2}d$ are defined, and each such $d'$ satisfies $\lvert\widetilde U(d')\rvert\le1$ and $\lvert\widetilde U(d')\rvert\le C_{\rm dom}\mathcal{V}(d')$. By Lemma~\ref{lem:shor_arra}(ii), $G_\Phi(\mathcal{T}_xd)=F_2(x;d)-\frac12F_1(x;d)-\frac12F_1(x+1;d)$; the shortfalls of $d$ of lengths $1$ and $2$ lie in $[0,1]$ and are Borel functions of $x$, as $d$ is admissible, and $0\le\Phi\le1$ on $[0,1]$ (Lemma~\ref{lem:pro_the_cap}); so $x\mapsto G_\Phi(\mathcal{T}_xd)$ is Borel with values in $[-1,1]$, and $G_{\rm av}(d)$ is defined with $\lvert G_{\rm av}(d)\rvert\le1$. Hence $Z(d)$ is defined and $\lvert Z(d)\rvert\le1+1+\frac12+\frac12=3$.

By the definition of $Z$, as $\frac12\overline{\widetilde U}(d)+\frac12\overline{\widetilde U}(\mathcal{T}_1d)=\frac14\sum_{y\in\{0,\frac12,1,\frac32\}}\widetilde U(\mathcal{T}_yd)$,
\[
\bigl\lvert Z(d)-G_{\rm av}(d)\bigr\rvert\le\lvert\widetilde U(\mathcal{D}d)\rvert+\frac14\sum_{y\in\{0,\frac12,1,\frac32\}}\lvert\widetilde U(\mathcal{T}_yd)\rvert\le C_{\rm dom}\Bigl(\mathcal{V}(\mathcal{D}d)+\frac14\sum_{y\in\{0,\frac12,1,\frac32\}}\mathcal{V}(\mathcal{T}_yd)\Bigr).
\]
By Lemma~\ref{lem:shor_arra}(ii), $V(I)(\mathcal{D}d)=V(2I)(d)$ and $V(I)(\mathcal{T}_yd)=V(I+y)(d)$; substituting these identities in the definition of $\mathcal{V}$ gives the inequality. The intervals $[j/32,j/32+2]$ with $-256\le j\le320$ lie in $[-8,12]$, and the intervals $[j/64+y,j/64+y+1]$ lie in $[-4,\frac{15}2]$; their endpoints and midpoints lie in $\frac1{64}\mathbb{Z}$; and the coefficients sum to $\frac{577}{64}C_{\rm dom}\bigl(1+4\cdot\frac14\bigr)=\frac{577}{32}C_{\rm dom}<20C_{\rm dom}$.

(ii) Let $d\in\mathscr{S}$. By Lemma~\ref{lem:geo_sla_pro}(iii), $\mathcal{D}d$ and the translates of $d$ are slack profiles, exceptional exactly when $d$ is, and slack profiles are admissible (Lemma~\ref{lem:shor_arra}(iv)); so (i) applies, and by Lemma~\ref{lem:ext_the_ave}(c) every extended corrector in $Z(d)$ is an averaged corrector. If $d=d_\sigma$ is exceptional, these correctors vanish and $G_\Phi(\mathcal{T}_xd)=\Phi(1-\sigma)(1-\frac12-\frac12)=0$, so $Z(d)=0$. Let $d$ be not exceptional. By Lemma~\ref{lem:cor_ine_sla}(i),
\[
\begin{gathered}
\widetilde U(\mathcal{D}d)=2\int_0^{1/2}U_2(2x;d)\,dx=\int_0^1U_2(x;d)\,dx,\\
\overline{\widetilde U}(d)=\int_0^1U_1(x;d)\,dx,\qquad\overline{\widetilde U}(\mathcal{T}_1d)=\int_1^2U_1(x;d)\,dx.
\end{gathered}
\]
Since $F_t-U_t=\nu_t\mu_t$ with $\nu_t$ independent of $x$ (Definition~\ref{def:the_pro_cor}), and $\mu_2(x)=\frac12\bigl(\mu_1(x)+\mu_1(x+1)\bigr)$ (Lemma~\ref{lem:geo_sla_pro}(iv)),
\[
Z(d)=\int_0^1\Bigl(\nu_2\mu_2(x;d)-\tfrac12\nu_1\mu_1(x;d)-\tfrac12\nu_1\mu_1(x+1;d)\Bigr)dx=(\nu_2-\nu_1)\int_0^1\mu_2(x;d)\,dx.
\]
Finally $\nu_2\le\nu_1$: this is clear unless $d$ is two-ray, and then $\nu_2-\nu_1=-\bigl(\psi(2/L)-\psi(1/L)\bigr)\bigl(1-\alpha_{\rm cap}\Delta/L\bigr)\le0$, because $\psi$ is nondecreasing and $\alpha_{\rm cap}\Delta/L\le\alpha_{\rm cap}C_{\rm cr}=\frac1{20}$ (Lemma~\ref{lem:geo_sla_pro}(iv)). Hence $Z(d)\le0$.

(iii) Let $\eta=\lVert d-d''\rVert_{[-16,16]}$, which we may assume finite. For $0\le x\le1$ the intervals $[x,x+2]$, $[x,x+1]$ and $[x+1,x+2]$ lie in $[0,3]$; their shortfalls lie in $[0,1]$ for both arrays and change by at most $\eta/2$, $\eta$ and $\eta$ (Lemma~\ref{lem:shor_arra}(iii)); as $\lvert\Phi'\rvert\le2^{20}$ on $[0,1]$ (Lemma~\ref{lem:pro_the_cap}), $G_\Phi(\mathcal{T}_xd)$ and hence $G_{\rm av}(d)$ change by at most $2^{20}(\frac12+\frac12+\frac12)\eta<2^{21}\eta$. By Lemma~\ref{lem:shor_arra}(ii), $\lVert\mathcal{D}d-\mathcal{D}d''\rVert_{[-8,8]}=\frac12\eta$ and $\lVert\mathcal{T}_yd-\mathcal{T}_yd''\rVert_{[-8,8]}=\lVert d-d''\rVert_{[y-8,y+8]}\le\eta$ for $0\le y\le\frac32$. By Lemma~\ref{lem:ext_the_ave}(b), the three correction terms of $Z$ change by at most $2^{152}\bigl(\frac12+\frac12+\frac12\bigr)\eta=3\cdot2^{151}\eta$. Hence $\lvert Z(d)-Z(d'')\rvert\le(2^{21}+3\cdot2^{151})\eta<2^{153}\eta$.

(iv) By Lemma~\ref{lem:dil_and_tim}(vi), $g_n$ and $\mathcal{D}g_n=g_{2n}$ are admissible, so $Z(g_n)$ is defined by (i). Let $0\le x\le1$ and $j=\lfloor nx\rfloor$. By Lemma~\ref{lem:shor_arra}(ii) and Lemma~\ref{lem:dil_and_tim}(iv),(i), the shortfalls of $g_n$ on $[x,x+t]$ with $t\in\{1,2\}$ are those of $\mathcal{T}_xg_n$ on $[0,t]$, which agree with those of $\mathcal{T}_{j/n}g_n$; so $h_1(x;g_n)=h_g(j,j+n)$, and $G_\Phi(\mathcal{T}_xg_n)=G_\Phi(\mathcal{T}_{j/n}g_n)=G_\Phi(g_n)\circ f_k^j$ by Lemma~\ref{lem:dil_and_tim}(iii). Integrating over $x$, these functions being constant on each interval $[j/n,(j+1)/n)$, gives the two finite averages. Every $g(a,b)$ is a continuous function of the initial point (Lemma~\ref{lem:dil_and_tim}(v)), and $G_\Phi(\mathcal{T}_{j/n}g_n)$ is a continuous function of $g(j,j+n)$, $g(j+n,j+2n)$ and $g(j,j+2n)$; so $F_{\rm av}(g_n)$ and $G_{\rm av}(g_n)$ are continuous. By Lemma~\ref{lem:dil_and_tim}(v), the restrictions of $g_n$ and of $g_{2n}$ to $[-r,r]^2$ depend continuously on the initial point for the uniform distance, for every $r>0$; since $\lVert\mathcal{T}_xd-\mathcal{T}_xd''\rVert_{[-8,8]}=\lVert d-d''\rVert_{[x-8,x+8]}$ for all arrays $d,d''$ (Lemma~\ref{lem:shor_arra}(ii)), the same holds for the restrictions of $\mathcal{T}_xg_n$ to $[-8,8]^2$. These arrays are admissible (Lemma~\ref{lem:dil_and_tim}(vi)), so by Lemma~\ref{lem:ext_the_ave}(b) the extended correctors $\widetilde U(\mathcal{T}_xg_n)$ and $\widetilde U(g_{2n})$ are continuous functions of the initial point, and so is $Z(g_n)$.

Explicit values: $\frac{577}{32}<20$; $2^{21}+3\cdot2^{151}=2^{152.58}<2^{153}$.
\end{proof}

\subsection{Approximation by slack profiles} \label{subsec:app_sla_pro}

This subsection proves that rescaled growth arrays are uniformly close to slack profiles outside an exceptional set whose mass is summable over the scales (Lemma~\ref{lem:sum_exc_mas:1}); the defect bound uses this on the body. The first two lemmas are deterministic. They concern a path semimetric, a semimetric on $\{0,\dots,R\}$ with steps $d(i,i+1)$ at most one, whose nonaffine markers, where it bends or changes step length, are joined to its nonunit markers, where it is not at nearly unit speed, beyond a fixed range only by slow secants. If it is $0$-hyperbolic, Lemma~\ref{lem:sha_hyp_pat} gives its shape; if it is $0$-hyperbolic at the basepoint $0$ and has no near-total cancellation, Lemma~\ref{lem:app_fin_tre} finds a slack profile uniformly close to it.

\begin{definition}
 Let $R\ge2$ be an integer and let $d$ be a \textbf{path semimetric}, a semimetric on $\{0,\dots,R\}$ (Definition~\ref{def:gro_hyp_sem}) whose \textbf{steps} $l_i=d(i,i+1)$ satisfy $l_i\le1$ for $0\le i<R$; the Gromov products $(i-1|i+1)_i$ below are those of $d$. Let $\zeta>0$.
\begin{itemize}
\item The \textbf{nonaffine markers} of $d$ at tolerance $\zeta$ are the elements of
\[
P=\bigl\{i\in\{1,\dots,R-1\}:\ (i-1|i+1)_i>\zeta\ \text{or}\ \lvert l_i-l_{i-1}\rvert>\zeta\bigr\}.
\]
\item The \textbf{nonunit markers} of $d$ at tolerance $\zeta$ are the indices $i\in\{0,\dots,R\}$ that have an incident step shorter than $1-\zeta/2$ (that is, $i\ge1$ and $l_{i-1}<1-\zeta/2$, or $i<R$ and $l_i<1-\zeta/2$), or that satisfy $1\le i\le R-1$ and $(i-1|i+1)_i>\zeta/2$; their set is $Q$.
\item For $c\in[1/2,1)$ and $B\ge1$, $d$ has \textbf{slow marked pairs at range $B$} for $c$, at tolerance $\zeta$, if
\[
d(i,j)\le c\lvert i-j\rvert\qquad\text{whenever }i\in P,\ j\in Q,\ \lvert i-j\rvert\ge B.
\]
\end{itemize}
\label{def:fin_tre_pat:1}
\end{definition}

\begin{lemma}
 Let $R\ge2$ be an integer and let $d$ be a $0$-hyperbolic semimetric on $\{0,\dots,R\}$ whose steps $l_i=d(i,i+1)$ satisfy $l_i\le1$ for $0\le i<R$. Let $c\in[1/2,1)$, $0<\gamma\le(1-c)/4$, $0<\zeta\le\gamma/(4R)$ and $B\ge1$. Let $P$ and $Q$ be the markers of $d$ at tolerance $\zeta$, and suppose that $d$ has slow marked pairs at range $B$ for $c$.
\begin{enumerate}
\item[(a)] If $P=\varnothing$, there is $\sigma\in[0,1]$ with $\bigl\lvert d(i,j)-\sigma\lvert i-j\rvert\bigr\rvert\le2R^2\zeta$ for $0\le i,j\le R$.
\item[(b)] Suppose $P\ne\varnothing$, and let $a=\min P$ and $b=\max P$. Put $a_*=0$ if $l_i\le c+\gamma$ for all $0\le i<a$, and $a_*=a$ otherwise; put $b_*=R$ if $l_i\le c+\gamma$ for all $b\le i<R$, and $b_*=b$ otherwise. Then
\[
d(i,j)\le(c+\gamma)\lvert i-j\rvert+3B\qquad(a_*\le i,j\le b_*).
\]
If $a_*=a$, then $d(i,j)\ge\lvert i-j\rvert-B-1$ for $0\le i,j\le a$. If $b_*=b$, then $d(i,j)\ge\lvert i-j\rvert-B-1$ for $b\le i,j\le R$.
\end{enumerate}
\label{lem:sha_hyp_pat}
\end{lemma}

\begin{proof}
The proof uses only the four-point condition. Between consecutive nonaffine markers the Gromov products are tiny and the steps nearly constant, so distances there are sums of steps. Such a stretch of speed above $c$ that ends at a nonaffine marker contains no nonunit marker at distance $B$ or more from that marker, by the slow pairs; hence it is shorter than $B$ when both its ends are nonaffine markers, and it has nearly unit speed when it is an end of the path.

Throughout, $\kappa_i=(i-1|i+1)_i$ for $1\le i\le R-1$. Gromov products of $d$ are nonnegative, and $d(i,j)\le\lvert i-j\rvert$ as $l_i\le1$, both by the triangle inequality.

\textbf{Straight stretches.} We first show that along a stretch on which every step exceeds the Gromov products at its two ends, distances are sums of steps minus twice the interior Gromov products. Let $0\le i_0<i_1\le R$ be integers, put $\kappa'_r=\kappa_r$ for $i_0<r<i_1$ and $\kappa'_{i_0}=\kappa'_{i_1}=0$, and suppose that $l_r>\kappa'_r+\kappa'_{r+1}$ for $i_0\le r<i_1$. We claim that
\begin{equation}
d(i,j)=\sum_{i\le r<j}l_r-2\sum_{i<r<j}\kappa_r\qquad(i_0\le i<j\le i_1).\tag{$\ast$}
\end{equation}
Since Gromov products are nonnegative, the hypothesis for $[i_0,i_1]$ implies the hypothesis for $[i,i_1]$ whenever $i_0\le i<i_1$; so it suffices to prove $(\ast)$ for $i=i_0$, by induction on $j$. For $j=i_0+1$, $(\ast)$ is the definition of $l_{i_0}$. For $j=i_0+2$ it is the identity $2\kappa_{i_0+1}=l_{i_0}+l_{i_0+1}-d(i_0,i_0+2)$, which is the definition of the Gromov product. Let $i_0+3\le j\le i_1$, suppose that $(\ast)$ holds at $j-1$ and at $j-2$, and put $w=j-1$; then $i_0<j-2$ and $w<i_1$. By $(\ast)$ at $j-1$ and at $j-2$, $d(i_0,w)-d(i_0,j-2)=l_{j-2}-2\kappa_{j-2}$, hence
\[
(i_0|j-2)_w=\tfrac12\bigl(d(i_0,w)+l_{j-2}-d(i_0,j-2)\bigr)=l_{j-2}-\kappa_{j-2}>\kappa_{j-1}=(j-2|j)_w,
\]
the inequality being the hypothesis at $r=j-2$. By $0$-hyperbolicity at $w$, $(i_0|j)_w\ge\min\{(i_0|j-2)_w,(j-2|j)_w\}=\kappa_{j-1}$. Again by $0$-hyperbolicity at $w$, $\kappa_{j-1}=(j-2|j)_w\ge\min\{(j-2|i_0)_w,(i_0|j)_w\}$; as $(j-2|i_0)_w>\kappa_{j-1}$, this minimum is $(i_0|j)_w$, so $(i_0|j)_w\le\kappa_{j-1}$. Hence $(i_0|j)_w=\kappa_{j-1}$, that is, $d(i_0,j)=d(i_0,w)+l_{j-1}-2\kappa_{j-1}$, which together with $(\ast)$ at $j-1$ is $(\ast)$ at $j$.

\textbf{Unmarked intervals.} We next show that the steps are nearly constant between consecutive nonaffine markers, and that a fast stretch there has speed above $c+\gamma/2$. Call an integer interval $[i_0,i_1]$ \emph{unmarked} if $i_0<i_1$ are consecutive points of $P\cup\{0,R\}$, and call it \emph{fast} if moreover $l_r>c+\gamma$ for some $i_0\le r<i_1$. On an unmarked interval, $\kappa_r\le\zeta$ and $\lvert l_r-l_{r-1}\rvert\le\zeta$ for $i_0<r<i_1$, because such $r$ do not lie in $P$; so any two of the steps $l_r$, $i_0\le r<i_1$, differ by at most $R\zeta\le\gamma/4$. On a fast unmarked interval every step therefore exceeds $c+3\gamma/4$, which exceeds $2\zeta\ge\kappa'_r+\kappa'_{r+1}$ (as $c\ge1/2$ and $\zeta\le\gamma/8\le1/64$); so $(\ast)$ holds on it, and since $2\zeta\le\gamma/4$,
\begin{equation}
d(i,j)>\bigl(c+\tfrac34\gamma\bigr)(j-i)-2\zeta(j-i)\ge\bigl(c+\tfrac12\gamma\bigr)(j-i)\qquad(i_0\le i<j\le i_1).\tag{1}
\end{equation}
On an unmarked interval that is not fast, every step is at most $c+\gamma$, and $d(i,j)\le(c+\gamma)(j-i)$ for $i_0\le i<j\le i_1$ by the triangle inequality.

\textbf{The affine case.} We prove (a). Suppose $P=\varnothing$, so that $[0,R]$ is unmarked. If every step is at most $2R\zeta$, then $d(i,j)\le2R\zeta\lvert i-j\rvert\le2R^2\zeta$ by the triangle inequality, and $\sigma=0$ works. Otherwise some step exceeds $2R\zeta$, so every step exceeds $R\zeta\ge2\zeta\ge\kappa'_r+\kappa'_{r+1}$, and $(\ast)$ holds on $[0,R]$. Put $\sigma=l_0\in[0,1]$. For $0\le i<j\le R$, $(\ast)$ gives
\[
\bigl\lvert d(i,j)-\sigma(j-i)\bigr\rvert\le\sum_{i\le r<j}\lvert l_r-l_0\rvert+2\sum_{i<r<j}\kappa_r\le R\cdot R\zeta+2R\zeta\le2R^2\zeta,
\]
as $2R\le R^2$. By symmetry, and as $d(i,i)=0$, this holds for all $0\le i,j\le R$. This proves (a).

From now on $P\ne\varnothing$, and $a,b,a_*,b_*$ are as in (b). Since $0<a\le b<R$, the intervals $[0,a]$ and $[b,R]$ are unmarked, and $a_*=a$ exactly when $[0,a]$ is fast, $b_*=b$ exactly when $[b,R]$ is fast.

\textbf{Markers.} We show that $P\subset Q$ and that a fast unmarked interval with both endpoints in $P$ is shorter than $B$. Let $i\in P$. If $\kappa_i>\zeta$, then $\kappa_i>\zeta/2$ and $i\in Q$. Otherwise $\lvert l_i-l_{i-1}\rvert>\zeta$; two steps in $[1-\zeta/2,1]$ differ by at most $\zeta/2<\zeta$, so $l_{i-1}<1-\zeta/2$ or $l_i<1-\zeta/2$, and $i\in Q$. Now let $[i_0,i_1]$ be a fast unmarked interval with $i_0,i_1\in P$. If $i_1-i_0\ge B$, the slow pair $(i_0,i_1)$, with $i_0\in P$ and $i_1\in P\subset Q$, would give $d(i_0,i_1)\le c(i_1-i_0)$, against (1). Hence $i_1-i_0<B$.

\textbf{The core.} We prove the upper bound of (b). For $0\le i\le j\le R$ put $E(i,j)=d(i,j)-(c+\gamma)(j-i)$; by the triangle inequality, $E(i,l)\le E(i,j)+E(j,l)$ for $i\le j\le l$, and $E(i,i)=0$. Two bounds on $E$ come first.
\begin{enumerate}
\item[(i)] \emph{If $i<j$ and $[i,j]$ is contained in an unmarked interval $[i_0,i_1]$ that is contained in $[a_*,b_*]$, then $E(i,j)\le B$.} If $[i_0,i_1]$ is not fast, $E(i,j)\le0$ by the paragraph on unmarked intervals. Let it be fast. Then $i_0\ne0$: the only unmarked interval with left endpoint $0$ is $[0,a]$, and if it is contained in $[a_*,b_*]$, then $a_*=0$ and $[0,a]$ is not fast. Likewise $i_1\ne R$. So $i_0,i_1\in P$, hence $i_1-i_0<B$ by the paragraph on markers, and $E(i,j)\le d(i,j)\le j-i<B$.
\item[(ii)] \emph{If $i\le j$ lie in $P$, then $E(i,j)\le B$.} If $j-i\ge B$, the slow pair $(i,j)$, with $j\in P\subset Q$, gives $E(i,j)\le-\gamma(j-i)\le0$; otherwise $E(i,j)\le d(i,j)\le j-i<B$.
\end{enumerate}
We also note: \emph{if $i<j$ and no point of $P$ lies strictly between $i$ and $j$, then $[i,j]$ is contained in an unmarked interval contained in $[a_*,b_*]$, provided $a_*\le i$ and $j\le b_*$.} Indeed, let $i_0$ be the largest point of $P\cup\{0,R\}$ with $i_0\le i$ and $i_1$ the least with $i_1\ge j$. No point of $P\cup\{0,R\}$ lies strictly between $i_0$ and $i_1$: there is none in $(i_0,i]$ by the choice of $i_0$, none in $[j,i_1)$ by the choice of $i_1$, and none in $(i,j)$, since $P$ has none there and $0\le i<j\le R$. So $[i_0,i_1]$ is unmarked and contains $[i,j]$. If $a_*=a$, then $a\le i$ and $a\in P$, so $a\le i_0$; if $a_*=0$, then $0\le i_0$ trivially. Likewise $i_1\le b_*$.

Now let $a_*\le i<j\le b_*$. If $P$ has no point in $[i,j]$, then $E(i,j)\le B$ by the note and (i). Otherwise let $i'$ and $j'$ be the least and the largest point of $P$ in $[i,j]$. If $i<i'$, then $P$ has no point in $[i,i')$, so $E(i,i')\le B$ by the note and (i); if $i=i'$, then $E(i,i')=0$. Likewise $E(j',j)\le B$, as $P$ has no point in $(j',j]$. By (ii), $E(i',j')\le B$. Hence
\[
E(i,j)\le E(i,i')+E(i',j')+E(j',j)\le3B,
\]
that is, $d(i,j)\le(c+\gamma)(j-i)+3B$. As $d$ is symmetric and $d(i,i)=0$, the upper bound of (b) holds for all $a_*\le i,j\le b_*$.

\textbf{Fast exteriors.} We prove the lower bounds of (b). Let $a_*=a$, so that $[0,a]$ is a fast unmarked interval. No $j\le a-B$ lies in $Q$: otherwise the slow pair $(a,j)$, with $a\in P$, would give $d(j,a)\le c(a-j)$, against (1). Hence, for $0\le r\le a-B$, the index $r$ is not a nonunit marker, and its incident step $l_r$ is at least $1-\zeta/2$. So fewer than $B$ of the steps $l_r$ with $0\le r<a$ are shorter than $1-\zeta/2$, namely at most those with $a-B<r<a$. Let $0\le i<j\le a$. By $(\ast)$ on $[0,a]$ and $\kappa_r\le\zeta$,
\[
d(i,j)\ge(1-\zeta/2)(j-i-B)-2\zeta(j-i)\ge(j-i)-B-\tfrac52\zeta R\ge(j-i)-B-1,
\]
since $\tfrac52\zeta R\le\tfrac58\gamma<1$; the first inequality holds trivially when $j-i<B$, its right-hand side being negative. As $d$ is symmetric and $d(i,i)=0$, the bound holds for all $0\le i,j\le a$. The case $b_*=b$ is the same after the reflection $i\mapsto R-i$ of the indices, which maps $d$ to a semimetric with the same properties, $P$ and $Q$ to the markers of the reflected semimetric, $b$ to its least nonaffine marker and $[b,R]$ to its first unmarked interval.
\end{proof}

\begin{lemma}
 Let $R\ge2$ be an integer, and let $d$ be a semimetric on $\{0,\dots,R\}$ that is $0$-hyperbolic at the basepoint $0$ and satisfies $d(i,i+1)\le1$ for $0\le i<R$. Let $B\ge1$ and $0<\zeta\le\beta_{\rm br}/(64R)$, and put $H=\lceil100B/\beta_{\rm br}\rceil$. Suppose that $d$ has slow marked pairs at range $B$ for $c=1-\beta_{\rm br}$, at tolerance $\zeta$, and that $d$ has no $\delta_{\rm can}/2$-cancellation (Definition~\ref{def:can_the_gro}, with $Y=\{0,\dots,R\}$). Then some $d_*\in\mathscr{S}$ satisfies
\[
\bigl\lvert d(i,j)-d_*(i,j)\bigr\rvert\le6H+2R^2\zeta\qquad(0\le i,j\le R).
\]
\label{lem:app_fin_tre}
\end{lemma}

\begin{proof}
 We realise $d$ in a real tree and replace the path through the images of $0,\dots,R$ by a slack profile. The slow core is coarse-grained at scale $H$ and its time interval is stretched by $H$ at each end, which leaves speed slack; each fast exterior is replaced by a unit geodesic ray. Every condition on a two-ray profile is then automatic except the deficit bound. An excessive deficit means that the two exteriors fold back over each other far beyond the core, and this fold is a cancellation of $d$.

\textbf{The tree and the shape.} After identifying points at distance zero, $(\{0,\dots,R\},d)$ embeds isometrically in a real tree $X$ \cite[Lemma 2.13]{bib:bestvina1997}; write $d_X$ for the metric of $X$ and $p_i$ for the image of $i$, so that $d_X(p_i,p_j)=d(i,j)$. A real tree is $0$-hyperbolic at every basepoint \cite[Section 2.2, Example 2.11]{bib:bestvina1997}, so $d$ is $0$-hyperbolic. Put $c=1-\beta_{\rm br}$ and $\gamma=\beta_{\rm br}/16$; then $c\in[1/2,1)$, $0<\gamma\le(1-c)/4$ and $\zeta\le\gamma/(4R)$, so Lemma~\ref{lem:sha_hyp_pat} applies to $d$ with $c,\gamma,\zeta,B$. Let $P$ be the set of nonaffine markers of $d$ at tolerance $\zeta$.

If $P=\varnothing$, part (a) of Lemma~\ref{lem:sha_hyp_pat} gives $\sigma\in[0,1]$ with $\lvert d(i,j)-\sigma\lvert i-j\rvert\rvert\le2R^2\zeta$ for $0\le i,j\le R$, and $d_*=d_\sigma$ is an affine slack profile: $d_\sigma(x,y)\le\lvert x-y\rvert$, and $(\mathbb{R},d_\sigma)$ is $0$-hyperbolic, all its Gromov products vanishing if $\sigma=0$, and it being isometric to the real line, a real tree, if $\sigma>0$.

Let now $P\ne\varnothing$, and let $a,b,a_*,b_*$ be as in part (b) of Lemma~\ref{lem:sha_hyp_pat}. Put
\[
L_0=b_*-a_*,\qquad a'=a_*-H,\qquad b'=b_*+H,\qquad L_*=b'-a'=L_0+2H.
\]
We call the left side \emph{fast} if $a_*=a$ and the right side fast if $b_*=b$. Let $p\colon[0,R]\to X$ run on $[i,i+1]$ along the geodesic segment from $p_i$ to $p_{i+1}$ at constant speed $d(i,i+1)\le1$ (a real tree is a geodesic space, \cite[Section 2, Definitions 2.1--2.2]{bib:bestvina1997}); then $p(i)=p_i$, and $p$ is $1$-Lipschitz. Let $X'$ be the real tree obtained from $X$ by gluing two distinct copies of $[0,\infty)$, the left half-line by its endpoint to $p_0$ and the right half-line to $p_R$ (two real trees glued at one point form a real tree). We write $d_X$ also for the metric of $X'$: it extends that of $X$, and the distance from a point $u\in X$ to the point at distance $\xi$ from $p_0$ on the left half-line is $d_X(u,p_0)+\xi$, likewise for the right half-line and $p_R$, and the points at distance $\xi$ from $p_0$ on the left half-line and at distance $\xi'$ from $p_R$ on the right half-line are at distance $\xi+d_X(p_0,p_R)+\xi'$. We now build a continuous path $\widetilde p\colon\mathbb{R}\to X'$ and put $d_*(x,y)=d_X(\widetilde p(x),\widetilde p(y))$ for $x,y\in\mathbb{R}$.

\textbf{The core.} We define $\widetilde p$ on $[a',b']$, with $\widetilde p(a')=p_{a_*}$ and $\widetilde p(b')=p_{b_*}$, so that $\widetilde p$ is $(1-b_{\rm sl})$-Lipschitz on $[a',b']$ and
\begin{equation}
d_X\bigl(\widetilde p(x),p(x)\bigr)<3H\qquad(x\in[a',b']\cap[0,R]).\tag{$\dagger$}
\end{equation}
Suppose first that $L_0\ge H$. Partition $[a_*,b_*]$ into consecutive intervals with integer endpoints and lengths between $H$ and $2H-1$ (possible as $L_0\ge H$: take $\lfloor L_0/H\rfloor$ pieces of length $H$ and distribute the remainder, which is below $H$, among them), and let $\widehat p\colon[a_*,b_*]\to X$ agree with $p$ at the endpoints of these pieces and run along the geodesic segment between them at constant speed on each piece. For a piece $[y,y']$, part (b) of Lemma~\ref{lem:sha_hyp_pat} gives $d(y,y')\le(c+\gamma)(y'-y)+3B$, so the speed of $\widehat p$ on it is at most $c+\gamma+3B/H<1-b_{\rm sl}$; hence $\widehat p$ is $(1-b_{\rm sl})$-Lipschitz. For $x$ in a piece $[y,y']$, since $p$ and $\widehat p$ are $1$-Lipschitz and agree at $y$ and $y'$, $d_X(\widehat p(x),p(x))\le2\min(x-y,y'-x)\le y'-y<2H$. For $x\in[a',b']$ put $\bar x=a_*+(x-a')L_0/L_*$ and $\widetilde p(x)=\widehat p(\bar x)$. The map $x\mapsto\bar x$ is affine and increasing from $[a',b']$ onto $[a_*,b_*]$, with slope $L_0/L_*<1$ and $\lvert\bar x-x\rvert\le H$ (its value is $a'+H$ at $a'$ and $b'-H$ at $b'$). So $\widetilde p$ is $(1-b_{\rm sl})$-Lipschitz on $[a',b']$, and for $x\in[a',b']\cap[0,R]$,
\[
d_X\bigl(\widetilde p(x),p(x)\bigr)\le d_X\bigl(\widehat p(\bar x),p(\bar x)\bigr)+d_X\bigl(p(\bar x),p(x)\bigr)<2H+H.
\]
Suppose next that $L_0<H$. Let $\widetilde p$ run on $[a',b']$ along the geodesic segment from $p_{a_*}$ to $p_{b_*}$ at constant speed $d(a_*,b_*)/L_*\le L_0/L_*<1/3<1-b_{\rm sl}$. For $x\in[a',b']\cap[0,R]$, since $x-a'\le L_*<3H$ and $\lvert x-a_*\rvert\le L_0+H<2H$,
\[
d_X\bigl(\widetilde p(x),p(x)\bigr)\le d_X\bigl(\widetilde p(x),\widetilde p(a')\bigr)+d_X\bigl(p_{a_*},p(x)\bigr)\le\tfrac13(x-a')+\lvert x-a_*\rvert<H+2H.
\]
In both cases ($\dagger$) holds and $\widetilde p$ is $(1-b_{\rm sl})$-Lipschitz on $[a',b']$.

\textbf{The sides.} We extend $\widetilde p$ to $(-\infty,a']$ and to $[b',\infty)$. If the left side is not fast, put $\widetilde p(x)=p_0$ for $x\le a'$ (then $a_*=0$, so this agrees with $\widetilde p(a')=p_{a_*}=p_0$). If it is fast, put $\ell_-=d(0,a)$, let $\rho_-\colon[0,\infty)\to X'$ run at unit speed from $p_a$ along the geodesic segment to $p_0$ and then along the left half-line, which is a geodesic ray by the distance formula on $X'$, and put $\widetilde p(x)=\rho_-(a'-x)$ for $x\le a'$. We show that in this case
\begin{equation}
d_X\bigl(\widetilde p(x),p(x)\bigr)<2H\qquad(0\le x\le a').\tag{$\ddagger$}
\end{equation}
By part (b) of Lemma~\ref{lem:sha_hyp_pat}, $d(i,a)\ge a-i-B-1$ for $0\le i\le a$, in particular $\ell_-\ge a-B-1$. With $d(i,0)\le i$ this gives, for the Gromov products of $d_X$ at $p_a$,
\begin{equation}
2(p_i|p_0)_{p_a}=d(a,i)+\ell_--d(i,0)\ge2(a-i)-2B-2\qquad(0\le i\le a).\tag{$\sharp$}
\end{equation}
As $\rho_-$ is a geodesic ray with $\rho_-(0)=p_a$ and $\rho_-(\ell_-)=p_0$, $(p_0|\rho_-(\xi))_{p_a}=\frac12\bigl(\ell_-+\xi-\lvert\ell_--\xi\rvert\bigr)=\min(\xi,\ell_-)$ for $\xi\ge0$. By $0$-hyperbolicity of the real tree $X'$ at $p_a$ \cite[Section 2.2, Example 2.11]{bib:bestvina1997},
\[
\bigl(p_i\big|\rho_-(a-i)\bigr)_{p_a}\ge\min\Bigl\{(p_i|p_0)_{p_a},\bigl(p_0\big|\rho_-(a-i)\bigr)_{p_a}\Bigr\}\ge a-i-B-1,
\]
since $\min(a-i,\ell_-)\ge a-i-B-1$. Hence, as $d(a,i)\le a-i$,
\[
d_X\bigl(p_i,\rho_-(a-i)\bigr)=d(a,i)+(a-i)-2\bigl(p_i\big|\rho_-(a-i)\bigr)_{p_a}\le2B+2.
\]
For $0\le x\le a'$ let $i$ be an integer nearest to $x$; then $0\le i\le a'$, $\lvert(a'-x)-(a-i)\rvert=\lvert H+x-i\rvert\le H+\frac12$, and
\[
\begin{aligned}
d_X\bigl(\widetilde p(x),p(x)\bigr)&\le d_X\bigl(\rho_-(a'-x),\rho_-(a-i)\bigr)+d_X\bigl(\rho_-(a-i),p_i\bigr)+d_X\bigl(p_i,p(x)\bigr)\\
&\le\bigl(H+\tfrac12\bigr)+(2B+2)+\tfrac12,
\end{aligned}
\]
which is below $2H$ because $H\ge100B/\beta_{\rm br}>2B+3$. The right side is treated alike: if it is not fast, $\widetilde p(y)=p_R$ for $y\ge b'$; if it is fast, $\ell_+=d(b,R)$, $\rho_+$ runs at unit speed from $p_b$ to $p_R$ and then along the right half-line, $\widetilde p(y)=\rho_+(y-b')$ for $y\ge b'$, and $d_X(\widetilde p(y),p(y))<2H$ for $b'\le y\le R$.

\textbf{The profile.} The path $\widetilde p$ is continuous and $1$-Lipschitz on each of $(-\infty,a']$, $[a',b']$ and $[b',\infty)$, hence on $\mathbb{R}$. So $d_*$ is an array with $d_*(x,y)\le\lvert x-y\rvert$ for all $x,y\in\mathbb{R}$, and $(\mathbb{R},d_*)$ is $0$-hyperbolic, its Gromov products being those of the real tree $X'$ at the image points. On a fast side, $\widetilde p$ is a unit-speed geodesic ray, so that side, $(-\infty,a']$ or $[b',\infty)$, is a unit geodesic ray for $d_*$. Let $I_0$ be the union of $[a',b']$ with the sides that are not fast. For $x<y$ in $I_0$ put $x^\circ=\min(\max(x,a'),b')$ and $y^\circ=\min(\max(y,a'),b')$; then $\widetilde p(x^\circ)=\widetilde p(x)$ and $\widetilde p(y^\circ)=\widetilde p(y)$, since $\widetilde p$ is constant on a side that is not fast, and $a'\le x^\circ\le y^\circ\le b'$ with $y^\circ-x^\circ\le y-x$, so $d_*(x,y)=d_*(x^\circ,y^\circ)\le(1-b_{\rm sl})(y-x)$. Thus $I_0$ is slow at speed $1-b_{\rm sl}$. Hence $d_*$ is a slow profile if no side is fast ($I_0=\mathbb{R}$); a one-ray profile if exactly one side is fast ($\mathbb{R}$ is the union of that side, a unit geodesic ray, and the closed half-line $I_0$, whose interior is the complementary half-line); and, if both sides are fast, $(-\infty,a']$ and $[b',\infty)$ are unit geodesic rays and $[a',b']$, of length $L_*$, is slow, so $d_*$ is a two-ray profile with slow length $L_*$ provided the deficit of $d_*$ on every interval is at most $C_{\rm cr}L_*$ (Definition~\ref{def:slac_prof}); this is proved below. Every point of $[0,R]$ lies in $[a',b']$ or in a fast side, since a side that is not fast lies outside $[0,R]$ (for instance $a'=-H$ if $a_*=0$). So ($\dagger$), ($\ddagger$) and its right-hand analogue give, for $0\le i,j\le R$,
\[
\bigl\lvert d(i,j)-d_*(i,j)\bigr\rvert\le d_X\bigl(p_i,\widetilde p(i)\bigr)+d_X\bigl(p_j,\widetilde p(j)\bigr)<6H.
\]

\textbf{The crossing deficit.} Let both sides be fast, so that $a_*=a$, $b_*=b$ and $L_0=b-a$. We prove that the deficit of $d_*$ on every interval is at most $C_{\rm cr}L_*$. For $x\le a'-\ell_-$ and $y\ge b'+\ell_+$, $\widetilde p(x)$ is the point at distance $a'-x-\ell_-$ from $p_0$ on the left half-line and $\widetilde p(y)$ the point at distance $y-b'-\ell_+$ from $p_R$ on the right half-line, so that
\[
d_*(x,y)=(a'-x-\ell_-)+d(0,R)+(y-b'-\ell_+),
\]
and the deficit $(y-x)-d_*(x,y)$ of $d_*$ on $[x,y]$ equals $\Delta_*=L_*+\ell_-+\ell_+-d(0,R)$. The deficit of an interval does not decrease when the interval grows: for $x'<x<y$, by the triangle inequality and $d_*(x',x)\le x-x'$,
\[
\bigl((y-x')-d_*(x',y)\bigr)-\bigl((y-x)-d_*(x,y)\bigr)\ge(x-x')-d_*(x',x)\ge0,
\]
and likewise at the right end. Every interval lies in some $[x,y]$ with $x\le a'-\ell_-$ and $y\ge b'+\ell_+$, so it suffices to show that $\Delta_*\le C_{\rm cr}L_*$.

Put $D_0=d(a,b)\le L_0$ and $\lambda=\frac12\bigl(\ell_-+\ell_+-D_0-d(0,R)\bigr)$, so that $\Delta_*=L_*+D_0+2\lambda$. By the triangle inequality $\ell_+=d(b,R)\le D_0+\ell_-+d(0,R)$, hence $\lambda\le\ell_-$; likewise $\lambda\le\ell_+$. As $D_0\le L_*$, it suffices to show that $2\lambda\le(C_{\rm cr}-2)L_*$. Suppose, to the contrary, that $2\lambda>(C_{\rm cr}-2)L_*$, and put
\[
\Gamma=3L_0+2B+3,\qquad n_1=\lceil2\Gamma/\delta_{\rm can}\rceil.
\]
Since $6H>2B+3$, $\Gamma\le3L_*$; and since $L_*\ge2$, $n_1+1\le6L_*/\delta_{\rm can}+2\le7L_*/\delta_{\rm can}<(C_{\rm cr}-2)L_*/2<\lambda$. So $n_1\le\lambda-1$, and as $\lambda\le\ell_-\le a$ and $\lambda\le\ell_+\le R-b$, we get $a-n_1\ge1$ and $b+n_1+1\le R$. Let
\[
r=a+\lceil L_0/2\rceil,\qquad t=n_1+\lceil L_0/2\rceil,\qquad n_1'=n_1+2\lceil L_0/2\rceil-L_0\in\{n_1,n_1+1\},
\]
so that $r-t=a-n_1$ and $r+t=b+n_1'$ lie in $\{0,\dots,R\}$, and $a\le r\le b$.

We bound three Gromov products of $d$ at the basepoint $a$; they are those of $d_X$ at $p_a$ of the corresponding points $p_i$. First, by ($\sharp$) with $i=a-n_1$, $(a-n_1|0)_a\ge n_1-B-1$. Second, since $d(a,R)\ge d(b,R)-D_0=\ell_+-D_0$,
\[
2(0|R)_a=\ell_-+d(a,R)-d(0,R)\ge2\lambda>2n_1.
\]
Third, for $b\le j\le R$, part (b) of Lemma~\ref{lem:sha_hyp_pat} and $d(j,R)\le R-j$ give $2(j|R)_b=d(b,j)+\ell_+-d(j,R)\ge2(j-b)-2B-2$; and since $d(a,u)\ge d(b,u)-D_0$ for every $u$, $2(j|R)_a\ge2(j|R)_b-2D_0$. With $j=b+n_1'$ this gives $(b+n_1'|R)_a\ge n_1-B-1-D_0$. Two applications of $0$-hyperbolicity at $a$ now give
\[
(a-n_1|b+n_1')_a\ge\min\bigl\{(a-n_1|0)_a,\,(0|R)_a,\,(R|b+n_1')_a\bigr\}\ge n_1-B-1-D_0.
\]
Hence, as $d(a,a-n_1)\le n_1$ and $d(a,b+n_1')\le D_0+n_1'$,
\[
d(r-t,r+t)=d(a,a-n_1)+d(a,b+n_1')-2(a-n_1|b+n_1')_a\le3D_0+2B+2+(n_1'-n_1)\le\Gamma.
\]
On the wings, part (b) of Lemma~\ref{lem:sha_hyp_pat} and steps of length at most one give
\[
d(r,r-t)\ge d(a,a-n_1)-d(a,r)\ge(n_1-B-1)-\lceil L_0/2\rceil\ge t-\Gamma,
\]
\[
d(r,r+t)\ge d(b,b+n_1')-d(b,r)\ge(n_1'-B-1)-\lfloor L_0/2\rfloor\ge t-\Gamma,
\]
because $2\lceil L_0/2\rceil\le L_0+1$, $\lceil L_0/2\rceil+\lfloor L_0/2\rfloor=L_0$, $n_1'\ge n_1$ and $\Gamma\ge L_0+B+2$. As $t\ge n_1\ge2\Gamma/\delta_{\rm can}$, we have $\Gamma\le(\delta_{\rm can}/2)t$; so both wings are at least $(1-\delta_{\rm can}/2)t$ and $d(r-t,r+t)\le(\delta_{\rm can}/2)t$. That is, $d$ has a $\delta_{\rm can}/2$-cancellation of length $t$ at $r$ (Definition~\ref{def:can_the_gro}), contrary to the hypothesis. Hence $2\lambda\le(C_{\rm cr}-2)L_*$ and $\Delta_*\le C_{\rm cr}L_*$, and $d_*$ is a two-ray slack profile.

In every case $d_*\in\mathscr{S}$ and $\lvert d(i,j)-d_*(i,j)\rvert\le\max(6H,2R^2\zeta)$ for $0\le i,j\le R$.

Explicit values: $c+\gamma+3B/H\le1-\beta_{\rm br}+\beta_{\rm br}/16+3\beta_{\rm br}/100=1-0.9075\beta_{\rm br}<1-0.75\beta_{\rm br}=1-b_{\rm sl}$; $H\ge100B/\beta_{\rm br}=240000B>2B+3$; $7/\delta_{\rm can}<(C_{\rm cr}-2)/2$, as $C_{\rm cr}=10^6/\delta_{\rm can}$, with a factor about $7\cdot10^4$ to spare.
\end{proof}

The next lemma applies Lemma~\ref{lem:app_fin_tre} to a rescaled growth array sampled on a fine grid, after replacing the samples by a nearby semimetric that is $0$-hyperbolic at the first grid point, as in Gromov's approximation of finite hyperbolic spaces by trees \cite{locbib:GhysHarpe1990}. The two quantitative inputs of the paper enter the approximation here: the scarcity of cancellations Theorem~\ref{thm:exp_sca_can} and the bridge product estimate Theorem~\ref{thm:bri_pro_est:2} bound the mass of the exceptional set.

\begin{lemma}
 Let $M\ge M_{\rm an}$, and let $n_0<N$ be of the form $2^j$ with $n_0\ge n_{\rm bal}$. For every $n=2^j$ with $n_0\le n<N$ there is a Borel set $\mathrm{Exc}_n\subset\mathbb{T}^2$ such that, for every initial point outside $\mathrm{Exc}_n$, some $d'\in\mathscr{S}$ satisfies
\[
\lVert g_n-d'\rVert_{[-2048,2048]}\le2^{-5572};
\]
and, summing over $n=2^j$ with $n_0\le n<N$,
\[
\sum_n\mathbb{P}(\mathrm{Exc}_n)\le2^{45767}e_N^2+2^{-4000}e_N+\frac{2^{-26000}}{n_0}.
\]
\label{lem:sum_exc_mas:1}
\end{lemma}

\begin{proof}
 We sample $g_n$ on a fine grid of $[-2048,2048]$, approximate the grid by a semimetric that is $0$-hyperbolic at a basepoint, and apply Lemma~\ref{lem:app_fin_tre} to it. That lemma fails only through a forbidden cancellation or a slow-pair violation joining a nonaffine marker to a nonunit one. The scarcity of cancellations (Theorem~\ref{thm:exp_sca_can}) controls the first kind of failure; the bridge product estimate (Theorem~\ref{thm:bri_pro_est:2}), applied to the summable weights of Lemma~\ref{lem:sum_mid_can}, controls the second.

\textbf{Parameters.} Throughout, $s=2^{-5600}$, $R=2^{5612}$, $\Lambda=\{0,\dots,R\}$, $x_i=-2048+is$ for $i\in\Lambda$ (so that $x_R=2048$), $B=64$, $H=\lceil100B/\beta_{\rm br}\rceil=15360000$, $\zeta=2^{-11500}$ and $c=1-\beta_{\rm br}$. Fix such an $n$, so that $n_0\le n<N$, and put $\rho_n=\lceil\log_2R\rceil/n=5612/n$. Since $n\ge n_0\ge n_{\rm bal}=2^{49000}$, the numbers $n_0s$, $ns$ and $nx_i$ are integers and $ns=2^{-5600}n$. The proof uses the following numerical facts, named by their roles.
\begin{itemize}
\item Tree error: $\rho_n\le\zeta s/100$, which holds as soon as $n_0\ge2^{17120}$.
\item Bridge length: $Bsn_0\ge2^{43406}>1000+100\log(8/\zeta)$.
\item Bridge window: $\kappa_{\rm br}B=16/5\ge1$.
\item Tolerance: $\zeta\le\min\bigl(\delta_{\rm can},\beta_{\rm br}/(64R)\bigr)$.
\item Cancellation length: $\gamma_{\rm can}sn_0>2^{43378}$.
\end{itemize}

\textbf{The comparison semimetric.} We first approximate the grid by a semimetric that is $0$-hyperbolic at the basepoint $0$. Write $d(i,j)=g_n(x_i,x_j)=g(nx_i,nx_j)/n$ for $i,j\in\Lambda$ (Lemma~\ref{lem:dil_and_tim}(i)). It is a semimetric on $\Lambda$ by Lemma~\ref{lem:nor_and_gro}(ii), and its Gromov product of $i,j$ at $w$ is that of $g$ of $nx_i,nx_j$ at $nx_w$, divided by $n$; so by Lemma~\ref{lem:hyp_the_gro}, $(\Lambda,d)$ is $\delta$-hyperbolic at every basepoint (Definition~\ref{def:gro_hyp_sem}), with $\delta=(\log_M2)/(2n)$. Write $(i|j)_0$ for the Gromov products of $d$ at the basepoint $0$ and $\lvert i\rvert_0=d(0,i)$, put $\beta(i,i)=\lvert i\rvert_0$, and for $i\ne j$ let $\beta(i,j)$ be the largest, over chains $i=i_0,i_1,\dots,i_p=j$ of distinct points of $\Lambda$, of $\min_{1\le l\le p}(i_{l-1}|i_l)_0$. Then:
\begin{enumerate}
\item[(a)] $(i|j)_0\le\beta(i,j)\le\min(\lvert i\rvert_0,\lvert j\rvert_0)$, since a one-step chain is admitted and every Gromov product satisfies $(u|u')_0\le\min(\lvert u\rvert_0,\lvert u'\rvert_0)$;
\item[(b)] $\beta(i,j')\ge\min\{\beta(i,j),\beta(j,j')\}$: concatenate two optimal chains and delete the loop between two visits of a repeated point, which does not decrease the minimum;
\item[(c)] along a chain $i_0,\dots,i_p$ with $p\le2^\ell$, $(i_0|i_p)_0\ge\min_l(i_{l-1}|i_l)_0-\ell\delta$, by induction on $\ell$, splitting the chain at a middle point and applying $\delta$-hyperbolicity once.
\end{enumerate}
A chain of distinct points has at most $R\le2^{\lceil\log_2R\rceil}$ steps, so (c) gives $\beta(i,j)\le(i|j)_0+\lceil\log_2R\rceil\delta$. Put $d^{\rm tr}(i,j)=\lvert i\rvert_0+\lvert j\rvert_0-2\beta(i,j)$. By (a) and (c), and as $\log_M2\le1$,
\[
d^{\rm tr}\le d\le d^{\rm tr}+2\lceil\log_2R\rceil\delta\le d^{\rm tr}+\rho_n.
\]
By (a) and (b), $d^{\rm tr}$ is a semimetric on $\Lambda$ with $d^{\rm tr}(0,i)=\lvert i\rvert_0$ (a chain from $0$ starts with a product $(0|i_1)_0=0$, so $\beta(0,i)=0$); its Gromov products at $0$ are the $\beta(i,j)$, so it is $0$-hyperbolic at $0$ by (b); and its triangle inequality reads $\beta(i,j)+\beta(j,j')\le\lvert j\rvert_0+\beta(i,j')$, which follows from (a) and (b). This is Gromov's approximation of finite hyperbolic spaces by trees (\cite{locbib:GhysHarpe1990}), here written out. Put $\widehat d=d^{\rm tr}/s$, a semimetric on $\Lambda$ that is $0$-hyperbolic at $0$; its steps satisfy $\widehat d(i,i+1)\le d(i,i+1)/s\le1$, because $g(nx_i,nx_{i+1})\le ns$ by Lemma~\ref{lem:nor_and_gro}(ii). Each $\beta(i,j)$ is a finite maximum of minima of continuous functions of the initial point (Lemma~\ref{lem:dil_and_tim}(v)), so every set defined below by finitely many inequalities between such functions is Borel.

\textbf{The exceptional set.} Let $P$ and $Q$ be the markers of $\widehat d$ at tolerance $\zeta$ (Definition~\ref{def:fin_tre_pat:1}). Let $\mathrm{Exc}_n$ be the set of initial points at which $\widehat d$ does not have slow marked pairs at range $B$ for $c$, or has a $\delta_{\rm can}/2$-cancellation (Definition~\ref{def:can_the_gro}, on the set $\Lambda$); it is Borel. Outside $\mathrm{Exc}_n$, Lemma~\ref{lem:app_fin_tre} applies to $\widehat d$ (by the tolerance bound, $\zeta\le\beta_{\rm br}/(64R)$) and gives $d_*\in\mathscr{S}$ with $\lvert\widehat d(i,j)-d_*(i,j)\rvert\le6H+2R^2\zeta$; hence, for $i,j\in\Lambda$,
\[
\bigl\lvert g_n(x_i,x_j)-sd_*(i,j)\bigr\rvert\le\bigl(d-d^{\rm tr}\bigr)(i,j)+s\bigl\lvert\widehat d(i,j)-d_*(i,j)\bigr\rvert\le\rho_n+s\bigl(6H+2R^2\zeta\bigr).
\]
Put $d'(x,y)=sd_*\bigl((x+2048)/s,(y+2048)/s\bigr)$. It is obtained from $d_*$ by $5600$ dilations, which give $(x,y)\mapsto sd_*(x/s,y/s)$, followed by the translation $\mathcal{T}_{2048}$ (Definition~\ref{def:arr_mid_can}); so $d'\in\mathscr{S}$ by Lemma~\ref{lem:geo_sla_pro}(iii), and $d'(x_i,x_j)=sd_*(i,j)$. Let $x,y\in[-2048,2048]$, and choose grid points $x_i,x_j$ with $\lvert x-x_i\rvert\le s$ and $\lvert y-x_j\rvert\le s$. As $g_n$ is a semimetric with $g_n(u,u')\le\lvert u-u'\rvert+1/n$ (Lemma~\ref{lem:dil_and_tim}(i)), $\lvert g_n(x,y)-g_n(x_i,x_j)\rvert\le g_n(x,x_i)+g_n(y,x_j)\le2s+2/n$; and as $d'(u,u')\le\lvert u-u'\rvert$ for all $u,u'$, Lemma~\ref{lem:shor_arra}(iv) gives $\lvert d'(x,y)-d'(x_i,x_j)\rvert\le2s$. Hence
\[
\lVert g_n-d'\rVert_{[-2048,2048]}\le\rho_n+\frac2n+s\bigl(6H+2R^2\zeta+4\bigr)\le\frac{5614}{n_0}+s\bigl(6H+2R^2\zeta+4\bigr)<2^{-5572}.
\]
It remains to bound $\sum_n\mathbb{P}(\mathrm{Exc}_n)$.

\textbf{The weights at the markers.} We show that a nonaffine marker carries a weight at least $\zeta^2/2$, and a nonunit marker a shortfall above $\zeta/8$. For $1\le i\le R-1$ let $I=[x_{i-1},x_{i+1}]$, with halves $I_l=[x_{i-1},x_i]$ and $I_r=[x_i,x_{i+1}]$, and put
\[
\begin{gathered}
W_i=(J+V)(I)(g_n)=F\bigl(h(I)(g_n),h(I_l)(g_n),h(I_r)(g_n)\bigr),\\
F(u_1,u_2,u_3)=u_1-\tfrac12(u_2+u_3)+(u_2-u_3)^2,
\end{gathered}
\]
the second form by Lemma~\ref{lem:shor_arra}(i). By Definition~\ref{def:arr_mid_can}, $J(I)(g_n)$ is the Gromov product $(i-1|i+1)_i$ of $d$ divided by $s$, and $h(I_l)(g_n)-h(I_r)(g_n)=\bigl(d(i,i+1)-d(i-1,i)\bigr)/s$. Since $0\le d-d^{\rm tr}\le\rho_n$ and $\rho_n/s\le\zeta/100$, these differ from the Gromov product $(i-1|i+1)_i$ of $\widehat d$ and from the step difference $\widehat d(i,i+1)-\widehat d(i-1,i)$ by at most $3\zeta/200$ and $\zeta/100$ respectively. So if $i\in P$, then $J(I)(g_n)\ge0.98\zeta$ or $V(I)(g_n)\ge(0.98\zeta)^2$, and as $J(I)$ and $V(I)$ are nonnegative (Lemma~\ref{lem:shor_arra}(i)), $W_i\ge\zeta^2/2$ in either case. For $j\in\Lambda$ let $Y_j$ be the largest of the shortfalls of $g_n$ on those of the intervals $[x_{j-1},x_j]$, $[x_j,x_{j+1}]$, $[x_{j-1},x_{j+1}]$ whose endpoints are grid points. If $j\in Q$, then $Y_j>\zeta/8$. Indeed, if a step of $\widehat d$ at $j$, say $\widehat d(j,j+1)$, is below $1-\zeta/2$, then the shortfall of $[x_j,x_{j+1}]$ is $1-d(j,j+1)/s>\zeta/2-\zeta/100$. If $1\le j\le R-1$ and the Gromov product $(j-1|j+1)_j$ of $\widehat d$ exceeds $\zeta/2$, then $\widehat d(j-1,j+1)<2-\zeta$, as the steps are at most one, and the shortfall of $[x_{j-1},x_{j+1}]$ is $1-d(j-1,j+1)/(2s)>\zeta/2-\zeta/200$.

By Lemma~\ref{lem:dil_and_tim}(i), the shortfall of $g_n$ on an interval with grid endpoints $x<x'$ is $h_g(nx,nx')$. So $W_i$ and $Y_j$ are shortfall observables (Definition~\ref{def:shor_obse}). The observable $W_i\ge0$ has Lipschitz constant at most $5$, since $\lvert\partial_1F\rvert+\lvert\partial_2F\rvert+\lvert\partial_3F\rvert=1+\max\bigl(1,4\lvert u_2-u_3\rvert\bigr)\le5$ on $[0,1]^3$; the observable $Y_j$, a maximum, has Lipschitz constant $1$. The indices occurring in $W_i$ lie within $ns$ of $nx_i$, and those occurring in $Y_j$ within $ns$ of $nx_j$.

\textbf{Slow-pair violations.} At an initial point where $\widehat d$ does not have slow marked pairs at range $B$ for $c$, there are $i\in P$ and $j\in Q$ with $\lvert i-j\rvert\ge B$ and $\widehat d(i,j)>c\lvert i-j\rvert$. For an ordered pair $(i,j)$ of points of $\Lambda$ with $1\le i\le R-1$ and $\lvert i-j\rvert\ge B$, let $\mathrm{Viol}_{i,j}$ be the Borel set of initial points at which $i\in P$, $\widehat d(i,j)>c\lvert i-j\rvert$ and $Y_j>\zeta/8$; by the previous paragraph, every such failure lies in one of these sets. Let $a<b$ be the integers with $\{a,b\}=\{nx_i,nx_j\}$; then $b-a=ns\lvert i-j\rvert\ge Bsn\ge Bsn_0$. On $\mathrm{Viol}_{i,j}$, since $d\ge d^{\rm tr}=s\widehat d$, we have $g(a,b)=nd(i,j)>(1-\beta_{\rm br})(b-a)$, that is $h_g(a,b)<\beta_{\rm br}$, and $W_i\ge\zeta^2/2$. By the bridge length, $b-a\ge1000+100\log(8/\zeta)$, and by the bridge window, $ns\le\kappa_{\rm br}(b-a)$. So Theorem~\ref{thm:bri_pro_est:2}, with threshold $\zeta/8$, the weight $W_i$ at one end and the observable $Y_j$ at the other, gives
\[
\mathbb{P}(\mathrm{Viol}_{i,j})\le\frac2{\zeta^2}\,\mathbb{E}\Bigl[\mathbf{1}_{\{h_g(a,b)\le\beta_{\rm br}\}}W_i\,\mathbf{1}_{\{Y_j>\zeta/8\}}\Bigr]\le\frac{2C_{\rm br}}{\zeta^2}\Bigl(\mathbb{E}\,W_i+M^{-\sigma_{\rm br}Bsn}\Bigr)\,\mathbb{P}(Y_j>\zeta/16).
\]
Each shortfall entering $Y_j$ is $h(\Theta)(g_n)$ with $\Theta=[x,x+s']$, $x$ a grid point and $s'\in\{s,2s\}$. As $nx$ is an integer, Lemma~\ref{lem:dil_and_tim}(iii) shows that its mean equals the mean of $h([0,s'])(g_n)=h_g(0,ns')$, which is $e_{ns}$ or $e_{2ns}$; these are at most $e_N$ by Lemma~\ref{lem:exa_dou_ide}, because $N$ is obtained from $ns$ and from $2ns$ by doublings. So $\mathbb{E}\,Y_j\le3e_N$, and Markov's inequality gives $\mathbb{P}(Y_j>\zeta/16)\le48e_N/\zeta$. The interval $[x_{i-1},x_{i+1}]$ has length $2s=2^{-5599}$, at most $2$, and is aligned at scale $n_0$, since its endpoints and midpoint lie in $s\mathbb{Z}$ and $n_0s$ is an integer; so Lemma~\ref{lem:sum_mid_can} gives $\sum_n\mathbb{E}\,W_i\le9e_N$. Moreover $\sum_nM^{-\sigma_{\rm br}Bsn}\le2M^{-\sigma_{\rm br}Bsn_0}$, a sum over the scales $n$, each twice the previous one, whose first term is at most $\frac12$, as $M\ge M_{\rm an}$ and $\sigma_{\rm br}Bsn_0\ge100$. Summing over the scales, for each ordered pair,
\[
\sum_n\mathbb{P}(\mathrm{Viol}_{i,j})\le\frac{96C_{\rm br}}{\zeta^3}e_N\bigl(9e_N+2M^{-\sigma_{\rm br}Bsn_0}\bigr)\le\frac{864C_{\rm br}}{\zeta^3}e_N\bigl(e_N+M^{-\sigma_{\rm br}Bsn_0}\bigr).
\]
There are at most $(R+1)^2$ ordered pairs, and $864C_{\rm br}(R+1)^2\zeta^{-3}<2^{45767}$. Finally $2^{45767}M^{-\sigma_{\rm br}Bsn_0}\le2^{-4000}$, as $\sigma_{\rm br}Bsn_0\ge2^{43402}$ and $\log_2M\ge1$. So the slow-pair violations contribute at most $2^{45767}e_N^2+2^{-4000}e_N$ to $\sum_n\mathbb{P}(\mathrm{Exc}_n)$.

\textbf{Forbidden cancellations.} Suppose that $\widehat d$ has a $\delta_{\rm can}/2$-cancellation of length $t$ at $r$, where $r-t,r,r+t\in\Lambda$, so that $t\in\{1,\dots,R\}$. Then $d^{\rm tr}$ has wings at least $(1-\delta_{\rm can}/2)ts$ and span at most $(\delta_{\rm can}/2)ts$ there. Since $d^{\rm tr}\le d\le d^{\rm tr}+\rho_n$ and $\rho_n\le\zeta s/100\le\delta_{\rm can}ts/100$ (tree error and tolerance), the wings of $d$ are at least $(1-\delta_{\rm can})ts$ and its span is at most $\delta_{\rm can}ts$. Multiplying by $n$, the growth semimetric $g$ has a $\delta_{\rm can}$-cancellation of length $nts$ at the time $nx_r$, so the initial point lies in $\mathrm{Can}_{\delta_{\rm can}}(nx_r,nts)$, whose measure is at most $C_{\rm can}\exp(-\gamma_{\rm can}nts)\le C_{\rm can}\exp(-\gamma_{\rm can}sn)$ by Theorem~\ref{thm:exp_sca_can}. There are at most $(R+1)^2$ pairs $(r,t)$, and $\sum_n\exp(-\gamma_{\rm can}sn)\le2\exp(-\gamma_{\rm can}sn_0)$, the first term being at most $\frac12$ by the cancellation length. Since $\exp(y)\ge y^2/2$ for $y\ge0$, and $2C_{\rm can}<2^{16}$, the cancellations contribute at most
\[
2C_{\rm can}(R+1)^2\exp(-\gamma_{\rm can}sn_0)\le\frac{2^{17}(R+1)^2}{(\gamma_{\rm can}s)^2\,n_0^2}<\frac{2^{22484}}{n_0^2}\le\frac{2^{-26000}}{n_0},
\]
the last step as $n_0\ge2^{49000}$. Adding the two contributions proves the bound on $\sum_n\mathbb{P}(\mathrm{Exc}_n)$.

Explicit values: $H=15360000$; profile error: $6Hs=2^{-5573.54}$, $2R^2\zeta s=2^{-5875}$, $4s=2^{-5598}$, $5614/n_0\le2^{-48987}$, total below $2^{-5573.5}$; the tree error needs $n_0\ge2^{17119.1}$; $Bsn_0\ge2^{43406}$ against $1000+100\log(8/\zeta)<8\cdot10^5$; $\sigma_{\rm br}Bsn_0\ge2^{43402.6}$; $\gamma_{\rm can}sn_0\ge2^{43378.8}$; $864C_{\rm br}=2^{41.76}$ and $2^{42}(R+1)^2\zeta^{-3}=2^{45766.0}$; $2C_{\rm can}=40000$ and $2^{17}(R+1)^2(\gamma_{\rm can}s)^{-2}=2^{22483.39}$; $1-3/200\ge0.98$, $1-1/100\ge0.98$, $0.98^2\ge1/2$, $1/2-1/100>1/8$.
\end{proof}

\subsection{The defect bound} \label{subsec:the_bal_def}

This subsection proves the defect bound Proposition~\ref{prop:def_the_ext:1}, the remaining input of the balance. At each scale its proof separates the body, the set of initial points at which some midpoint cancellation or imbalance of a fixed dyadic family of intervals is not tiny, from its complement. Most of the work lies off the body, where the arrays are close to affine arrays (Lemma~\ref{lem:aff_app_off}) and Lemma~\ref{lem:bal_def_off} bounds the defect at slow, near-unit and middle speeds.

\begin{lemma}
 Let $\mathcal{I}$ be the family of the intervals $[-2048+i\,2^{12-l},\,-2048+(i+1)\,2^{12-l}]$ with integers $0\le l<1112$ and $0\le i<2^l$, the internal nodes of the binary subdivision of $[-2048,2048]$ down to length $2^{-1100}$. Let $n=2^j\ge2^{1100}$, and consider an initial point at which $J(I)(g_n)<2^{-2300}$ and $V(I)(g_n)<2^{-2300}$ for every $I\in\mathcal{I}$. Let $\sigma$ be the speed of the secant $[-2048,2048]$ of $g_n$. Then $0\le\sigma\le1$ and
\[
\bigl\lvert g_n(x,y)-\sigma\lvert x-y\rvert\bigr\rvert\le2^{-1020}\qquad(x,y\in[-2048,2048]).
\]
\label{lem:aff_app_off}
\end{lemma}

\begin{proof}
Halving a node changes the speed little when its midpoint cancellation and its imbalance are small. Along the $1112$ levels of the subdivision this gives every leaf nearly the speed of the root, and the triangle inequality, applied along the leaves for an upper bound and across the root interval for a lower bound, controls every distance between grid points.

Write $d=g_n$ and call the points of $2^{-1100}\mathbb{Z}\cap[-2048,2048]$ grid points. As $n2^{-1100}$ is an integer, Lemma~\ref{lem:dil_and_tim}(i) gives $h([s,t])(d)\in[0,1]$ for grid points $s<t$, so the speed $\sigma([s,t])=1-h([s,t])(d)$ lies in $[0,1]$; in particular $\sigma=\sigma([-2048,2048])\in[0,1]$. Every $I\in\mathcal{I}$ has length $2^{12-l}\ge2^{-1099}$, so its endpoints and its midpoint are grid points.

We first compare the speed of a node with the speeds of its halves. By Lemma~\ref{lem:shor_arra}(i), for $I\in\mathcal{I}$ the speeds of both halves of $I$ differ from the speed of $I$ by at most $J(I)+\frac12\sqrt{V(I)}<2^{-2300}+2^{-1151}<2^{-1150}$.

We next compare the leaves with the root. The leaves are the intervals $\Theta=[s,s+2^{-1100}]$ with $s$ a grid point below $2048$; each is obtained from $[-2048,2048]$ by $1112$ successive halvings, through members of $\mathcal{I}$. Hence $\lvert\sigma(\Theta)-\sigma\rvert<\varepsilon$ for every leaf, where $\varepsilon=1112\cdot2^{-1150}$.

We then bound the distances between grid points. Let $s<t$ be grid points. The interval $[s,t]$ is the union of consecutive leaves, so the triangle inequality gives $d(s,t)\le(\sigma+\varepsilon)(t-s)$, and likewise $d(-2048,s)\le(\sigma+\varepsilon)(s+2048)$ and $d(t,2048)\le(\sigma+\varepsilon)(2048-t)$. The triangle inequality across the root interval gives
\[
4096\,\sigma=d(-2048,2048)\le d(-2048,s)+d(s,t)+d(t,2048)\le(\sigma+\varepsilon)\bigl(4096-(t-s)\bigr)+d(s,t),
\]
so $d(s,t)\ge\sigma(t-s)-4096\varepsilon$. Together, $\lvert d(s,t)-\sigma(t-s)\rvert\le4096\varepsilon<2^{-1127}$.

We finally pass to real times. For $x\in[-2048,2048]$ let $x'=2^{-1100}\lfloor2^{1100}x\rfloor$, a grid point with $0\le x-x'<2^{-1100}$, and define $y'$ likewise. By the triangle inequality and Lemma~\ref{lem:dil_and_tim}(i), $\lvert d(x,y)-d(x',y')\rvert\le d(x,x')+d(y,y')\le2(2^{-1100}+1/n)$, and $\bigl\lvert\sigma\lvert x-y\rvert-\sigma\lvert x'-y'\rvert\bigr\rvert\le2\cdot2^{-1100}$. As $n\ge2^{1100}$,
\[
\bigl\lvert d(x,y)-\sigma\lvert x-y\rvert\bigr\rvert<2^{-1127}+6\cdot2^{-1100}<2^{-1020}.
\]

Explicit values: $4096\cdot1112\cdot2^{-1150}=2^{-1127.88}$ and $2^{-1127}+6\cdot2^{-1100}<2^{-1097}$.
\end{proof}

\begin{lemma}
 Let $M\ge M_{\rm an}$, and let $n_0<N$ be of the form $2^j$ with $n_0\ge2^{1100}$. Let $\mathcal{I}$ be the family of intervals of Lemma~\ref{lem:aff_app_off}, and for every $n=2^j$ with $n_0\le n<N$ let $\mathrm{Off}_n\subset\mathbb{T}^2$ be the set of initial points at which $J(I)(g_n)<2^{-2300}$ and $V(I)(g_n)<2^{-2300}$ for every $I\in\mathcal{I}$. Then each $\mathrm{Off}_n$ is Borel, and, summing over $n=2^j$ with $n_0\le n<N$,
\[
\sum_n\mathbb{E}\bigl[Z(g_n)\,\mathbf{1}_{\mathrm{Off}_n}\bigr]\le2^{-846}e_N+2^{1188}e_N^2+\frac{\alpha_{\rm cap}}{n_0}.
\]
\label{lem:bal_def_off}
\end{lemma}

\begin{proof}
 Off the body the array is affine up to an error $\varepsilon=2^{-1020}$ on $[-2048,2048]$, at the speed $\sigma_n$ of the secant $[-2048,2048]$ (Lemma~\ref{lem:aff_app_off}), and we distinguish three ranges of $\sigma_n$. At slow speeds every capped shortfall and every side factor saturates, and the defect vanishes. At near-unit speeds every shortfall in the balance functional lies where the cap is linear, so the balance functional is $\alpha_{\rm cap}$ times a midpoint cancellation of length two, which truncation of the wings (Lemma~\ref{lem:gro_pro_und}) reduces to $\zeta'$ times a midpoint cancellation of length $2\zeta'$; the correction terms are imbalances whose halves have small shortfalls; both are summable with a factor $\zeta'$. At middle speeds the time interval $[0,1024n]$ of the orbit has small shortfall and the unit interval at its right end a shortfall bounded below; the bridge product estimate decouples the weights near time zero, which bound the defect, from that right-end shortfall, whose mean is $e_n$.

\textbf{Setup.} Within this proof sums over $n$ run over the scales $n$ of the statement, $n_0\le n<N$, and $\zeta'=2^{-1000}$, $\varepsilon=2^{-1020}$. Fix such an $n$; it is even and at least $2^{1100}$, so $n\zeta'$ and $n2^{-1100}$ are integers. Every $I\in\mathcal{I}$ has its endpoints and midpoint in $2^{-1100}\mathbb{Z}$, so $J(I)(g_n)$ and $V(I)(g_n)$ are continuous functions of finitely many values of $g$, and hence of the initial point (Lemma~\ref{lem:dil_and_tim}(i),(v)); thus $\mathrm{Off}_n$ is open, hence Borel. The same holds for the speed $\sigma_n=g_n(-2048,2048)/4096$ of the secant $[-2048,2048]$, so the sets
\[
\begin{aligned}
\mathrm{Slow}_n&=\mathrm{Off}_n\cap\bigl\{\sigma_n\le1-\tfrac34b_{\rm sl}\bigr\},\\
\mathrm{Mid}_n&=\mathrm{Off}_n\cap\bigl\{1-\tfrac34b_{\rm sl}<\sigma_n<1-\tfrac12\zeta'\bigr\},\\
\mathrm{Near}_n&=\mathrm{Off}_n\cap\bigl\{\sigma_n\ge1-\tfrac12\zeta'\bigr\}
\end{aligned}
\]
are Borel and partition $\mathrm{Off}_n$. Write $d=g_n$. At an initial point of $\mathrm{Off}_n$, Lemma~\ref{lem:aff_app_off} gives $\bigl\lvert d(x,y)-\sigma_n\lvert x-y\rvert\bigr\rvert\le\varepsilon$ for $x,y\in[-2048,2048]$, and hence
\begin{equation}
\bigl\lvert h(\Theta)(d)-(1-\sigma_n)\bigr\rvert\le\frac{\varepsilon}{\lvert\Theta\rvert}\qquad(\Theta\subset[-2048,2048]\text{ an interval}).\tag{$\ast$}
\end{equation}
By Lemma~\ref{lem:dil_and_tim}(vi), $g_n$ and $\mathcal{D}g_n$ are admissible; so, by Lemma~\ref{lem:pro_the_bal}(i),(iv), $Z(g_n)$ is defined and is a continuous function of the initial point, $\lvert Z(g_n)\rvert\le3$, and
\begin{equation}
\bigl\lvert Z(g_n)-G_{\rm av}(g_n)\bigr\rvert\le C_{\rm dom}\sum_Ic_I\,V(I)(g_n),\tag{$\ast\ast$}
\end{equation}
where $I$ runs over the fixed finite family of intervals of Lemma~\ref{lem:pro_the_bal}(i), which lie in $[-8,12]$, have length $1$ or $2$ and are aligned at scale $n_0$ (their endpoints and midpoints lie in $\frac1{64}\mathbb{Z}$), and the coefficients $c_I\ge0$ satisfy $\sum_Ic_I<20$. For intervals $I$ with endpoints and midpoint in $\frac1n\mathbb{Z}$, the functions $J(I)(g_n)$ and $V(I)(g_n)$ of the initial point are continuous and nonnegative (Lemma~\ref{lem:dil_and_tim}(i),(v), Lemma~\ref{lem:shor_arra}(i)); so all integrals below are defined.

\textbf{The slow range.} We show that $Z(g_n)=0$ on $\mathrm{Slow}_n$. Let the initial point lie in $\mathrm{Slow}_n$. By $(\ast)$, every interval $\Theta\subset[-4,8]$ of length $1$ or $2$ has $h(\Theta)(d)\ge1-\sigma_n-\varepsilon\ge\frac34b_{\rm sl}-\varepsilon>\frac12b_{\rm sl}$. For $0\le x\le1$, Lemma~\ref{lem:shor_arra}(ii) gives $G_\Phi(\mathcal{T}_xd)=F_2(x;d)-\frac12F_1(x;d)-\frac12F_1(x+1;d)$, whose three shortfalls are those of the intervals $[x,x+2]$, $[x,x+1]$ and $[x+1,x+2]$ in $[0,3]$; as $\Phi=1$ on $[b_{\rm sl}/2,\infty)$ (Lemma~\ref{lem:pro_the_cap}), $G_\Phi(\mathcal{T}_xd)=1-\frac12-\frac12=0$, and so $G_{\rm av}(d)=0$. The side factor of an array $d'$ reads the shortfalls of $d'$ on $[-2,-1]$ and $[3,4]$ (Definition~\ref{def:ext_cor_and}). For $d'=\mathcal{D}d$ these are the shortfalls of $d$ on $[-4,-2]$ and $[6,8]$, and for $d'=\mathcal{T}_yd$ with $y\in\{0,\frac12,1,\frac32\}$ those on $[y-2,y-1]$ and $[y+3,y+4]$ (Lemma~\ref{lem:shor_arra}(ii)); these intervals lie in $[-4,8]$ and have length $1$ or $2$, so all these shortfalls are at least $b_{\rm sl}/4$. The five arrays $\mathcal{D}d$ and $\mathcal{T}_yd$ are admissible (Lemma~\ref{lem:dil_and_tim}(vi)), so their extended correctors vanish by Lemma~\ref{lem:ext_the_ave}(a). By the definition of $Z$, $Z(g_n)=0$.

\textbf{The near-unit range.} We show that
\[
\sum_n\mathbb{E}\bigl[Z(g_n)\,\mathbf{1}_{\mathrm{Near}_n}\bigr]\le\bigl(160\,C_{\rm dom}+\alpha_{\rm cap}\bigr)\zeta'e_N+\frac{\alpha_{\rm cap}}{n_0}.
\]
Let the initial point lie in $\mathrm{Near}_n$. For an interval $\Theta\subset[-2048,2048]$ of length at least $\frac12$, $(\ast)$ gives
\[
1-\sigma_n-\frac{\varepsilon}{\lvert\Theta\rvert}\le h(\Theta)(d)\le1-\sigma_n+\frac{\varepsilon}{\lvert\Theta\rvert}\le\frac{\zeta'}2+2\varepsilon<\zeta'<u_{\rm lo}.
\]

We first bound the correction terms. The halves of each interval $I$ in $(\ast\ast)$ have length at least $\frac12$, so their shortfalls are below $\zeta'$ at our initial point; hence $V(I)(g_n)\mathbf{1}_{\mathrm{Near}_n}\le\mathbf{1}_{\{\max(h(I_l)(g_n),h(I_r)(g_n))\le\zeta'\}}V(I)(g_n)$ everywhere. By $(\ast\ast)$ and Lemma~\ref{lem:sum_mid_can} with threshold $\zeta'$, for intervals of length at most $2$,
\[
\sum_n\mathbb{E}\Bigl[\mathbf{1}_{\mathrm{Near}_n}\bigl(Z(g_n)-G_{\rm av}(g_n)\bigr)\Bigr]\le C_{\rm dom}\sum_Ic_I\cdot8\zeta'e_N\le160\,C_{\rm dom}\zeta'e_N.
\]

We next bound the balance term. For $0\le x\le1$, the shortfalls of $[x,x+1]$, $[x+1,x+2]$ and $[x,x+2]$ for $d$ are below $u_{\rm lo}$, where $\Phi(u)=\alpha_{\rm cap}u$ (Definition~\ref{def:the_cap_and}). So Lemma~\ref{lem:shor_arra}(i),(ii) give $G_\Phi(\mathcal{T}_xd)=\alpha_{\rm cap}J([x,x+2])(d)$. Let $j=\lfloor nx\rfloor$; by Lemma~\ref{lem:shor_arra}(ii) and Lemma~\ref{lem:dil_and_tim}(iv), $J([x,x+2])(d)=J([j/n,j/n+2])(d)$. Put $i=j+n$ and $s=\zeta'n$, an integer with $1\le s\le n$. By Lemma~\ref{lem:dil_and_tim}(i), the Gromov product of $g$ satisfies $2(i-n|i+n)_i=2n\,J([j/n,j/n+2])(d)$, and the wing deficit of Lemma~\ref{lem:gro_pro_und} is $D=n\max\bigl(h([j/n,j/n+1])(d),\,h([j/n+1,j/n+2])(d)\bigr)$. The shortfall bounds and Lemma~\ref{lem:shor_arra}(i) give $D\le n(\zeta'/2+\varepsilon)$ and $J([j/n,j/n+2])(d)\le(1-\sigma_n+\varepsilon/2)-(1-\sigma_n-\varepsilon)=\frac32\varepsilon$, so
\[
2s-2D\ge n(\zeta'-2\varepsilon)>3n\varepsilon\ge2(i-n|i+n)_i,
\]
because $\zeta'>5\varepsilon$. Thus the minimum in \eqref{eq:trunc_gromov} is $2(i-n|i+n)_i$, and \eqref{eq:trunc_gromov} gives $(i-n|i+n)_i\le(i-s|i+s)_i+\frac12$. Since $(i-s|i+s)_i=s\,J([i/n-\zeta',i/n+\zeta'])(d)$ by Lemma~\ref{lem:dil_and_tim}(i), this reads
\[
J([j/n,j/n+2])(d)\le\zeta'J\bigl([j/n+1-\zeta',\,j/n+1+\zeta']\bigr)(d)+\frac1{2n}.
\]
Midpoint cancellations are nonnegative (Lemma~\ref{lem:shor_arra}(i)). Hence on $\mathrm{Near}_n$,
\[
0\le G_{\rm av}(g_n)\le\alpha_{\rm cap}\int_0^1\Bigl(\zeta'J\bigl([\tfrac{\lfloor nx\rfloor}n+1-\zeta',\,\tfrac{\lfloor nx\rfloor}n+1+\zeta']\bigr)(g_n)+\frac1{2n}\Bigr)dx,
\]
and the right side is nonnegative everywhere. By Lemma~\ref{lem:shor_arra}(ii), Lemma~\ref{lem:dil_and_tim}(iii) and Fubini--Tonelli, the integral of the right side over the initial point is $\alpha_{\rm cap}\bigl(\zeta'\,\mathbb{E}\,J([1-\zeta',1+\zeta'])(g_n)+1/(2n)\bigr)$. The interval $[1-\zeta',1+\zeta']$ has length $2\zeta'=2^{-999}$ and is aligned at scale $n_0$, as $n_0\zeta'$ is an integer; so Lemma~\ref{lem:sum_mid_can} and $\sum_n1/(2n)\le1/n_0$ give
\[
\sum_n\mathbb{E}\bigl[\mathbf{1}_{\mathrm{Near}_n}G_{\rm av}(g_n)\bigr]\le\alpha_{\rm cap}\Bigl(\zeta'e_N+\frac1{n_0}\Bigr).
\]
Adding the two bounds gives the claim.

\textbf{The middle range.} We show that
\[
\sum_n\mathbb{E}\bigl[\lvert Z(g_n)\rvert\,\mathbf{1}_{\mathrm{Mid}_n}\bigr]\le\bigl(2^{40}+20\,C_{\rm dom}\bigr)\frac{72\,C_{\rm br}}{\zeta'}\,e_N\bigl(e_N+M^{-100n_0}\bigr).
\]
Put $Y=h_g(1024n,1025n)$ and
\[
\mathrm{Br}_n=\bigl\{h_g(0,1024n)\le\beta_{\rm br}\bigr\}\cap\bigl\{Y>\zeta'/4\bigr\}.
\]
We first show $\mathrm{Mid}_n\subset\mathrm{Br}_n$. Let the initial point lie in $\mathrm{Mid}_n$. By Lemma~\ref{lem:dil_and_tim}(i) and $(\ast)$, $h_g(0,1024n)=h([0,1024])(d)\le1-\sigma_n+\varepsilon/1024<\frac34b_{\rm sl}+\varepsilon<\beta_{\rm br}$, since $\frac34b_{\rm sl}=\frac9{16}\beta_{\rm br}$ by the values of Convention~\ref{conv:not_and_con}; and $Y=h([1024,1025])(d)\ge1-\sigma_n-\varepsilon>\zeta'/2-\varepsilon>\zeta'/4$.

We next bound the defect by weights. Let $I$ be an interval whose shortfall and the shortfalls of whose halves lie in $[0,1]$, and put $\bar h=\frac12(h(I_l)+h(I_r))$, so that $h(I)-\bar h=J(I)\ge0$ (Lemma~\ref{lem:shor_arra}(i)) and $h(I_l)-\bar h=-(h(I_r)-\bar h)=\frac12(h(I_l)-h(I_r))$. Taylor's formula at $\bar h$ with $\lvert\Phi'\rvert\le2^{20}$ and $\lvert\Phi''\rvert\le2^{40}$ on $[0,1]$ (Lemma~\ref{lem:pro_the_cap}) gives, the first-order terms cancelling,
\[
\begin{aligned}
\bigl\lvert\Phi(h(I))-\tfrac12\Phi(h(I_l))-\tfrac12\Phi(h(I_r))\bigr\rvert&\le\bigl\lvert\Phi(h(I))-\Phi(\bar h)\bigr\rvert+\tfrac12\cdot2^{40}\cdot\frac{V(I)}4\\
&\le2^{20}J(I)+2^{37}V(I)\le2^{40}(J+V)(I).
\end{aligned}
\]
For $0\le x\le1$ let $j=\lfloor nx\rfloor$. By Lemma~\ref{lem:shor_arra}(ii) and Lemma~\ref{lem:dil_and_tim}(iv), $G_\Phi(\mathcal{T}_xd)=G_\Phi(\mathcal{T}_{j/n}d)$ is the left side of this display for $I=[j/n,j/n+2]$, whose shortfall and the shortfalls of whose halves lie in $[0,1]$ (Lemma~\ref{lem:dil_and_tim}(i)). With $W_x=(J+V)([j/n,j/n+2])(g_n)$, therefore $\lvert G_{\rm av}(g_n)\rvert\le2^{40}\int_0^1W_x\,dx$, and by $(\ast\ast)$
\[
\lvert Z(g_n)\rvert\le2^{40}\int_0^1W_x\,dx+C_{\rm dom}\sum_Ic_I\,V(I)(g_n).
\]

We now apply the bridge to each weight. By Lemma~\ref{lem:dil_and_tim}(i) and Lemma~\ref{lem:shor_arra}(i), $W_x$ is $\mathcal{F}\bigl(h_g(j,j+2n),h_g(j,j+n),h_g(j+n,j+2n)\bigr)$ with $j=\lfloor nx\rfloor\in[0,n]$ and $\mathcal{F}(u_1,u_2,u_3)=u_1-\frac12(u_2+u_3)+(u_2-u_3)^2$, whose restriction to $[0,1]^3$ has Lipschitz constant $5$ for the maximum norm; and each $V(I)(g_n)$ is $(u_1-u_2)^2$ at two shortfalls of $g$ with indices in $[-8n,12n]$, of Lipschitz constant $4$ on $[0,1]^2$. So every weight $W$ among the $W_x$ and the $V(I)(g_n)$ is a nonnegative shortfall observable with Lipschitz constant at most $6$ (Definition~\ref{def:shor_obse}), all its indices lie in $[-8n,12n]$, within $1024\kappa_{\rm br}n=51.2\,n$ of $0$, and the indices of $Y$ lie within $n$ of $1024n$. Apply Theorem~\ref{thm:bri_pro_est:2} with $a=0$, $b=1024n$ and threshold $\zeta'/4\le1$, which is allowed as $1024n\ge1000+100\log(4/\zeta')$: since $W\ge0$ and $\mathrm{Mid}_n\subset\mathrm{Br}_n$,
\[
\mathbb{E}\bigl[\mathbf{1}_{\mathrm{Mid}_n}W\bigr]\le\mathbb{E}\bigl[\mathbf{1}_{\{h_g(0,1024n)\le\beta_{\rm br}\}}W\,\mathbf{1}_{\{Y>\zeta'/4\}}\bigr]\le C_{\rm br}\bigl(\mathbb{E}W+M^{-1024\sigma_{\rm br}n}\bigr)\,\mathbb{P}(Y>\zeta'/8).
\]
By Lemma~\ref{lem:inv_are_und}(i),(iii), $\mathbb{E}Y=e_n$, and $e_n\le e_N$ by Lemma~\ref{lem:exa_dou_ide}; so Markov's inequality gives $\mathbb{P}(Y>\zeta'/8)\le8e_N/\zeta'$. By Lemma~\ref{lem:shor_arra}(ii) and Lemma~\ref{lem:dil_and_tim}(iii), $\mathbb{E}W_x=\mathbb{E}\,(J+V)([0,2])(g_n)$, and Lemma~\ref{lem:sum_mid_can} (intervals of length $1$ or $2$, aligned at scale $n_0$) gives $\sum_n\mathbb{E}W\le9e_N$ for each weight. As $M\ge M_{\rm an}$ and $n_0\ge1000$, $M^{-1024\sigma_{\rm br}n_0}\le\frac12$, so $\sum_nM^{-1024\sigma_{\rm br}n}\le2M^{-1024\sigma_{\rm br}n_0}\le2M^{-100n_0}$, because $1024\sigma_{\rm br}\ge100$. Hence, for each weight,
\[
\sum_n\mathbb{E}\bigl[\mathbf{1}_{\mathrm{Mid}_n}W\bigr]\le\frac{8C_{\rm br}}{\zeta'}\,e_N\bigl(9e_N+2M^{-100n_0}\bigr)\le\frac{72\,C_{\rm br}}{\zeta'}\,e_N\bigl(e_N+M^{-100n_0}\bigr);
\]
for the weights $W_x$ this holds for each $x$, and Fubini--Tonelli integrates it over $0\le x\le1$. Summing with the coefficients $2^{40}$ and $C_{\rm dom}c_I$, whose total is below $2^{40}+20C_{\rm dom}$, gives the claim.

\textbf{Addition.} As $\mathrm{Slow}_n$, $\mathrm{Mid}_n$ and $\mathrm{Near}_n$ partition $\mathrm{Off}_n$, and $(2^{40}+20C_{\rm dom})\cdot72C_{\rm br}/\zeta'<2^{1188}$, the three ranges give
\[
\sum_n\mathbb{E}\bigl[Z(g_n)\,\mathbf{1}_{\mathrm{Off}_n}\bigr]\le\bigl(160\,C_{\rm dom}+\alpha_{\rm cap}\bigr)\zeta'e_N+\frac{\alpha_{\rm cap}}{n_0}+2^{1188}e_N^2+2^{1188}M^{-100n_0}e_N.
\]
Here $(160C_{\rm dom}+\alpha_{\rm cap})\zeta'<2^{-847}$, and $2^{1188}M^{-100n_0}\le2^{1188}\exp(-100n_0)<2^{-1000}$, as $\log M\ge1$ and $n_0\ge2^{1100}$. Since $2^{-847}+2^{-1000}<2^{-846}$, the lemma follows.

Explicit values: $\zeta'<u_{\rm lo}$; $\zeta'>5\varepsilon$; $\zeta'/2+2\varepsilon<\zeta'$; $\frac34b_{\rm sl}-\varepsilon>\frac12b_{\rm sl}\ge\frac14b_{\rm sl}$; $\frac34b_{\rm sl}+\varepsilon<\beta_{\rm br}$; $\zeta'/2-\varepsilon>\zeta'/4$; $1000+100\log(4/\zeta')<70454\le1024n_0$; $1024\kappa_{\rm br}=51.2$ and $1024\sigma_{\rm br}=102.4$; $(160C_{\rm dom}+\alpha_{\rm cap})\zeta'=2^{-847.68}$; $(2^{40}+20C_{\rm dom})\cdot72C_{\rm br}/\zeta'=2^{1187.49}$.
\end{proof}

The defect bound assembles these estimates with the approximation by slack profiles on the body.

\begin{prop}
 Let $M\ge M_{\rm an}$, and let $n_0<N$ be of the form $2^j$ with $n_0\ge n_{\rm bal}$. Then, summing over $n=2^j$ with $n_0\le n<N$,
\[
\sum_n\mathbb{E}\,Z(g_n)\le2^{-780}e_N+2^{45770}e_N^2+\frac2{n_0}.
\]
\label{prop:def_the_ext:1}
\end{prop}

\begin{proof}
 At each scale we split the initial points into the \emph{body}, where some midpoint cancellation or imbalance of the dyadic family $\mathcal{I}$ of Lemma~\ref{lem:aff_app_off} is not tiny, and its complement $\mathrm{Off}_n$. Off the body, Lemma~\ref{lem:bal_def_off} bounds the defect. The body has small mass summed over the scales, by Markov's inequality and the summability of midpoint cancellations and imbalances; on it the arrays are uniformly close to slack profiles outside an exceptional set of small mass (Lemma~\ref{lem:sum_exc_mas:1}), and the balance defect, which is nonpositive on slack profiles and Lipschitz (Lemma~\ref{lem:pro_the_bal}), is small there.

\textbf{The numbers.} Within this proof sums over $n$ run over $n=2^j$ with $n_0\le n<N$. The proof uses the family $\mathcal{I}$ of Lemma~\ref{lem:aff_app_off}, which has fewer than $2^{1112}$ members, each of length $2^{12-l}$ with $0\le l<1112$ and with endpoints and midpoint in $2^{-1100}\mathbb{Z}$; the body threshold $2^{-2300}$ of Lemma~\ref{lem:bal_def_off}; eleven doublings, from $N$ to $2^{11}N$, which cover the lengths of the members of $\mathcal{I}$; the profile error $2^{-5572}$ of Lemma~\ref{lem:sum_exc_mas:1}; and the Lipschitz constant $2^{153}$ of the balance defect (Lemma~\ref{lem:pro_the_bal}(iii)). The hypothesis $n_0\ge n_{\rm bal}=2^{49000}$ is used only to apply Lemma~\ref{lem:sum_exc_mas:1}, whose statement requires it; every other step needs only $n_0\ge2^{1100}$. Every $n$ is even, so $g_n$ and $\mathcal{D}g_n$ are admissible (Lemma~\ref{lem:dil_and_tim}(vi)), and $Z(g_n)$ is defined, is a continuous function of the initial point and satisfies $\lvert Z(g_n)\rvert\le3$ (Lemma~\ref{lem:pro_the_bal}(i),(iv)). Let $\mathrm{Off}_n$ be the Borel set of Lemma~\ref{lem:bal_def_off}, and $\mathrm{Body}_n=\mathbb{T}^2\setminus\mathrm{Off}_n$.

\textbf{The mass of the body.} By Corollary~\ref{cor:dou_ine_for} (as $M\ge M_{\rm an}$) and induction on $r$, $e_{2^rN}+2^{-r}/N\le D_{\rm dbl}^r(e_N+1/N)$ for every integer $r\ge0$; at $r=11$, with $D_{\rm dbl}=2^{28}$, $e_{2^{11}N}\le2^{308}(e_N+1/N)$. Each $I\in\mathcal{I}$ has length $p=2^{12-l}\le2^{12}$ and is aligned at scale $n_0$, as $n_02^{-1100}$ is an integer; so Lemma~\ref{lem:sum_mid_can} and Lemma~\ref{lem:exa_dou_ide} give
\[
\sum_n\mathbb{E}\,\bigl(J(I)+V(I)\bigr)(g_n)\le9e_{pN/2}\le9e_{2^{11}N}\le9\cdot2^{308}\Bigl(e_N+\frac1N\Bigr).
\]
On $\mathrm{Body}_n$ some $J(I)(g_n)$ or $V(I)(g_n)$ with $I\in\mathcal{I}$ is at least $2^{-2300}$, and all of them are nonnegative (Lemma~\ref{lem:shor_arra}(i)); so $\mathbf{1}_{\mathrm{Body}_n}\le2^{2300}\sum_{I\in\mathcal{I}}\bigl(J(I)+V(I)\bigr)(g_n)$, and
\[
\sum_n\mathbb{P}(\mathrm{Body}_n)\le2^{2300}\cdot2^{1112}\cdot9\cdot2^{308}\Bigl(e_N+\frac1N\Bigr)<2^{3725}\Bigl(e_N+\frac1{n_0}\Bigr).
\]

\textbf{The body.} Let $\mathrm{Exc}_n$ be the Borel sets of Lemma~\ref{lem:sum_exc_mas:1}, which applies as $M\ge M_{\rm an}$ and $n_0\ge n_{\rm bal}$. Let the initial point lie in $\mathrm{Body}_n\setminus\mathrm{Exc}_n$. Then some $d'\in\mathscr{S}$ satisfies $\lVert g_n-d'\rVert_{[-2048,2048]}\le2^{-5572}$, and in particular $\lVert g_n-d'\rVert_{[-16,16]}\le2^{-5572}$. By Lemma~\ref{lem:pro_the_bal}(ii), $d'$ satisfies the hypothesis of Lemma~\ref{lem:pro_the_bal}(i) and $Z(d')\le0$; so Lemma~\ref{lem:pro_the_bal}(iii), applied to $g_n$ and $d'$, gives
\[
Z(g_n)\le Z(d')+2^{153}\cdot2^{-5572}\le2^{-5419}.
\]
As $Z(g_n)\le3$ on $\mathrm{Exc}_n$, Lemma~\ref{lem:sum_exc_mas:1} and the mass of the body give
\[
\begin{aligned}
\sum_n\mathbb{E}\bigl[Z(g_n)\,\mathbf{1}_{\mathrm{Body}_n}\bigr]&\le3\sum_n\mathbb{P}(\mathrm{Exc}_n)+2^{-5419}\sum_n\mathbb{P}(\mathrm{Body}_n)\\
&\le3\cdot2^{45767}e_N^2+3\cdot2^{-4000}e_N+\frac{3\cdot2^{-26000}}{n_0}+2^{-1694}\Bigl(e_N+\frac1{n_0}\Bigr).
\end{aligned}
\]

\textbf{Off the body.} By Lemma~\ref{lem:bal_def_off}, which applies as $M\ge M_{\rm an}$ and $n_0\ge2^{1100}$,
\[
\sum_n\mathbb{E}\bigl[Z(g_n)\,\mathbf{1}_{\mathrm{Off}_n}\bigr]\le2^{-846}e_N+2^{1188}e_N^2+\frac{\alpha_{\rm cap}}{n_0}.
\]

\textbf{Addition.} Since $\mathrm{Body}_n$ and $\mathrm{Off}_n$ partition $\mathbb{T}^2$, $\sum_n\mathbb{E}\,Z(g_n)$ is at most the sum of the last two bounds. The coefficient of $e_N$ is $2^{-1694}+3\cdot2^{-4000}+2^{-846}<2^{-845}<2^{-780}$, that of $e_N^2$ is $3\cdot2^{45767}+2^{1188}<2^{45770}$, and that of $1/n_0$ is $2^{-1694}+3\cdot2^{-26000}+\alpha_{\rm cap}<2$.

Explicit values: $9\cdot2^{2300+1112+308}=2^{3723.17}$; $2^{153-5572}=2^{-5419}$ and $2^{-5419+3725}=2^{-1694}$; the coefficient of $e_N$ is just above $2^{-846}$; the coefficient of $e_N^2$ is $2^{45768.58}$; the coefficient of $1/n_0$ is below $2^{-45.4}$.
\end{proof}

\section{The bootstrap and the proof of the main results} \label{sec:the_boo_and}

This section completes the proof: it derives the uniform bound Proposition~\ref{prop:uni_bou_the:2} from the bound at short times, the doubling inequality and the balance across scales, and then proves the main results.

\begin{prop}
 Let $M\ge M_{\rm ent}$. Then
\[
e_n\le\frac{C_{\rm ent}}{\log M}
\]
for every $n=2^j$.
\label{prop:uni_bou_the:2}
\end{prop}

\begin{proof}
We argue at the least $n=2^j$ at which the mean shortfall exceeds the level $C_{\rm ent}/\log M$: the shortfall starts below this level at short times, one doubling cannot carry it far above the level, and the balance across scales prevents it from creeping across.

Put $\theta=C_{\rm ent}/\log M$, and let $n_0$ be the least $2^j$ with $n_0\ge\log M$, so that $\log M\le n_0<2\log M$. By the values of Convention~\ref{conv:not_and_con}, $\log M\ge\log M_{\rm ent}\ge n_{\rm bal}$ and $M_{\rm ent}\ge M_{\rm an}$; hence $n_0\ge n_{\rm bal}$, $M\ge M_{\rm an}$ and $1/n_0\le1/\log M=\theta/C_{\rm ent}$.

We first show that $e_n<3/\log M<\theta$ for every $n=2^j\le n_0$. With $K=M-4$ (Convention~\ref{conv:sta_map_par}) one has $K\ge M/2\ge2^{16}$. Since $2x\le\sqrt{\exp(x)/2}$ for $x\ge6$ and $\log M\ge6$, every such $n$ satisfies $1\le n\le n_0<2\log M\le\sqrt{M/2}\le\sqrt K$, and Lemma~\ref{lem:mea_sho_sho} gives $e_n<3/\log M$; also $3<C_{\rm ent}$.

Suppose that some $n=2^j$ has $e_n>\theta$, and let $N$ be the least one. By the previous paragraph $N>n_0$, so $N/2\ge n_0$ is of the form $2^j$ and $e_{N/2}\le\theta$. We next show that one doubling cannot carry the shortfall far above the level. Corollary~\ref{cor:dou_ine_for} at $n=N/2$, together with $2/N\le1/n_0\le\theta/C_{\rm ent}$, gives
\[
e_N<e_N+\frac1N\le D_{\rm dbl}\Bigl(e_{N/2}+\frac2N\Bigr)\le D_{\rm dbl}\Bigl(1+\frac1{C_{\rm ent}}\Bigr)\theta<2^{29}\theta
\]
with the values of Convention~\ref{conv:not_and_con}. Hence $e_N<2^{29}C_{\rm ent}/\log M\le2^{29}C_{\rm ent}/\log M_{\rm ent}\le\zeta_{\rm bal}$.

We finally show that the shortfall cannot creep across the level. Since $M\ge M_{\rm an}$, $n_{\rm bal}\le n_0<N$ and $e_N\le\zeta_{\rm bal}$, Proposition~\ref{prop:bal_bet_dya} applies at the scales $n_0<N$; with $e_{n_0}<3/\log M$ and $1/n_0\le1/\log M$ it gives
\[
e_N\le C_{\rm bal}\Bigl(e_{n_0}+\frac1{n_0}\Bigr)<\frac{4C_{\rm bal}}{\log M}<\frac{C_{\rm ent}}{\log M}=\theta,
\]
which contradicts $e_N>\theta$. Hence $e_n\le\theta$ for every $n=2^j$.

Explicit values: $\log M_{\rm ent}=2^{50000}$ exceeds $n_{\rm bal}=2^{49000}$ and $\log M_{\rm an}=10^{12}$; $2x\le\sqrt{\exp(x)/2}$ holds at $x=6$, and the difference $\frac x2-\frac12\log2-\log(2x)$ of the logarithms increases for $x\ge2$; $D_{\rm dbl}(1+1/C_{\rm ent})=2^{28}(1+2^{-430})<2^{29}$; $2^{29}C_{\rm ent}/\log M_{\rm ent}=2^{-49541}<\zeta_{\rm bal}=2^{-47000}$; $4C_{\rm bal}=2^{422}<C_{\rm ent}=2^{430}$.
\end{proof}

With Proposition~\ref{prop:met_ent_fro}, the uniform bound gives the main theorem. The passage from $n=2^j$ to the limit $e_\infty$ uses only that this limit exists over all integers, which Proposition~\ref{prop:met_ent_fro} proves.

\begin{thm}
 Let $k\ge\exp(2^{50000})$. Then
\[
0<\log M-2^{430}\le h_m(f_k)\le\log M.
\]
\label{thm:met_ent_the:1}
\end{thm}

\begin{proof}
By Proposition~\ref{prop:met_ent_fro}, $h_m(f_k)=(1-e_\infty)\log M$, where the limit $e_\infty=\lim_{n\to\infty}e_n$ over all integers $n\ge1$ exists and lies in $[0,1]$; in particular $h_m(f_k)\le\log M$. Since this limit exists, it equals the limit along $n=2^j$: $e_\infty=\lim_{j\to\infty}e_{2^j}$. By the values of Convention~\ref{conv:not_and_con}, $M_{\rm ent}=\exp(2^{50000})$ and $C_{\rm ent}=2^{430}$, and $M=2\pi k+4>k\ge M_{\rm ent}$; so Proposition~\ref{prop:uni_bou_the:2} gives $e_{2^j}\le C_{\rm ent}/\log M$ for every $j\ge0$. Hence $e_\infty\le C_{\rm ent}/\log M$ and $h_m(f_k)\ge\log M-C_{\rm ent}$. Finally $\log M>\log M_{\rm ent}=2^{50000}>C_{\rm ent}$, which proves strict positivity.
\end{proof}

\begin{corollary}
 Let $k\ge\exp(2^{50000})$. Then
\[
\bigl\lvert h_m(f_k)-\log k\bigr\rvert\le2^{430},
\]
and for every $c\in(0,1)$,
\[
m\{z\in\mathbb{T}^2:\lambda_+(z)<c\log M\}\le\frac{2^{430}}{(1-c)\log M}.
\]
Consequently, for every fixed $c\in(0,1)$, $m\{z:\lambda_+(z)\ge c\log k\}\to1$ as $k\to\infty$.
\label{cor:ent_asy_and:2}
\end{corollary}

\begin{proof}
By Theorem~\ref{thm:met_ent_the:1}, $\log M-C_{\rm ent}\le h_m(f_k)\le\log M$, where $C_{\rm ent}=2^{430}$ (Convention~\ref{conv:not_and_con}). Since $\log k<\log M=\log k+\log(2\pi+4/k)\le\log k+\log(2\pi+1)$ and $\log(2\pi+1)<2<C_{\rm ent}$, this gives $\log k-C_{\rm ent}<h_m(f_k)\le\log k+C_{\rm ent}$, the first assertion.

By Proposition~\ref{prop:met_ent_fro}, $0\le\lambda_+\le\log M$ almost everywhere and $\int\lambda_+\,dm=h_m(f_k)$. Hence Theorem~\ref{thm:met_ent_the:1} gives
\[
\int\Bigl(1-\frac{\lambda_+}{\log M}\Bigr)dm=1-\frac{h_m(f_k)}{\log M}\le\frac{C_{\rm ent}}{\log M}.
\]
The integrand is nonnegative and exceeds $1-c$ on $\{\lambda_+<c\log M\}$, so Markov's inequality proves the area bound. Since $\log M>\log k$, the set $\{\lambda_+<c\log k\}$ is contained in $\{\lambda_+<c\log M\}$, and the bound tends to zero as $k\to\infty$.
\end{proof}

Letting $c\to0$ in Corollary~\ref{cor:ent_asy_and:2} gives $m\{z:\lambda_+(z)=0\}\le C_{\rm ent}/\log M$: the stochastic sea has area at least $1-C_{\rm ent}/\log M$. Both the entropy deficit $C_{\rm ent}$ and the threshold $\exp(2^{50000})$ are set by the balance across scales: the deficit must exceed $4C_{\rm bal}$, and the threshold must make $n_0\ge n_{\rm bal}$ and, after one doubling, $2^{459}/\log M\le\zeta_{\rm bal}$.

\section{Scarcity of cancellations} \label{sec:scar_canc}

This appendix proves Theorem~\ref{thm:exp_sca_can}, one of the two quantitative inputs of the paper, restated before the doubling inequality Corollary~\ref{cor:dou_ine_for}: if $M\ge M_{\rm an}$ and $\delta=\delta_{\rm can}$, then for every $i\in\mathbb Z$ and every integer $n\ge1$
\[
\mathbb P\bigl(\mathrm{Can}_\delta(i,n)\bigr)\le C_{\rm can}\exp(-\gamma_{\rm can}n).
\]
The proof is a nested covering argument in one coordinate $\theta$ of the central state. A cancellation produces many matched pairs of hyperbolic times (first subsection); at each of them $\theta$ is confined, through a cubic relation (second subsection), to a set with small covers (third subsection); and the covers nest (last subsection). The covering theorem \ref{thm:sma_cov_the} is the analytic core of the appendix and its hardest part. The appendix uses only Section~\ref{sec:tra_mat_and} and the definition of cancellations (Definition~\ref{def:can_the_gro}).

\subsection{Endpoint-small solutions and hyperbolic times} \label{subsec:end_sma_sol}

The Wronskian of two solutions of the path recurrence is constant, as for any Jacobi operator \cite{locbib:Teschl2000}, and the recurrence at the central time is the centre equation that links the two outward paths (Lemma~\ref{lem:wro_and_cen:2}). A cancellation gives a linearized solution, normalized at the centre, that is small at both ends (Lemma~\ref{lem:can_giv_sol}). Along each outward path the variation sequence then grows, and the solution decays, at nearly full rate, with total deficit $O(\delta n)$ (Lemma~\ref{lem:tot_def_bou}). By the counting form of the Pliss lemma (Lemma~\ref{lem:pli_typ_cou}), all but $O(\delta n+1)$ indices are hyperbolic times (Lemma~\ref{lem:cou_non_hyp}), and these pair across the two paths (Lemma~\ref{lem:man_mat_pai}). At each hyperbolic time the solution has an explicit boundary expansion (Lemma~\ref{lem:bou_val_for}), which is the input of the next subsection.

\begin{definition}
 Fix a \textbf{central time} $\iota\in\{0,-1\}$ and a position $q_\iota$, and let the \textbf{free coordinate} $\theta=q_{\iota-1}$, the other coordinate of the state at time $\iota$, vary; the orbit is the one whose state at time $\iota$ is $(q_\iota,\theta)$, and $q_{\iota+1}=\phi(q_\iota)-\theta$ modulo one. All objects below are functions of $\theta$ for fixed $q_\iota$, computed along real lifts.

The \textbf{outward paths} are
\[
q^-_j=q_{\iota-j},\qquad q^+_j=q_{\iota+j}\qquad(j\ge0).
\]
When only one path is involved its sign is omitted: $q_j$, $v_j=v(q_j)$. The \textbf{path recurrence} of a path is the relation $w_{j+1}=v_jw_j-w_{j-1}$ $(j\ge1)$ for real sequences $(w_j)_{j\ge0}$. The \textbf{variation sequence} $t$ and the \textbf{complementary sequence} $\widetilde t$ of a path are the solutions of the path recurrence with
\begin{equation}\label{eq:cs_variation}
t_0=0,\quad t_1=1,\qquad\widetilde t_0=1,\quad\widetilde t_1=0.
\end{equation}

The \textbf{central pair} of a linearized solution $(u_i)$ is $(u_0,u_{-1})$. A \textbf{normalized solution} is a linearized solution $(u_i)$ along the orbit with $u_\iota=1$; on the paths $u^-_j=u_{\iota-j}$ and $u^+_j=u_{\iota+j}$, so $u^\pm_0=1$, and the sign is omitted as above. When the sign is omitted on a single path, $u_j$ is the path value and $u_0=1$; otherwise the indices of $(u_i)$ are orbit times, as in the central pair $(u_0,u_{-1})$ and in the endpoint pairs below, and for $\iota=-1$ the orbit value $u_0$ equals $u^+_1$.

For $n\ge1$ the \textbf{path lengths} are $N_+=n-\iota$ and $N_-=n+\iota+1$ (one is $n$, the other $n+1$), so that $(u^+_{N_+},u^+_{N_+-1})=(u_n,u_{n-1})$ and $(u^-_{N_-},u^-_{N_--1})=(u_{-n-1},u_{-n})$; when only one path is involved, its length is written $N$.

For $j\ge1$ with $t_j\ne0$, the \textbf{boundary value} is $R_j=-\widetilde t_j/t_j$.
\label{def:out_pat_and}
\end{definition}

\begin{definition}
 The \textbf{pair-growth logarithms} of a path, for a normalized solution $u$, are
\[
\tau_j=\log_M\max(\lvert t_j\rvert,\lvert t_{j-1}\rvert),\qquad\upsilon_j=\log_M\max(\lvert u_j\rvert,\lvert u_{j-1}\rvert)\qquad(j\ge1),
\]
defined whenever the pair is nonzero; with a sign, $\tau^\pm_j$ and $\upsilon^\pm_j$ refer to the plus or the minus path.

Given a normalized solution and one of its outward paths, an index $j\ge2$ is a \textbf{hyperbolic time} for that solution if
\begin{equation}\label{eq:cs_hyperbolic_time}
\begin{aligned}
\tau_j-\tau_l&\ge c_{\rm hyp}(j-l)&&(1\le l<j),\\
\tau_{j+1}-\tau_j&\ge c_{\rm hyp},\\
\upsilon_j-\upsilon_{j+1}&\ge c_{\rm hyp},&\qquad\upsilon_{j+1}-\upsilon_{j+2}&\ge c_{\rm hyp}.
\end{aligned}
\end{equation}
The first line is a hyperbolic-time condition for the sequence $\tau$ in the manner of \cite[Definition 2.6]{locbib:AlvesBonattiViana2000}; the other three are one-step inequalities, one on the variation and two on the solution. An index $j\ge2$ that is not a hyperbolic time is \textbf{non-hyperbolic}.
\label{def:pai_gro_log}
\end{definition}

\begin{lemma}
 Fix $\iota$ and $q_\iota$.
\begin{enumerate}
\item[(i)] On either path and for every real lift, the positions satisfy $q_{j+1}=\phi(q_j)-q_{j-1}$ and the restriction of every linearized solution satisfies the path recurrence $(j\ge1)$. For local real lifts of the positions, $\partial_\theta q^-_j=t^-_j$ and $\partial_\theta q^+_j=-t^+_j$ for all $j\ge0$.
\item[(ii)] On either path, every solution $u$ of the path recurrence with $u_0=1$ satisfies
\begin{equation}\label{eq:cs_wronskian}
u_{j-1}t_j-u_jt_{j-1}=1\qquad(j\ge1).
\end{equation}
\item[(iii)] The outward restrictions of a normalized solution satisfy
\begin{equation}\label{eq:cs_centre_equation}
u^-_1+u^+_1=v(q_\iota).
\end{equation}
Conversely, two solutions of the path recurrences with value one at index zero satisfying \eqref{eq:cs_centre_equation} assemble into a normalized solution across the centre.
\item[(iv)] On either path the pairs $(t_j,t_{j-1})$ and, for a normalized solution, $(u_j,u_{j-1})$, $j\ge1$, are nonzero, and
\begin{equation}\label{eq:cs_log_increments}
\lvert\tau_{j+1}-\tau_j\rvert\le1,\qquad\lvert\upsilon_{j+1}-\upsilon_j\rvert\le1 .
\end{equation}
\end{enumerate}
\label{lem:wro_and_cen:2}
\end{lemma}

\begin{proof}
All four items rest on the fact that the outward positions, their $\theta$-derivatives and all outward linearized solutions satisfy one second-order recurrence of determinant one.

(i) On the plus path, \eqref{eq:orbit_recurrence} at time $\iota+j$ reads $q^+_{j+1}=\phi(q^+_j)-q^+_{j-1}$; on the minus path, the same relation at time $\iota-j$ reads $q_{\iota-j+1}=\phi(q_{\iota-j})-q_{\iota-j-1}$, that is $q^-_{j+1}=\phi(q^-_j)-q^-_{j-1}$. Both hold on real lifts. In the same way the linearized relation $u_{i+1}=v(q_i)u_i-u_{i-1}$ of \eqref{eq:linearized_matrix}, read at the time $\iota\pm j$, is the path recurrence for $u^\pm$. The position $q^\pm_0=q_\iota$ is fixed, $q^-_1=\theta$ and $q^+_1=\phi(q_\iota)-\theta$. Differentiating the recurrence of the positions along a lift shows, since $\phi'=v$, that $\partial_\theta q^\pm_j$ solves the path recurrence, with initial data $(0,1)$ on the minus path and $(0,-1)$ on the plus path; by \eqref{eq:cs_variation} these are $t^-_j$ and $-t^+_j$. The sign difference comes solely from the initial positions at index one; it does not change the recurrence or its Wronskian identity.

(ii) The left side of \eqref{eq:cs_wronskian} equals $u_0t_1-u_1t_0=1$ at $j=1$, and the recurrence preserves it, since each step matrix has determinant one; this is the constancy of the Wronskian of a Jacobi difference equation, see \cite[Section 1.1, equations (1.20) and (1.21)]{locbib:Teschl2000}, in the case of off-diagonal coefficients $-1$.

(iii) The linearized recurrence at time $\iota$ reads $u_{\iota+1}=v(q_\iota)u_\iota-u_{\iota-1}$, and $u_\iota=1$, $u_{\iota+1}=u^+_1$, $u_{\iota-1}=u^-_1$. Conversely, the two outward solutions agree at the centre, and \eqref{eq:cs_centre_equation} is exactly the recurrence at time $\iota$; at all other times the recurrence holds on one of the paths, by (i).

(iv) The pairs are images of the nonzero pairs $(t_1,t_0)=(1,0)$ and $(u_1,u_0)$, with $u_0=1$, under invertible step matrices, so they are nonzero and the logarithms are defined. One step of the recurrence and its inverse have norm at most $M$ in the operator norm induced by the maximum norm (Lemma~\ref{lem:nor_and_gro}(i)), so the maximum norm of a pair changes by a factor in $[M^{-1},M]$ per step; this is \eqref{eq:cs_log_increments}.
\end{proof}

\begin{definition}
 Let $n\ge1$ and $\delta>0$. An \textbf{endpoint-small configuration} of length $n$ and defect $\delta$ is an orbit together with a central time $\iota\in\{0,-1\}$ and a normalized solution $(u_i)$ such that
\begin{equation}\label{eq:cs_small_solution}
\lVert(u_0,u_{-1})\rVert\le\sqrt2,\qquad\lVert(u_n,u_{n-1})\rVert,\ \lVert(u_{-n},u_{-n-1})\rVert\le3M^{-(1-3\delta)n}.
\end{equation}
Its two outward paths are those of Definition~\ref{def:out_pat_and} for this $\iota$; on each, the two small endpoint entries sit at the indices $N-1$ and $N$.
\label{def:end_sma_con}
\end{definition}

\begin{lemma}
 Let $\delta>0$, $n\ge1$ and $z\in\mathrm{Can}_\delta(0,n)$. Then there are a centre $\iota\in\{0,-1\}$ and a normalized solution along the orbit of $z$ which, together with this orbit, form an endpoint-small configuration of length $n$ and defect $\delta$.
\label{lem:can_giv_sol}
\end{lemma}

\begin{proof}
The solution through the least expanded input direction of $T_{0,-n}$ is small at the left end; since the full product barely grows, this direction is nearly the least expanded direction of $T_{0,n}$ as well, so the solution is small at both ends.

By Definition~\ref{def:can_the_gro} and \eqref{eq:growth_orbit} (Definition~\ref{def:gro_sem_orb}), the transfer matrices of Definition~\ref{def:orb_sta_and} satisfy
\[
\min\bigl(\log_M\lVert T_{0,-n}\rVert,\ \log_M\lVert T_{0,n}\rVert\bigr)\ge(1-\delta)n,\qquad\log_M\lVert T_{-n,n}\rVert\le\delta n .
\]
Write
\[
\sigma_-=\lVert T_{0,-n}\rVert,\qquad\sigma_+=\lVert T_{0,n}\rVert,
\]
and let $a_-$, $a_+$ be unit most expanded input directions for these two products. Let $b$ be a unit vector perpendicular to $a_-$. Since transfer matrices have determinant one (Lemma~\ref{lem:nor_and_gro}(i)), the singular values of either product are its norm and the reciprocal of its norm, and the singular input directions are orthogonal. Hence
\[
\lVert T_{0,-n}b\rVert=\sigma_-^{-1},\qquad\lVert T_{0,n}b\rVert\ge\sigma_+\lvert a_-\wedge a_+\rvert,
\]
the second because the component of $b$ along $a_+$ has absolute value $\lvert a_-\wedge a_+\rvert$ (Convention~\ref{conv:wed_pro_and}). Applying the full product $T_{-n,n}$ to the vector $T_{0,-n}b$ gives $T_{0,n}b$ (Lemma~\ref{lem:nor_and_gro}(i)), so
\[
\lVert T_{-n,n}\rVert\ge\sigma_-\sigma_+\lvert a_-\wedge a_+\rvert,\qquad\text{hence}\qquad\lvert a_-\wedge a_+\rvert\le M^{\delta n-2(1-\delta)n}=M^{-(2-3\delta)n}.
\]
Resolving $b$ in the singular input directions of $T_{0,n}$ also gives the upper bound
\[
\lVert T_{0,n}b\rVert\le\sigma_+^{-1}+\sigma_+\lvert a_-\wedge a_+\rvert\le2M^{-(1-3\delta)n},
\]
using $\sigma_+^{-1}\le M^{-(1-\delta)n}$ and $\sigma_+\le M^n$ (Lemma~\ref{lem:nor_and_gro}(i)). The left endpoint has the smaller bound $\sigma_-^{-1}\le M^{-(1-\delta)n}$.

Let $(u_i)$ be the linearized solution with central pair $(u_0,u_{-1})=b$ (Definition~\ref{def:out_pat_and}); its endpoint pairs are $(u_n,u_{n-1})=T_{0,n}b$ and $(u_{-n},u_{-n-1})=T_{0,-n}b$, so it is small at both endpoints. Since $\lVert b\rVert=1$, one of $u_0,u_{-1}$ has absolute value at least $1/\sqrt2$; choose that time $\iota\in\{0,-1\}$ and divide the whole solution by $u_\iota$. The normalized solution has $u_\iota=1$, central pair of norm at most $\sqrt2$, and endpoint pairs of norm at most $\sqrt2\cdot2M^{-(1-3\delta)n}\le3M^{-(1-3\delta)n}$, which is \eqref{eq:cs_small_solution}.

This construction is used only as an existence statement; no measurable selection of a singular direction is required.
\end{proof}

\begin{lemma}
 For every $\delta>0$, $n\ge1$ and every endpoint-small configuration of length $n$ and defect $\delta$, the following hold on each of its two outward paths:
\begin{equation}\label{eq:cs_total_deficit}
\sum_{j=1}^{N-1}\bigl[1-(\tau_{j+1}-\tau_j)\bigr]\le3\delta n+\log6,\qquad\sum_{j=1}^{N-1}\bigl[1-(\upsilon_j-\upsilon_{j+1})\bigr]\le3\delta n+\log3,
\end{equation}
every summand is nonnegative, and $\tau_N\ge(1-3\delta)n-\log_M6$.
\label{lem:tot_def_bou}
\end{lemma}

\begin{proof}
The Wronskian forces the variation sequence to be large where the solution is small, so near-full growth of $\tau$ and near-full decay of $\upsilon$ over the whole path leave only an $O(\delta n+1)$ total deficit.

By Lemma~\ref{lem:wro_and_cen:2}(i), the restriction of the normalized solution to the path solves the path recurrence, with $u_0=1$. At index $N$, \eqref{eq:cs_wronskian} (Lemma~\ref{lem:wro_and_cen:2}(ii)) and the endpoint bound in \eqref{eq:cs_small_solution} imply
\[
1\le2\max(\lvert u_N\rvert,\lvert u_{N-1}\rvert)\max(\lvert t_N\rvert,\lvert t_{N-1}\rvert),
\]
and hence
\[
\tau_N\ge(1-3\delta)n-\log_M6,\qquad\upsilon_N\le-(1-3\delta)n+\log_M3 .
\]
Since $\tau_1=0$, $\upsilon_1\ge0$ (because $u_0=1$) and $N-1\le n$, the total deficits of the sequences $\tau$ and $-\upsilon$ over the path satisfy \eqref{eq:cs_total_deficit}; here $\log_Mx\le\log x$ for $x\ge1$, since $M>4$ (Convention~\ref{conv:sta_map_par}). Every summand is nonnegative by \eqref{eq:cs_log_increments} (Lemma~\ref{lem:wro_and_cen:2}(iv)).
\end{proof}

\begin{lemma}
 For every $\delta>0$, $n\ge1$ and every endpoint-small configuration of length $n$ and defect $\delta$, the number of indices in $\{2,\dots,N-2\}$ of either outward path that are not hyperbolic times is at most $1200(\delta n+1)$.
\label{lem:cou_non_hyp}
\end{lemma}

\begin{proof}
Put $c=c_{\rm hyp}$ (Convention~\ref{conv:not_and_con}). Each non-hyperbolic index costs a definite part of the total deficit of Lemma~\ref{lem:tot_def_bou}: failures of the first line of \eqref{eq:cs_hyperbolic_time} are counted by the Pliss lemma, failures of the one-step lines directly.

Fix one of the two paths, of length $N$, and put
\[
\Gamma_\tau=\sum_{j=1}^{N-1}\bigl[1-(\tau_{j+1}-\tau_j)\bigr],\qquad\Gamma_\upsilon=\sum_{j=1}^{N-1}\bigl[1-(\upsilon_j-\upsilon_{j+1})\bigr].
\]
By \eqref{eq:cs_total_deficit} (Lemma~\ref{lem:tot_def_bou}), $\Gamma_\tau\le3\delta n+\log6$ and $\Gamma_\upsilon\le3\delta n+\log3$, every summand being nonnegative; the one-step increments of $\tau$ and of $-\upsilon$ are at most one by \eqref{eq:cs_log_increments} (Lemma~\ref{lem:wro_and_cen:2}(iv)).

We first count the failures of the first line. Put $x_i=\tau_{i+1}$ for $0\le i\le N-1$; this sequence has increments at most one and total deficit $\Gamma_\tau$. An index $j$ fails $\tau_j-\tau_l\ge c(j-l)$ for some $1\le l<j$ exactly when $x_{j-1}-x_{l'}<c(j-1-l')$ for some $0\le l'<j-1$. By the Pliss lemma \ref{lem:pli_typ_cou} there are at most $\Gamma_\tau/(1-c)$ such indices.

If $\tau_{j+1}-\tau_j<c$, the summand $1-(\tau_{j+1}-\tau_j)$ of $\Gamma_\tau$ exceeds $1-c$; so at most $\Gamma_\tau/(1-c)$ indices fail the second line.

If $\upsilon_j-\upsilon_{j+1}<c$, the corresponding summand of $\Gamma_\upsilon$ exceeds $1-c$, so at most $\Gamma_\upsilon/(1-c)$ increments of $-\upsilon$ are small; a small increment at position $j$ makes the third line fail at the indices $j$ and $j-1$ only, so at most $2\Gamma_\upsilon/(1-c)$ indices fail the third line.

The three counts add up to at most
\[
\frac{2\Gamma_\tau+2\Gamma_\upsilon}{1-c}\le\frac{12\delta n+2\log6+2\log3}{1-c}\le\frac{12(\delta n+1)}{1-c}=1200(\delta n+1),
\]
since $2\log6+2\log3<5.79\le12$ and $1-c_{\rm hyp}=1/100$.
\end{proof}

Hyperbolic times on the two paths are matched when their pair-growth logarithms differ by at most $3/5$; a cancellation produces at least $n/20$ such pairs, in increasing order on both paths.

\begin{definition}
 Fix $\iota$ and $q_\iota$. A pair $(j_+,j_-)$ of indices at least $2$ on the plus and the minus path is \textbf{matched} at $\theta$ if both are hyperbolic times at $\theta$ for one common normalized solution and
\begin{equation}\label{eq:cs_matching}
\bigl\lvert\tau^+_{j_+}-\tau^-_{j_-}\bigr\rvert\le\frac35 .
\end{equation}
\label{def:mat_pai_hyp:2}
\end{definition}

\begin{lemma}
 Let $0<\delta\le\delta_{\rm can}$ and $n\ge1/\delta_{\rm can}$.
\begin{enumerate}
\item[(i)] Every endpoint-small configuration of length $n$ and defect $\delta$ has at least $n/20$ pairs $(j^+_l,j^-_l)$ of indices $j^\pm_l\in\{2,\dots,N_\pm-2\}$, hyperbolic times on the respective paths for the configuration's normalized solution, such that both sequences $(j^+_l)_l$ and $(j^-_l)_l$ increase strictly and
\[
\bigl\lvert\tau^+_{j^+_l}-\tau^-_{j^-_l}\bigr\rvert\le\frac12\qquad\text{for every }l;
\]
in particular each pair is matched at the configuration's free coordinate $\theta$.
\item[(ii)] The strictly increasing patterns $\bigl((j^+_l,j^-_l)\bigr)_{l=1}^{\lceil n/20\rceil}$ of $\lceil n/20\rceil$ index pairs with entries in $\{2,\dots,n\}$ number at most $\exp(2n/5)$.
\end{enumerate}
\label{lem:man_mat_pai}
\end{lemma}

\begin{proof}
We take levels two apart for the plus logarithms, pair the first plus index reaching each level with a minus index whose logarithm lies within $\frac12$ of it, and discard the levels spoiled by non-hyperbolic indices; first passage makes the pairs increase, and each non-hyperbolic index spoils at most one level. Put $c=c_{\rm hyp}$, so that $\delta_{\rm can}=10^{-3}(1-c)/24$ (Convention~\ref{conv:not_and_con}).

We first bound the non-hyperbolic indices. By Lemma~\ref{lem:cou_non_hyp}, each path has at most $1200(\delta n+1)$ non-hyperbolic indices in $\{2,\dots,N-2\}$. With the values of Convention~\ref{conv:not_and_con}, $1200\delta_{\rm can}=10^{-3}/2$, and $1200\le10^{-3}n/2$ for $n\ge1/\delta_{\rm can}$; so each path has at most $10^{-3}n$ non-hyperbolic indices.

We next locate the ends of the paths. On each path $\tau_1=0$, so $\tau_2\le1$ by \eqref{eq:cs_log_increments} (Lemma~\ref{lem:wro_and_cen:2}(iv)); and, by Lemma~\ref{lem:tot_def_bou}, $\log_M6\le\log6$ (as $M>4$, Convention~\ref{conv:sta_map_par}) and \eqref{eq:cs_log_increments},
\[
\tau_{N-2}\ge\tau_N-2\ge(1-3\delta)n-\log6-2>\frac{3n}4+1 .
\]

We now choose levels and partners. For $0\le l\le l_0=\lfloor n/4\rfloor-2$ put $a_l=\lceil n/4\rceil+2l$, so that $n/4\le a_l\le3n/4-2$. Let $j^+_l$ be the least plus index with $\tau^+_{j^+_l}\ge a_l$. Then $a_l\le\tau^+_{j^+_l}<a_l+1$, and $3\le j^+_l\le N_+-2$ by the bounds on the ends. Let $j^-$ be the least minus index with $\tau^-_{j^-}\ge\tau^+_{j^+_l}$. Then $3\le j^-\le N_--2$, again by the bounds on the ends, and $\tau^-_{j^--1}<\tau^+_{j^+_l}\le\tau^-_{j^-}\le\tau^-_{j^--1}+1$; so one of $j^--1,j^-$, called $j^-_l$, satisfies $\lvert\tau^-_{j^-_l}-\tau^+_{j^+_l}\rvert\le\frac12$.

We show that the pairs increase. Let $l<l'$. Then $\tau^+_i<a_l$ for $i<j^+_l$ and $\tau^+_{j^+_l}<a_l+1<a_{l'}$, so $j^+_l<j^+_{l'}$ by first passage. Also $\tau^-_{j^-_l}\in[a_l-\frac12,a_l+\frac32)$, so $\tau^-_{j^-_l}<a_l+\frac32\le a_{l'}-\frac12\le\tau^-_{j^-_{l'}}$; when $j^-_l$ and $j^-_{l'}$ are hyperbolic times, the first line of \eqref{eq:cs_hyperbolic_time} at the larger of the two indices forces $j^-_l<j^-_{l'}$.

We discard the levels spoiled by non-hyperbolic indices. If $j^+_l=j$, then $a_l\in(\max_{i<j}\tau^+_i,\tau^+_j]$, an interval of length at most one. If $j^-_l=j'$, then $a_l\in(\tau^-_{j'}-\frac32,\tau^-_{j'}+\frac12]$, an interval of length two. Since the levels are two apart, each such interval contains at most one level. So, by the first step, at most $2\cdot10^{-3}n$ levels have $j^+_l$ or $j^-_l$ non-hyperbolic. The remaining levels number at least $n/4-2-2\cdot10^{-3}n\ge n/20$, since $n\ge1/\delta_{\rm can}$, and their pairs $(j^+_l,j^-_l)$ are hyperbolic times, strictly increasing on both paths, with matching error at most $\frac12\le\frac35$. In particular each pair is matched at $\theta$ (Definition~\ref{def:mat_pai_hyp:2}). This proves (i).

We count the patterns. A strictly increasing pattern with exactly $p$ pairs is determined by the $p$-subsets of $\{2,\dots,n\}$ formed by its plus and by its minus indices, paired in order; there are $\binom{n-1}{p}^2$ such choices. Since $p!\ge(p/\exp(1))^p$, $\binom{n-1}{p}\le n^p/p!\le(\exp(1)n/p)^p$. For $p=\lceil n/20\rceil$ one has $n/p\le20$ and $2p<n/10+2$, hence
\[
\binom{n-1}{p}^2\le\bigl(20\exp(1)\bigr)^{2p}\le\bigl(20\exp(1)\bigr)^{n/10+2}\le\exp(2n/5),
\]
the last inequality because $(n/10+2)\log(20\exp(1))\le2n/5$ for $n\ge18726$, while $1/\delta_{\rm can}\ge18726$ by the values of Convention~\ref{conv:not_and_con}. This proves (ii).
\end{proof}

At a hyperbolic time the value $u_1$ of the solution is the boundary value $R_j$ plus an explicit reciprocal term, up to an error $O(M^{-3c_{\rm hyp}}/L_j^2)$; this expansion is the source of the cubic relation.

\begin{lemma}
 The following holds for every normalized solution $u$, either outward path, and every index $j$ that is a hyperbolic time for $u$ on that path; put $L_j=\lvert t_j\rvert$.
\begin{enumerate}
\item[(a)] The variation sequence is dominated at $j$:
\begin{equation}\label{eq:cs_dominated_variation}
L_j=M^{\tau_j},\qquad\lvert t_l\rvert\le L_jM^{-c_{\rm hyp}(j-l)}\quad(1\le l<j),\qquad\lvert t_{j+1}\rvert\ge M^{c_{\rm hyp}}L_j .
\end{equation}
\item[(b)] The solution decays across $j$:
\begin{equation}\label{eq:cs_decaying_solution}
\lvert u_{j+1}\rvert\le M^{-c_{\rm hyp}}\lvert u_j\rvert\le M^{-2c_{\rm hyp}}\lvert u_{j-1}\rvert,\qquad\lvert u_{j-1}\rvert\le\frac2{L_j}.
\end{equation}
\item[(c)] The boundary values $R_j$ and $R_{j+1}$ exist, $R_{j+1}-R_j=1/(t_jt_{j+1})$, $u_1-R_{j+1}=u_{j+1}/t_{j+1}$, and
\begin{equation}\label{eq:cs_boundary_value}
u_1=R_{j+1}+\frac{u_{j+1}}{t_{j+1}}=R_j+\frac1{t_j^2\,(v_j-t_{j-1}/t_j)}+\varepsilon_j,\qquad\lvert\varepsilon_j\rvert\le\frac{2M^{-3c_{\rm hyp}}}{L_j^2}.
\end{equation}
\end{enumerate}
\label{lem:bou_val_for}
\end{lemma}

\begin{proof}
Put $c=c_{\rm hyp}$. The hyperbolic-time condition says that $t$ grows and $u$ decays at rate at least $c$ around $j$; the Wronskian then pins $u_{j-1}$ to size $1/L_j$, and the boundary value at index $j+1$ differs from $u_1$ by the small quantity $u_{j+1}/t_{j+1}$. By Lemma~\ref{lem:wro_and_cen:2}(i), the restriction of $u$ to the path solves the path recurrence, and $u_0=1$.

(a) For $1\le l<j$, the first line of \eqref{eq:cs_hyperbolic_time} gives $\max(\lvert t_l\rvert,\lvert t_{l-1}\rvert)=M^{\tau_l}\le M^{\tau_j-c(j-l)}$. With $l=j-1$ this shows $\lvert t_{j-1}\rvert<M^{\tau_j}$, so $t_j$ is the dominant entry of its pair and $L_j=M^{\tau_j}$; hence $\lvert t_l\rvert\le L_jM^{-c(j-l)}$. The second line gives $\max(\lvert t_{j+1}\rvert,\lvert t_j\rvert)\ge M^cL_j>L_j$, so $t_{j+1}$ is dominant at the next pair and $\lvert t_{j+1}\rvert\ge M^cL_j$.

(b) The inequality $\upsilon_{j+1}-\upsilon_{j+2}\ge c$ forces $\lvert u_j\rvert>\lvert u_{j+1}\rvert$ (otherwise the maximum $\lvert u_{j+1}\rvert$ would exceed itself by the factor $M^c$), and then $\lvert u_j\rvert\ge M^c\max(\lvert u_{j+2}\rvert,\lvert u_{j+1}\rvert)$, so $\lvert u_{j+1}\rvert\le M^{-c}\lvert u_j\rvert$. In the same way $\upsilon_j-\upsilon_{j+1}\ge c$ and $\max(\lvert u_{j+1}\rvert,\lvert u_j\rvert)=\lvert u_j\rvert$ give $\lvert u_j\rvert\le M^{-c}\lvert u_{j-1}\rvert$. Using \eqref{eq:cs_wronskian} (Lemma~\ref{lem:wro_and_cen:2}(ii)) and $\lvert t_{j-1}/t_j\rvert\le M^{-c}$ from (a),
\[
\lvert u_{j-1}\rvert L_j=\lvert1+u_jt_{j-1}\rvert\le1+M^{-2c}\lvert u_{j-1}\rvert L_j ,
\]
so $\lvert u_{j-1}\rvert\le L_j^{-1}(1-M^{-2c})^{-1}\le2/L_j$, because $M>4$ (Convention~\ref{conv:sta_map_par}) gives $M^{-2c}<1/2$.

(c) By (a), $t_j\ne0$ and $t_{j+1}\ne0$, so $R_j$ and $R_{j+1}$ exist. By linearity and \eqref{eq:cs_variation}, the sequence $w=\widetilde t+R_jt$ solves the path recurrence with values $w_0=1$, $w_1=R_j$ and $w_j=\widetilde t_j+R_jt_j=0$. Its value at index $j-1$ is $1/t_j$, by \eqref{eq:cs_wronskian} with $w_j=0$; its value at index $j+1$ is therefore $v_j\cdot0-1/t_j=-1/t_j$. Since $w_{j+1}=\widetilde t_{j+1}+R_jt_{j+1}$ and $R_{j+1}=-\widetilde t_{j+1}/t_{j+1}$,
\[
R_{j+1}-R_j=-\frac{w_{j+1}}{t_{j+1}}=\frac1{t_jt_{j+1}} .
\]
Every solution $u$ of the path recurrence with $u_0=1$ equals $\widetilde t+u_1t$, since both have the same values at the indices $0$ and $1$; so $u_{j+1}=\widetilde t_{j+1}+u_1t_{j+1}=(u_1-R_{j+1})t_{j+1}$, that is $u_1-R_{j+1}=u_{j+1}/t_{j+1}$. Since $t_{j+1}=v_jt_j-t_{j-1}=t_j(v_j-t_{j-1}/t_j)$,
\[
u_1=R_{j+1}+\frac{u_{j+1}}{t_{j+1}}=R_j+\frac1{t_j^2(v_j-t_{j-1}/t_j)}+\frac{u_{j+1}}{t_{j+1}},
\]
and by (a) and (b)
\[
\Bigl\lvert\frac{u_{j+1}}{t_{j+1}}\Bigr\rvert\le\frac{M^{-2c}\lvert u_{j-1}\rvert}{M^cL_j}\le\frac{2M^{-3c}}{L_j^2}.
\]
This is \eqref{eq:cs_boundary_value}.
\end{proof}

\subsection{Local tests and the cubic relation} \label{subsec:loc_tes_and}

A local test is attached to a pair of hyperbolic times matched at a value $\theta_*$ of the free coordinate. Its base scales are the sizes of the variation sequences at the two times, and its successful set records the rescaled offsets $h\in[-1,1]$ for which the two indices are again hyperbolic times, for a common normalized solution, at $\theta=\theta_*+h/L$, where $L$ is the larger base scale. The cancellation estimate needs only that successful sets are small. The main lemma of this subsection, Lemma~\ref{lem:cub_rel_suc:1}, says that at every successful point the boundary cubic $P$, the cubic Taylor polynomial at $h=0$ of $KL^2\bigl(v(q_\iota)-R_{j_1}-R_{j_2}\bigr)$, equals, up to a small error, the sum $\Phi_1+\Phi_2$ of two reciprocal cosine terms in the rescaled phases, the positions at the two hyperbolic times as functions of $h$. Derivative bounds along a rescaled path (Lemma~\ref{lem:der_bou_res}) show that the boundary values have small fourth derivatives (Lemma~\ref{lem:fou_der_the:1}) and the rescaled phases are nearly affine (Lemma~\ref{lem:pha_bou_loc}), and the expansions of Lemma~\ref{lem:bou_val_for} on the two paths are added through the centre equation.

\begin{definition}
 A \textbf{local test} is a tuple $(\iota,q_\iota,\theta_*,j_+,j_-)$ whose pair $(j_+,j_-)$ is matched at the \textbf{base coordinate} $\theta_*$. Its \textbf{base scales} are
\[
L_+=\lvert t^+_{j_+}(\theta_*)\rvert,\qquad L_-=\lvert t^-_{j_-}(\theta_*)\rvert,\qquad L=\max(L_+,L_-).
\]
Its \textbf{successful set} $\mathcal S\subset[-1,1]$ is the set of $h\in[-1,1]$ such that, at $\theta=\theta_*+h/L$ (on the real lift), $j_+$ and $j_-$ are hyperbolic times for some common normalized solution, which may depend on $h$; elements of $\mathcal S$ are \textbf{successful points}.
\label{def:loc_tes_and}
\end{definition}

\begin{definition}
 Fix $\iota$, $q_\iota$, a value $\theta_*$ of the free coordinate (the base coordinate), an outward path and an index $j\ge1$ on it with $t_j(\theta_*)\ne0$. Put
\[
L_j=\lvert t_j(\theta_*)\rvert,\qquad\xi=L_j(\theta-\theta_*),\qquad D=\frac{d}{d\xi},
\]
and regard the lifted positions $q_l$ $(0\le l\le j)$ and $t_l$ as functions of $\xi$. The \textbf{quotient} is $E=t_{j-1}/t_j$, wherever $t_j\ne0$.
\label{def:sin_pat_res}
\end{definition}

\begin{lemma}
 Let $M\ge\exp(100)$. Fix $\iota$, $q_\iota$, a base coordinate $\theta_*$, an outward path, and an index $j\ge2$ on it that is a hyperbolic time at $\theta_*$ for some normalized solution, and put $\rho_l=M^{-c_{\rm hyp}(j-l)}$. Then for $\lvert\xi\rvert\le1$ and $0\le l\le j$:
\begin{enumerate}
\item[(i)] $\lvert q_l(\xi)-q_l(0)\rvert\le2\rho_l$ and $\lvert Dq_l(\xi)\rvert\le2\rho_l$;
\item[(ii)] $\lvert D^pq_l(\xi)\rvert\le10^p\rho_l^pM^{1-pc_{\rm hyp}}$ for $2\le p\le5$;
\item[(iii)] $\frac12L_j\le\lvert t_j(\xi)\rvert\le\frac32L_j$, and $\lvert t_l(\xi)\rvert\le2L_j\rho_l$ for $l<j$;
\item[(iv)] the quotient $E=t_{j-1}/t_j=Dq_{j-1}/Dq_j$ satisfies $\lvert E\rvert\le4M^{-c_{\rm hyp}}$ and $\lvert D^pE\rvert\le M^{-c_{\rm hyp}}$ for $1\le p\le4$.
\end{enumerate}
\label{lem:der_bou_res}
\end{lemma}

\begin{proof}
Put $c=c_{\rm hyp}$ and $\gamma=M^{1-2c}$; since $M\ge\exp(100)$, $\gamma\le\exp(-98)$. For a smooth function $F$ we use Faà di Bruno's formula over set partitions, $D^p[F(q)]=\sum_\mathcal{P} F^{(\lvert\mathcal{P}\rvert)}(q)\prod_{\beta\in\mathcal{P}}D^{\lvert\beta\rvert}q$, the sum running over the partitions $\mathcal{P}$ of $\{1,\dots,p\}$, with $\lvert\mathcal{P}\rvert$ the number of blocks and $\lvert\beta\rvert$ the size of a block $\beta$; and, in (iv), its form $D^i[F(q)]=\sum_{s=1}^iF^{(s)}(q)B_{i,s}(Dq,\dots,D^{i-s+1}q)$ with the partial Bell polynomials $B_{i,s}$.

We first control displacements and first derivatives. By Lemma~\ref{lem:bou_val_for}(a) at $\theta_*$, the base variation obeys $\lvert t_l(0)\rvert\le L_j\rho_l$, and $L_j^{-1}\le\rho_1$. By Lemma~\ref{lem:der_bou_the}, $\lvert\phi^{(s)}\rvert\le(2\pi)^{s-1}M$ for $s\ge2$. Put $\Delta_l=q_l(\xi)-q_l(0)$. Taylor's formula for the recurrence of the positions (Lemma~\ref{lem:wro_and_cen:2}(i)) and propagation by the base step matrices, of norm at most $M$ per step in the maximum norm (Lemma~\ref{lem:nor_and_gro}(i)), give
\[
\lvert\Delta_l\rvert\le\rho_l+\pi\sum_{i=1}^{l-1}M^{l-i}\Delta_i^2 .
\]
Induction proves $\lvert\Delta_l\rvert\le2\rho_l$, since
\[
\sum_{i<l}M^{l-i}\rho_i^2\le\rho_l^2\frac{\gamma}{1-\gamma},\qquad\frac{4\pi\gamma}{1-\gamma}\le1 .
\]
Comparing the varied first derivative with the base variation then gives
\[
\lvert Dq_l\rvert\le\rho_l+2\pi\sum_{i<l}M^{l-i}\lvert\Delta_i\rvert\lvert Dq_i\rvert\le\rho_l+8\pi\rho_l^2\frac{\gamma}{1-\gamma}\le2\rho_l .
\]
In each induction the forcing depends only on indices below $l$; the initial data at indices zero and one have the required bounds. This proves (i).

We next prove (ii), by induction on $p$. For $p\ge2$, differentiation of the recurrence of the positions gives
\[
D^pq_{l+1}=v(q_l)D^pq_l-D^pq_{l-1}+F_{p,l},
\]
where $F_{p,l}$ is the sum of the terms of Faà di Bruno's formula for $D^p[\phi(q_l)]$ over the partitions with at least two blocks. Assume (ii) for the orders below $p$. The partition into singletons contributes at most $(2\pi)^{p-1}M(2\rho_l)^p$, by (i). Every other partition with at least two blocks has a block $\beta$ of size $\lvert\beta\rvert\ge2$, whose factor is at most $10^{\lvert\beta\rvert}\rho_l^{\lvert\beta\rvert}M^{1-\lvert\beta\rvert c}\le10^{\lvert\beta\rvert}\gamma\rho_l^{\lvert\beta\rvert}$, while each block of size one contributes at most $2\rho_l$; so such a partition contributes at most $10^5\gamma(2\pi)^{p-1}M(2\rho_l)^p$. A set of at most five elements has at most $52$ partitions, and $52\cdot10^5\gamma\le10^{-30}$; hence
\[
\lvert F_{p,l}\rvert\le(1+10^{-30})(2\pi)^{p-1}M2^p\rho_l^p .
\]
The initial higher derivatives vanish. Propagating the forcing by the varied step matrices, of norm at most $M$ per step (Lemma~\ref{lem:nor_and_gro}(i)), and using $\rho_i^p=\rho_l^pM^{-pc(l-i)}$, the geometric series gives
\[
\lvert D^pq_l\rvert\le(1+10^{-30})(2\pi)^{p-1}2^p\rho_l^p\sum_{i<l}M^{(l-i)(1-pc)}\le2(1+10^{-29})(4\pi)^{p-1}\rho_l^pM^{1-pc}\le10^p\rho_l^pM^{1-pc},
\]
since $M^{1-pc}\le\gamma$ and $2(4\pi)^{p-1}\le10^p/2$ for $2\le p\le5$. This proves (ii).

For (iii), by Lemma~\ref{lem:wro_and_cen:2}(i), $t_l(\xi)=\pm L_jDq_l(\xi)$, with one sign on the path, and $\lvert Dq_j(0)\rvert=1$. By (ii) with $p=2$, $\lvert Dq_j(\xi)-Dq_j(0)\rvert\le100\gamma<\frac12$. Together with (i), this proves (iii).

For (iv), write $Y_1=Dq_{j-1}$ and $Y_0=Dq_j$, so that $E=Y_1/Y_0$. By (i) and (iii), $\lvert Y_1\rvert\le2M^{-c}$ and $\frac12\le\lvert Y_0\rvert\le\frac32$; by (ii), with $\rho_j=1$ and $\rho_{j-1}=M^{-c}$, for $1\le i\le4$
\[
\lvert D^iY_0\rvert\le10^{i+1}\gamma,\qquad\lvert D^iY_1\rvert\le10^{i+1}M^{1-2(i+1)c}\le10^{i+1}\gamma M^{-c}.
\]
Faà di Bruno's formula for $Y_0^{-1}$, with $\lvert Y_0^{-s-1}\rvert\le2^{s+1}$ and $\gamma\le1$, gives $\lvert D^i(Y_0^{-1})\rvert\le\gamma\Omega_i$ for $1\le i\le4$, where
\[
\Omega_i=\sum_{s=1}^is!\,2^{s+1}B_{i,s}\bigl(10^2,10^3,\dots,10^{i-s+2}\bigr).
\]
Leibniz' rule then gives $\lvert D^pE\rvert\le\gamma M^{-c}C_p$ for $1\le p\le4$, with
\[
C_p=2\cdot10^{p+1}+2\Omega_p+\sum_{i=1}^{p-1}\binom pi10^{p-i+1}\Omega_i .
\]
Since $C_p<10^{12}$ for $1\le p\le4$ and $10^{12}\gamma\le1$, $\lvert D^pE\rvert\le M^{-c}$. Finally $\lvert E\rvert\le\lvert Y_1\rvert/\lvert Y_0\rvert\le4M^{-c}$. This proves (iv).

Explicit values. $\Omega_1,\dots,\Omega_4$ are $400$, $164000$, $1.0084\cdot10^8$ and $8.26724\cdot10^{10}$, and $C_1,\dots,C_4$ are $1000$, $410000$, $2.521\cdot10^8$ and $2.06681\cdot10^{11}$.
\end{proof}

\begin{lemma}
 Let $M\ge\exp(100)$. Fix $\iota$, $q_\iota$, a base coordinate $\theta_*$, an outward path, and an index $j\ge2$ on it that is a hyperbolic time at $\theta_*$ for some normalized solution, and let $L_j$, $\xi$ and $D$ be as in Definition~\ref{def:sin_pat_res}. Then $R_j$ exists and is smooth on $\lvert\xi\rvert\le1$, and
\[
L_j^2\lvert D^4R_j\rvert\le10^9M^{-72/25}.
\]
\label{lem:fou_der_the:1}
\end{lemma}

\begin{proof}
Put $c=c_{\rm hyp}$ (Convention~\ref{conv:not_and_con}) and $\gamma=M^{1-2c}\le\exp(-98)$. We write the boundary value as an entry of a Green matrix. By Lemma~\ref{lem:der_bou_res}(iii), $\lvert t_j\rvert\ge L_j/2$ on $\lvert\xi\rvert\le1$. Let $J$ be the matrix indexed by $0<l,p<j$ with diagonal entries $J_{pp}=v(q_p)$, entries $-1$ next to the diagonal and $0$ elsewhere. The sequence $t$ of Lemma~\ref{lem:gre_mat_fin} with $(a,b)=(0,j)$ is the variation sequence \eqref{eq:cs_variation}, so that lemma gives $\det J=t_j\ne0$, a unique solution $U$ of the path recurrence on $\{0,\dots,j\}$ with $U_0=1$ and $U_j=0$, and
\[
G_{lp}=(J^{-1})_{lp}=t_{\min(l,p)}U_{\max(l,p)}\qquad(0<l,p<j).
\]
By linearity and \eqref{eq:cs_variation}, $\widetilde t+R_jt$ solves the path recurrence with values $1$, $R_j$ and $\widetilde t_j+R_jt_j=0$ at the indices $0$, $1$, $j$; hence $U=\widetilde t+R_jt$ by uniqueness, $R_j=U_1=G_{11}$, and $G_{1a}=G_{a1}=U_a$. All entries of $J$ are smooth in $\xi$, hence so are $G$ and $R_j$.

We next bound the entries of $G$. Put $\ell_a=j-a$ for $0<a<j$. By \eqref{eq:cs_wronskian} (Lemma~\ref{lem:wro_and_cen:2}(ii)) with $U_j=0$, $(U_j,U_{j-1})=(0,1/t_j)$. Backward propagation from this pair, by steps of norm at most $M$ in the maximum norm (Lemma~\ref{lem:nor_and_gro}(i)), and $\lvert t_l\rvert\le2L_j\rho_l$ with $\rho_l=M^{-c(j-l)}$ (Lemma~\ref{lem:der_bou_res}(iii)) give
\[
L_j\lvert U_a\rvert\le2M^{\ell_a-1},\qquad\lvert G_{ab}\rvert\le4M^{\min(\ell_a,\ell_b)-1-c\max(\ell_a,\ell_b)}.
\]

We bound the derivatives of the diagonal. By Faà di Bruno's formula over set partitions, Lemma~\ref{lem:der_bou_res}(i), (ii) and Lemma~\ref{lem:der_bou_the}, for $1\le p\le4$,
\[
\lvert D^pv(q_a)\rvert\le V_pM^{1-cp\ell_a},\qquad V_p=(4\pi)^p(1+10^{-30}):
\]
the partition into singletons contributes at most $(2\pi)^pM(2M^{-c\ell_a})^p=(4\pi)^pM^{1-cp\ell_a}$, each of the other at most fourteen partitions has a block $\beta$ of size $\lvert\beta\rvert\ge2$ and contributes at most $10^4\gamma$ times this, and $15\cdot10^4\gamma\le10^{-30}$.

We now differentiate the inverse. Induction from $DG=-G(DJ)G$ gives
\[
D^4G=\sum_{\substack{p_1+\dots+p_s=4\\ p_i\ge1}}(-1)^s\frac{4!}{p_1!\cdots p_s!}\,GJ^{(p_1)}G\cdots J^{(p_s)}G .
\]
The matrices $J^{(p_i)}$ are diagonal, with entries $D^{p_i}v(q_a)$. At the $(1,1)$ entry, write each factor $J^{(p_i)}$ as a sum over its sites $a_i$, and put $\ell_i=j-a_i\ge1$. The preceding estimates bound the absolute value of each site term, multiplied by $L_j^2$, by
\[
4^s\prod_{i=1}^sV_{p_i}\,M^H,\qquad H=-1+\ell_1+\ell_s+\sum_{i=1}^{s-1}\bigl\{\min(\ell_i,\ell_{i+1})-c\max(\ell_i,\ell_{i+1})\bigr\}-c\sum_ip_i\ell_i .
\]
Using $\min(x,y)-c\max(x,y)=-\lvert x-y\rvert+(1-c)\max(x,y)$, the total variation $\sum_i\lvert\ell_i-\ell_{i+1}\rvert$ is at least $\ell_1+\ell_s-2\min_i\ell_i$, the sum of adjacent maxima is at most $2\sum_ip_i\ell_i$, and $4\min_i\ell_i\le\sum_ip_i\ell_i$. Consequently
\[
H\le-1+\Bigl(\frac52-3c\Bigr)\sum_ip_i\ell_i=-1-\frac{47}{100}\sum_ip_i\ell_i .
\]
Summing over the sites, the site sum of one composition is at most
\[
4^s\prod_iV_{p_i}\,M^{-72/25}\prod_i\bigl(1-M^{-47p_i/100}\bigr)^{-1},
\]
and each factor $(1-M^{-47p_i/100})^{-1}$ is at most $1+10^{-20}$, since $M^{-47/100}\le\exp(-47)$. Summing over the compositions of $4$, and since there are at most four factors of each kind,
\[
\begin{aligned}
L_j^2\lvert D^4R_j\rvert&\le(1+10^{-19})(4\pi)^4\Bigl(\sum_{\substack{p_1+\dots+p_s=4\\ p_i\ge1}}\frac{4!\,4^s}{p_1!\cdots p_s!}\Bigr)M^{-72/25}\\
&=(1+10^{-19})\,8676\,(4\pi)^4M^{-72/25}<2.2\cdot10^8M^{-72/25},
\end{aligned}
\]
which proves the bound.

Explicit values. The composition sum is $8676=4!\,[x^4]\bigl(1-4(\exp(x)-1)\bigr)^{-1}$, and $8676(4\pi)^4<2.17\cdot10^8$.
\end{proof}

\begin{definition}
 Let $(\iota,q_\iota,\theta_*,j_+,j_-)$ be a local test. Put $h=L(\theta-\theta_*)$, and number the two paths $i=1,2$ so that the first has the larger base scale (the plus path if $L_+=L_-$); $j_i$ is the selected index on path $i$, so that $\lvert t_{j_1}(\theta_*)\rvert=L$ on path $1$. The \textbf{rescaled phases} and \textbf{rescaled quotients} are
\[
Q_i(h)=q_{j_i}\ \text{on path }i,\qquad E_i(h)=\frac{t_{j_i-1}}{t_{j_i}}\ \text{on path }i,
\]
with $b_i=Q_i'(0)$, $d_i=Q_i''(0)$, $\omega_i=Q_i(0)$ and the \textbf{slope ratio} $r=\lvert b_2\rvert$. The \textbf{reciprocal terms} are, on their domains,
\[
\Phi_i(h)=\frac1{(Q_i'(h))^2\bigl(\cos(2\pi Q_i(h))+(2-E_i(h))/K\bigr)}\qquad(i=1,2).
\]
The \textbf{boundary function} and the \textbf{boundary cubic} are
\[
\Lambda(h)=KL^2\bigl(v(q_\iota)-R_{j_1}(\theta)-R_{j_2}(\theta)\bigr),\qquad P(h)=\sum_{p=0}^{3}\frac{\Lambda^{(p)}(0)}{p!}h^p,
\]
where $R_{j_i}$ is the boundary value on path $i$ and $\theta=\theta_*+h/L$. The letter $\Phi$ without index denotes another object, the cap of the balance argument; the two never occur together.
\label{def:res_pha_and}
\end{definition}

\begin{lemma}
 Let $M\ge\exp(100)$, and write $\lVert f\rVert_{C^1([-1,1])}=\sup_{[-1,1]}\lvert f\rvert+\sup_{[-1,1]}\lvert f'\rvert$. For every local test, $i=1,2$ and $h\in[-1,1]$:
\begin{align*}
Q_i(h)&=\omega_i+b_ih+\tfrac12d_ih^2+\mathcal E_i(h),&\lVert\mathcal E_i\rVert_{C^1([-1,1])}&\le10^3\lvert b_i\rvert^3M^{1-3c_{\rm hyp}},\\
\lvert d_i\rvert&\le100\lvert b_i\rvert^2M^{1-2c_{\rm hyp}},&\lvert b_1\rvert&=1,\quad\lvert b_2\rvert=r\in[M^{-3/5},1],\\
\lvert Q_i^{(p)}(h)\rvert&\le10^p\lvert b_i\rvert^pM^{1-pc_{\rm hyp}}\quad(2\le p\le5),&\tfrac12\lvert b_i\rvert&\le\lvert Q_i'(h)\rvert\le\tfrac32\lvert b_i\rvert,\\
\lvert\partial_h^pE_i(h)\rvert&\le4M^{-c_{\rm hyp}}\lvert b_i\rvert^p\quad(0\le p\le4).&&
\end{align*}
\label{lem:pha_bou_loc}
\end{lemma}

\begin{proof}
Put $c=c_{\rm hyp}$. We read the bounds of Lemma~\ref{lem:der_bou_res} on each path in the variable $h$. On path $i$, the single-path variable of Definition~\ref{def:sin_pat_res} with $j=j_i$ is $\xi_i=L_{j_i}(\theta-\theta_*)=(L_{j_i}/L)h$, where $L_{j_i}=\lvert t_{j_i}(\theta_*)\rvert$; and $Q_i(h)=q_{j_i}(\xi_i)$, $E_i(h)=E(\xi_i)$. Write $D_i=d/d\xi_i$. By Lemma~\ref{lem:wro_and_cen:2}(i), $\lvert D_iq_{j_i}(0)\rvert=\lvert t_{j_i}(\theta_*)\rvert/L_{j_i}=1$, so $\lvert b_i\rvert=\lvert Q_i'(0)\rvert=L_{j_i}/L$ and $\partial_h=\lvert b_i\rvert D_i$. By the numbering $\lvert b_1\rvert=1$, and $\lvert h\rvert\le1$ gives $\lvert\xi_i\rvert\le\lvert b_i\rvert\le1$.

Since $j_i$ is a hyperbolic time at $\theta_*$ on path $i$, Lemma~\ref{lem:der_bou_res} applies on $\lvert\xi_i\rvert\le1$ with $j=j_i$ and $\rho_{j_i}=1$. It gives $\lvert Q_i^{(p)}\rvert=\lvert b_i\rvert^p\lvert D_i^pq_{j_i}\rvert\le10^p\lvert b_i\rvert^pM^{1-pc}$ for $2\le p\le5$; $\lvert Q_i'\rvert=\lvert b_i\rvert\lvert D_iq_{j_i}\rvert\in[\frac12\lvert b_i\rvert,\frac32\lvert b_i\rvert]$, by its item (iii) and $\lvert t_{j_i}\rvert=L_{j_i}\lvert D_iq_{j_i}\rvert$; $\lvert d_i\rvert=\lvert Q_i''(0)\rvert\le100\lvert b_i\rvert^2M^{1-2c}$; $\lvert E_i\rvert\le4M^{-c}$; and $\lvert\partial_h^pE_i\rvert=\lvert b_i\rvert^p\lvert D_i^pE\rvert\le M^{-c}\lvert b_i\rvert^p$ for $1\le p\le4$.

The remainder $\mathcal E_i(h)=Q_i(h)-\omega_i-b_ih-\frac12d_ih^2$ satisfies $\lvert\mathcal E_i(h)\rvert\le\frac{\lvert h\rvert^3}6\sup\lvert Q_i'''\rvert$ and $\lvert\mathcal E_i'(h)\rvert\le\frac{h^2}2\sup\lvert Q_i'''\rvert$, so $\lVert\mathcal E_i\rVert_{C^1([-1,1])}\le\frac23\cdot10^3\lvert b_i\rvert^3M^{1-3c}\le10^3\lvert b_i\rvert^3M^{1-3c}$.

Finally, at $\theta_*$ both indices are hyperbolic times, so $L_\pm=M^{\tau^\pm_{j_\pm}}$ by \eqref{eq:cs_dominated_variation} (Lemma~\ref{lem:bou_val_for}(a)), and the matching condition \eqref{eq:cs_matching} gives $r=L_{j_2}/L=M^{-\lvert\tau^+_{j_+}-\tau^-_{j_-}\rvert}\ge M^{-3/5}$; also $r\le1$ by the numbering. The value $\frac35$ of the matching tolerance enters the argument only through this lower bound on $r$; the matched pairs used in the proof of the scarcity of cancellations satisfy the matching condition with $\frac12$.
\end{proof}

\begin{lemma}
 Let $M\ge\exp(100)$. For every local test and every successful point $h$, the reciprocal terms $\Phi_1(h)$ and $\Phi_2(h)$ are defined, and
\[
P(h)=\Phi_1(h)+\Phi_2(h)+\varepsilon(h),\qquad\lvert\varepsilon(h)\rvert\le10^8\bigl(M^{-47/25}+r^{-2}M^{1-3c_{\rm hyp}}\bigr).
\]
\label{lem:cub_rel_suc:1}
\end{lemma}

\begin{proof}
The relation is the centre equation in rescaled form: at a successful point the boundary-value formula writes each outward value $u^{(i)}_1$ as the boundary value $R_{j_i}$ plus a reciprocal term, the centre equation adds the two, and after rescaling the difference between $v(q_\iota)$ and the two boundary values is the boundary function $\Lambda$, which the fourth-derivative bound replaces by its cubic Taylor polynomial $P$. Put $c=c_{\rm hyp}$. As in the proof of Lemma~\ref{lem:pha_bou_loc}, on path $i$ put $L_{j_i}=\lvert t_{j_i}(\theta_*)\rvert=\lvert b_i\rvert L$, $\xi_i=(L_{j_i}/L)h$ and $D_i=d/d\xi_i$, so that $\partial_h=\lvert b_i\rvert D_i$.

We first derive the relation at a successful point. Let $h\in\mathcal S$ and $\theta=\theta_*+h/L$. There is a normalized solution $u$ for which $j_1$ and $j_2$ are hyperbolic times on their paths at $\theta$; let $u^{(i)}$ be its restriction to path $i$. By Lemma~\ref{lem:bou_val_for}(c) on each path,
\[
u_1^{(i)}=R_{j_i}(\theta)+\frac1{t_{j_i}(\theta)^2\bigl(v(Q_i(h))-E_i(h)\bigr)}+\varepsilon_{j_i},\qquad\lvert\varepsilon_{j_i}\rvert\le\frac{2M^{-3c}}{t_{j_i}(\theta)^2},
\]
where $v_{j_i}=v(q_{j_i})=v(Q_i(h))$ and $t_{j_i-1}/t_{j_i}=E_i(h)$. Here $v_{j_i}-E_i=t_{j_i+1}/t_{j_i}\ne0$, because $\lvert t_{j_i+1}\rvert\ge M^c\lvert t_{j_i}\rvert>0$ at a hyperbolic time by \eqref{eq:cs_dominated_variation}, and $Q_i'(h)\ne0$ by Lemma~\ref{lem:pha_bou_loc}; since $v(x)=2+K\cos(2\pi x)$ (Convention~\ref{conv:sta_map_par}), $v(Q_i)-E_i=K\bigl(\cos(2\pi Q_i)+(2-E_i)/K\bigr)$, so the denominators of $\Phi_1(h)$ and $\Phi_2(h)$ are nonzero. The centre equation \eqref{eq:cs_centre_equation} (Lemma~\ref{lem:wro_and_cen:2}(iii)) gives $u_1^{(1)}+u_1^{(2)}=v(q_\iota)$. Multiply by $KL^2$. Since $\lvert t_{j_i}(\theta)\rvert=L_{j_i}\lvert D_iq_{j_i}(\xi_i)\rvert=L\lvert Q_i'(h)\rvert$ (Lemma~\ref{lem:wro_and_cen:2}(i)),
\[
\Lambda(h)=\Phi_1(h)+\Phi_2(h)+KL^2(\varepsilon_{j_1}+\varepsilon_{j_2}),
\]
and, by $\lvert Q_i'\rvert\ge\lvert b_i\rvert/2$ and $K<M$,
\[
KL^2\lvert\varepsilon_{j_i}\rvert\le\frac{2KL^2M^{-3c}}{L^2(Q_i'(h))^2}\le\frac{8KM^{-3c}}{b_i^2}\le8\lvert b_i\rvert^{-2}M^{1-3c},
\]
so $KL^2\lvert\varepsilon_{j_1}+\varepsilon_{j_2}\rvert\le8(1+r^{-2})M^{1-3c}\le16r^{-2}M^{1-3c}$.

We next replace $\Lambda$ by its cubic Taylor polynomial. Each $R_{j_i}$ is smooth in $\xi_i$ on $\lvert\xi_i\rvert\le1$, with $L_{j_i}^2\lvert D_i^4R_{j_i}\rvert\le10^9M^{-72/25}$, by Lemma~\ref{lem:fou_der_the:1}; so $\Lambda$ is smooth on $[-1,1]$, and
\[
KL^2\lvert\partial_h^4R_{j_i}\rvert=KL^2\lvert b_i\rvert^4\lvert D_i^4R_{j_i}\rvert\le\frac{KL^2\lvert b_i\rvert^4\,10^9M^{-72/25}}{b_i^2L^2}=10^9K\lvert b_i\rvert^2M^{-72/25}\le10^9M^{-47/25}.
\]
Hence $\sup_{[-1,1]}\lvert\Lambda^{(4)}\rvert\le2\cdot10^9M^{-47/25}$, and Taylor's formula at $0$ gives, for $\lvert h\rvert\le1$,
\[
\lvert\Lambda(h)-P(h)\rvert\le\frac{h^4}{24}\sup\lvert\Lambda^{(4)}\rvert\le\frac{10^9}{12}M^{-47/25}.
\]
Therefore $\varepsilon(h)=P(h)-\Phi_1(h)-\Phi_2(h)=\bigl(P(h)-\Lambda(h)\bigr)+KL^2(\varepsilon_{j_1}+\varepsilon_{j_2})$ satisfies
\[
\lvert\varepsilon(h)\rvert\le\frac{10^9}{12}M^{-47/25}+16r^{-2}M^{1-3c}\le10^8\bigl(M^{-47/25}+r^{-2}M^{1-3c}\bigr),
\]
since $10^9/12<10^8$.
\end{proof}

\subsection{Small covers of successful sets} \label{subsec:sma_cov_suc}

Every successful set is covered by boundedly many intervals of total length at most $\eta$ once $M$ is large in terms of $\eta$ (Theorem~\ref{thm:sma_cov_the}); for $\eta=\eta_{\rm cov}$ it suffices that $M\ge\exp(1.1\cdot10^{11})$. The proof distinguishes three cases. At a small slope ratio between the two phases (Theorem~\ref{thm:sma_cov_suc}), the fourth derivative of the slow reciprocal term is negligible, dominant or nearly constant (Lemma~\ref{lem:fou_der_the:2}), and a sublevel estimate under a lower bound on a derivative covers the successful set. At a near-opposite alignment of the two phases, the two reciprocal terms nearly cancel, and the relation is compared with a normalized model (Lemma~\ref{lem:the_bou_cub}) whose sublevel sets have small covers (Theorem~\ref{thm:sma_cov_for}). In the generic case (Theorem~\ref{thm:sma_cov_for:1}), clearing the cosine denominators turns the relation into a sublevel condition for an entire function that is not small on $[-1,1]$. The tools, collected first, are a Remez-type inequality \cite{locbib:Remez1936} (Lemma~\ref{lem:rem_ine_and}), used in all three cases; a sublevel estimate under a lower bound on a derivative (Lemma~\ref{lem:sub_set_est}), a variant, with an exceptional set and a count of components, of the one-dimensional sublevel set estimate \cite{locbib:CarberyChristWright1999}; and a sublevel estimate for bounded entire functions (Lemma~\ref{lem:sub_set_ent}). The cover constants are explicit, but a first reading needs only that they are finite.

\begin{lemma}
 Let $P$ be a real polynomial of degree at most $d$ and let $E\subset[-1,1]$ be measurable with measure $s>0$. Then
\[
\sup_{[-1,1]}\lvert P\rvert\le(24/s)^d\sup_E\lvert P\rvert .
\]
For every real cubic $P$, $\lVert P\rVert_{\rm coeff}\le24\sup_{[-1,1]}\lvert P\rvert$. For every $\zeta\ge0$, the sublevel set $\{x\in[-1,1]:\lvert P(x)\rvert\le\zeta\}$ of a real polynomial $P$ of degree at most $d$ is a union of at most $d+1$ intervals. If $P$ is a real cubic with $\lVert P\rVert_{\rm coeff}>H$ for some $H>0$, then for every $\zeta>0$ the set $\{h\in[-1,1]:\lvert P(h)\rvert\le\zeta\}$ is a union of at most four intervals of total length less than $24(24\zeta/H)^{1/3}$.
\label{lem:rem_ine_and}
\end{lemma}

\begin{proof}
The first assertion follows from Remez's inequality \cite{locbib:Remez1936} (see also \cite[Section 5.1]{locbib:BorweinErdelyi1995}): $\sup_{[-1,1]}\lvert P\rvert$ is at most the value at $4/s-1$ of the Chebyshev polynomial of degree $d$, times $\sup_E\lvert P\rvert$. Since that Chebyshev polynomial is at most $(2x)^d$ at every $x\ge1$, the right side is at most $(8/s)^d\sup_E\lvert P\rvert\le(24/s)^d\sup_E\lvert P\rvert$.

For the coefficient bound, interpolate a cubic at $-1,-1/3,1/3,1$. Each denominator in the Lagrange basis has magnitude at least $16/27$, and each coefficient of its cubic numerator has magnitude at most three. Thus every coefficient of $P$ is bounded by $4\cdot3\cdot(27/16)\sup_{[-1,1]}\lvert P\rvert<24\sup_{[-1,1]}\lvert P\rvert$.

For the sublevel sets, a nonconstant real polynomial of degree at most $d$ has at most $d-1$ critical points. On each of the resulting at most $d$ monotonicity intervals, the inverse image of $[-\zeta,\zeta]$ is an interval. So the sublevel set has at most $d$ components when $P$ is nonconstant, and at most one when $P$ is constant; both are within the asserted count $d+1$.

For the last assertion, the set is a union of at most four intervals by the preceding count with $d=3$; let $s$ be its measure. If $s>0$, the first two assertions give $H<\lVert P\rVert_{\rm coeff}\le24\sup_{[-1,1]}\lvert P\rvert\le24(24/s)^3\zeta$, that is $s^3<24^4\zeta/H$, so $s<24(24\zeta/H)^{1/3}$. If $s=0$ there is nothing to prove.
\end{proof}

\begin{lemma}
Let $p\ge1$, let $I\subset\mathbb{R}$ be an interval, let $\lambda>0$, and let $G\in C^p(I)$ satisfy $\lvert G^{(p)}(x)\rvert\ge\lambda$ for every $x\in I$. Then for every $\zeta>0$ the set $\{x\in I:\lvert G(x)\rvert\le\zeta\}$ has at most $p+1$ connected components and Lebesgue measure
\[
\bigl\lvert\{x\in I:\lvert G(x)\rvert\le\zeta\}\bigr\rvert\le\bigl(2^{p+1}-2\bigr)\Bigl(\frac\zeta\lambda\Bigr)^{1/p}.
\]
More generally, if $\lvert G^{(p)}\rvert\ge\lambda$ holds on $I\setminus E$, where $E\subset I$ is a union of at most $N_E$ intervals and $G\in C^p(I\setminus E)$, then $\{x\in I:\lvert G(x)\rvert\le\zeta\}$ is covered by at most $N_E+(N_E+1)(p+1)$ intervals and
\[
\bigl\lvert\{x\in I:\lvert G(x)\rvert\le\zeta\}\bigr\rvert\le\lvert E\rvert+(N_E+1)\bigl(2^{p+1}-2\bigr)\Bigl(\frac\zeta\lambda\Bigr)^{1/p}.
\]
\label{lem:sub_set_est}
\end{lemma}

\begin{proof}
The measure bound is a variant of the one-dimensional sublevel set estimate \cite[Proposition 2.1]{locbib:CarberyChristWright1999}, which is stated under $u^{(p)}\ge1$ with the constant $2\exp(1)((p+1)!)^{1/p}$; the proof below is self-contained and gives the constant $2^{p+1}-2$, the component count and the version with exceptional intervals.

For $p=1$, continuity and nonvanishing of $G'$ make $G$ strictly monotone, and the inverse image of $[-\zeta,\zeta]$ has length at most $2\zeta/\lambda=(2^2-2)(\zeta/\lambda)$.

For $p\ge2$ we prove the measure bound by induction on $p$. The derivative $H=G^{(p-1)}$ is strictly monotone, so for $\mu>0$ the set $\{\lvert H\rvert<\mu\}$ is an interval of length at most $2\mu/\lambda$. Its complement consists of at most two intervals on which $\lvert G^{(p-1)}\rvert\ge\mu$. The induction hypothesis therefore gives
\[
\bigl\lvert\{\lvert G\rvert\le\zeta\}\bigr\rvert\le\frac{2\mu}\lambda+2\bigl(2^p-2\bigr)\Bigl(\frac\zeta\mu\Bigr)^{1/(p-1)}.
\]
Take $\mu=\lambda^{(p-1)/p}\zeta^{1/p}$; then $\mu/\lambda$ and $(\zeta/\mu)^{1/(p-1)}$ both equal $(\zeta/\lambda)^{1/p}$, and $2+2(2^p-2)=2^{p+1}-2$.

For the component bound, repeated Rolle's theorem shows that $G'$ has at most $p-1$ distinct zeros: $p$ such zeros would force a zero of $G^{(p)}$. Cutting at these zeros partitions $I$ into at most $p$ intervals on whose interiors $G$ is strictly monotone. The inverse image of $[-\zeta,\zeta]$ on each is an interval. Their union has at most $p$ components, some possibly single points, which is within the asserted bound $p+1$.

Finally remove the at most $N_E$ exceptional intervals. Their complement has at most $N_E+1$ interval components. Apply the preceding bounds on each component and add the exceptional intervals themselves. This proves the covering and measure estimates in the general statement.
\end{proof}

\begin{lemma}
Let $\eta\in(0,1]$ and $R=192/\eta$. Let $F$ be an entire function, real on $\mathbb{R}$, with $\lvert F\rvert\le B^{\mathbb C}$ on the closed disc $\{h\in\mathbb{C}:\lvert h\rvert\le R\}$ and $\sup_{[-1,1]}\lvert F\rvert\ge u$, where $0<u\le B^{\mathbb C}$. Put
\[
D=\Bigl\lceil\frac{\log(1000B^{\mathbb C}/u)}{\log4}\Bigr\rceil,\qquad s=\frac u{1000}\Bigl(\frac\eta{48}\Bigr)^D.
\]
Then $\{h\in[-1,1]:\lvert F(h)\rvert\le s\}$ is contained in at most $D+1$ intervals of total length less than $\eta/2$.
\label{lem:sub_set_ent}
\end{lemma}

\begin{proof}
Let $Q$ be the Taylor polynomial of $F$ of degree $D$ at $0$; it is real, since $F$ is real on $\mathbb{R}$. We first show that $Q$ is close to $F$ on $[-1,1]$. Cauchy's estimates give $\lvert F^{(p)}(0)\rvert/p!\le B^{\mathbb C}R^{-p}$, so on $[-1,1]$
\[
\lvert F-Q\rvert\le\sum_{p>D}B^{\mathbb C}R^{-p}=\frac{B^{\mathbb C}R^{-D}}{R-1}\le\frac u{1000}\Bigl(\frac4R\Bigr)^D=s,
\]
since $R-1\ge1$ and $4^D\ge1000B^{\mathbb C}/u$; here $4/R=\eta/48$. Hence $\sup_{[-1,1]}\lvert Q\rvert\ge u-s\ge u/2$, and the set of the statement lies in $E=\{h\in[-1,1]:\lvert Q(h)\rvert\le2s\}$.

We next show that $\lvert E\rvert<\eta/2$. Otherwise the Remez inequality of Lemma~\ref{lem:rem_ine_and}, applied to $Q$ on $E$, gives $\sup_{[-1,1]}\lvert Q\rvert\le(48/\eta)^D\,2s=2u/1000<u/2$, a contradiction. Finally, $E$ is a union of at most $D+1$ intervals, by the same lemma.
\end{proof}

The next definition fixes the affine phases, used in all three cases, and the normalized model of the opposite case.

\begin{definition}
 For $\omega,b\in\mathbb{R}$, the \textbf{affine phase} of $(\omega,b)$ is $a(h)=\omega+bh$. Let $\omega\in\mathbb{R}$ and $b\in\{-1,1\}$, let $p$ be a real quadratic polynomial and let $\beta\in[0,1]$, with
\[
\max(\lVert p\rVert_{\rm coeff},\beta)=1 .
\]
With $a$ the affine phase of $(\omega,b)$, the \textbf{normalized opposite-cosine model} is
\[
\mathcal R_{\omega,b,p,\beta}(h)=\frac{2p'(h)}{b^3}\sec(2\pi a(h))-\frac{2\pi p(h)}{b^2}\sin(2\pi a(h))\sec^2(2\pi a(h))-\frac{4\beta}{b^2}\sec^2(2\pi a(h)),
\]
defined where $\cos(2\pi a(h))\ne0$.
\label{def:nor_opp_cos:1}
\end{definition}

At a small slope ratio the second phase is slow. Its reciprocal term is written through the slow denominator, and the fourth derivative of the fast term is compared with that of the secant of the affine fast phase.

\begin{definition}
 For a local test, the \textbf{affine fast phase} is the affine phase $a$ of $(\omega_1,b_1)$ (Definition~\ref{def:nor_opp_cos:1}), that is $a(h)=\omega_1+b_1h$. The \textbf{slow denominator} and the \textbf{affine fast fourth derivative} are
\[
S(h)=\Bigl(\frac{Q_2'(h)}r\Bigr)^2\Bigl(\cos(2\pi Q_2(h))+\frac{2-E_2(h)}K\Bigr),
\]
\[
G(h)=(2\pi)^4\Bigl(\frac{24}{\cos^5(2\pi a(h))}-\frac{20}{\cos^3(2\pi a(h))}+\frac1{\cos(2\pi a(h))}\Bigr).
\]
The letter $G$ with a subscript denotes functionals of arrays in the balance argument; the two meanings never occur together.
\label{def:slo_den_and}
\end{definition}

\begin{lemma}
 Let $M\ge\exp(100)$, and consider a local test with slope ratio $r$. Then $\Phi_2=(r^2S)^{-1}$ wherever $\Phi_2$ is defined, $G$ is the fourth derivative of $\sec(2\pi a)$ away from its poles, and
\[
\lvert S\rvert\le2,\qquad\lvert S'\rvert\le10r,\qquad\lvert S^{(j)}\rvert\le10^{20}r^j\quad(2\le j\le4).
\]
At every point with $\lvert S\rvert\le1/100$, one has $5r\le\lvert S'\rvert\le7r$. For $\zeta\in(0,1/100]$ with $10^{20}M^{1-2c_{\rm hyp}}\le\zeta/2$,
\[
\lvert\Phi_1^{(4)}-G\rvert\le10^{100}\zeta^{-20}M^{1-2c_{\rm hyp}}\quad\text{on }\{\lvert\cos(2\pi a)\rvert\ge\zeta\}.
\]
Let $0<w_0,d_0\le10^{-110}$ and suppose $r\le w_0/100$. If $\lvert S(0)\rvert\ge w_0$, then $S$ has no zero and
\[
\lvert\Phi_2^{(4)}\rvert\le10^{100}r^2w_0^{-5}.
\]
If $\lvert S(0)\rvert<w_0$, put $w=\max(\lvert S(0)\rvert,r)$. The function $S$ is strictly monotone and has at most one zero. On its nonzero set the following alternatives hold:
\begin{enumerate}
\item If $r/w>d_0$, then
\[
\lvert\Phi_2^{(4)}\rvert\ge\frac{d_0^5}{100r^3}.
\]
\item If $r/w\le d_0$, then $S$ has no zero and, with $\alpha=24S'(0)^4/(r^2S(0)^5)$,
\[
\lvert\Phi_2^{(4)}-\alpha\rvert\le10^{100}(w_0+d_0)\lvert\alpha\rvert.
\]
\end{enumerate}
\label{lem:fou_der_the:2}
\end{lemma}

\begin{proof}
Put $c=c_{\rm hyp}$. The two identities hold by the definitions: $r^2S=(Q_2')^2\bigl(\cos(2\pi Q_2)+(2-E_2)/K\bigr)$ is the denominator of $\Phi_2$; and, since $\lvert b_1\rvert=1$ (Lemma~\ref{lem:pha_bou_loc}), the fourth derivative of $h\mapsto\sec(2\pi a(h))$ is $(2\pi)^4$ times $\sec^{(4)}=24\sec^5-20\sec^3+\sec$ evaluated at $2\pi a(h)$, which is $G$.

We next bound the slow denominator and the derivatives of the fast term. Write $\gamma=M^{1-2c}$, $\bar Y_2=Q_2'/r$ and $\sigma_i=(2-E_i)/K$. The bounds of Lemma~\ref{lem:pha_bou_loc} imply
\[
\lvert\bar Y_2(0)\rvert=1,\qquad\lvert\bar Y_2-\bar Y_2(0)\rvert\le100r\gamma,\qquad\lvert\bar Y_2^{(j)}\rvert\le10^5r^j\gamma\quad(1\le j\le4),
\]
\[
\lvert\sigma_i\rvert\le6/M,\qquad\lvert\sigma_i^{(j)}\rvert\le3\cdot10^6\lvert b_i\rvert^jM^{-1-c}\quad(1\le j\le4).
\]
Here $K=M-4\ge M/2$. Since $M\ge\exp(100)$, all terms involving $\gamma$ in these bounds are smaller than $10^{-29}$ after multiplication by $10^8$. The product rule applied to $S=\bar Y_2^2(\cos(2\pi Q_2)+\sigma_2)$ gives the displayed derivative bounds: each derivative of order $j\le4$ contributes a factor $r^j$, and the remaining coefficients are bounded by $10^{20}$. More explicitly the principal term of $S'$ is
\[
-2\pi r\bar Y_2^3\sin(2\pi Q_2);
\]
the other terms have absolute value at most $10^8r\gamma$. Thus $\lvert S'\rvert\le10r$. If $\lvert S\rvert\le1/100$, then $\lvert\cos(2\pi Q_2)\rvert\le11/1000$, since $\lvert\bar Y_2^2-1\rvert<10^{-30}$ and $\lvert\sigma_2\rvert<10^{-30}$. Hence the principal term has magnitude between $6r$ and $6.4r$, giving $5r\le\lvert S'\rvert\le7r$.

We justify the fast derivative estimate by a finite derivative norm. On any component of $\{\lvert\cos(2\pi a)\rvert\ge\zeta\}$ use
\[
\lVert f\rVert_4=\sum_{j=0}^4\frac{\sup\lvert f^{(j)}\rvert}{j!},
\]
which is submultiplicative. Put $Y_1=Q_1'$ and $U=Y_1^2(\cos(2\pi Q_1)+\sigma_1)$, so that $\Phi_1=U^{-1}$. The same bounds, Taylor's formula at zero, and the chain rule through order four yield
\[
\lVert U-\cos(2\pi a)\rVert_4\le10^{20}\gamma,\qquad\lVert U\rVert_4,\ \lVert\cos(2\pi a)\rVert_4\le10^5.
\]
For these bounds it suffices to use $\lVert Q_1-\omega_1-b_1h\rVert_4\le10^6\gamma$ and $\lVert Y_1-b_1\rVert_4\le10^6\gamma$ (in fact at most $900\gamma$ and $6600\gamma$), and $\lVert\cos(2\pi Q_1)-\cos(2\pi a)\rVert_4\le10^{15}\gamma$; the last follows by differentiating
\[
-2\pi(Q_1-\omega_1-b_1h)\int_0^1\sin\bigl(2\pi(\omega_1+b_1h+s(Q_1-\omega_1-b_1h))\bigr)\,ds .
\]
The hypothesis on $\zeta$ gives $\lvert U\rvert\ge\zeta/2$. The reciprocal derivative formula through order four gives
\[
\lVert U^{-1}\rVert_4,\ \lVert(\cos(2\pi a))^{-1}\rVert_4\le10^{25}\zeta^{-5}.
\]
Indeed the reciprocal Taylor coefficients at a point are bounded by $(2/\zeta)\sum_{j=0}^4(2\cdot10^5/\zeta)^j$, which is smaller than the stated bound. Using $U^{-1}-(\cos(2\pi a))^{-1}=(\cos(2\pi a)-U)U^{-1}(\cos(2\pi a))^{-1}$ and multiplying the resulting norm estimate by $4!$ proves the bound $10^{100}\zeta^{-20}\gamma$ for $\lvert\Phi_1^{(4)}-G\rvert$, since $G$ is the fourth derivative of $(\cos(2\pi a))^{-1}$.

We now use the explicit fourth derivative of a reciprocal:
\[
\Phi_2^{(4)}=\frac1{r^2}\Bigl(\frac{24(S')^4}{S^5}-\frac{36(S')^2S''}{S^4}+\frac{6(S'')^2+8S'S'''}{S^3}-\frac{S''''}{S^2}\Bigr).
\]
If $\lvert S(0)\rvert\ge w_0$, the bound $\lvert S'\rvert\le10r$ gives $\lvert S\rvert\ge w_0-10r\ge w_0/2$. Inserting this and the preceding derivative bounds into the formula proves $\lvert\Phi_2^{(4)}\rvert\le10^{100}r^2w_0^{-5}$.

Suppose $\lvert S(0)\rvert<w_0$. Then $\lvert S(h)\rvert\le w_0+10r<2w_0<1/100$ on the whole interval. The lower bound for $\lvert S'\rvert$ and continuity make $S$ strictly monotone. On the nonzero set, comparison with the first term in the reciprocal formula gives
\[
\Bigl\lvert\Phi_2^{(4)}-\frac{24(S')^4}{r^2S^5}\Bigr\rvert\le10^{70}\lvert S\rvert\,\Bigl\lvert\frac{24(S')^4}{r^2S^5}\Bigr\rvert .
\]
For example the relative size of the second term is at most $(3/2)\lvert S\rvert\lvert S''\rvert/\lvert S'\rvert^2\le10^{20}\lvert S\rvert$; the other relative sizes are bounded by $10^{41}\lvert S\rvert^2$ and $10^{21}\lvert S\rvert^3$. These are dominated by the stated bound since $\lvert S\rvert<1$. Since $\lvert S\rvert\le11w$, $\lvert S'\rvert\ge5r$ and $10^{70}\lvert S\rvert<1/2$, it follows that
\[
\lvert\Phi_2^{(4)}\rvert\ge\frac{12\cdot5^4r^2}{(11w)^5}\ge\frac{r^2}{100w^5}.
\]
If $r/w>d_0$ this proves the first alternative.

For the second alternative, $d_0<1$ implies $w=\lvert S(0)\rvert>0$. We have
\[
\Bigl\lvert\frac{S(h)}{S(0)}-1\Bigr\rvert\le10d_0,\qquad\Bigl\lvert\frac{S'(h)}{S'(0)}-1\Bigr\rvert\le10^{20}r/5\le10^{20}d_0 .
\]
Thus $S$ has no zero. The mean value theorem applied to $(x,y)\mapsto x^4y^{-5}$ on the box where both variables are within $1/2$ of one bounds the relative change of the leading reciprocal term by $10^{25}d_0$. The other terms in the reciprocal formula contribute at most $10^{72}w_0$ times that constant leading term. Both are bounded by $10^{100}(w_0+d_0)$, which proves the assertion.
\end{proof}

\begin{thm}
 For every $\eta\in(0,1]$ there are $r_\eta\in(0,1)$ and $M^{\rm sl}_\eta\ge\exp(100)$ such that, for $M\ge M^{\rm sl}_\eta$, the successful set of every local test with slope ratio $r\le r_\eta$ is contained in a union of at most $500$ intervals of total length at most $\eta$. For $\eta=\eta_{\rm cov}$ one may take $r_\eta=10^{-773}$ and $M^{\rm sl}_\eta=10^{794}$.
\label{thm:sma_cov_suc}
\end{thm}

\begin{proof}
We fix the constants. Put $\zeta=\eta/100$, $B=10^5\zeta^{-5}$,
\[
\lambda=\frac{(\eta/10^6)^{10}}{100B},\qquad w_0=d_0=\frac\lambda{10^{106}},
\]
\[
r_\eta=\min\Bigl\{\frac{w_0}{100},\Bigl(\frac{\lambda w_0^5}{10^{104}}\Bigr)^{1/2},\Bigl(\frac{d_0^5}{10^4B}\Bigr)^{1/3},10^{-106}\Bigr\},
\]
\[
M^{\rm sl}_\eta=\max\Bigl\{\exp(100),\Bigl(\frac{10^{105}\zeta^{-20}}\lambda\Bigr)^2,\Bigl(\frac{10^{40}}{\lambda\eta^4}\Bigr)^2\Bigr\},
\]
and let $M\ge M^{\rm sl}_\eta$ and $r\le r_\eta$. Then $w_0=d_0=\eta^{15}10^{-183}\le10^{-110}$ and $r\le w_0/100$, so Lemma~\ref{lem:fou_der_the:2} applies with these $w_0,d_0$ and with $\zeta$. We use the slow denominator $S$, the affine fast phase $a$ and the affine fast fourth derivative $G$ of Definition~\ref{def:slo_den_and}, and the reciprocal terms $\Phi_i$ of Definition~\ref{def:res_pha_and}.

We first bound a model fourth-derivative band uniformly in the additive constant. On the pole-safe set $\lvert\cos(2\pi a)\rvert\ge\zeta$,
\[
5(2\pi)^4\le\lvert G\rvert\le B .
\]
For every $\alpha\in\mathbb{R}$, the band
\[
\{\lvert\cos(2\pi a)\rvert\ge\zeta,\ \lvert G+\alpha\rvert\le4\lambda(1+\lvert\alpha\rvert)\}
\]
is contained in at most $60$ intervals of total length at most $\eta/1000$. If $\lvert\alpha\rvert>2B$, the band is empty: $\lvert G+\alpha\rvert\ge\lvert\alpha\rvert/2>(1+\lvert\alpha\rvert)/4$, while $4\lambda<1/4$. Otherwise clear the fifth cosine power and use the degree-five polynomial
\[
H(y)=24-20y^2+y^4+\frac\alpha{(2\pi)^4}y^5 .
\]
The band is contained in the pullback under $y=\cos(2\pi a)$ of $\{y\in[-1,1]:\lvert H(y)\rvert\le12B\lambda\}$. Since $H(0)=24$, the Remez inequality of Lemma~\ref{lem:rem_ine_and} bounds the measure $s$ of this latter set by
\[
s\le24(B\lambda)^{1/5}.
\]
Its number of interval components is at most six. The cosine has at most five monotonicity intervals on $[-1,1]$, since $\lvert b_1\rvert=1$ (Lemma~\ref{lem:pha_bou_loc}). Its pullback therefore lies in at most $6\cdot5=30$ intervals, so the count of $60$ intervals stated for the band holds with room to spare. For any measurable subset of $[-1,1]$ of measure $s$, its pullback has measure at most $8\sqrt s$. Indeed the change of variables on the two periods gives density $2/(\pi\sqrt{1-y^2})$, whose integral over a set of measure $s$ is bounded by splitting at zero and integrating the endpoint majorants $(1-\lvert y\rvert)^{-1/2}$. Consequently the model band has measure at most
\[
40(B\lambda)^{1/10}=\frac{40}{100^{1/10}}\,\frac\eta{10^6}<\frac\eta{1000}.
\]

We next establish a fourth-derivative lower bound for the actual reciprocal sum. Remove the fast pole strip $\{\lvert\cos(2\pi a)\rvert<\zeta\}$, which is covered by at most five intervals of total length at most $3\zeta$. By the onset $M\ge M^{\rm sl}_\eta$ and Lemma~\ref{lem:fou_der_the:2}, the first denominator has no zero on the remaining set, and
\[
\lvert\Phi_1^{(4)}-G\rvert\le\lambda/10^5 .
\]
If $\lvert S(0)\rvert\ge w_0$, that lemma and the choice of $r_\eta$ give
\[
\lvert\Phi_2^{(4)}\rvert\le10^{100}r^2w_0^{-5}\le\lambda/10^4 .
\]
The lower bound for $\lvert G\rvert$ proves $\lvert(\Phi_1+\Phi_2)^{(4)}\rvert\ge\lambda$ there.

Suppose $\lvert S(0)\rvert<w_0$ and put $w=\max(\lvert S(0)\rvert,r)$. If $r/w>d_0$, the first alternative of that lemma gives
\[
\lvert\Phi_2^{(4)}\rvert\ge d_0^5/(100r^3)\ge100B .
\]
Thus the same lower bound holds on the nonzero set of $S$; by that lemma $S$ is strictly monotone, so it has at most one zero, which we remove. If $r/w\le d_0$, put $\alpha=24S'(0)^4/(r^2S(0)^5)$. The denominator has no zero, and
\[
\lvert\Phi_2^{(4)}-\alpha\rvert\le10^{100}(w_0+d_0)\lvert\alpha\rvert=\frac{2\lambda}{10^6}\lvert\alpha\rvert .
\]
Outside the model band just estimated, it follows that $\lvert(\Phi_1+\Phi_2)^{(4)}\rvert\ge\lambda$. In all cases we have removed a union $\mathrm{Exc}$ of at most $66$ intervals with
\[
\lvert\mathrm{Exc}\rvert\le3\zeta+\eta/1000,
\]
and on its complement the actual reciprocal sum is smooth and has fourth derivative of magnitude at least $\lambda$.

We finally apply the sublevel estimate to the successful set. The cubic relation of Lemma~\ref{lem:cub_rel_suc:1} and $r\ge M^{-3/5}$ (Lemma~\ref{lem:pha_bou_loc}) give at every successful point
\[
\lvert P-\Phi_1-\Phi_2\rvert\le10^8\bigl(M^{-47/25}+r^{-2}M^{1-3c_{\rm hyp}}\bigr)\le2\cdot10^8M^{11/5-3c_{\rm hyp}}=:\varepsilon_M .
\]
The fourth derivative of the cubic $P$ vanishes. By Lemma~\ref{lem:sub_set_est} with $p=4$, applied on the at most $67$ complementary intervals, the successful set is covered by at most
\[
66+67\cdot5=401<500
\]
intervals of total length at most
\[
3\zeta+\eta/1000+67\cdot30\,(\varepsilon_M/\lambda)^{1/4}.
\]
The onset gives $\varepsilon_M\le2\cdot10^8M^{-1/2}\le2\cdot10^{-32}\lambda\eta^4$, so this total is less than $0.032\eta<\eta$. Isolated removed zeros can be covered by degenerate intervals, or by arbitrarily short intervals using this strict margin. This proves the assertion.

Explicit values. For $\eta=\eta_{\rm cov}=10^{-5}$: $\zeta=10^{-7}$, $B=10^{40}$, $\lambda=10^{-152}$ and $w_0=d_0=10^{-258}$. The four terms of the minimum defining $r_\eta$ are $10^{-260}$, $10^{-773}$, $10^{-1334/3}$ and $10^{-106}$, and the three terms of the maximum defining $M^{\rm sl}_\eta$ are $\exp(100)<10^{44}$, $10^{794}$ and $10^{424}$; so $r_\eta=10^{-773}$ and $M^{\rm sl}_\eta=10^{794}$.
\end{proof}

When the two phases are nearly opposite, the two reciprocal cosine terms nearly cancel, and the boundary cubic is compared with the normalized model:

\begin{definition}
 For a local test, choose $\kappa\in\{-1,1\}$ and $j\in\mathbb{Z}$, and put
\[
\widetilde Q_2=\kappa Q_2-\bigl(j+\tfrac12\bigr),\qquad\widetilde\omega_2=\widetilde Q_2(0),\qquad\widetilde b_2=\widetilde Q_2'(0),\qquad\widetilde d_2=\widetilde Q_2''(0).
\]
This choice is an \textbf{opposite phase alignment}. Its \textbf{alignment polynomial} and \textbf{size} are
\[
p(h)=(\widetilde\omega_2-\omega_1)+(\widetilde b_2-b_1)h+\tfrac12(\widetilde d_2-d_1)h^2,\qquad\ell=\max(\lVert p\rVert_{\rm coeff},K^{-1}).
\]
Its normalized parameters are $\bar p=p/\ell$ and $\beta_*=(K\ell)^{-1}$; when $b_1\in\{-1,1\}$, its model is $\mathcal R_{\omega_1,b_1,\bar p,\beta_*}$ of Definition~\ref{def:nor_opp_cos:1}.
\label{def:opp_pha_ali:1}
\end{definition}

\begin{lemma}
 Let $\zeta\in(0,1/100]$ and
\[
M\ge\max\bigl(\exp(100),(10^6/\zeta)^2\bigr).
\]
Consider an opposite phase alignment of a local test with $\ell\le\zeta/1000$, and let $a(h)=\omega_1+b_1h$. For every successful point $h\in[-1,1]$ with $\lvert\cos(2\pi a(h))\rvert\ge\zeta$,
\[
\Bigl\lvert\frac{P(h)}\ell-\mathcal R_{\omega_1,b_1,\bar p,\beta_*}(h)\Bigr\rvert\le10^{15}\zeta^{-10}\bigl(\ell+M^{-22/25}\bigr).
\]
\label{lem:the_bou_cub}
\end{lemma}

\begin{proof}
Put $c=c_{\rm hyp}$ (Convention~\ref{conv:not_and_con}). We keep the actual phase before replacing it by the affine phase. Put $y=Q_1'$, and let $\widetilde Q_2=\kappa Q_2-(j+\frac12)$ be the aligned second phase. By Lemma~\ref{lem:pha_bou_loc}, $\lvert b_1\rvert=1$, and $\max(\lVert\bar p\rVert_{\rm coeff},\beta_*)=\max(\lVert p\rVert_{\rm coeff},K^{-1})/\ell=1$, so the model is defined. Since $\mathcal E_2$ of Lemma~\ref{lem:pha_bou_loc} vanishes to second order at $0$, $\widetilde\omega_2=\kappa\omega_2-j-\frac12$, $\widetilde b_2=\kappa b_2$ and $\widetilde d_2=\kappa d_2$, and the second-order Taylor forms of $Q_1$ and $Q_2$ in that lemma give
\[
\widetilde Q_2=Q_1+p+\mathcal E,\qquad\mathcal E=\kappa\mathcal E_2-\mathcal E_1,\qquad\lVert\mathcal E\rVert_{C^1([-1,1])}\le2000M^{1-3c},
\]
\[
\max(\lvert Q_1-a\rvert,\lvert y-b_1\rvert)\le500M^{1-2c}.
\]
Since $\ell\ge K^{-1}>1/M$ and $2000M^{1-3c}\le1/M$ for $M\ge\exp(100)$, the error $\mathcal E$ and its derivative are at most $\ell$. Also $\lvert p\rvert,\lvert p'\rvert\le3\ell$ on $[-1,1]$, so the second pair $(\widetilde Q_2,\widetilde Q_2')$ differs from the first pair $(Q_1,y)$ by at most $4\ell$ in each coordinate. The alignment gives $r=\lvert b_2\rvert=\lvert\widetilde b_2\rvert\ge1-\ell>1/2$.

At the point under consideration, and throughout the line segments between the pairs $(a,b_1)$, $(Q_1,y)$ and $(Q_1+p+\mathcal E,y+p'+\mathcal E')$, the phase cosine has magnitude at least $\zeta/4$ and the derivative coordinate has magnitude at least $1/2$. Indeed $500M^{1-2c}\le\zeta/2000$ and $4\ell\le\zeta/250$, so every pair on these segments lies within $9\zeta/2000<\zeta/100$ of $(a,b_1)$ in each coordinate; since $x\mapsto\cos(2\pi x)$ is $2\pi$-Lipschitz, the phase cosine differs from $\cos(2\pi a(h))$ by less than $2\pi\zeta/100<\zeta/10$, and $\lvert b_1\rvert=1$. Set
\[
f(x,y)=y^{-2}\sec(2\pi x),\qquad B(x,y)=y^{-2}\sec^2(2\pi x).
\]
Direct differentiation on this domain shows that the sum of the absolute values of the first partial derivatives of $f$ is at most $10^3\zeta^{-2}$, that of the second partial derivatives of $f$ at most $10^5\zeta^{-3}$, and that of the first partial derivatives of $B$ at most $10^4\zeta^{-3}$. For example $f_{xx}=(2\pi)^2y^{-2}\bigl(\sec(2\pi x)\tan^2(2\pi x)+\sec^3(2\pi x)\bigr)$; the other derivatives have smaller powers of $\zeta^{-1}$. These estimates justify all Taylor bounds below with the maximum norm on pairs.

Write $\sigma_i=(2-E_i)K^{-1}$. The bounds of Lemma~\ref{lem:pha_bou_loc} imply
\[
\lvert\sigma_i\rvert\le3\ell,\qquad\lvert\sigma_i-2K^{-1}\rvert\le4K^{-1}M^{-c}.
\]
Since $\cos(2\pi Q_2)=-\cos(2\pi\widetilde Q_2)$ and $(Q_2')^2=(\widetilde Q_2')^2$, the reciprocal terms of Definition~\ref{def:res_pha_and} satisfy
\[
\Phi_1+\Phi_2=f(Q_1,y)-f(Q_1+p+\mathcal E,y+p'+\mathcal E')-\sigma_1B(Q_1,y)-\sigma_2B(Q_1+p+\mathcal E,y+p'+\mathcal E')+\mathcal E_{\rm rec},
\]
with
\[
\lvert\mathcal E_{\rm rec}\rvert\le10^4\ell^2\zeta^{-3}.
\]
This follows from the exact remainder $1/(X+\sigma)-1/X+\sigma/X^2=\sigma^2/[X^2(X+\sigma)]$ and its version with $X-\sigma$; the denominators have magnitudes at least $\zeta/4$ and $\zeta/8$.

Taylor's formula for the difference of the two $f$ terms, followed by replacement of its derivatives at $(Q_1,y)$ by their values at $(a,b_1)$, gives
\[
f(Q_1,y)-f(Q_1+p+\mathcal E,y+p'+\mathcal E')=-f_x(a,b_1)p-f_y(a,b_1)p'+\mathcal E_f,
\]
where
\[
\lvert\mathcal E_f\rvert\le2\cdot10^6\zeta^{-2}M^{1-3c}+8\cdot10^5\zeta^{-3}\ell^2+1.5\cdot10^8\zeta^{-3}\ell M^{1-2c}.
\]
The three terms bound respectively the phase remainder, the quadratic Taylor remainder, and the derivative replacement.

Replacing $\sigma_i$ by $2K^{-1}$ in the shift terms costs at most $512\zeta^{-2}\ell M^{-c}$. Replacing their two $B$ factors by $B(a,b_1)$ costs at most
\[
2\cdot10^7\zeta^{-3}\ell M^{1-2c}+8\cdot10^4\zeta^{-3}\ell^2,
\]
by the gradient bound for $B$, since the sum of the two pair displacements from $(a,b_1)$ is at most $1000M^{1-2c}+4\ell$.

The resulting main term is
\[
-f_x(a,b_1)p-f_y(a,b_1)p'-4K^{-1}B(a,b_1),
\]
which, divided by $\ell$, is exactly the model $\mathcal R_{\omega_1,b_1,\bar p,\beta_*}$ of Definition~\ref{def:nor_opp_cos:1}. Finally, by Lemma~\ref{lem:cub_rel_suc:1}, $P=\Phi_1+\Phi_2+\varepsilon$ at the successful point, and $\varepsilon$ divided by $\ell\ge1/M$ is at most $10^8M^{-22/25}+4\cdot10^8M^{2-3c}$, since $r^{-2}<4$. Divide the preceding remainder bounds by $\ell$, use $\ell^{-1}<M$, $\zeta\le1$ and $M\ge1$, and sum them; every power of $M$ that occurs is at most $M^{-22/25}$, and the total is at most $7\cdot10^8\zeta^{-3}(\ell+M^{-22/25})\le10^{15}\zeta^{-10}(\ell+M^{-22/25})$, as claimed.
\end{proof}

\begin{thm}
 Let $\eta\in(0,1]$ and $\zeta_\eta=\eta/100$. There are $s_\eta>0$ and an integer $d_\eta\ge1$ with the following property. For every normalized opposite-cosine model $\mathcal R=\mathcal R_{\omega,b,p,\beta}$, with affine phase $a$, and every real cubic $P$, the set
\[
\{h\in[-1,1]:\lvert\cos(2\pi a(h))\rvert\le\zeta_\eta\}\cup\{h\in[-1,1]:\lvert\cos(2\pi a(h))\rvert>\zeta_\eta,\ \lvert P(h)-\mathcal R(h)\rvert\le s_\eta\}
\]
is contained in a union of at most $d_\eta+6$ intervals of total length at most $\eta$. For $\eta=\eta_{\rm cov}$ one may take $d_\eta=2\cdot10^8$ and $s_\eta=\exp(-3.3\cdot10^9)$.
\label{thm:sma_cov_for}
\end{thm}

\begin{proof}
We fix the constants. Put $B_\eta=100\zeta_\eta^{-2}$, $H_\eta=10^9\eta^{-3}B_\eta$, $R=192/\eta$ and $B^{\mathbb C}_\eta=40H_\eta(1+R)^3\exp(4\pi R)$, and
\[
d_\eta=\Bigl\lceil\frac{\log(8000B^{\mathbb C}_\eta)}{\log4}\Bigr\rceil,\qquad s_\eta=\frac1{8000}\Bigl(\frac\eta{48}\Bigr)^{d_\eta};
\]
these are the $D$ and $s$ of Lemma~\ref{lem:sub_set_ent} for $u=1/8$ and $B^{\mathbb C}=B^{\mathbb C}_\eta$, and $d_\eta\ge3$. The conclusion remains true for a larger count and a smaller tolerance, which, with the explicit values below, gives the last sentence of the statement.

We first treat cubics with large coefficients. On the pole-safe set $\lvert\cos(2\pi a)\rvert>\zeta_\eta$ we have $\lvert p\rvert\le3$, $\lvert p'\rvert\le3$ and $\beta\le1$, so
\[
\lvert\mathcal R\rvert\le6/\zeta_\eta+(6\pi+4)/\zeta_\eta^2\le30/\zeta_\eta^2 .
\]
Thus $\lvert P-\mathcal R\rvert\le s_\eta$ implies $\lvert P\rvert\le B_\eta$ there. If $\lVert P\rVert_{\rm coeff}>H_\eta$, then by Lemma~\ref{lem:rem_ine_and} the set $\{h\in[-1,1]:\lvert P(h)\rvert\le B_\eta\}$ is a union of at most four intervals of total length less than $24(24B_\eta/H_\eta)^{1/3}=24(24\cdot10^{-9})^{1/3}\eta<\eta/4$.

Suppose now that $\lVert P\rVert_{\rm coeff}\le H_\eta$. We clear the cosine denominator: put
\[
F(h)=b^2P(h)\cos^2(2\pi a(h))-\frac2bp'(h)\cos(2\pi a(h))+2\pi p(h)\sin(2\pi a(h))+4\beta .
\]
This is entire, real on $\mathbb{R}$, and $F=b^2\cos^2(2\pi a)\,(P-\mathcal R)$ away from the poles. We prove a lower bound independent of $P$. Choose four consecutive zeros $x_0,\dots,x_3$ of $\cos(2\pi a)$ in $[-1,1]$, at spacing $1/2$, and let $u_F=\sup_{[-1,1]}\lvert F\rvert$. At these zeros
\[
F(x_i)=2\pi p(x_i)\sin(2\pi a(x_i))+4\beta .
\]
The sine signs alternate. The third finite difference of the quadratic $p$ is zero, so these four equations give $32\beta\le8u_F$, and hence $\beta\le u_F/4$. They then imply $\lvert p(x_i)\rvert\le u_F/\pi$. Interpolating $p$ at any three of these zeros bounds every coefficient by $24u_F/\pi<8u_F$: the Lagrange denominators have magnitude at least $1/4$ and the numerator coefficients at most two. The normalization $\max(\lVert p\rVert_{\rm coeff},\beta)=1$ therefore yields $u_F\ge1/8$.

For complex $h$ with $\lvert h\rvert\le R$, the coefficient bounds and $\lvert\cos(2\pi a(h))\rvert,\lvert\sin(2\pi a(h))\rvert\le\exp(2\pi R)$ give
\[
\lvert F(h)\rvert\le40H_\eta(1+R)^3\exp(4\pi R)=B^{\mathbb C}_\eta .
\]
By Lemma~\ref{lem:sub_set_ent} with $u=1/8$, the set $\{h\in[-1,1]:\lvert F(h)\rvert\le s_\eta\}$ is contained in at most $d_\eta+1$ intervals of total length less than $\eta/2$. At every point of the second set of the statement, $\lvert F\rvert\le\lvert P-\mathcal R\rvert\le s_\eta$, since $b^2=1$; so these intervals cover that set.

Finally, the pole strip $\lvert\cos(2\pi a)\rvert\le\zeta_\eta$ has at most five interval components in $[-1,1]$. Because $\lvert b\rvert=1$, its measure is $4\arcsin(\zeta_\eta)/\pi\le3\zeta_\eta$. Adding those intervals gives at most $d_\eta+6$ intervals, with total length less than $\eta/2+3\eta/100<\eta$.

Explicit values. For $\eta=\eta_{\rm cov}=10^{-5}$: $\zeta_\eta=10^{-7}$, $B_\eta=10^{16}$, $H_\eta=10^{40}$, $R=1.92\cdot10^7$ and $\log B^{\mathbb C}_\eta<2.5\cdot10^8$; hence, with $\log4>1.3$, $d_\eta<2\cdot10^8$, and $\log(1/s_\eta)\le\log8000+2\cdot10^8\log(4.8\cdot10^6)<3.3\cdot10^9$.
\end{proof}

In the generic case the two phases are neither slow nor nearly opposite. Clearing the cosine denominators turns the relation into a sublevel condition for an entire function whose supremum on $[-1,1]$ is bounded below, and Lemma~\ref{lem:sub_set_ent} gives the cover.

\begin{thm}
 Let $\eta,r_0\in(0,1]$, put $\zeta_0=\eta r_0/1000$, and let $0<\nu\le10^{-29}\eta^5r_0^4$. Put
\[
R=\frac{192}\eta,\qquad u=\frac{r_0^2\nu}{100},\qquad B^{\mathbb C}=10^{18}\eta^{-5}r_0^{-4}(1+R)^3\exp(4\pi R),
\]
\[
D=\Bigl\lceil\frac{\log(1000B^{\mathbb C}/u)}{\log4}\Bigr\rceil,\qquad s=\frac u{1000}\Bigl(\frac\eta{48}\Bigr)^D.
\]
Let $r\in[r_0,1]$, $\omega_1,\omega_2\in\mathbb{R}$, $b_1\in\{-1,1\}$ and $b_2\in\{-r,r\}$; put $a_i(h)=\omega_i+b_ih$; and let $\kappa\in\{-1,1\}$ satisfy $\kappa b_2=b_1r$. Suppose there is no integer $j$ with both $\lvert\kappa\omega_2-\omega_1-(j+\frac12)\rvert\le2\nu$ and $1-r\le\nu$. Then, for every real cubic $P$, the set of $h\in[-1,1]$ with $\lvert\cos(2\pi a_1(h))\rvert\ge\zeta_0$, $\lvert\cos(2\pi a_2(h))\rvert\ge\zeta_0$ and
\[
\Bigl\lvert P(h)-\frac1{\cos(2\pi a_1(h))}-\frac1{r^2\cos(2\pi a_2(h))}\Bigr\rvert\le s
\]
lies in at most $D+1$ intervals of total length less than $\eta/2$.
\label{thm:sma_cov_for:1}
\end{thm}

\begin{proof}
Put $B_0=100r_0^{-2}\zeta_0^{-2}=10^8\eta^{-2}r_0^{-4}$ and $H_0=10^9\eta^{-3}B_0=10^{17}\eta^{-5}r_0^{-4}$, so that $\nu\le(10^{12}H_0)^{-1}$, $B^{\mathbb C}=10H_0(1+R)^3\exp(4\pi R)$ and $D\ge3$. Write $y_i(h)=\cos(2\pi a_i(h))$.

We first treat cubics with large coefficients. At every point of the set in the statement, $\lvert y_1\rvert\ge\zeta_0$, $r^2\lvert y_2\rvert\ge r_0^2\zeta_0$ and $s\le1$, so $\lvert P\rvert\le1+\zeta_0^{-1}(1+r_0^{-2})\le B_0$. If $\lVert P\rVert_{\rm coeff}>H_0$, then by Lemma~\ref{lem:rem_ine_and} these points lie in at most four intervals of total length less than $24(24B_0/H_0)^{1/3}<\eta/4$; since $D\ge3$, this proves the assertion in that case.

Assume henceforth that the coefficients of $P$ have magnitude at most $H_0$. Clear the denominators: the entire function
\[
F(h)=r^2P(h)y_1(h)y_2(h)-r^2y_2(h)-y_1(h)
\]
is real on $\mathbb{R}$ and equals $r^2y_1y_2\bigl(P-1/y_1-1/(r^2y_2)\bigr)$ away from the zeros of $y_1y_2$. We prove that $\sup_{[-1,1]}\lvert F\rvert\ge u$. Choose consecutive zeros $x_0$ and $x_1=x_0+1/2$ of $y_1$, with $x_0\in[-1/2,0]$. If $\lvert y_2(x_i)\rvert>\nu$ at either zero, then $\lvert F(x_i)\rvert=r^2\lvert y_2(x_i)\rvert\ge r_0^2\nu>u$. Otherwise, since $2\pi a_2(x_1)=2\pi a_2(x_0)+\pi b_2$,
\[
\lvert\sin(\pi b_2)\rvert\le\lvert y_2(x_0)\rvert+\lvert y_2(x_1)\rvert\le2\nu .
\]
Since $\sin(\pi r)\ge2\min(r,1-r)$ and $r\ge r_0>\nu$, this implies $1-r\le\nu$.

At $x_0$, the phases $a_1(x_0)$ and $\kappa a_2(x_0)$ lie within $\arcsin(\nu)/(2\pi)\le\nu$ of the set $\frac14+\frac12\mathbb{Z}$ of zeros of $\cos(2\pi\,\cdot)$, the first exactly. Since $\kappa a_2-a_1=\kappa\omega_2-\omega_1+b_1(r-1)h$ and $\lvert x_0\rvert\le\frac12$, it follows that for some integer $i$
\[
\lvert\kappa\omega_2-\omega_1-i/2\rvert\le2\nu .
\]
If $i$ were odd, then $i/2=j+\frac12$ with $j=(i-1)/2$, and this and $1-r\le\nu$ would contradict the hypothesis. Thus $i$ is even. Write $\kappa a_2=a_1+i/2+\varepsilon$, with $\lvert\varepsilon\rvert\le3\nu$ on $[-1,1]$. Then $y_2=\cos(2\pi(a_1+\varepsilon))$, and at $x_0$, where $y_1=0$,
\[
y_2'(x_0)=r\cos(2\pi\varepsilon(x_0))\,y_1'(x_0).
\]
Since $\lvert P(x_0)\rvert\le4H_0$ and $\lvert y_2(x_0)\rvert\le\nu$,
\[
\lvert F'(x_0)\rvert=\lvert y_1'(x_0)\rvert\,\bigl\lvert r^2P(x_0)y_2(x_0)-r^3\cos(2\pi\varepsilon(x_0))-1\bigr\rvert\ge5,
\]
because $\lvert y_1'(x_0)\rvert=2\pi$, $4H_0\nu\le4\cdot10^{-12}$ and $\cos(2\pi\varepsilon(x_0))>0$. The product rule bounds $\lvert F''\rvert\le10^8H_0$ on $[-1,1]$. Put $h_1=(10^9H_0)^{-1}$. Taylor's formula, with $x_0+h_1\in[-1,1]$, yields
\[
\lvert F(x_0+h_1)-F(x_0)\rvert\ge5h_1-\tfrac12\cdot10^8H_0h_1^2\ge4h_1 .
\]
It follows that $\sup_{[-1,1]}\lvert F\rvert\ge2h_1>u$.

On the complex disc $\lvert h\rvert\le R$, the coefficient bound and $\lvert y_i(h)\rvert\le\exp(2\pi R)$ give
\[
\lvert F(h)\rvert\le10H_0(1+R)^3\exp(4\pi R)=B^{\mathbb C},
\]
and $u\le B^{\mathbb C}$. By Lemma~\ref{lem:sub_set_ent}, the set $\{h\in[-1,1]:\lvert F(h)\rvert\le s\}$ lies in at most $D+1$ intervals of total length less than $\eta/2$. At every point of the set in the statement, $\lvert F\rvert=r^2\lvert y_1y_2\rvert\,\bigl\lvert P-1/y_1-1/(r^2y_2)\bigr\rvert\le s$, which proves the assertion.
\end{proof}

The three cases together cover every local test:

\begin{thm}
 For every $\eta\in(0,1]$ there are $M^{\rm cov}_\eta<\infty$ and an integer $N^{\rm cov}_\eta$ such that, for $M\ge M^{\rm cov}_\eta$, the successful set of every local test is contained in a union of at most $N^{\rm cov}_\eta$ intervals of total length at most $\eta$. For $\eta=\eta_{\rm cov}$ one may take $M^{\rm cov}_\eta=\exp(1.1\cdot10^{11})$ and $N^{\rm cov}_\eta=4\cdot10^9$.
\label{thm:sma_cov_the}
\end{thm}

\begin{proof}
 At every successful point the boundary cubic $P$ nearly equals the sum $\Phi_1+\Phi_2$ of the two reciprocal terms (Lemma~\ref{lem:cub_rel_suc:1}), and the successful set is covered by showing that a cubic can stay this close to the reciprocal sum only on a small set. The argument splits into three cases according to the two rescaled phases: small slope ratios, covered by Theorem~\ref{thm:sma_cov_suc}; a near-opposite alignment of the two phases, which reduces through Lemma~\ref{lem:the_bou_cub} to the model covers of Theorem~\ref{thm:sma_cov_for}; and the generic case, which is Theorem~\ref{thm:sma_cov_for:1}.

We fix the constants. Put $c=c_{\rm hyp}$. Let $r_\eta$ and $M^{\rm sl}_\eta$ be as in Theorem~\ref{thm:sma_cov_suc}, and $\zeta_\eta=\eta/100$, $s_\eta$ and $d_\eta$ as in Theorem~\ref{thm:sma_cov_for}. Put $r_*=r_\eta$,
\[
\lambda=\frac{s_\eta\zeta_\eta^{10}}{4\cdot10^{15}},\qquad\zeta_*=\frac{\eta r_*}{1000},\qquad\nu=\min\Bigl\{\frac\lambda{100},10^{-29}\eta^5r_*^4\Bigr\},\qquad C'_*=10^{20}r_*^{-10}\zeta_*^{-20},
\]
and let $D_*$ and $s_*$ be the $D$ and $s$ of Theorem~\ref{thm:sma_cov_for:1} for $(\eta,r_0,\nu)=(\eta,r_*,\nu)$, whose $\zeta_0$ is $\zeta_*$. Put
\[
M^{\rm cov}_\eta=\max\Bigl\{M^{\rm sl}_\eta,\Bigl(\frac{10^6}{\zeta_\eta}\Bigr)^2,\Bigl(\frac{10^5}\lambda\Bigr)^2,\Bigl(\frac{10^6}{\zeta_*}\Bigr)^2,\Bigl(\frac{C'_*}{s_*}\Bigr)^2\Bigr\},\qquad N^{\rm cov}_\eta=\max\{500,d_\eta+6,D_*+11\},
\]
and let $M\ge M^{\rm cov}_\eta$; in particular $M\ge\exp(100)$. The conclusion remains true for a larger threshold and a larger count.

If $r\le r_*$, Theorem~\ref{thm:sma_cov_suc} proves the assertion.

Suppose next that the test admits an opposite phase alignment (Definition~\ref{def:opp_pha_ali:1}) of size $\ell\le\lambda$. Since $\lambda\le\zeta_\eta/1000$ and $M\ge\max(\exp(100),(10^6/\zeta_\eta)^2)$, Lemma~\ref{lem:the_bou_cub} applies with $\zeta=\zeta_\eta$, and its error is at most $s_\eta$:
\[
10^{15}\zeta_\eta^{-10}(\ell+M^{-22/25})\le2\cdot10^{15}\zeta_\eta^{-10}\lambda=s_\eta/2,
\]
since $M^{-22/25}\le M^{-1/2}\le\lambda$. Apply Theorem~\ref{thm:sma_cov_for} to the cubic $P/\ell$ and the model $\mathcal R_{\omega_1,b_1,\bar p,\beta_*}$, whose affine phase is $a(h)=\omega_1+b_1h$. It covers every successful point outside its fast pole strip, and that strip as well, by at most $d_\eta+6$ intervals of total length at most $\eta$.

It remains to consider $r>r_*$ and tests with no opposite phase alignment of size at most $\lambda$. Put $a_i(h)=\omega_i+b_ih$. Remove the two affine pole strips $\{\lvert\cos(2\pi a_i)\rvert<\zeta_*\}$. Each is covered by at most five intervals; using monotonicity near its cosine zeros, the sum of their lengths is at most $10\zeta_*/r_*=\eta/100$. By Lemma~\ref{lem:pha_bou_loc}, $\lvert Q_i-a_i\rvert,\lvert Q_i'-b_i\rvert\le500M^{1-2c}$ on $[-1,1]$; the cosine changes by at most $2\pi$ times the phase error, the shift $(2-E_i)/K$ is at most $6/M$, and $\lvert Q_i'\rvert\ge\lvert b_i\rvert/2$. The onset $(10^6/\zeta_*)^2$ makes all these perturbations smaller than $\zeta_*/100$, so outside the strips the actual denominators are nonzero, and subtracting the two reciprocals gives
\[
\bigl\lvert\Phi_i-b_i^{-2}\cos(2\pi a_i)^{-1}\bigr\rvert\le10^{10}r_*^{-4}\zeta_*^{-4}M^{1-2c}.
\]
At successful pole-safe points the cubic relation of Lemma~\ref{lem:cub_rel_suc:1} and $r\ge M^{-3/5}$ therefore give
\[
\Bigl\lvert P-\frac1{\cos(2\pi a_1)}-\frac1{r^2\cos(2\pi a_2)}\Bigr\rvert\le10^{12}r_*^{-4}\zeta_*^{-4}M^{11/5-3c}\le C'_*M^{-1/2}\le s_* .
\]

We check the hypothesis of Theorem~\ref{thm:sma_cov_for:1}. Here $\lvert b_1\rvert=1$ and $r=\lvert b_2\rvert\in(r_*,1]$ (Lemma~\ref{lem:pha_bou_loc}); let $\kappa\in\{-1,1\}$ satisfy $\kappa b_2=b_1r$. Suppose some integer $j$ had $\lvert\kappa\omega_2-\omega_1-(j+\frac12)\rvert\le2\nu$ and $1-r\le\nu$. The opposite phase alignment with this $\kappa$ and this $j$ would have alignment polynomial with constant coefficient $\kappa\omega_2-\omega_1-(j+\frac12)$, of magnitude at most $2\nu\le\lambda$, linear coefficient $\kappa b_2-b_1=b_1(r-1)$, of magnitude at most $\nu\le\lambda$, and quadratic coefficient $\frac12(\kappa d_2-d_1)$, of magnitude at most $100M^{1-2c}\le\lambda$; and $K^{-1}\le\lambda$, by the onset $(10^5/\lambda)^2$. So its size would be at most $\lambda$, which this case excludes. Hence Theorem~\ref{thm:sma_cov_for:1} applies to the boundary cubic $P$, with $r_0=r_*$, and covers the successful pole-safe points by at most $D_*+1$ intervals of total length less than $\eta/2$.

We count. Adding the at most ten pole-strip intervals gives at most $D_*+11$ intervals, of total length less than $\eta/100+\eta/2<\eta$. With the first two cases, the count is at most $N^{\rm cov}_\eta$, and the theorem is proved.

Explicit values. For $\eta=\eta_{\rm cov}=10^{-5}$, take the values $r_\eta=10^{-773}$ and $M^{\rm sl}_\eta=10^{794}$ of Theorem~\ref{thm:sma_cov_suc} and $d_\eta=2\cdot10^8$ and $s_\eta=\exp(-3.3\cdot10^9)$ of Theorem~\ref{thm:sma_cov_for}. Then $\log(1/\lambda)<3.4\cdot10^9$, $\zeta_*=10^{-781}$, $\nu=\lambda/100$ with $\log(1/\nu)<3.5\cdot10^9$, the $u=r_*^2\nu/100$ of Theorem~\ref{thm:sma_cov_for:1} has $\log(1/u)<3.6\cdot10^9$, its $B^{\mathbb C}$ has $\log B^{\mathbb C}<2.5\cdot10^8$, and $D_*<3\cdot10^9$, $\log(1/s_*)<5.2\cdot10^{10}$ and $\log C'_*<54000$. Hence $\log M^{\rm cov}_\eta<2(5.2\cdot10^{10}+54000)<1.1\cdot10^{11}$, the other terms of the maximum being smaller, and $N^{\rm cov}_\eta<3\cdot10^9+11<4\cdot10^9$.
\end{proof}

\subsection{Nesting the covers} \label{subsec:nes_the_cov}

The last step assembles the three subsections into the scarcity of cancellations. Fix the central time, the position there and a pattern of matched indices; a value of $\theta$ is admissible if the hyperbolic times of the pattern occur at it. Near an admissible value, the admissible values, rescaled by the expansion at the current pair, are successful points of a local test. The expansion grows by at least the factor $M^{c_{\rm hyp}}$ from one pair to the next, so the neighbourhoods used at the next pair are negligible at this scale and stay close to the successful set. A Vitali selection and the covers of Theorem~\ref{thm:sma_cov_the} then shrink a neighbourhood of the admissible set by the factor $10^{-4}$ at each pair, which beats the at most $\exp(2n/5)$ patterns, and Fubini--Tonelli in the state coordinates bounds the measure of the cancellation set. The theorem is used in the proofs of the doubling inequality Corollary~\ref{cor:dou_ine_for} and of the approximation by slack profiles Lemma~\ref{lem:sum_exc_mas:1}.

\begin{thm}
 Let $M\ge M_{\rm an}$. For every $i\in\mathbb{Z}$ and every integer $n\ge1$,
\[
\mathbb{P}\bigl(\mathrm{Can}_{\delta_{\rm can}}(i,n)\bigr)\le C_{\rm can}\exp(-\gamma_{\rm can}n).
\]
\label{thm:exp_sca_can}
\end{thm}

\begin{proof}
We count index patterns against a contraction of the set of admissible values of the free coordinate $\theta=q_{\iota-1}$ (Definition~\ref{def:out_pat_and}) from one matched pair to the next, as in the parameter exclusion of \cite{locbib:BenedicksCarleson1991}.

Put $c=c_{\rm hyp}$ and $\delta=\delta_{\rm can}$ (Convention~\ref{conv:not_and_con}). We first read off the local cover. By Theorem~\ref{thm:sma_cov_the} with $\eta=\eta_{\rm cov}$, since $M\ge M_{\rm an}\ge\exp(1.1\cdot10^{11})$ (Convention~\ref{conv:not_and_con}), the successful set of every local test (Definition~\ref{def:loc_tes_and}) lies in at most $4\cdot10^9$ intervals of total length at most $\eta_{\rm cov}$. So its $\rho$-neighbourhood, with $\rho=M^{-c}/4$, has measure at most $\eta_{\rm cov}+2\cdot10^9M^{-c}\le2\eta_{\rm cov}$, since $M\ge M_{\rm an}$; and $10\eta_{\rm cov}=10^{-4}$.

We reduce to $i=0$ and to large $n$. By Lemma~\ref{lem:inv_are_und}(iii), $\mathbf 1_{\mathrm{Can}_\delta(i,n)}=\mathbf 1_{\mathrm{Can}_\delta(0,n)}\circ f_k^i$, and $f_k^i$ preserves $m$ by Lemma~\ref{lem:inv_are_und}(i); so it suffices to take $i=0$. For $n<1/\delta_{\rm can}$ the claimed bound exceeds $C_{\rm can}/\exp(1)>1$, since $\gamma_{\rm can}=\delta_{\rm can}$; this is why the rate $\gamma_{\rm can}$ is chosen equal to $\delta_{\rm can}$, below the rate $0.06$ that the last step of this proof gives. Let $n\ge1/\delta_{\rm can}$ and $z\in\mathrm{Can}_\delta(0,n)$. By Lemma~\ref{lem:can_giv_sol} there are $\iota\in\{0,-1\}$ and a normalized solution that together with the orbit of $z$ form an endpoint-small configuration of length $n$ and defect $\delta$ (Definition~\ref{def:end_sma_con}). By Lemma~\ref{lem:man_mat_pai}, this solution has $p=\lceil n/20\rceil$ pairs $(j^+_l,j^-_l)$ of hyperbolic times (the first $p$ of those of the lemma, with entries in $\{2,\dots,n\}$, since $j^\pm_l\le N_\pm-2\le n-1$ by Definition~\ref{def:out_pat_and}), strictly increasing on both paths and matched at $\theta=q_{\iota-1}$ (Definition~\ref{def:mat_pai_hyp:2}); and the strictly increasing patterns of $p$ index pairs number at most $\exp(2n/5)$.

We next define the admissible sets. Fix $\iota$, a pattern and $q_\iota$. Let $\mathrm{Adm}$ be the set of $\theta$ for which one normalized solution with $\lVert(u_0,u_{-1})\rVert\le\sqrt2$ has hyperbolic times (Definition~\ref{def:pai_gro_log}) at all indices of the pattern, with every pair satisfying \eqref{eq:cs_matching}. Parametrize the central pair (Definition~\ref{def:out_pat_and}) by its entry other than $u_\iota=1$; this entry lies in $[-1,1]$. All quantities in these conditions are continuous in this entry and in $(q_\iota,\theta)$; the pairs that enter logarithms are nonzero (Lemma~\ref{lem:wro_and_cen:2}(iv)); and the conditions are closed. Hence the set of $(q_\iota,\theta)$ with $\theta\in\mathrm{Adm}$ is the projection of a compact set. It is therefore compact, and so is each $\mathrm{Adm}$.

For $\theta\in\mathrm{Adm}$ put
\[
s_l(\theta)=\max\bigl(\lvert t^+_{j^+_l}(\theta)\rvert,\lvert t^-_{j^-_l}(\theta)\rvert\bigr)^{-1},\qquad V_l=\bigcup_{\theta\in\mathrm{Adm}}B\bigl(\theta,s_l(\theta)/10\bigr),
\]
where $B(\theta,s)$ is the open arc of radius $s$ about $\theta$. The bound $\lvert t_{l'}\rvert\le L_jM^{-c(j-l')}$ of \eqref{eq:cs_dominated_variation} (Lemma~\ref{lem:bou_val_for}(a)), applied at the later index on each path, gives $s_{l+1}\le M^{-c}s_l$; and $s_1\le M^{-c}$, since the same bound at $l'=1$ gives $1=\lvert t_1\rvert\le L_jM^{-c(j-1)}$ at the indices $j=j^\pm_1\ge2$.

We show that $\lvert V_{l+1}\rvert\le10^{-4}\lvert V_l\rvert$. By the Vitali covering lemma (see \cite[Theorem~2.1]{locbib:Mattila1995} and its proof, applied to the closed arcs), there are disjoint balls $B_\alpha=B(\theta_\alpha,s_\alpha/10)$ of the family defining $V_l$, at most countably many, with $s_\alpha=s_l(\theta_\alpha)$, such that every ball $B(\theta,s_l(\theta)/10)$ of the family meets some $B_\alpha$ with $s_l(\theta)\le2s_\alpha$; then $\lvert\theta-\theta_\alpha\rvert\le3s_\alpha/10$. Since $\theta_\alpha\in\mathrm{Adm}$, $(\iota,q_\iota,\theta_\alpha,j^+_l,j^-_l)$ is a local test with $L=1/s_\alpha$ (Definition~\ref{def:loc_tes_and}). Every $\theta\in\mathrm{Adm}$ attached to $B_\alpha$ gives the successful point $h=(\theta-\theta_\alpha)/s_\alpha\in[-1,1]$ (on the lift; $s_\alpha\le M^{-c}$). The next ball of $\theta$ has rescaled radius $s_{l+1}(\theta)/(10s_\alpha)\le M^{-c}/5<\rho$. So the next balls of the points attached to $B_\alpha$ have a union of measure at most $2\eta_{\rm cov}s_\alpha=10\eta_{\rm cov}\lvert B_\alpha\rvert=10^{-4}\lvert B_\alpha\rvert$. Every point of $\mathrm{Adm}$ is attached to some $B_\alpha$, and the $B_\alpha$ are disjoint subsets of $V_l$; summing over $\alpha$ proves the claim.

We conclude. Since $\lvert V_1\rvert\le1$, $\lvert\mathrm{Adm}\rvert\le\lvert V_p\rvert\le(10^{-4})^{p-1}\le10^4\exp(-n\log(10^4)/20)$. For this $\iota$, the union over the patterns of the compact sets above contains the states $(q_\iota,q_{\iota-1})$ of every $z\in\mathrm{Can}_\delta(0,n)$. By Lemma~\ref{lem:inv_are_und}(ii) the state map carries $m$ to $m$. Fubini's theorem and a sum over the patterns and the two values of $\iota$ then give
\[
\mathbb{P}\bigl(\mathrm{Can}_\delta(0,n)\bigr)\le2\cdot10^4\exp\Bigl(-\Bigl(\frac{\log10^4}{20}-\frac25\Bigr)n\Bigr)\le C_{\rm can}\exp(-\gamma_{\rm can}n),
\]
since $2\cdot10^4=C_{\rm can}$ and $\log(10^4)/20-2/5>0.06>\gamma_{\rm can}$, by the values of Convention~\ref{conv:not_and_con}.
\end{proof}

\section{Notation and constants} \label{sec:not_and_con}

\begin{convention}
\textbf{Symbols.} The following symbols have one meaning throughout: $k,K,M,\phi,v$ (Convention~\ref{conv:sta_map_par}), with $v_i=v(q_i)$; $f_k$, $h_m$, $\lambda_+$ and the area $m$ (Definition~\ref{def:sta_map_ent}); $\mathbb{P}$, $\mathbb{E}$ (Convention~\ref{conv:mea_and_int}); the initial point $z$ of an orbit, $q_i$, $\Pi$, $A_i$, $T_{a,b}$ (Definition~\ref{def:orb_sta_and}); $g$, $e_n$, $e^\Psi_n$ (Definition~\ref{def:gro_sem_orb}); $h_d$ (Definition~\ref{def:sho_sem_the}); $(x|y)_w$ (Definition~\ref{def:gro_hyp_sem}); $\mathrm{Can}_\delta(i,n)$ (Definition~\ref{def:can_the_gro}); and the constants listed below. A symbol that any other definition fixes keeps its meaning in every definition and result that cites that definition, directly or through the definitions it cites; a variable over which a definition only quantifies, and a symbol introduced in a statement or a proof, are local to that definition, or to that statement and its proof. The coordinates $(x_i,y_i)$ of $f_k^i(z)$ (Definition~\ref{def:orb_sta_and}) serve only to define the positions and the states; elsewhere $x$ and $y$, with or without index, are real variables. A bare $T$ may denote a generic matrix of $\mathrm{SL}(2,\mathbb{R})$. Besides the entropy $h_m(f_k)$ and the shortfalls $h_d(s,t)$, the letter $h$ is used for the shortfalls $h(I)$, $h_t(x)$ of arrays (Definition~\ref{def:arr_mid_can}) and, in the rescaling of local tests, for the rescaled coordinate (Definition~\ref{def:res_pha_and}). A letter $t$ with an integer index is a solution of a three-term recurrence that vanishes at its first index and equals one at the next, such as the variation sequences $t^\pm_j$ (Definition~\ref{def:out_pat_and}) and the solutions $t^\pm_p$ (Definition~\ref{def:pre_gro_eli}), whose indices run backwards for $t^-_p$; without an index, $t$ is a real parameter. Euler's number is written $\exp(1)$, and $\mathrm{i}$ is the imaginary unit.

\textbf{Devices.} $\log$ is the natural logarithm, $\log_M=\log/\log M$, and $\log_2$ is written explicitly. $\lVert\cdot\rVert$ is the Euclidean norm of a vector and the operator norm of a matrix; other norms carry a subscript: $\lVert\cdot\rVert_\infty$ is the maximum norm on $\mathbb{R}^2$ and its operator norm, $\lVert\cdot\rVert_{C^1}$ the $C^1$ norm on an interval, $\lVert\cdot\rVert_{\rm coeff}$ the largest absolute coefficient of a polynomial, $\lVert\cdot\rVert_*$ a weighted norm (proof of Lemma~\ref{lem:sta_gra_und:1}), and $\lVert d-d''\rVert_\Theta$ the uniform distance of two arrays on an interval $\Theta$ (Definition~\ref{def:arr_mid_can}); $\sup_E\lvert F\rvert$ is the supremum of $\lvert F\rvert$ on a set $E$. $\lvert\cdot\rvert$ is absolute value, length or Lebesgue measure. A Gromov product is always written with its basepoint, $(x|y)_w$. Calligraphic letters denote families, index sets, operators (the translations $\mathcal{T}_x$ and the dilation $\mathcal{D}$ of arrays, Definition~\ref{def:arr_mid_can}) and auxiliary functions; $\mathscr{S}$ is the family of slack profiles (Definition~\ref{def:slac_prof}), and $\mathbf{1}$ an indicator. A tilde marks a modified version, a bar an average or a normalized version, a hat a transformed version, a star a base or selected element, and a superscript $\pm$ an orientation. Sets defined by conditions are roman words with their parameters as subscripts, as $\mathrm{Can}_\delta(i,n)$ and $\mathrm{Exc}$.

\textbf{Constants.} A constant written with a roman subscript is fixed once for the whole paper, with the value listed below; a threshold on $M$ is written $M_{\rm name}$, or $M^{\rm name}_\eta$ when it depends on a parameter $\eta$, and a constant that depends on a parameter carries it as a subscript, as $r_\eta$. Unless a dependence is stated, constants in estimates for large $M$ depend neither on the orbit, nor on the times, nor on $M$. In a proof, $C$ may denote an absolute constant (independent of the orbit, the times, $M$ and the parameters of the statement) that changes from line to line, if the proof says so; such a $C$ never appears in a statement. A proof that fixes numbers ends with a paragraph of explicit values.
\begin{description}
\item[$M_{\rm ent}=\exp(2^{50000})$, $C_{\rm ent}=2^{430}$] parameter threshold and entropy deficit of the main results (Proposition~\ref{prop:uni_bou_the:2}).
\item[$M_{\rm an}=\exp(10^{12})$] parameter threshold of the doubling inequality and of the cancellation and bridge estimates (Corollary~\ref{cor:dou_ine_for}, Theorem~\ref{thm:exp_sca_can}, Theorem~\ref{thm:bri_pro_est:2}).
\item[$C_{\rm bal}=2^{420}$, $\zeta_{\rm bal}=2^{-47000}$, $n_{\rm bal}=2^{49000}$] constant, shortfall threshold and initial-scale threshold of the balance across scales (Proposition~\ref{prop:bal_bet_dya}).
\item[$\delta_{\rm can}=\gamma_{\rm can}=1/2400000$, $C_{\rm can}=20000$] defect, rate and prefactor of the scarcity of cancellations (Theorem~\ref{thm:exp_sca_can}); $\delta_{\rm can}=10^{-3}(1-c_{\rm hyp})/24$.
\item[$c_{\rm hyp}=99/100$] hyperbolic-time rate (Definition~\ref{def:pai_gro_log}).
\item[$\eta_{\rm cov}=10^{-5}$] cover length for the successful sets of local tests (Theorem~\ref{thm:sma_cov_the}).
\item[$D_{\rm dbl}=2^{28}$] doubling factor (Corollary~\ref{cor:dou_ine_for}).
\item[$c_{\rm pre}=0.98$, $w_{\rm pre}=0.8$, $\alpha_{\rm pre}=3c_{\rm pre}-2=0.94$] prefix-growth rate, contraction exponent and direction exponent of the bridge argument (Definition~\ref{def:pre_gro_eli}, Lemma~\ref{lem:sta_gra_und:1}, Lemma~\ref{lem:mos_con_dir:1}).
\item[$\beta_{\rm br}=1/2400$, $\kappa_{\rm br}=1/20$, $\sigma_{\rm br}=1/10$, $C_{\rm br}=2^{32}$] shortfall allowance, window, rate and prefactor of the bridge product estimate (Theorem~\ref{thm:bri_pro_est:2}).
\item[$b_{\rm sl}=3\beta_{\rm br}/4=1/3200$, $C_{\rm cr}=10^6/\delta_{\rm can}$] slack margin and crossing constant of slack profiles (Definition~\ref{def:slac_prof}).
\item[$u_{\rm lo}=b_{\rm sl}/16$, $\alpha_{\rm cap}=1/(20C_{\rm cr})$] linear threshold and slope of the cap (Definition~\ref{def:the_cap_and}).
\item[$c_{\rm sav}=\alpha_{\rm cap}b_{\rm sl}u_{\rm lo}/C_{\rm cr}=\alpha_{\rm cap}b_{\rm sl}^2/(16C_{\rm cr})$] terminal saving (Lemma~\ref{lem:cor_ine_sla}).
\item[$A_{\rm sl}=4(C_{\rm cr}+1)/u_{\rm lo}$] onset ratio of the slow weight of the profile corrector (Definition~\ref{def:the_pro_cor}).
\item[$C_{\rm env}=2^{23}$, $C_{\rm dom}=2^{145}$] Lipschitz constant of the envelope of the averaged corrector (Definition~\ref{def:ext_cor_and}, Lemma~\ref{lem:mod_the_pro}), and constant of imbalance domination (Lemma~\ref{lem:spa_dom_the}).
\end{description}
\label{conv:not_and_con}
\end{convention}

\end{document}